\documentclass[12pt,a4paper,onecolumn]{article}
\usepackage{textcomp}
\usepackage{tikz}
\usepackage{amssymb}
\usepackage{caption}
\usepackage{color}
\usepackage{verbatim}
\usepackage{fancybox}
\usepackage{titlesec} 
\renewcommand\thesection{\arabic{section}} 
\usepackage{esvect}
\usepackage{cite}
\usepackage{setspace}
\usepackage{mathtools}
\usepackage{amsmath}
\usepackage{mathrsfs}
\usetikzlibrary{positioning,fit,calc}
\usepackage[a4paper,left=1in,top=1in,right=1in,bottom=1in,nohead]{geometry} 
\usepackage[hidelinks]{hyperref}
\usepackage{array}
\usepackage{bm}
\usepackage{multirow}
\usepackage{hhline}
\usepackage{amsthm}
\usepackage[bottom]{footmisc} 
\usepackage{makecell}
\numberwithin{equation}{section}
\usepackage{ragged2e}

\newcommand*\lineexample[1]{%
  #1~(\tikz[baseline=-\the\fontdimen22\textfont2]\draw[#1, line width=0.3mm](0,0)--(1.3em,0);)%
}

\def\correspondingauthor{\footnote{Corresponding author. Email: williewong088@gmail.com.}}

\tikzset{block/.style={draw,thick,text width=2cm,minimum height=1cm,align=center},
         line/.style={-latex}}
\newcolumntype{P}[1]{>{\centering\arraybackslash}m{#1}} 

\titleformat{\section}[block]{\large\scshape\bfseries}{\thesection.}{1em}{} 
\titleformat{\subsection}[block]{\bfseries}{\thesubsection.}{1em}{} 

\newtheorem{thm}{Theorem}[section]
\newtheorem{ppn}[thm]{Proposition}
\newtheorem{cor}[thm]{Corollary}
\newtheorem{lem}[thm]{Lemma}

\theoremstyle{definition}
\newtheorem{defn}[thm]{Definition}
\newtheorem{dis}[thm]{Discussion}
\newtheorem{rmk}[thm]{Remark}
\newtheorem{eg}[thm]{Example}

\newtheorem{prob}[thm]{Problem}

\begin{document}
\pagenumbering{arabic}
\begin{center}
    \textbf{\Large Application of some techniques in Sperner Theory: Optimal orientations of vertex-multiplications\\ of trees with diameter 4}
\vspace{0.1 in} 
    \\{\large H.W. Willie Wong\correspondingauthor{}, E.G. Tay}
\vspace{0.1 in} 
\\National Institute of Education\\Nanyang Technological University\\Singapore
\end{center}

\begin{abstract}
\noindent Koh and Tay proved a fundamental classification of $G$ vertex-multiplications into three classes $\mathscr{C}_0, \mathscr{C}_1$ and $\mathscr{C}_2$. They also showed that any vertex-multiplication of a tree with diameter at least 3 does not belong to the class $\mathscr{C}_2$. Of interest, $G$ vertex-multiplications are extensions of complete $n$-partite graphs and Gutin characterised complete bipartite graphs with orientation number 3 (or 4 resp.) via an ingenious use of Sperner's theorem. In this paper, we investigate vertex-multiplications of trees with diameter $4$ in $\mathscr{C}_0$ (or $\mathscr{C}_1$) and exhibit its intricate connections with problems in Sperner Theory, thereby extending Gutin's approach. Let $s$ denote the vertex-multiplication of the central vertex. We almost completely characterise the case of even $s$ and give a complete characterisation for the case of odd $s\ge 3$.
\end{abstract}
\section{Introduction}
Let $G$ be a graph with vertex set $V(G)$ and edge set $E(G)$. In this paper, we consider only graphs with no loops or parallel edges. For any vertices $v,x\in V(G)$, the $\textit{distance}$ from $v$ to $x$, $d_G(v,x)$, is defined as the length of a shortest path from $v$ to $x$. For $v\in V(G)$, its $\textit{eccentricity}$ $e_G(v)$ is defined as $e_G(v)=\max\{d_G(v,x) \mid x\in V(G)\}$. The $\textit{diameter}$ of $G$, denoted by $d(G)$, is defined as $d(G)=\max\{e_G(v)\mid v\in V(G)\}$. The above notions are defined similarly for a digraph $D$; and we refer the reader to \cite{BJJ GG} for any undefined terminology. For a digraph $D$, a vertex $x$ is said to be \textit{reachable} from another vertex $v$ if $d_D(v,x)<\infty$. The \textit{outset} and \textit{inset} of a vertex $v\in V(D)$ are defined to be $O_D(v)=\{x\in V(D)\mid v\rightarrow x\}$ and $I_D(v)=\{y\in V(D)\mid y\rightarrow v\}$ respectively. The \textit{outdegree} $\deg^+_D(v)$ and \textit{indegree} $\deg^-_D(v)$ of a vertex $v\in V(D)$ are defined by $\deg^+_D(v)=|O_D(v)|$ and $\deg^-_D(v)=|I_D(v)|$ respectively. If there is no ambiguity, we shall omit the subscript for the above notation.
\noindent\par An $\textit{orientation}$ $D$ of a graph $G$ is a digraph obtained from $G$ by assigning a direction to every edge $e\in E(G)$. An orientation $D$ of $G$ is said to be \textit{strong} if every two vertices in $V(D)$ are mutually reachable. An edge $e\in E(G)$ is a \textit{bridge} if $G-e$ is disconnected. Robbins' One-way Street Theorem \cite{RHE} states that a connected graph $G$ has a strong orientation if and only if $G$ is bridgeless.
\indent\par Given a connected and bridgeless graph $G$, let $\mathscr{D}(G)$ be the family of strong orientations of $G$. The $\textit{orientation number}$ of $G$ is defined as 
\begin{align*}
\bar{d}(G)=\min\{d(D)\mid D\in \mathscr{D}(G)\}.
\end{align*}
Any orientation $D$ in $\mathscr{D}(G)$ with $d(D)=\bar{d}(G)$ is called an \textit{optimal orientation} of $G$. The general problem of finding the orientation number of a connected and bridgeless graph is very difficult. Moreover, Chv{\'a}tal and Thomassen \cite{CV TC} proved that it is NP-hard to determine whether a graph admits an orientation of diameter 2. Hence, it is natural to focus on special classes of graphs. The orientation number was evaluated for various classes of graphs, such as the complete graphs \cite{BF TR,MSB,PJ 1} and complete bipartite graphs \cite{GG 1,SL}. Of interest, Gutin ingeniously made use of a celebrated result in combinatorics, Sperner's theorem (see Theorem \ref{thmC6.2.1}), to determine a characterisation of complete bipartite graphs with orientation number 3 (or 4 resp.).
\begin{thm} ({\v S}olt\'es \cite{SL} and Gutin \cite{GG 1})
\label{thmC6.1.1} For $q\ge p \ge 2$,
\begin{equation}
 \bar{d}(K(p,q))=\left\{
  \begin{array}{@{}ll@{}}
    3, & \text{if}\ q\le{{p}\choose{\lfloor{p/2}\rfloor}}, \nonumber\\
    4, & \text{if}\ q>{{p}\choose{\lfloor{p/2}\rfloor}}. \nonumber
  \end{array}\right.
\end{equation}
\end{thm}
In 2000, Koh and Tay \cite{KKM TEG 8} studied the orientation numbers of a family of graphs known as the $G$ vertex-multiplications. They extended the results on complete $n$-partite graphs. Let $G$ be a given connected graph with vertex set $V(G)=\{v_1,v_2,\ldots, v_n\}$. For any sequence of $n$ positive integers $(s_i)$, a $G$ \textit{vertex-multiplication} (also known as an \textit{extension} of $G$ in \cite{BJJ GG}), denoted by $G(s_1, s_2,\ldots, s_n)$, is the graph with vertex set $V^*=\bigcup_{i=1}^n{V_i}$ and edge set $E^*$, where $V_i$'s are pairwise disjoint sets with $|V_i|=s_i$, for $i=1,2,\ldots,n$; and for any $u,v\in V^*$, $uv\in E^*$ if and only if $u\in V_i$ and $v\in V_j$ for some $i,j\in \{1,2,\ldots, n\}$ with $i\neq j$ such that $v_i v_j\in E(G)$. For instance, if $G\cong K_n$, then the graph $G(s_1, s_2,\ldots, s_n)$ is a complete $n$-partite graph with partite sizes $s_1, s_2,\ldots, s_n$. Also, we say $G$ is a \textit{parent graph} of a graph $H$ if $H\cong G(s_1, s_2,\ldots, s_n)$ for some sequence $(s_i)$ of positive integers.
\noindent\par For $i=1,2,\ldots, n$, we denote the $x$-th vertex in $V_i$ by $(x,v_i)$, i.e., $V_i=\{(x,v_i)\mid x=1,2,\ldots,s_i\}$. Hence, two vertices $(x,v_i)$ and $(y,v_j)$ in $V^*$ are adjacent in $G(s_1,s_2,\ldots, s_n)$ if and only if $i\neq j$ and $v_i v_j\in E(G)$. For convenience, we write $G^{(s)}$ in place of $G(s,s,\ldots,s)$ for any positive integer $s$, and it is understood that the number of $s$'s is equal to the order of $G$, $n$. Thus, $G^{(1)}$ is simply the graph $G$ itself.
\indent\par The $G$ vertex-multiplications are a natural generalisation of complete multipartite graphs. Optimal orientations minimising the diameter can also be used to solve a variant of the Gossip Problem on a graph $G$. The Gossip Problem attributed to Boyd by Hajnal et al. \cite{HA MEC SE} is stated as follows:
\begin{justify}
``There are $n$ ladies, and each one of them knows an item of scandal which is not known to any of the others. They communicate by telephone, and whenever two ladies make a call, they pass on to each other, as much scandal as they know at that time.  How many calls are needed before all ladies know all the scandal?"
\end{justify}
\indent\par The Problem has been the source of many papers that have studied the spread of information by telephone calls, conference calls, letters and computer networks. One can imagine a network of people modelled by a $G$ vertex-multiplication where the parent graph is $G$ and persons within a partite set are not allowed to communicate directly with each other, for perhaps secrecy or disease containment reasons.
\indent\par The following theorem by Koh and Tay \cite{KKM TEG 8} provides a fundamental classification on $G$ vertex-multiplications.
\begin{thm} (Koh and Tay \cite{KKM TEG 8}) \label{thmC6.1.2} Let $G$ be a connected graph of order $n\ge 3$. If $s_i\ge 2$ for $i=1,2,\ldots, n$, then $d(G)\le \bar{d}(G(s_1,s_2,\ldots,s_n))\le d(G)+2$.
\end{thm}
In view of Theorem \ref{thmC6.1.2}, all graphs of the form $G(s_1,s_2,\ldots, s_n)$, with $s_i\ge 2$ for all $i=1,2,\ldots, n$, can be classified into three classes $\mathscr{C}_j$, where 
\begin{align*}
\mathscr{C}_j=\{G(s_1,s_2,\ldots, s_n)\mid\bar{d}(G(s_1,s_2,\ldots, s_n))=d(G)+j\},
\end{align*}
for $j=0,1,2$. Henceforth, we assume $s_i\ge 2$ for $i=1,2,\ldots, n$. The following lemma was found useful in proving Theorem \ref{thmC6.1.2}.
\begin{lem}(Koh and Tay \cite{KKM TEG 8}) \label{lemC6.1.3} Let $s_i,t_i$ be integers such that $s_i\le t_i$ for $i=1,2,\ldots, n$. If the graph $G(s_1,s_2,\ldots, s_n)$ admits an orientation $F$ in which every vertex $v$ lies on a cycle of length not exceeding $m$, then $\bar{d}(G(t_1, t_2,\ldots, t_n))\le \max\{m, d(F)\}$.
\end{lem}
To discuss further, we need some notation. In this paper, let $T_4$ (or simply $T$ unless stated otherwise) be a tree of diameter $4$ with vertex set $V(T_4)=\{v_1,v_2,\ldots, v_n\}$. We further denote by $\mathtt{c}$, the unique central vertex of $T_4$, i.e., $e_{T_4}(\mathtt{c})=2$, and the neighbours of $\mathtt{c}$ by $[i]$, i.e., $N_{T_4}(\mathtt{c})=\{[i]\mid i=1,2,\ldots, \deg_{T_4}(\mathtt{c}) \}$. For each $i =1,2,\ldots, \deg_{T_4}(\mathtt{c})$, we further denote the neighbours of $[i]$, excluding $\mathtt{c}$, by $[\alpha, i]$, i.e., $N_{T_4}([i])-\{\mathtt{c}\}=\{[\alpha, i] \mid \alpha=1,2,\ldots, \deg_{T_4}([i])-1\}$. Figure \ref{figC6.1.1} illustrates the use of this notation.
\begin{center}
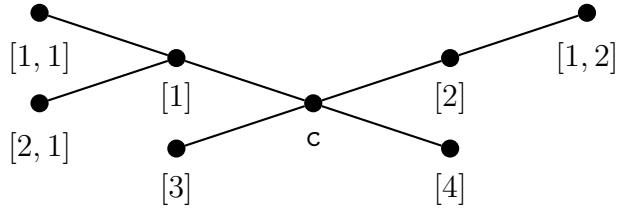

\tikzstyle{every node}=[circle, draw, fill=black!100,
                       inner sep=0pt, minimum width=6pt]
\begin{tikzpicture}[thick,scale=0.6]%
\draw(0,0)node[label={[yshift=-0.2cm]270:{$\mathtt{c}$}}](u){};

\draw(-6,2)node[label={[yshift=0cm] 270:{$[1,1]$}}](11){};
\draw(-6,0)node[label={[yshift=0cm] 270:{$[2,1]$}}](12){};
\draw(6,2)node[label={[yshift=0cm] 270:{$[1,2]$}}](21){};

\draw(-3,1)node[label={[yshift=-0.1cm] 270:{$[1]$}}](1){};
\draw(-3,-1)node[label={[yshift=-0.1cm] 270:{$[3]$}}](3){};
\draw(3,1)node[label={[yshift=-0.1cm] 270:{$[2]$}}](2){};
\draw(3,-1)node[label={[yshift=-0.1cm] 270:{$[4]$}}](4){};

\draw(u)--(1);
\draw(u)--(2);
\draw(u)--(3);
\draw(u)--(4);

\draw(1)--(11);
\draw(1)--(12);
\draw(2)--(21);
\end{tikzpicture}
{\captionof{figure}{Labelling vertices in a $T_4$}\label{figC6.1.1}}
\end{center}
\indent\par From here onwards, let $\mathcal{T}=T_4(s_1,s_2,\ldots, s_n)$ be a vertex-multiplication of a tree $T_4$. In $\mathcal{T}$, the integer $s_i$ corresponds to the vertex $v_i$, $i\neq n$, while $s_n=s$ corresponds to $\mathtt{c}$. We will loosely use the two denotations of a vertex, for example, if $v_i=[j]$, then $s_i=s_{[j]}$. Also, if $X\subseteq \mathbb{N}_k=\{1,2,\ldots,k\}$, where $k\in\mathbb{Z}^+$, and $v\in V(T)$, then set $(X,v)=\{(i,v)\mid i\in X\}$. In particular, $(\mathbb{N}_s,\mathtt{c})=\{(1,\mathtt{c}),(2,\mathtt{c})\ldots, (s,\mathtt{c})\}$. For any set $\lambda\subseteq (\mathbb{N}_s,\mathtt{c})$, $\bar{\lambda}=(\mathbb{N}_s,\mathtt{c})-\lambda$ denotes its complement set.
\indent\par A vertex $v$ in a graph $G$ is a \textit{leaf} if $deg_G(v)=1$. For a given $T_4$, set $E=\{[i]\mid [i]$ is an leaf in $T_4\}$. For a given $\mathcal{T}$ of $T_4$, set $\mathcal{T}(A_j)=\{[i]\mid s_{[i]}=j,\ 1\le i\le \deg_T(\mathtt{c})$, $[i]\not\in E\}$, where $j$ is a positive integer. If there is no ambiguity, we will use $A_j$ instead of $\mathcal{T}(A_j)$. Similarly, $A_{\le j}$ and $A_{\ge j}$ denote the corresponding sets, when the condition $s_{[i]}=j$ is replaced by $s_{[i]}\le j$ and $s_{[i]}\ge j$ respectively.  For example, if $T_4$ is as given in Figure \ref{figC6.1.1}, then $E=\{[3], [4]\}$; furthermore, if $s_i=2$ for all $i=1,2,\ldots, n$, in $\mathcal{T}$, then $A_2=\{[1],[2]\}$.
\noindent\par Theorem \ref{thmC6.1.2} was generalised to digraphs by Gutin et al. \cite{GG KKM TEG AY}. Ng and Koh \cite{NKL KKM} and Wong and Tay \cite{WHW TEG 8A} investigated vertex-multiplications of cycles and Cartesian products of graphs respectively. Koh and Tay \cite{KKM TEG 11} studied vertex-multiplications of trees. Since trees with diameter at most 2 are parent graphs of complete bipartite graphs and are completely solved, Koh and Tay considered trees of diameter at least 3. They proved that vertex-multiplications of trees with diameter 3, 4 or 5 does not belong to the class $\mathscr{C}_2$ and those with diameter at least 6 belong to the class $\mathscr{C}_0$.
\begin{thm}(Koh and Tay \cite{KKM TEG 11})\label{thm:CC2.2.6} If $T$ is a tree of order $n$ and $d(T)=3, 4$ or $5$, then $T(s_1,s_2,\ldots, s_n)\in \mathscr{C}_0\cup \mathscr{C}_1$.
\end{thm}
\begin{thm}(Koh and Tay \cite{KKM TEG 11})\label{thm:CC2.2.7}
If $T$ is a tree of order $n$ and $d(T)\ge 6$, then $T(s_1,s_2,\ldots, s_n)\in \mathscr{C}_0$.
\end{thm}
Wong and Tay \cite{WHW TEG 3} proved a characterisation for vertex-multiplications of trees with diameter $5$ in $\mathscr{C}_0$ and $\mathscr{C}_1$.
\begin{thm}(Wong and Tay \cite{WHW TEG 3})
Let $T$ be a tree of diameter $5$ with central vertices $\mathtt{c}_1$ and $\mathtt{c}_2$, and $s_i$ corresponding to $\mathtt{c}_i$ for $i=1,2$. Furthermore for $i=1,2$, denote $E'_i=\{u\mid u\in N_T(\mathtt{c}_i)-\{\mathtt{c}_{3-i}\}, u$ is not an leaf in $T\}$ and $m_i=\min\{s_{u}\mid u\in E'_i \}$.
\\(a) If $s_{1}\ge 3$, or $s_{2}\ge 3$, or $m_1, m_2\ge 4$, then $T(s_1,s_2,\ldots, s_n)\in \mathscr{C}_0$.
\\(b) Suppose $s_1=s_2=2$ and $2\le m_{1}\le 3$ or $2\le m_{2}\le 3$. Then, $T(s_1,s_2\ldots, s_n)\in \mathscr{C}_0$ if and only if $|E'_j|=1 \text{ for some } j=1,2$.
\end{thm}
Koh and Tay \cite{KKM TEG 11} obtained some results regarding membership in $\mathscr{C}_0$ or $\mathscr{C}_1$ for vertex-multiplications of trees with diameter 4.
\begin{thm}\label{thmC6.1.6} (Koh and Tay \cite{KKM TEG 11}) For a given $T_4$,
\\(a) if $\deg_{T_4}(\mathtt{c})=2$, then $\mathcal{T}\in \mathscr{C}_0$.
\\(b) if $\deg_{T_4}(\mathtt{c})\ge3$, then $T_4^{(2)}\in \mathscr{C}_1$.
\end{thm}
In this paper, we further investigate vertex-multiplications of trees with diameter $4$ and almost completely classify them as $\mathscr{C}_0$ or $\mathscr{C}_1$. The techniques required here exhibit intricate connections with problems in Sperner Theory. In Section 2, we provide the main tools, which comprise well-known results from Sperner Theory and structural properties of optimal orientations of a $\mathcal{T}$. Section 3 focuses primarily on the case where $s$ is even and the findings are summarised in Theorem \ref{thmC6.1.7}. In Section 4, we prove a complete characterisation of vertex-multiplications of trees with odd $s\ge 3$, namely Theorem \ref{thmC6.1.8}.

\begin{thm}\label{thmC6.1.7}
Let $T_4$ be a tree of diameter $4$ with the only central vertex $\mathtt{c}$. Suppose $s$ is even for a $\mathcal{T}$. Then,
\end{thm}
(a) For $s=2$:
\begin{spacing}{1.23}
\centering
\begin{tabular}{| P{2cm} | P{1cm} | m{9.2cm}| P{2.3cm}|}
\hline
$\bm{|A_2\cup A_3|}$ & $\bm{|A_{\ge 4}|}$ & $\bm{\mathcal{T}\in\mathscr{C}_0\iff\ldots}$ & \textbf{Proposition}\\
\hline
0 & $\ge 2$ & Always true. & \ref{ppnC6.3.4} \\ 
\hline\hline
$\ge 1$ & $\ge 0$ & $\deg_T(\mathtt{c})=2$. & \ref{ppnC6.3.2} \\ 
\hline
\end{tabular}
{\captionof{table}{Summary for $\mathcal{T}=T_4(s_1, s_2, \ldots, s_n)$, where $s=2$.}}
\end{spacing}
(b) For $s\ge 4$:
\begin{spacing}{1.23}
\centering
\begin{tabular}{| P{0.7cm} | P{0.7cm} | P{1cm} | m{9.4cm}| P{2.3cm}|}
\hline
$\bm{|A_2|}$ & $\bm{|A_3|}$ & $\bm{|A_{\ge 4}|}$ & $\bm{\mathcal{T}\in\mathscr{C}_0\iff\ldots}$ & \textbf{Proposition}\\
\hline
0 & 0 & $\ge 2$ & Always true. & \ref{ppnC6.3.4}\protect\footnotemark\\
\hline\hline
0 & $\ge 1$ & $\ge 0$ & $|A_3|\le {{s}\choose{s/2}}+{{s}\choose{(s/2)+1}}-2$. &  \ref{ppnC6.3.9} \\ 
\hline\hline
$\ge 2$ & 0 & 0 & 
\makecell[l]{(i) $|A_2|\le{{s}\choose{\lceil{s/2}\rceil}}-1$, if $|A_2|< \deg_T(\mathtt{c})$,
\\ (ii) $|A_2|\le{{s}\choose{\lceil{s/2}\rceil}}$, if $|A_2|=\deg_T(\mathtt{c})$.}
& \ref{ppnC6.3.5}\protect\footnotemark\\ 
\hline\hline
$\ge 1$ & 0 & $\ge 1$ & 
\makecell[l]{(i) $|A_2|\le {{s}\choose{s/2}}-2$, if $|A_{\ge 4}|\ge 2$ or $|A_{\ge 2}|<\deg_T(\mathtt{c})$,
\\ (ii) $|A_2|\le{{s}\choose{s/2}}-1$, otherwise.}
& \ref{ppnC6.3.10}\\
\hline
\hline
$\ge 1$ & 1 & 0 & 
\makecell[l]{(i) $|A_2|\le{{s}\choose{s/2}}-2$ if $|A_{\ge 2}|< \deg_T(\mathtt{c})$,
\\(ii) $|A_2|\le{{s}\choose{s/2}}-1$, if $|A_{\ge 2}|=\deg_T(\mathtt{c})$.}
 & \ref{ppnC6.3.11}\\
\hline
\end{tabular}
\addtocounter{footnote}{-1}
\footnotetext{holds for all integers $s\ge 2$.}
\addtocounter{footnote}{1}
\footnotetext{holds for all integers $s\ge 3$.}
\newpage
The final case is incomplete and excludes the case $|A_2|\ge 1$, $|A_3|=1$ and $|A_{\ge 4}|=0$.
\begin{tabular}{| P{0.7cm} | P{0.7cm} | P{1cm} | m{9.4cm}| P{2.3cm}|}
\hline
$\bm{|A_2|}$ & $\bm{|A_3|}$ & $\bm{|A_{\ge 4}|}$ & $\bm{\mathcal{T}\in\mathscr{C}_0}$ & \textbf{Proposition}\\
\hline
$\ge 1$ & $\ge 1$ & $\ge 0$ & 
\makecell[l]{(a) $\mathcal{T}\in \mathscr{C}_0\Rightarrow 2|A_2|+|A_3|\le {{s}\choose{s/2}}+{{s}\choose{(s/2)+1}}$
\\$-\kappa_{s,\frac{s}{2}}^*(k)$ for some $k\le |A_2|+|A_3|$.
\\(b) There exists some 
\\$|A_2|+1\le k\le \min\{|A_2|+|A_3|, {{s}\choose{s/2}}-1\}$ such that 
\\$2|A_2|+|A_3|\le {{s}\choose{s/2}}+{{s}\choose{(s/2)+1}}-\kappa_{s,\frac{s}{2}}(k)-3$
\\$\Rightarrow \mathcal{T}\in \mathscr{C}_0$.
\\Note: $\kappa_{s,\frac{s}{2}}(\cdot)$ and $\kappa_{s,\frac{s}{2}}^*(\cdot)$ will be defined later.}
& \ref{ppnC6.3.12}\\
\hline
\end{tabular}
{\captionof{table}{Summary for $\mathcal{T}=T_4(s_1, s_2, \ldots, s_n)$, where $s\ge 4$ is even.}}
\end{spacing}

\begin{thm}\label{thmC6.1.8}
Let $T_4$ be a tree of diameter $4$ with the only central vertex $\mathtt{c}$. Suppose $s\ge 3$ is odd for a $\mathcal{T}$. Then,
\end{thm}
\begin{spacing}{1.23}
\begin{centering}
\begin{tabular}{| P{0.8cm} | P{0.8cm} | P{1.1cm} | m{8.2cm}| P{2.4cm}|}
\hline
$\bm{|A_2|}$ & $\bm{|A_3|}$ & $\bm{|A_{\ge 4}|}$ & $\bm{\mathcal{T}\in\mathscr{C}_0\iff \ldots}$ & \textbf{Proposition}\\
\hline
$0$ & $\ge 1$ & $\ge 0$ & $|A_3|\le 2{{s}\choose{\lceil{s/2}\rceil}}-2$. & \ref{ppnC6.4.1}\\ 
\hline\hline
$0$ & $0$ & $\ge 2$ & Always true. & \ref{ppnC6.3.4}\protect\footnotemark\\ 
\hline\hline
$\ge 2$ & $0$ & $0$ & 
\makecell[l]{(i) $|A_2|\le{{s}\choose{\lceil{s/2}\rceil}}-1$, if $|A_2|< \deg_T(\mathtt{c})$,
\\ (ii) $|A_2|\le{{s}\choose{\lceil{s/2}\rceil}}$, if $|A_2|=\deg_T(\mathtt{c})$.}
& \ref{ppnC6.3.5}\protect\footnotemark\\ 
\hline\hline
$\ge 1$ & $0$ & $\ge 1$ & $|A_2|\le {{s}\choose{\lceil{s/2}\rceil}}-1$. & \ref{ppnC6.4.11}\\
\hline\hline
$\ge 1$ & $1$ & $0$ & $|A_2|\le {{s}\choose{\lceil{s/2}\rceil}}-1$. & \multirow{2}{*}{\ref{ppnC6.4.3}} \\
\hhline{----~}
$\ge 1$ & $\ge 2$ & $0$ &
\makecell[l]{(i) $2|A_2|+|A_3|\le 2{{s}\choose{\lceil{s/2}\rceil}}-2$, or
\\(ii) $2|A_2|+|A_3|= 2{{s}\choose{\lceil{s/2}\rceil}}-1$, 
\\$|A_2|\ge \lceil\frac{s}{2}\rceil \lfloor\frac{s}{2}\rfloor$ and $s\ge 5$.} &
\\
\hline\hline
 $\ge 1$ & $1$ & $\ge 1$ & $|A_2|\le {{s}\choose{\lceil{s/2}\rceil}}-2$. & \multirow{2}{*}{\ref{ppnC6.4.13}} \\
\hhline{----~}
$\ge 1$ & $\ge 2$ & $\ge 1$ & $2|A_2|+|A_3|\le 2{{s}\choose{\lceil{s/2}\rceil}}-2$. &\\
\hline
\end{tabular}
{\captionof{table}{Summary for $\mathcal{T}=T_4(s_1, s_2, \ldots, s_n)$, where $s\ge 3$ is odd.}}
\end{centering}
\end{spacing}
\addtocounter{footnote}{-1}
\footnotetext{holds for all integers $s\ge 2$.}
\addtocounter{footnote}{1}
\footnotetext{holds for all integers $s\ge 3$.}

\noindent\par As we shall see in the proofs of Theorems \ref{thmC6.1.7} and \ref{thmC6.1.8}, it is a key insight to partition $N_T(\mathtt{c})$ into 4 types of vertices, $A_2$, $A_3$, $A_{\ge 4}$ and $E$. Their sizes will then determine the equivalent conditions of an optimal orientation (except possibly Proposition \ref{ppnC6.3.12}). We shall consider cases dependent on these 4 sets. The lack of conformity in the equivalent conditions across all cases gives a compelling indication that the case distinctions are required.
\section{Preliminaries}
Our overarching approach is to reduce the investigation of optimal orientations of tree vertex-multiplication graphs to variants of problems in Sperner Theory; particularly concerning cross-intersecting antichains. The change in perspective grants us leverage on the following useful results in Sperner Theory.
\indent\par For any $n\in \mathbb{Z}^+$, let $\mathbb{N}_n=\{1,2,\ldots, n\}$ and $2^{\mathbb{N}_n}$ denote the power set of $\mathbb{N}_n$. For any integer $k$, $0\le k\le n$, ${{\mathbb{N}_n}\choose{k}}$ denotes the collection of all $k$-subsets (i.e., subsets of cardinality $k$) of $\mathbb{N}_n$. Two families $\mathscr{A},\mathscr{B}\subseteq 2^{\mathbb{N}_n}$ are said to be \textit{cross-intersecting} if $A\cap B\neq\emptyset$ for all $A\in\mathscr{A}$ and all $B\in\mathscr{B}$. Two subsets $X$ and $Y$ of $\mathbb{N}_n$ are said to be \textit{independent} if $X\not\subseteq Y$ and $Y\not\subseteq X$. An \textit{antichain} or \textit{Sperner family} $\mathscr{A}$ on $\mathbb{N}_n$ is a collection of pairwise independent subsets of $\mathbb{N}_n$, i.e., for all $X,Y\in \mathscr{A}$, $X\not\subseteq Y$.
\begin{thm}(Sperner \cite{SE})\label{thmC6.2.1}
For any $n\in \mathbb{Z}^+$, if $\mathscr{A}$ is an antichain on $\mathbb{N}_n$, then $|\mathscr{A}|\le {{n}\choose{\lfloor{n/2}\rfloor}}$. Furthermore, equality holds if and only if all members in $\mathscr{A}$ have the same size, ${\lfloor\frac{n}{2}\rfloor}$ or ${\lceil{\frac{n}{2}}\rceil}$.
\end{thm}
Lih's theorem \cite{LKW} provides the maximum size of an antichain with each member intersecting a fixed set and Griggs \cite{GJR} determined all such maximum-sized antichains.
\begin{thm}(Lih \cite{LKW})
Let $n\in \mathbb{Z}^+$ and $Y\subseteq \mathbb{N}_n$. If $\mathscr{A}$ is an antichain on $\mathbb{N}_n$ such that $A\cap Y\neq \emptyset$ for all $A\in \mathscr{A}$, then 
\begin{align*}
|\mathscr{A}|\le {{n}\choose{\lceil{n/2}\rceil}}-{{n-|Y|}\choose{\lceil{n/2}\rceil}}.
\end{align*}
\end{thm}
\begin{thm}(Griggs \cite{GJR}) Let $n\in \mathbb{Z}^+$ and $Y\subseteq \mathbb{N}_n$. If $\mathscr{A}$ is an antichain on $\mathbb{N}_n$ such that $A\cap Y\neq \emptyset$ for all $A\in \mathscr{A}$ and $|\mathscr{A}|={{n}\choose{\lceil{n/2}\rceil}}-{{n-|Y|}\choose{\lceil{n/2}\rceil}}$, then $\mathscr{A}$ consists of exactly one of the following:
\\(i) $\lceil\frac{n}{2}\rceil$-sets, or
\\(ii) $\frac{n-1}{2}$-sets for odd $n$ and $|Y|\ge \frac{n+3}{2}$, or
\\(iii) $\frac{n+2}{2}$-sets for even $n$ and $|Y|=1$.
\end{thm}
Given two cross-intersecting antichains $\mathscr{A}$ and $\mathscr{B}$ on $\mathbb{N}_n$, Ou \cite{OY}, Frankl and Wong \cite{FP WHW} and Wong and Tay \cite{WHW TEG 5} independently derived an upper bound for $|\mathscr{A}|+|\mathscr{B}|$. Furthermore, Wong and Tay \cite{WHW TEG 4} determined all extremal and almost-extremal cross-intersecting antichains for $\mathscr{A}$ and $\mathscr{B}$.
\begin{thm}(Ou \cite{OY}, Frankl and Wong \cite{FP WHW} and Wong and Tay \cite{WHW TEG 5})\label{thmC6.2.4} Let $\mathscr{A}$ and $\mathscr{B}$ be two cross-intersecting antichains on $\mathbb{N}_n$, where $n\in\mathbb{Z}^+$ and $n\ge 3$. Then, 
\begin{align*}
|\mathscr{A}|+|\mathscr{B}|\le {{n}\choose{\lfloor{(n+1)/2}\rfloor}}+{{n}\choose{\lceil{(n+1)/2}\rceil}}
\end{align*}
Furthermore, equality holds if and only if $\{\mathscr{A},\mathscr{B}\}=\{{{\mathbb{N}_n}\choose{\lfloor{(n+1)/2}\rfloor}}, {{\mathbb{N}_n}\choose{\lceil{(n+1)/2}\rceil}}\}$.
\end{thm}

\begin{thm}(Wong and Tay \cite{WHW TEG 4})\label{thmC6.2.5}
Let $\mathscr{A}$ and $\mathscr{B}$ be two cross-intersecting antichains on $\mathbb{N}_n$, where $n\ge 3$ is an odd integer and $|\mathscr{A}|\ge |\mathscr{B}|$. Then, $|\mathscr{A}|+|\mathscr{B}|=2{{n}\choose{\lceil{n/2}\rceil}}-1$ if and only if $\mathscr{A}={{\mathbb{N}_n}\choose{\lceil{n/2}\rceil}}$, $\mathscr{B}\subset{{\mathbb{N}_n}\choose{\lceil{n/2}\rceil}}$ and $|\mathscr{B}|={{n}\choose{\lceil{n/2}\rceil}}-1$.
\end{thm}

\begin{thm}(Wong and Tay \cite{WHW TEG 4})\label{thmC6.2.6} Let $\mathscr{A}$ and $\mathscr{B}$ be two cross-intersecting antichains on $\mathbb{N}_n$, where $n\ge 4$ is an even integer and $|\mathscr{A}|\ge |\mathscr{B}|$. Then, $|\mathscr{A}|+|\mathscr{B}|={{n}\choose{n/2}}+{{n}\choose{(n/2)+1}}-1$ if and only if
\\(i) $\mathscr{A}={{\mathbb{N}_n}\choose{n/2}}$, $\mathscr{B}\subset{{\mathbb{N}_n}\choose{(n/2)+1}}$ and $|\mathscr{B}|={{n}\choose{(n/2)+1}}-1$, or
\\(ii) $\mathscr{A}\subset{{\mathbb{N}_n}\choose{n/2}}$, $|\mathscr{A}|={{n}\choose{n/2}}-1$, and $\mathscr{B}={{\mathbb{N}_n}\choose{(n/2)+1}}$.
\end{thm}
Kruskal-Katona Theorem (KKT) is closely related to the \textit{squashed order} of the $k$-sets. The squash relations $\le_s$ and $<_s$ are defined as follows. For $A,B\in {{\mathbb{N}_n}\choose{k}}$, $A\le_s B$ if the largest element of the symmetric difference $(A-B)\cup (B-A)$ is in $B$. Furthermore, denote $A <_s B$ if $A\le_s B$ and $A\neq B$. For example, the $3$-subsets of $\mathbb{N}_5$ in squashed order are: $\bm{123} <_s \bm{124} <_s \bm{134} <_s \bm{234} <_s \bm{125} <_s \bm{135} <_s \bm{235} <_s \bm{145} <_s \bm{245} <_s \bm{345}$. Here, we omit the braces and write $\bm{abc}$ to represent the set $\{a,b,c\}$, if there is no ambiguity. We shall denote the collections of the first $m$ and last $m$ $k$-subsets of $\mathbb{N}_n$ in squashed order by $F_{n,k}(m)$ and $L_{n,k}(m)$ respectively.
\indent\par For a family $\mathscr{A}\subseteq {{\mathbb{N}_n}\choose{k}}$, the \textit{shadow} and \textit{shade} of $\mathscr{A}$ are defined as
\begin{align*}
&\Delta \mathscr{A}=\{ X\subseteq \mathbb{N}_n\mid |X|=k-1, X\subset Y\text{ for some } Y\in \mathscr{A}\},\text{ if }k>0,\text{ and }\\
&\nabla \mathscr{A}=\{ X\subseteq \mathbb{N}_n\mid |X|=k+1, Y\subset X\text{ for some }Y\in \mathscr{A}\},\text{ if }k<n 
\end{align*}
respectively.
\indent\par Then, KKT says that the shadow of a family $\mathscr{A}$ of $k$-sets has size at least that of the shadow of the first $|\mathscr{A}|$ $k$-sets in squashed order.
\begin{thm}(Kruskal \cite{KJB}, Katona \cite{GOHK} and Clements and Lindstr{\"o}m \cite{CGF LB})
Let $\mathscr{A}$ be a collection of $k$-sets of $\mathbb{N}_n$ and suppose the $k$-binomial representation of $|\mathscr{A}|$ is $|\mathscr{A}|={{a_k}\choose{k}}+{{a_{k-1}}\choose{k-1}}+\ldots+{{a_t}\choose{t}}$, where $a_k>a_{k-1}>\ldots>a_t\ge t\ge 1$. Then,
\begin{align*}
&|\Delta \mathscr{A}|\ge|\Delta F_{n,k}(|\mathscr{A}|)|={{a_k}\choose{k-1}}+{{a_{k-1}}\choose{k-2}}+\ldots+{{a_t}\choose{t-1}}.
\end{align*}
\end{thm}
\noindent\par By considering the complements of sets in $F_{n,k}(m)$, the next lemma can be proved.
\begin{lem} \label{lemC6.2.8}
For any integer $0\le m\le {{n}\choose{k}}$, $|\Delta F_{n,k}(m)|=|\nabla L_{n,n-k}(m)|$.
\end{lem}
\begin{defn}(Wong and Tay \cite{WHW TEG 5})
Let $n,r$ and $m$ be integers such that $0\le m \le {{n}\choose{r}}$. Define
\begin{align*}
&\kappa_{n,r}(m)=|\Delta F_{n,r}(m)|-m \text{ and }\kappa_{n,r}^*(m)=\min\limits_{0\le j\le m}\kappa_{n,r}(j).
\end{align*}
\end{defn}
We remark that $\kappa_{n,\frac{n}{2}}(m)=|\nabla L_{n,\frac{n}{2}}(m)|-m$ by Lemma \ref{lemC6.2.8}. Using KKT, Wong and Tay \cite{WHW TEG 5} derived an upper bound for cross-intersecting antichains with at most $k$ disjoint pairs.
\begin{thm}(Wong and Tay \cite{WHW TEG 5})\label{thmC6.2.10}
Let $n\ge 4$ be an even integer and $\mathscr{A}$ and $\mathscr{B}$ be two antichains on $\mathbb{N}_n$. Suppose there exist orderings of the elements $A_1,A_2,\ldots, A_{|\mathscr{A}|}$ in $\mathscr{A}$, and $B_1,B_2,\ldots, B_{|\mathscr{B}|}$ in $\mathscr{B}$, and some integer $k\le \min\{|\mathscr{A}|, |\mathscr{B}|\}$, such that $A_i\cap B_j=\emptyset$ only if $i=j\le k$. Then,  
\begin{align*}
|\mathscr{A}|+|\mathscr{B}|\le {{n}\choose{n/2}}+{{n}\choose{(n/2)+1}}-\kappa_{n,\frac{n}{2}}^*(k),
\end{align*}
where $\kappa_{n,\frac{n}{2}}^*(k)=0$ if $k<1+\sum\limits_{i=1}^{n/2}{{2i-1}\choose{i}}$ and $\kappa_{n,\frac{n}{2}}^*(k)<0$ otherwise.
Furthermore, equality holds if 
\\(i) $k<1+\sum\limits_{i=1}^{n/2}{{2i-1}\choose{i}}$, $\mathscr{A}={{\mathbb{N}_n}\choose{n/2}}$ and $\mathscr{B}={{\mathbb{N}_n}\choose{(n/2)+1}}$, or 
\\(ii) $k\ge 1+\sum\limits_{i=1}^{n/2}{{2i-1}\choose{i}}$, $\mathscr{A}={{\mathbb{N}_n}\choose{n/2}}$ and $\mathscr{B}=L_{n,\frac{n}{2}}(m)\cup {{\mathbb{N}_n}\choose{(n/2)+1}}-\nabla L_{n,\frac{n}{2}}(m)$, where $0<m\le k$ is an integer such that $\kappa_{n,\frac{n}{2}}^*(k)=\kappa_{n,\frac{n}{2}}(m)$.
\end{thm}
Lastly, we prove here some key properties of an optimal orientation of a vertex-multiplication graph $\mathcal{T}$ in $\mathscr{C}_0$. Let $D$ be an orientation of $\mathcal{T}$. If $v_p$ and $v_q$, $1\le p, q\le n$ and $p\neq q$, are adjacent vertices in the parent graph $G$, then for each $i$, $1\le i \le s_p$, we denote by $O_D^{v_q}((i,v_p))=\{(j,v_q)\mid (i,v_p) \rightarrow (j,v_q), 1\le j \le s_q\}$ and $I_D^{v_q}((i,v_p))=\{(j,v_q)\mid (j,v_q)\rightarrow (i,v_p), 1\le j \le s_q\}$. If there is no ambiguity, we shall omit the subscript $D$ for the above notation.
\indent\par The next lemma is important but easy to verify.
\begin{lem}\textbf{(Duality)}
Let $D$ be an orientation of a graph $G$. Let $\tilde{D}$ be the orientation of $G$ such that $uv\in A(\tilde{D})$ if and only if $vu \in A(D)$. Then, $d(\tilde{D})=d(D)$.
\end{lem}
\begin{lem}\label{lemC6.2.12}
Let $D$ be an orientation of a $\mathcal{T}$ where $d(D)=4$. Then, $d_D((p,[\alpha,i]),(q,[j]))$ $=d_D((q,[j]),(p,[\alpha,i]))=3$ for all $1\le i, j \le \deg_T(\mathtt{c})$, $i\neq j$, $1\le \alpha \le \deg_T([i])-1$, $1\le p\le s_{[\alpha,i]}$ and $1\le q\le s_{[j]}$.
\end{lem}
\noindent\textit{Proof}: Note that $3=d_T([\alpha,i],[j])\le d_D((p,[\alpha,i]),(q,[j]))\le d(D)=4$. Since there is no $[\alpha, i]-[j]$ path of even length in $T$, there is no $(p,[\alpha, i])-(q,[j])$ path of even length in $\mathcal{T}$, in particular, no path of length 4. Hence, $d_D((p,[\alpha,i]),(q,[j]))=3$. Similarly, $d_D((q,[j]),(p,[\alpha,i]))=3$ may be proved.
\qed

\indent\par Since we are going to use this fact repeatedly, we state the following obvious lemma.
\begin{lem}\label{lemC6.2.13}
Let $D$ be an orientation of a $\mathcal{T}$ where $d(D)=4$. For $1\le i \le \deg_T(\mathtt{c})$,
\\(a) if $s_{[i]}=2$, then for $1\le p \le s_{[\alpha,i]}$ and $1\le \alpha \le \deg_T([i])-1$, either $(2,[i])\rightarrow (p,[\alpha,i])\rightarrow (1,[i])$ or $(1,[i])\rightarrow (p,[\alpha,i])\rightarrow (2,[i])$.
\\(b) if $s_{[i]}=3$, then for $1\le p \le s_{[\alpha,i]}$ and $1\le \alpha \le \deg_T([i])-1$, either $|O((p,[\alpha,i]))|=1$ or $|I((p,[\alpha,i]))|=1$.
\end{lem}
\noindent\textit{Proof}: Both statements follow from the fact that $O((p,[\alpha,i]))\neq\emptyset$ and $I((p,[\alpha,i]))\neq\emptyset$ for all $p=1,2,\ldots,s_{[\alpha,i]}$ so that $D$ is a strong orientation.
\qed

\begin{eg}
To help the reader understand the following lemmas and the proof of Proposition \ref{ppnC6.3.9}, we use the orientation $D$ shown in Figures \ref{figC6.3.7} and \ref{figC6.3.8} (see pages 27-28) for this example. It will be shown later that $d(D)=4$.
\\(a) Observe that $O((1,[1,i]))=\{(3,[i])\}$ for $i=5,6$, and $O^\mathtt{c}((3,[5]))=\{(1,\mathtt{c}), (2,\mathtt{c}), (3,\mathtt{c})\}$ and $O^\mathtt{c}((3,[6]))=\{(2,\mathtt{c}), (3,\mathtt{c}), (4,\mathtt{c})\}$ are independent.
\\(b) Note that $O((1,[1,1]))=\{(1,[1]), (2,[1])\}$, $O^\mathtt{c}((1,[1]))=\{(1,\mathtt{c}), (2,\mathtt{c})\}$, $O^\mathtt{c}((2,[1])) =\{(3,\mathtt{c}), (4,\mathtt{c})\}$. It is easy to check that $O^\mathtt{c}((3,[5])) \not\subseteq O^\mathtt{c}((p,[1]))$ for $p=1,2$, and $O^\mathtt{c}((1,[1]))\cup O^\mathtt{c}((2,[1])) \not\subseteq O^\mathtt{c}((3,[5]))$.
\\ In Lemmas \ref{lemC6.2.15}-\ref{lemC6.2.18}, we prove that these observations and their duals hold generally.
\end{eg}

\begin{lem}\label{lemC6.2.15}
Let $D$ be an orientation of a $\mathcal{T}$ where $d(D)=4$ and $1\le i, j\le \deg_T(\mathtt{c})$, $i\neq j$, $1\le \alpha \le \deg_T([i])-1$, $1\le \beta \le \deg_T([j])-1$. Suppose $O^\mathtt{c}(u_i)=O^\mathtt{c}(v_i)$ for any $u_i,v_i\in O((1,[\alpha,i]))$ and $O^\mathtt{c}(u_j)=O^\mathtt{c}(v_j)$ for any $u_j,v_j\in O((1,[\beta,j]))$, then $O^\mathtt{c}(w_i)$ and $O^\mathtt{c}(w_j)$ are independent for any $w_i\in O((1,[\alpha,i]))$ and $w_j\in O((1,[\beta,j]))$.
\end{lem}
\noindent\textit{Proof}: By Lemma \ref{lemC6.2.12}, $d_D((1,[\alpha,i]), w_j)=3$. Hence, it follows that $d_D(w_i, w_j)=2$ and $O^\mathtt{c}(w_i)\not\subseteq O^\mathtt{c}(w_j)$. A similar argument shows $O^\mathtt{c}(w_j)\not\subseteq O^\mathtt{c}(w_i)$.
\qed

\begin{lem}\label{lemC6.2.16}
Let $D$ be an orientation of a $\mathcal{T}$ where $d(D)=4$ and $1\le i, j\le \deg_T(\mathtt{c})$, $i\neq j$, $1\le \alpha \le \deg_T([i])-1$, $1\le \beta \le \deg_T([j])-1$. Suppose $I^\mathtt{c}(u_i)=I^\mathtt{c}(v_i)$ for any $u_i,v_i\in I((1,[\alpha,i]))$ and $I^\mathtt{c}(u_j)=I^\mathtt{c}(v_j)$ for any $u_j,v_j\in I((1,[\beta,j]))$, then $I^\mathtt{c}(w_i)$ and $I^\mathtt{c}(w_j)$ are independent for any $w_i\in I((1,[\alpha,i]))$ and $w_j\in I((1,[\beta,j]))$.
\end{lem}
\noindent\textit{Proof}: This lemma follows from Lemma \ref{lemC6.2.15} and the Duality Lemma.
\qed

\begin{lem}\label{lemC6.2.17}
Let $D$ be an orientation of a $\mathcal{T}$ where $d(D)=4$. Suppose $O((1,[\alpha,i]))=\{(1,[i])\}$ and $O((1,[\beta,j]))=\{(1,[j]), (2,[j])\}$ for $1\le i, j\le \deg_T(\mathtt{c})$, $i\neq j$, $1\le \alpha \le \deg_T([i])-1$, and $1\le \beta \le \deg_T([j])-1$. Then, for each $p=1,2,\ldots, s_{[j]}$,
\\(a) $O^\mathtt{c}((1,[i]))\not\subseteq O^\mathtt{c}((p,[j]))$,
\\(b) $O^\mathtt{c}((1,[j]))\cup O^\mathtt{c}((2,[j]))\not\subseteq O^\mathtt{c}((1,[i]))$.
\end{lem}
\noindent\textit{Proof}: (a) can be proved similarly to Lemma \ref{lemC6.2.15}. By Lemma \ref{lemC6.2.12}, $d_D((1,[\beta,j]), (1,[i]))=3$, which implies $d_D((p,[j]), (1,[i]))=2$ for some $p=1,2$. Hence, (b) follows.
\qed

\begin{lem}\label{lemC6.2.18}
Let $D$ be an orientation of a $\mathcal{T}$ where $d(D)=4$. Suppose $I((1,[\alpha,i]))=\{(1,[i])\}$ and $I((1,[\beta,j]))=\{(1,[j]), (2,[j])\}$ for $1\le i, j\le \deg_T(\mathtt{c})$, $i\neq j$, $1\le \alpha \le \deg_T([i])-1$, and $1\le \beta \le \deg_T([j])-1$. Then, for each $p=1,2,\ldots, s_{[j]}$,
\\(a) $I^\mathtt{c}((1,[i]))\not\subseteq I^\mathtt{c}((p,[j]))$,
\\(b) $I^\mathtt{c}((1,[j]))\cup I^\mathtt{c}((2,[j]))\not\subseteq I^\mathtt{c}((1,[i]))$.
\end{lem}
\noindent\textit{Proof}: This lemma follows from Lemma \ref{lemC6.2.17} and the Duality Lemma.
\qed

\begin{lem}\label{lemC6.2.19}
Let $D$ be an orientation of a complete bipartite graph $K(p,q)$ with partite sets $V_1=\{1_1,1_2,\ldots, 1_p\}$ and $V_2=\{2_1,2_2,\ldots,2_q\}$, $q\ge p\ge 3$. Suppose further for $1\le i\le p$ that $\lambda_i\rightarrow 2_i \rightarrow \bar{\lambda}_i$, where $\lambda_i=\{1_i, 1_{i+1},\ldots, 1_{i+\lceil\frac{p}{2}\rceil-1}\}$. Then, $d_D(1_i,1_j)=2$ for any $1\le i,j\le p$, $i\neq j$.
\end{lem}
\noindent\textit{Proof}: Let $t\in \mathbb{N}_p$ such that $t\equiv i-\lceil{\frac{p}{2}}\rceil+1 \pmod{p}$. Since $1_i\rightarrow \{2_i,2_t\}$ and $V_1-\{1_i\}\subseteq O(2_i)\cup O(2_t)$, it follows that $d_{D}(1_i,1_j)=2$ for $i\neq j$.
\qed

\section{Proof of Theorem \ref{thmC6.1.7}}
In proving the ``only if" direction of the following propositions, we shall use a common setup forged with the following notions. For a $\mathcal{T}$, let $D$ be an orientation of $\mathcal{T}$ with $d(D)=4$. If $A_2\neq\emptyset$, then by Lemma \ref{lemC6.2.13}(a), we may assume without loss of generality in $D$ that
\begin{align}
(2,[i])\rightarrow (1,[1,i])\rightarrow (1,[i]) \text{ for any }[i]\in A_2. \label{eqC6.3.1}
\end{align}
Also, we let
\begin{align}
B^O_2=\{O^\mathtt{c}((1,[i]))\mid [i]\in A_2\} \text{ and }B^I_2=\{I^\mathtt{c}((2,[i]))\mid [i]\in A_2\}.\label{eqC6.3.2}
\end{align}
Note that $(1,[i])$ ($(2,[i])$ resp.) is effectively the only `outlet' (`inlet' resp.) for the vertex $(1,[1,i])$ if $[i]\in A_2$.
\indent\par Analogously, if $A_3\neq \emptyset$, then by Lemma \ref{lemC6.2.13}(b), we can partition $A_3$ into $A^O_3$ and $A^I_3$, where 
\begin{align}
\left. \begin{array}{@{}ll@{}}
&A^O_3=\{[i]\mid \forall \alpha, 1\le\alpha\le \deg_T([i])-1,\forall p, 1\le p\le s_{[\alpha,i]}, |O((p,[\alpha,i]))|=1\},\\
&A^I_3=\{[i]\mid \exists \alpha, 1\le\alpha\le \deg_T([i])-1, \exists p, 1\le p\le s_{[\alpha,i]}, |O((p,[\alpha,i]))|=2\}.
  \end{array}\right\}
\label{eqC6.3.3}
\end{align}
Without loss of generality, we assume in $D$ that 
\begin{align}
\left. \begin{array}{@{}ll@{}}
&\{(1,[i]), (2,[i])\}\rightarrow (1,[1,i])\rightarrow (3,[i])\text{ if }[i]\in A^O_3,\\
\text{and}&(3,[i])\rightarrow (1,[1,i])\rightarrow \{(1,[i]), (2,[i])\} \text{ if }[i]\in A^I_3.
  \end{array}\right\}
\label{eqC6.3.4}
\end{align}
We also let
\begin{align} 
B^O_3=\{O^\mathtt{c}((3,[i]))\mid [i]\in A^O_3\}, \text{ and } B^I_3=\{I^\mathtt{c}((3,[i]))\mid [i]\in A^I_3\}.\label{eqC6.3.5}
\end{align}
Note that $(3,[i])$ is effectively the only `outlet' (`inlet' resp.) for the vertex $(1,[1,i])$ if $[i]\in A^O_3$  ($A^I_3$ resp.). Furthermore, both $B^O_2\cup B^O_3$ and $B^I_2\cup B^I_3$ are antichains on $(\mathbb{N}_s, \mathtt{c})$ by Lemmas \ref{lemC6.2.15} and \ref{lemC6.2.16} respectively.
\begin{eg}
Let $D$ be the orientation shown in Figures \ref{figC6.3.7} and \ref{figC6.3.8} (see pages 27-28). Then, $A^I_3=\{[1], [2], [3], [4]\}$, $A^O_3=\{[5], [6]\}$, $B^I_3=\{\{(2,\mathtt{c}), (3,\mathtt{c})\}, \{(1,\mathtt{c}), (4,\mathtt{c})\},$ $\{(1,\mathtt{c}), (3,\mathtt{c})\}, \{(2,\mathtt{c}), (4,\mathtt{c})\}\}$, and $B^O_3=\{\{(1,\mathtt{c}), (2,\mathtt{c}), (3,\mathtt{c})\}, \{(2,\mathtt{c}), (3,\mathtt{c}), (4,\mathtt{c})\}\}$.
\end{eg}
\noindent\par As the problem differs for $s=2$ from $s\ge 4$, we consider them separately.
\begin{ppn} \label{ppnC6.3.2}
Suppose $s=2$ and $A_2\cup A_3\neq \emptyset$ for a $\mathcal{T}$. Then, $\mathcal{T}\in \mathscr{C}_0$ if and only if $\deg_T(\mathtt{c})=2$.
\end{ppn}
\noindent\textit{Proof}: $(\Rightarrow)$ Since $\mathcal{T}\in \mathscr{C}_0$, there exists an orientation $D$ of $\mathcal{T}$, where $d(D)=4$. As $A_2\cup A_3\neq\emptyset$, we assume (\ref{eqC6.3.1})-(\ref{eqC6.3.5}) here. From $d(T)=4$, it follows that $|A_{\ge 2}|\ge 2$. We shall consider two cases to show $|A_{\ge 2}|=2$ and $E=\emptyset$.
\\
\\Case 1. $|A_2\cup A^O_3|>0$.
\indent\par Let $[i^*]\in A_2\cup A^O_3$ and $\delta=1$ if $[i^*]\in A_2$, and $\delta=3$ if $[i^*]\in A^O_3$. For all $[i]\in N_T(\mathtt{c})-\{[i^*]\}$ and all $p=1,2,\ldots, s_{[i]}$, since $d_D((1,[1,i^*]), (p,[i]))=3=d_D((1, [1,i]), (\delta,[i^*]))$ by Lemma \ref{lemC6.2.12}, we may assume without loss of generality that $(2,\mathtt{c})\rightarrow (\delta,[i^*])\rightarrow(1,\mathtt{c})$, which implies $(1,\mathtt{c})\rightarrow (p,[i])\rightarrow (2,\mathtt{c})$. Now, if $[i],[j]\in A_{\ge 2}-\{[i^*]\}$, then $d_D((1,[1,i]), (1,[j]))>3$, a contradiction to Lemma \ref{lemC6.2.12}. Hence, $|A_{\ge 2}-\{[i^*]\}|\le 1$ and thus, $|A_{\ge 2}|=2$. If $E\neq\emptyset$, then a similar argument follows for $[i]\in A_{\ge 2}-\{[i^*]\}$ and $[j]\in E$.
\\
\\Case 2. $|A_2\cup A^O_3|=0$ and $|A^I_3|>0$.
\indent\par Then, $A^I_3$ behaves like $A^O_3$ in $\tilde{D}$. The result follows from Case 1 by the Duality Lemma.
\begin{rmk}\label{rmkC6.3.3}
We note the difference in the definition (\ref{eqC6.3.3}) of $A^O_3$ and $A^I_3$ respectively. For the argument, we actually needed only a partition $A^O_3$ and $A^I_3$ of $A_3$ satisfying 
\begin{align*}
&[i]\in A^O_3\implies \exists \alpha, 1\le\alpha\le \deg_T([i])-1, \exists p, 1\le p\le s_{[\alpha,i]}, |O((p,[\alpha,i]))|=1,\\
\text{and }&[i]\in A^I_3\implies \exists \alpha, 1\le\alpha\le \deg_T([i])-1, \exists p, 1\le p\le s_{[\alpha,i]}, |I((p,[\alpha,i]))|=1.
\end{align*}
If $A^O_3$ and $A^I_3$ were each defined using existential quantifiers instead, their intersection may be nonempty. We may arbitrarily include these elements in $A^O_3$ or $A^I_3$ (but not both) to get a partition. However, for the sake of a well-defined partition, we used (\ref{eqC6.3.3}). We emphasize that this does not affect the duality effect in the argument and shall repeatedly apply this.
\end{rmk}
\noindent $(\Leftarrow)$ This follows from Theorem \ref{thmC6.1.6}(a).
\qed

\indent\par For the following, note that Proposition \ref{ppnC6.3.4} holds for all integers $s\ge 2$ while Proposition \ref{ppnC6.3.5} and Corollary \ref{corC6.3.8} hold for all integers $s\ge 3$.
\newpage
\begin{ppn}\label{ppnC6.3.4}
If $s\ge 2$ and $A_2=A_3=\emptyset$ for a $\mathcal{T}$, then $\mathcal{T}\in \mathscr{C}_0$.
\end{ppn}
\noindent\textit{Proof}: Let $\mathcal{H}=T(t_1,t_2,\ldots,t_n)$ be the subgraph of $\mathcal{T}$, where $t_{[i]}=4$ for all $[i]\in \mathcal{T}(A_{\ge 4})$ and $t_v=2$ otherwise. We will use $A_j$ for $\mathcal{H}(A_j)$ for the remainder of this proof. Note that $A_j\neq \emptyset$ if and only if $j=4$. Define an orientation $D$ for $\mathcal{H}$ as follows.
\begin{align}
&\{(2,[i]), (3,[i])\} \rightarrow (1,[\alpha, i])\rightarrow \{(1,[i]), (4,[i])\} \rightarrow (2,[\alpha, i])\rightarrow \{(2,[i]), (3,[i])\},\label{eqC6.3.6}\\
&\text{ and }\{(1,[i]), (2,[i])\} \rightarrow (1,\mathtt{c}) \rightarrow \{(3,[i]), (4,[i])\}\rightarrow (2,\mathtt{c}) \rightarrow \{(1,[i]), (2,[i])\}\label{eqC6.3.7}
\end{align}
for all $[i]\in A_4$ and $1\le \alpha \le \deg_T([i])-1$.
\begin{align}
(2,\mathtt{c})\rightarrow \{(1,[j]), (2,[j])\} \rightarrow (1,\mathtt{c})\label{eqC6.3.8}
\end{align}
for all $[j]\in E$. (See Figure \ref{figC6.2.2} for $D$ when $s=2$.)

\indent\par As there are several cases in verifying the following claim, and similarly for the ``only if" direction of the following propositions, we categorise the cases as follows. Let $[i], [j]\in N_T(\mathtt{c})$, $1\le \alpha \le \deg_T([i])-1$ and $1\le \beta \le \deg_T([j])-1$. 
\\(1) For $v\in \{(p,[\alpha, i]), (q,[i])\mid p=1,2; q=1,2,3,4\}$ and $w\in\{(p,[\beta, j]), (q,[j])\mid p=1,2; q=1,2,3,4\}$, or $v\in\{(p,[\beta, j]), (q,[j])\mid p=1,2; q=1,2,3,4\}$ and $w\in \{(p,[\alpha, i]), (q,[i])\mid p=1,2; q=1,2,3,4\}$:
\\(i) $[i]=[j]\in A_2$; $[i]=[j]\in A_3$; $[i]=[j]\in A_4$;
\\(ii) $i\neq j$ and $[i],[j]\in A_2$; $i\neq j$ and $[i],[j]\in A_3$; $i\neq j$ and  $[i],[j]\in A_4$;
\\(iii) $[i]\in A_2$ and $[j]\in A_3$; $[i]\in A_2$ and $[j]\in A_4$; $[i]\in A_3$ and $[j]\in A_4$.
\\(Here, it is understood that the values of $p$ and $q$ are as whenever applicable.)
\\(2) For $v\in \{(p,[\alpha, i]), (q,[i])\mid p=1,2; q=1,2,3,4\}$ and $w\in \{(r,[k])\mid 1\le r\le 2\}$, or $v\in \{(r,[k])\mid 1\le r\le 2\}$ and $w\in \{(p,[\alpha, i]), (q,[i])\mid p=1,2; q=1,2,3,4\}$:
\\(i) $[i]\in A_2$ and $[k]\in E$; $[i]\in A_3$ and $[k]\in E$; $[i]\in A_4$ and $[k]\in E$;
\\(3) $v,w\in(\mathbb{N}_s, \mathtt{c})$.
\\(4) $v\in \{(p,[\alpha, i]), (q,[i])\mid p=1,2; q=1,2,3,4\}$, where $[i]\in N_T(\mathtt{c})$, and $w\in(\mathbb{N}_s, \mathtt{c})$, or $v\in(\mathbb{N}_s, \mathtt{c})$ and $w\in \{(p,[\alpha, i]), (q,[i])\mid p=1,2; q=1,2,3,4\}$, where $[i]\in N_T(\mathtt{c})$.
\\(5) $v\in \{(p,[i])\mid p=1,2\}$ and $w\in \{(q,[j])\mid q=1,2\}$, where $[i], [j]\in N_T(\mathtt{c})$.
\\((5) is required as we primarily focused on vertices $v, w$ with distance at least $3$ in $\mathcal{T}$ in (1) and (2).)
\\
\\Claim: For all $v,w\in V(D)$, $d_{D}(v,w)\le 4$.
\\
\\Case 1.1. $v,w \in \{(1,[\alpha, i]),(2,[\alpha, i])\mid 1\le \alpha \le \deg_T([i])-1\}$ for each $[i]\in A_4$.
\indent\par This is clear since, by (\ref{eqC6.3.6}), $(2,[i])\rightarrow (1,[\alpha_1, i])\rightarrow (1,[i])\rightarrow (2,[\alpha_2, i])\rightarrow (2,[i])$ is a directed $C_4$ for all $1\le \alpha_1\le \alpha_2\le \deg_T([i])-1$.
\\
\\Case 1.2. $v\in\{(1,[\alpha, i]), (2,[\alpha, i])\}$ and $w\in\{(1,[i]), (2,[i]), (3,[i]), (4,[i])\}$ for all $[i]\in A_4$ and $1\le \alpha \le \deg_T([i])-1$.
\indent\par By symmetry, it suffices to show for the case $v=(1,[\alpha, i])$. By (\ref{eqC6.3.6}), $(1,[\alpha, i])\rightarrow \{(1,[i]), (4,[i])\} \rightarrow (2,[\alpha, i])\rightarrow \{(2,[i]), (3,[i])\}$.
\\
\\Case 1.3. $v\in\{(1,[\alpha, i]), (2,[\alpha, i])\}$ and $w\in\{(1,\mathtt{c}), (2,\mathtt{c})\}$ for all $[i]\in A_4$ and $1\le \alpha \le \deg_T([i])-1$.
\indent\par By symmetry, it suffices to show for the case $v=(1,[\alpha, i])$. By (\ref{eqC6.3.6})-(\ref{eqC6.3.7}), $(1,[\alpha, i])\rightarrow (1,[i])\rightarrow (1,\mathtt{c})$ and $(1,[\alpha, i])\rightarrow (4,[i])\rightarrow (2,\mathtt{c})$.
\\
\\Case 1.4. $v\in\{(1,[\alpha, i]), (2,[\alpha, i])\}$ and $w\in\{(1,[j]), (2,[j]), (3,[j]), (4,[j])\}$ for all $[i], [j]\in A_4$, $i\neq j$, and all $1\le \alpha \le \deg_T([i])-1$.
\indent\par By symmetry, it suffices to show for the case $v=(1,[\alpha, i])$. By (\ref{eqC6.3.6})-(\ref{eqC6.3.7}), $(1,[\alpha, i])\rightarrow (1,[i])\rightarrow (1,\mathtt{c})\rightarrow \{(3,[j]), (4,[j])\}$ and $(1,[\alpha, i])\rightarrow (4,[i])\rightarrow (2,\mathtt{c})\rightarrow \{(1,[j]), (2,[j])\}$.
\\
\\Case 1.5. $v\in\{(1,[\alpha, i]), (2,[\alpha, i])\}$ and $w\in\{(1,[\beta,j]), (2,[\beta,j])\}$ for all $[i], [j]\in A_4$, $i\neq j$, all $1\le \alpha \le \deg_T([i])-1$ and $1\le \beta \le \deg_T([j])-1$.
\indent\par By symmetry, it suffices to show for the case $v=(1,[\alpha, i])$. By (\ref{eqC6.3.6})-(\ref{eqC6.3.7}), $(1,[\alpha, i])\rightarrow (1,[i])\rightarrow (1,\mathtt{c})\rightarrow \{(3,[j]), (4,[j])\}$ and while $(3,[j])\rightarrow (1,[\beta,j])$ and $(4,[j])\rightarrow (2,[\beta,j])$.
\\
\\Case 2.1. $v\in\{(1,[i]), (2,[i]), (3,[i]), (4,[i])\}$ and $w\in\{(1,[\alpha, i]), (2,[\alpha, i])\}$ for all $[i]\in A_4$ and $1\le \alpha \le \deg_T([i])-1$.
\indent\par By symmetry, it suffices to show for the case $v=(1,[i])$. By (\ref{eqC6.3.6}), $(1,[i])\rightarrow (2,[\alpha, i])\rightarrow (2,[i])\rightarrow (1,[\alpha, i])$ .
\\
\\Case 2.2. $v, w\in\{(1,[i]), (2,[i]), (3,[i]), (4,[i])\}$ for all $[i]\in A_4$ and $1\le \alpha \le \deg_T([i])-1$, or $v\in\{(1,[i]), (2,[i]), (3,[i]), (4,[i])\}$ and $w\in\{(1,\mathtt{c}), (2,\mathtt{c})\}$ for all $[i]\in A_4$.
\indent\par This follows directly from (\ref{eqC6.3.7}).
\\
\\Case 2.3. $v\in\{(1,[i]), (2,[i]), (3,[i]), (4,[i])\}$ and $w\in\{(1,[j]), (2,[j]), (3,[j]), (4,[j])\}$ for $[i], [j]\in A_4$ and $i\neq j$.
\indent\par By symmetry, it suffices to show for the case $v\in\{(1,[i]), (2,[i])\}$. By (\ref{eqC6.3.7}), $\{(1,[i]),$ $(2,[i])\} \rightarrow (1,\mathtt{c}) \rightarrow \{(3,[j]), (4,[j])\}\rightarrow (2,\mathtt{c})\rightarrow \{(1,[j]), (2,[j])\}$.
\\
\\Case 2.4. $v\in\{(1,[i]), (2,[i]), (3,[i]), (4,[i])\}$ and $w\in\{(1,[\beta,j]), (2,[\beta,j])\}$ for $[i], [j]\in A_4$, $i\neq j$, and $1\le \beta \le \deg_T([j])-1$.
\indent\par By symmetry, it suffices to show for the case $v\in\{(1,[i]), (2,[i])\}$. By (\ref{eqC6.3.6})-(\ref{eqC6.3.7}), $\{(1,[i]),$ $(2,[i])\} \rightarrow (1,\mathtt{c}) \rightarrow \{(3,[j]), (4,[j])\}$, while $(3,[j])\rightarrow (1,[\beta,j])$ and $(4,[j])\rightarrow (2,[\beta,j])$.
\\
\\Case 3.1. $v, w\in\{(1,\mathtt{c}), (2,\mathtt{c})\}$, or $v\in\{(1,\mathtt{c}), (2,\mathtt{c})\}, w\in\{(1,[i]), (2,[i]), (3,[i]), (4,[i])\}$ for all $[i]\in A_4$.
\indent\par This follows directly from (\ref{eqC6.3.7}).
\\
\\Case 3.2. $v\in\{(1,\mathtt{c}), (2,\mathtt{c})\}, w\in\{(1,[\alpha, i]), (2,[\alpha, i])\}$ for all $[i]\in A_4$ and $1\le \alpha \le \deg_T([i])-1$.
\indent\par By symmetry, it suffices to show for the case $v=(1,\mathtt{c})$. By (\ref{eqC6.3.6})-(\ref{eqC6.3.7}), $(1,\mathtt{c}) \rightarrow \{(3,[i]),$ $(4,[i])\}$, while $(3,[i])\rightarrow (1,[\alpha, i])$ and $(4,[i])\rightarrow (2,[\alpha, i])$.
\\
\\Case 4.1. $v\in\{(1,[\alpha, i]), (2,[\alpha, i])\}$ and $w\in\{(1,[j]), (2,[j])\}$ for all $[i]\in A_4$, $1\le \alpha \le \deg_T([i])-1$ and $[j]\in E$.
\indent\par By (\ref{eqC6.3.6})-(\ref{eqC6.3.8}), $(1,[\alpha, i])\rightarrow (4,[i])\rightarrow (2,\mathtt{c})$, $(2,[\alpha, i])\rightarrow (3,[i])\rightarrow (2,\mathtt{c})$, and $(2,\mathtt{c})\rightarrow \{(1,[j]),(2,[j])\}$.
\\
\\Case 4.2. $v\in\{(1,[j]), (2,[j])\}$ and $w\in \{(1,[\alpha, i]), (2,[\alpha, i])\}$ for all $[i]\in A_4$, $1\le \alpha \le \deg_T([i])-1$ and $[j]\in E$.
\indent\par By (\ref{eqC6.3.6})-(\ref{eqC6.3.8}), $\{(1,[j]),(2,[j])\}\rightarrow (1,\mathtt{c})\rightarrow \{(3,[i]),(4,[i])\}$, $(3,[i])\rightarrow (1,[\alpha, i])$, and $(4,[i])\rightarrow (2,[\alpha, i])$.
\\
\\Case 4.3. $v\in\{(1,[i]), (2,[i]), (3,[i]), (4,[i])\}$ and $w\in\{(1,[j]), (2,[j])\}$ for all $[i]\in A_4$, $1\le \alpha \le \deg_T([i])-1$ and $[j]\in E$.
\indent\par By (\ref{eqC6.3.7})-(\ref{eqC6.3.8}), $\{(1,[i]),(2,[i])\}\rightarrow (1,\mathtt{c})\rightarrow\{(3,[i]),(4,[i])\}\rightarrow (2,\mathtt{c})\rightarrow \{(1,[j]),(2,[j])\}$.
\\
\\Case 4.4. $v\in\{(1,[j]), (2,[j])\}$ and $w\in\{(1,[i]), (2,[i]), (3,[i]), (4,[i])\}$ for all $[i]\in A_4$, $1\le \alpha \le \deg_T([i])-1$ and $[j]\in E$.
\indent\par By (\ref{eqC6.3.7})-(\ref{eqC6.3.8}), $\{(1,[j]), (2,[j])\}\rightarrow (1,\mathtt{c})\rightarrow \{(3,[i]),(4,[i])\}\rightarrow (2,\mathtt{c})\rightarrow \{(1,[i]),(2,[i])\}$.
\\
\\Case 4.5. 
\\(i) $v\in\{(1,\mathtt{c}), (2,\mathtt{c})\}$ and $w\in\{(1,[j]), (2,[j])\}$ for all $[j]\in E$, or
\\(ii) $v\in\{(1,[j]), (2,[j])\}$ and $w\in\{(1,\mathtt{c}), (2,\mathtt{c})\}$ for all $[j]\in E$.
\indent\par By (\ref{eqC6.3.7})-(\ref{eqC6.3.8}), $\{(1,[j]),(2,[j])\}\rightarrow (1,\mathtt{c})\rightarrow (3,[k])\rightarrow (2,\mathtt{c}) \rightarrow \{(1,[j]),(2,[j])\}$ is a directed $C_4$ for any $[k]\in A_4$.
\\
\\Case 4.6. $v\in\{(1,[i]), (2,[i])\}$ and $w\in\{(1,[j]), (2,[j])\}$ for all $[i],[j]\in E$.
\indent\par By (\ref{eqC6.3.7})-(\ref{eqC6.3.8}), $\{(1,[i]),(2,[i])\}\rightarrow (1,\mathtt{c})\rightarrow (3,[k])\rightarrow (2,\mathtt{c}) \rightarrow \{(1,[j]),(2,[j])\}$ is a directed $C_4$ for any $[k]\in A_4$.

\indent\par Therefore, the claim follows. Since every vertex lies in a directed $C_4$ for $D$ and $d(D)=4$, $\bar{d}(\mathcal{T})\le \max \{4, d(D)\}$ by Lemma \ref{lemC6.1.3}. With $\bar{d}(\mathcal{T})\ge d(T)=4$, it follows that $\bar{d}(\mathcal{T})=4$.
\qed

\begin{center}
\begin{tikzpicture}[thick,scale=0.7]%
\draw(-6,2)node[circle, draw, fill=black!100, inner sep=0pt, minimum width=6pt, label={[] 180:{\small $(1,[1,1])$}}](1_11){};
\draw(-6,0)node[circle, draw, fill=black!100, inner sep=0pt, minimum width=6pt, label={[] 180:{\small $(2,[1,1])$}}](2_11){};

\draw(-6,-2)node[circle, draw, fill=black!100, inner sep=0pt, minimum width=6pt, label={[] 180:{\small $(1,[2,1])$}}](1_21){};
\draw(-6,-4)node[circle, draw, fill=black!100, inner sep=0pt, minimum width=6pt, label={[] 180:{\small $(2,[2,1])$}}](2_21){};

\draw(-3,2)node[circle, draw, fill=black!100, inner sep=0pt, minimum width=6pt, label={[yshift=0cm, xshift=0.3cm] 90:{\small $(1,[1])$}}](1_1){};
\draw(-3,0)node[circle, draw, fill=black!100, inner sep=0pt, minimum width=6pt, label={[yshift=0cm, xshift=0.3cm] 270:{\small $(2,[1])$}}](2_1){};
\draw(-3,-2)node[circle, draw, fill=black!100, inner sep=0pt, minimum width=6pt, label={[yshift=0cm, xshift=0.3cm] 270:{\small $(3,[1])$}}](3_1){};
\draw(-3,-4)node[circle, draw, fill=black!100, inner sep=0pt, minimum width=6pt, label={[yshift=0cm, xshift=-0.3cm] 270:{\small $(4,[1])$}}](4_1){};

\draw(0,0)node[circle, draw, inner sep=0pt, minimum width=3pt](1_u){\scriptsize $(1,\mathtt{c})$};
\draw(0,-2)node[circle, draw, inner sep=0pt, minimum width=3pt](2_u){\scriptsize $(2,\mathtt{c})$};

\draw(3,2)node[circle, draw, fill=black!100, inner sep=0pt, minimum width=6pt, label={[yshift=0cm, xshift=-0.3cm] 90:{\small $(1,[2])$}}](1_2){};
\draw(3,0)node[circle, draw, fill=black!100, inner sep=0pt, minimum width=6pt, label={[yshift=0cm, xshift=-0.3cm] 270:{\small $(2,[2])$}}](2_2){};
\draw(3,-2)node[circle, draw, fill=black!100, inner sep=0pt, minimum width=6pt, label={[yshift=0cm, xshift=-0.3cm] 270:{\small $(3,[2])$}}](3_2){};
\draw(3,-4)node[circle, draw, fill=black!100, inner sep=0pt, minimum width=6pt, label={[yshift=0cm, xshift=0.3cm] 270:{\small $(4,[2])$}}](4_2){};

\draw(6,0)node[circle, draw, fill=black!100, inner sep=0pt, minimum width=6pt, label={[] 0:{\small $(1,[1,2])$}}](1_12){};
\draw(6,-2)node[circle, draw, fill=black!100, inner sep=0pt, minimum width=6pt, label={[] 0:{\small $(2,[1,2])$}}](2_12){};

\draw(-3,-6)node[circle, draw, fill=black!100, inner sep=0pt, minimum width=6pt, label={[] 180:{\small $(1,[3])$}}](1_3){};
\draw(-3,-8)node[circle, draw, fill=black!100, inner sep=0pt, minimum width=6pt, label={[] 180:{\small $(2,[3])$}}](2_3){};

\draw(3,-6)node[circle, draw, fill=black!100, inner sep=0pt, minimum width=6pt, label={[] 0:{\small $(1,[4])$}}](1_4){};
\draw(3,-8)node[circle, draw, fill=black!100, inner sep=0pt, minimum width=6pt, label={[] 0:{\small $(2,[4])$}}](2_4){};

\draw[->, line width=0.3mm, >=latex, shorten <= 0.2cm, shorten >= 0.1cm](2_1)--(1_u);
\draw[->, line width=0.3mm, >=latex, shorten <= 0.2cm, shorten >= 0.1cm](1_1)--(1_u);

\draw[->, line width=0.3mm, >=latex, shorten <= 0.2cm, shorten >= 0.1cm](2_2)--(1_u);
\draw[->, line width=0.3mm, >=latex, shorten <= 0.2cm, shorten >= 0.1cm](1_2)--(1_u);

\draw[->, line width=0.3mm, >=latex, shorten <= 0.2cm, shorten >= 0.1cm](3_1)--(2_u);
\draw[->, line width=0.3mm, >=latex, shorten <= 0.2cm, shorten >= 0.1cm](4_1)--(2_u);

\draw[->, line width=0.3mm, >=latex, shorten <= 0.2cm, shorten >= 0.1cm](3_2)--(2_u);
\draw[->, line width=0.3mm, >=latex, shorten <= 0.2cm, shorten >= 0.1cm](4_2)--(2_u);

\draw[->, line width=0.3mm, >=latex, shorten <= 0.2cm, shorten >= 0.15cm](1_11)--(1_1);
\draw[->, line width=0.3mm, >=latex, shorten <= 0.2cm, shorten >= 0.15cm](1_11)--(4_1);
\draw[->, line width=0.3mm, >=latex, shorten <= 0.2cm, shorten >= 0.15cm](1_12)--(1_2);
\draw[->, line width=0.3mm, >=latex, shorten <= 0.2cm, shorten >= 0.15cm](1_12)--(4_2);

\draw[->, line width=0.3mm, >=latex, shorten <= 0.2cm, shorten >= 0.15cm](2_11)--(2_1);
\draw[->, line width=0.3mm, >=latex, shorten <= 0.2cm, shorten >= 0.15cm](2_11)--(3_1);
\draw[->, line width=0.3mm, >=latex, shorten <= 0.2cm, shorten >= 0.15cm](2_12)--(2_2);
\draw[->, line width=0.3mm, >=latex, shorten <= 0.2cm, shorten >= 0.15cm](2_12)--(3_2);

\draw[->, line width=0.3mm, >=latex, shorten <= 0.2cm, shorten >= 0.15cm](1_21)--(1_1);
\draw[->, line width=0.3mm, >=latex, shorten <= 0.2cm, shorten >= 0.15cm](1_21)--(4_1);
\draw[->, line width=0.3mm, >=latex, shorten <= 0.2cm, shorten >= 0.15cm](2_21)--(2_1);
\draw[->, line width=0.3mm, >=latex, shorten <= 0.2cm, shorten >= 0.15cm](2_21)--(3_1);

\draw[->, line width=0.3mm, >=latex, shorten <= 0.2cm, shorten >= 0.1cm](1_4)--(1_u);
\draw[->, line width=0.3mm, >=latex, shorten <= 0.2cm, shorten >= 0.1cm](2_4)--(1_u);
\draw[->, line width=0.3mm, >=latex, shorten <= 0.2cm, shorten >= 0.1cm](1_3)--(1_u);
\draw[->, line width=0.3mm, >=latex, shorten <= 0.2cm, shorten >= 0.1cm](2_3)--(1_u);

\draw[dashed,->, line width=0.3mm, >=latex, shorten <= 0.2cm, shorten >= 0.15cm](2_1)--(1_11);
\draw[dashed,->, line width=0.3mm, >=latex, shorten <= 0.2cm, shorten >= 0.15cm](1_1)--(2_11);

\draw[dashed,->, line width=0.3mm, >=latex, shorten <= 0.2cm, shorten >= 0.15cm](3_1)--(1_21);
\draw[dashed,->, line width=0.3mm, >=latex, shorten <= 0.2cm, shorten >= 0.15cm](4_1)--(2_21);

\draw[dashed,->, line width=0.3mm, >=latex, shorten <= 0.2cm, shorten >= 0.15cm](1_2)--(2_12);
\draw[dashed,->, line width=0.3mm, >=latex, shorten <= 0.2cm, shorten >= 0.15cm](2_2)--(1_12);

\draw[dashed,->, line width=0.3mm, >=latex, shorten <= 0.2cm, shorten >= 0.15cm](3_2)--(1_12);
\draw[dashed,->, line width=0.3mm, >=latex, shorten <= 0.2cm, shorten >= 0.15cm](4_2)--(2_12);

\draw[dashed,->, line width=0.3mm, >=latex, shorten <= 0.2cm, shorten >= 0.15cm](1_1)--(2_21);
\draw[dashed,->, line width=0.3mm, >=latex, shorten <= 0.2cm, shorten >= 0.15cm](2_1)--(1_21);
\draw[dashed,->, line width=0.3mm, >=latex, shorten <= 0.2cm, shorten >= 0.15cm](3_1)--(1_11);
\draw[dashed,->, line width=0.3mm, >=latex, shorten <= 0.2cm, shorten >= 0.15cm](4_1)--(2_11);
\end{tikzpicture}
\captionsetup{justification=centering}
{\captionof{figure}{Orientation $D$ for $\mathcal{H}$, where $A_4=\{[1],[2]\}$, $E=\{[3],[4]\}$.\\Note that the parent graph is the tree in Figure \ref{figC6.1.1}.}\label{figC6.2.2}
\captionsetup{justification=justified}
\caption*{Note: For clarity, the arcs directed from $(p,\mathtt{c})$ to $(q,[i])$ are omitted, while the arcs directed from $(q,[i])$ to $(r,[\alpha, i])$ are represented by \lineexample{dashed} lines. The same simplification is used for Figures \ref{figC6.3.3} to \ref{figC6.4.18}.}}
\end{center}

\begin{ppn}\label{ppnC6.3.5}
Suppose $s\ge 3$ and $A_{\ge 3}=\emptyset$ for a $\mathcal{T}$. Then,
\begin{align*}
\mathcal{T}\in \mathscr{C}_0\iff \left\{
  \begin{array}{@{}ll@{}}
    |A_2|\le{{s}\choose{\lceil{s/2}\rceil}}-1, & \text{if}\ |A_2|< \deg_T(\mathtt{c}), \\
    |A_2|\le{{s}\choose{\lceil{s/2}\rceil}}, & \text{if}\ |A_2|=\deg_T(\mathtt{c}). 
  \end{array}\right.
\end{align*}
\end{ppn}
\noindent\textit{Proof}: $(\Rightarrow)$ Since $\mathcal{T}\in \mathscr{C}_0$, there exists an orientation $D$ of $\mathcal{T}$, where $d(D)=4$. As $A_2\neq\emptyset$, we assume (\ref{eqC6.3.1})-(\ref{eqC6.3.2}) here. By Sperner's theorem, $|A_2|=|B^O_2|\le {{s}\choose{\lceil{s/2}\rceil}}$. So, we are done if $|A_2|=\deg_T(\mathtt{c})$. 
\indent\par Now, assume $|A_2|< \deg_T(\mathtt{c})$ and let $[i^*]\in E$. If $|O^\mathtt{c}((1,[i^*]))|\ge \lceil\frac{s}{2}\rceil$, then $d_D((1,[1,i]),$ $(1,[i^*])) =3$ implies $O^\mathtt{c}((1,[i]))\cap I^\mathtt{c}((1,[i^*]))\neq \emptyset$ for all $[i]\in A_2$. It follows from Lih's theorem that $|A_2|=|B^O_2|\le {{s}\choose{\lceil{s/2}\rceil}}-{{s-|I^\mathtt{c}((1,[i^*]))|}\choose{\lceil{s/2}\rceil}}\le {{s}\choose{\lceil{s/2}\rceil}}-{{\lceil{s/2}\rceil}\choose{\lceil{s/2}\rceil}}={{s}\choose{\lceil{s/2}\rceil}}-1$. If $|O^\mathtt{c}((1,[i^*]))|\le \lfloor\frac{s}{2}\rfloor$, then $d_D((1,[i^*]),(1,[1,i]))=3$ implies $I^\mathtt{c}((2,[i]))\cap O^\mathtt{c}((1,[i^*]))\neq \emptyset$ for all $[i]\in A_2$. It follows from Lih's theorem that $|A_2|=|B^I_2|\le {{s}\choose{\lceil{s/2}\rceil}}-{{s-|O^\mathtt{c}((1,[i^*]))|}\choose{\lceil{s/2}\rceil}}\le {{s}\choose{\lceil{s/2}\rceil}}-{{\lceil{s/2}\rceil}\choose{\lceil{s/2}\rceil}}={{s}\choose{\lceil{s/2}\rceil}}-1$.

\begin{rmk}
On account of the above part, it is intuitive to let $O^\mathtt{c}((1,[i]))=O^\mathtt{c}((2,[i]))$ and $|O^\mathtt{c}((1,[i]))|=\lfloor\frac{s}{2}\rfloor$ in constructing an optimal orientation $D$ of $\mathcal{T}$. Indeed, this is our plan if $|A_2|$ is big enough (i.e., $|A_2|\ge s$). However, there are some potential drawbacks of this approach if $|A_2|$ is small (i.e., $|A_2|<s$).  For instance, consider $s=5$ and $\deg_T(\mathtt{c})=|A_2|=2$. If we assigned $O^\mathtt{c}((p,[1]))=\{(1,\mathtt{c}),(2,\mathtt{c})\}$ and $O^\mathtt{c}((p,[2]))=\{(1,\mathtt{c}),(3,\mathtt{c})\}$ for $p=1,2$, then $\deg^+((1,\mathtt{c}))=0$ and $\deg^-((j,\mathtt{c}))=0$ for $j=4,5$. Consequently, $D$ will not be a strong orientation. Hence, we consider cases dependent on $|A_2|$ to circumvent this problem; namely, they are Cases 1 and 2 for small $|A_2|$, and Cases 3 and 4 for large $|A_2|$.
\end{rmk}

\noindent $(\Leftarrow)$ Without loss of generality, assume $A_2=\{[i]\mid i\in\mathbb{N}_{|A_2|}\}$. Thus, it is taken that $E=\{[i]\mid i\in\mathbb{N}_{\deg_T(\mathtt{c})}-\mathbb{N}_{|A_2|}\}$ if $|A_2|<\deg_T(\mathtt{c})$. Let $\mathcal{H}=T(t_1,t_2,\ldots,t_n)$ be the subgraph of $\mathcal{T}$, where $t_\mathtt{c}=s$ and $t_v=2$ for all $v\neq \mathtt{c}$. We will use $A_j$ for $\mathcal{H}(A_j)$ for the remainder of this proof.
\\
\\Case 1. $|A_2|=\deg_T(\mathtt{c})$ (i.e., $E=\emptyset$) and $|A_2|\le s$.
\indent\par Define an orientation $D_1$ for $\mathcal{H}$ as follows.
\begin{align}
(2,[i])\rightarrow (1,[\alpha, i])\rightarrow (1,[i])\rightarrow (2,[\alpha, i])\rightarrow (2,[i])\label{eqC6.3.9}
\end{align}
for all $[i]\in A_2$ and $1\le \alpha \le \deg_T([i])-1$.
\begin{align}
(\mathbb{N}_s,\mathtt{c})-\{(i,\mathtt{c})\}\rightarrow \{(1,[i]),(2,[i])\}\rightarrow (i,\mathtt{c})\label{eqC6.3.10}
\end{align}
for all $1\le i\le |A_2|-1$.
\begin{align}
(\mathbb{N}_s,\mathtt{c})-\{(k,\mathtt{c})\mid |A_2|\le k\le s\}\rightarrow\{(1,[|A_2|]),(2,[|A_2|])\}\nonumber
\\\rightarrow \{(k,\mathtt{c})\mid |A_2|\le k\le s\}. \label{eqC6.3.11}
\end{align}
(See Figure \ref{figC6.3.3} for $D_1$ when $s=5$.)
\\
\\Claim 1: For all $v,w\in V(D_1)$, $d_{D_1}(v,w)\le 4$.
\\
\\Case 1.1. $v,w \in \{(1,[\alpha, i]),(2,[\alpha, i]),(1,[i]),(2,[i])\}$ for each $[i]\in A_2$ and $1\le \alpha \le \deg_T([i])-1$.
\indent\par This is clear since (\ref{eqC6.3.9}) is a directed $C_4$.
\\
\\Case 1.2. For each $[i],[j]\in A_2$, $i \neq j$, each $1\le \alpha \le \deg_T([i])-1$, and each $1\le \beta \le \deg_T([j])-1$,
\\(i) $v=(p,[\alpha, i]), w=(q,[\beta,j])$ for $p,q=1,2$, or
\\(ii) $v=(p,[\alpha, i]), w=(q,[i])$ for $p,q=1,2$, or
\\(iii) $v=(p,[i]), w=(q,[\beta,j])$ for $p,q=1,2$.
\indent\par If $i\neq j$, then, by (\ref{eqC6.3.9})-(\ref{eqC6.3.11}), $(p,[\alpha, i])\rightarrow (p,[i])\rightarrow (i,\mathtt{c}) \rightarrow (1,[j]) \rightarrow (2,[\beta,j])$ and $(p,[\alpha, i])\rightarrow (p,[i])\rightarrow (i,\mathtt{c}) \rightarrow (2,[j]) \rightarrow (1,[\beta,j])$.
\\
\\Case 1.3. $v=(x_1,\mathtt{c})$ and $w=(x_2,\mathtt{c})$ for $x_1\neq x_2$ and $1\le x_1, x_2 \le s$.
\indent\par If $x_2<|A_2|$, then $(x_1,\mathtt{c}) \rightarrow (1,[x_2]) \rightarrow (x_2,\mathtt{c})$ by (\ref{eqC6.3.10}). If $x_1<|A_2|\le x_2\le s$, then $(x_1,\mathtt{c}) \rightarrow (1,[|A_2|]) \rightarrow (x_2,\mathtt{c})$ by (\ref{eqC6.3.10})-(\ref{eqC6.3.11}). If $|A_2|\le x_1, x_2\le s$, then $(x_1,\mathtt{c}) \rightarrow (1,[1])\rightarrow (1,\mathtt{c}) \rightarrow (1,[|A_2|]) \rightarrow (x_2,\mathtt{c})$ by (\ref{eqC6.3.10})-(\ref{eqC6.3.11}).
\\
\\Case 1.4. $v\in \{(1,[i]), (2,[i]), (1,[\alpha, i]), (2,[\alpha, i])\}$ for each $[i]\in A_2$, $1\le \alpha \le \deg_T([i])-1$, and $w=(j,\mathtt{c})$ for $1\le j\le s$.
\indent\par If $j=i$, or $i=|A_2|\le j\le s$, then $(p,[\alpha, i])\rightarrow (p,[i])\rightarrow (j,\mathtt{c})$ for $p=1,2$, by (\ref{eqC6.3.9})-(\ref{eqC6.3.11}). If $j\neq i$ and $j<|A_2|$, then $(p,[\alpha, i])\rightarrow (p,[i])\rightarrow (i,\mathtt{c}) \rightarrow (1,[j]) \rightarrow (j,\mathtt{c})$ for $p=1,2$, by (\ref{eqC6.3.9})-(\ref{eqC6.3.10}). If $i<|A_2|\le j\le s$, then $(p,[\alpha, i])\rightarrow (p,[i])\rightarrow (i,\mathtt{c}) \rightarrow (1,[|A_2|]) \rightarrow (j,\mathtt{c})$ for $p=1,2$, by (\ref{eqC6.3.9})-(\ref{eqC6.3.11}).
\\
\\Case 1.5. $v=(j,\mathtt{c})$ for each $1\le j\le s$, and $w\in \{(1,[i]), (2,[i]), (1,[\alpha, i]), (2,[\alpha, i])\}$ for each $[i]\in A_2$ and $1\le \alpha \le \deg_T([i])-1$.
\indent\par If $j<|A_2|$ and $j\neq i$, or $i<|A_2|\le j\le s$, then $(j,\mathtt{c})\rightarrow (p,[i]) \rightarrow (3-p,[\alpha, i])$ for $p=1,2$, by (\ref{eqC6.3.9})-(\ref{eqC6.3.11}). If $i=j<|A_2|$, then $(j,\mathtt{c})\rightarrow (1,[|A_2|])\rightarrow (|A_2|,\mathtt{c})\rightarrow (p,[i]) \rightarrow (3-p,[\alpha, i])$ for $p=1,2$, by (\ref{eqC6.3.9})-(\ref{eqC6.3.11}). If $i=|A_2|\le j\le s$, then $(j,\mathtt{c}) \rightarrow (1,[1])\rightarrow (1,\mathtt{c}) \rightarrow (p,[|A_2|]) \rightarrow (3-p,[\alpha,|A_2|])$, for $p=1,2$, by (\ref{eqC6.3.9})-(\ref{eqC6.3.11}).
\\
\\Case 1.6. $v=(p,[i])$ and $w=(q, [j])$, where $1\le p,q\le 2$, $i\neq j$, and $[i],[j]\in A_2$.
\noindent\par This follows from the fact that $|O^\mathtt{c}((p,[i]))|>0$, $|I^\mathtt{c}((q,[j]))|>0$, and $d_{D_1}((r_1,\mathtt{c}), (r_2,\mathtt{c}))$ $=2$ for any $r_1\neq r_2$ and $1\le r_1, r_2\le |A_2|$ by Case 1.3.
\\
\\Case 2. $|A_2|<\deg_T(\mathtt{c})$ (i.e., $E\neq\emptyset$) and $|A_2|<s$.
\indent\par Define an orientation $D_2$ for $\mathcal{H}$ as follows.
\begin{align}
(2,[i])\rightarrow (1,[\alpha, i])\rightarrow (1,[i])\rightarrow (2,[\alpha, i])\rightarrow (2,[i])\label{eqC6.3.12}
\end{align}
for all $[i]\in A_2$ and $1\le \alpha \le \deg_T([i])-1$.
\begin{align}
(\mathbb{N}_s,\mathtt{c})-\{(i,\mathtt{c})\}\rightarrow \{(1,[i]),(2,[i])\}\rightarrow (i,\mathtt{c})\label{eqC6.3.13}
\end{align}
for all $1\le i\le |A_2|$.
\begin{align}
&(\mathbb{N}_{|A_2|},\mathtt{c})\rightarrow \{(p,[i])\mid p=1,2; [i]\in E\}\nonumber\\
&\rightarrow \{(k,\mathtt{c})\mid |A_2|< k\le s\}\rightarrow \{(q,[j])\mid q=1,2;[j]\in A_2\}.\label{eqC6.3.14}
\end{align}
(See Figure \ref{figC6.3.4} for $D_2$ when $s=5$.)
\\
\\Claim 2: For all $v,w\in V(D_2)$, $d_{D_2}(v,w)\le 4$.
\noindent\par In view of the similarity between $D_1$ and $D_2$, it suffices to check the following.
\\
\\Case 2.1. For each $[i]\in A_2$, each $[j]\in E$, and each $1\le \alpha \le \deg_T([i])-1$,
\\(i) $v=(p,[\alpha, i]), w=(q,[j])$ for $p,q=1,2$, or
\\(ii) $v=(p,[i]), w=(q,[j])$ for $p,q=1,2$, or
\\(iii) $v=(q,[j]), w=(p,[\alpha, i])$ for $p,q=1,2$.
\indent\par By (\ref{eqC6.3.12})-(\ref{eqC6.3.14}), (i) and (ii) follow from $(p,[\alpha, i])\rightarrow (p,[i])\rightarrow (i,\mathtt{c}) \rightarrow \{(1,[j]),$ $(2,[j])\}$. Similarly for (iii), $\{(1,[j]), (2,[j])\} \rightarrow (s,\mathtt{c})\rightarrow (3-p,[i])\rightarrow (p,[\alpha, i])$.
\\
\\Case 2.2. $v=(x_1,\mathtt{c})$ and $w=(x_2,\mathtt{c})$ for $x_1\neq x_2$ and $1\le x_1, x_2 \le s$.
\indent\par If $x_2\le |A_2|$, then $(x_1,\mathtt{c}) \rightarrow (1,[x_2]) \rightarrow (x_2,\mathtt{c})$ by (\ref{eqC6.3.13}). If $x_1\le |A_2|$ and $|A_2|+1\le x_2\le s$, then $(x_1,\mathtt{c}) \rightarrow (1,[j]) \rightarrow (x_2,\mathtt{c})$ for any $[j]\in E$ by (\ref{eqC6.3.14}). If $|A_2|+1\le x_1, x_2 \le s$, then $(x_1,\mathtt{c}) \rightarrow (1,[1])\rightarrow (1,\mathtt{c}) \rightarrow (1,[j]) \rightarrow (x_2,\mathtt{c})$ for any $[j]\in E$ by (\ref{eqC6.3.13}) and (\ref{eqC6.3.14}).
\\
\\Case 2.3. $v\in \{(1,[i]), (2,[i])\}$ for each $[i]\in E$, and $w=(j,\mathtt{c})$ for $1\le j\le s$.
\indent\par For $1\le j\le |A_2|$, $\{(1,[i]),(2,[i])\}\rightarrow (s,\mathtt{c})\rightarrow (1,[j])\rightarrow (j,\mathtt{c})$ by (\ref{eqC6.3.13}) and (\ref{eqC6.3.14}). For $|A_2|+1\le j\le s$, $\{(1,[i]),(2,[i])\}\rightarrow (j,\mathtt{c})$ by (\ref{eqC6.3.14}). 
\\
\\Case 2.4. $v=(j,\mathtt{c})$ for each $1\le j\le s$, and $w\in \{(1,[i]), (2,[i])\}$ for each $[i]\in E$.
\noindent\par By (\ref{eqC6.3.14}), for any $1\le j\le |A_2|$, $(j,\mathtt{c})\rightarrow \{(1,[i]), (2,[i])\}$. For $|A_2|+1\le j\le s$, $(j,\mathtt{c})\rightarrow (1,[1])\rightarrow (1,\mathtt{c})\rightarrow \{(1,[i]),(2,[i])\}$ by (\ref{eqC6.3.13}) and (\ref{eqC6.3.14}).
\\
\\Case 2.5. $v=(p,[i])$ and $w=(q, [j])$, where $1\le p,q\le 2$, and $[i],[j]\in E$.
\indent\par Here, it is possible that $i=j$. Note that $\{(1,[i]),(2,[i])\}\rightarrow (s,\mathtt{c})\rightarrow (1,[1])\rightarrow (1,\mathtt{c})\rightarrow \{(1,[j]),(2,[j])\}$ by (\ref{eqC6.3.13}) and (\ref{eqC6.3.14}).
\\
\indent\par To settle Cases 3 and 4 (and forthcoming propositions), we require the following notation.

\begin{defn}\label{defnC6.3.7}
Set $\{\lambda_1,\lambda_2,\ldots, \lambda_{{s}\choose{\lceil s/2\rceil}}\}={{(\mathbb{N}_s,\mathtt{c})}\choose{\lceil s/2\rceil}}$, i.e., the set containing all $\lceil\frac{s}{2}\rceil$-subsets of $(\mathbb{N}_s,\mathtt{c})$. In particular, for $1\le i\le s$, let $\lambda_i=\{(i,\mathtt{c}), (i+1,\mathtt{c}), \ldots, (i+\lceil{\frac{s}{2}}\rceil-1, \mathtt{c})\}$, the sets containing $\lceil\frac{s}{2}\rceil$ vertices in consecutive (cyclic) order starting from $(i,\mathtt{c})$. For example, $\lambda_2=\{(2,\mathtt{c}),(3,\mathtt{c}),\ldots, (\lceil\frac{s}{2}\rceil+1,\mathtt{c})\}$. The denotation of the remaining $\lambda_i$'s can be arbitrary.
\end{defn}
\noindent Case 3. $|A_2|=\deg_T(\mathtt{c})$ (i.e., $E=\emptyset$) and $s< |A_2|\le {{s}\choose{\lceil{s/2}\rceil}}$. (If $s=3$, this case does not apply, and we refer to Case 1 instead.)
\indent\par Define an orientation $D_3$ for $\mathcal{H}$ as follows.
\begin{align}
&(2,[i])\rightarrow (1,[\alpha, i])\rightarrow (1,[i])\rightarrow (2,[\alpha, i])\rightarrow (2,[i]),\label{eqC6.3.15}\\
&\text{and }\lambda_i\rightarrow \{(1,[i]),(2,[i])\}\rightarrow \bar{\lambda}_i \label{eqC6.3.16}
\end{align}
for all $[i]\in A_2$ and $1\le \alpha \le \deg_T([i])-1$. We point out that the $\lceil\frac{s}{2}\rceil$-sets $\lambda_1,\lambda_2,\ldots,\lambda_{|A_2|}$ ($\lfloor\frac{s}{2}\rfloor$-sets $\bar{\lambda}_1,\bar{\lambda}_2,\ldots,\bar{\lambda}_{|A_2|}$ resp.) are used as `in-sets' (`out-sets' resp.) to construct $B^I_2$ ($B^O_2$ resp.). (See Figure \ref{figC6.3.5} for $D_3$ when $s=5$.)
\\
\\Claim 3: For all $v,w\in V(D_3)$, $d_{D_3}(v,w)\le 4$.
\\
\\Case 3.1. $v,w \in \{(1,[\alpha, i]),(2,[\alpha, i]),(1,[i]),(2,[i])\}$ for each $[i]\in A_2$ and $1\le \alpha \le \deg_T([i])-1$.
\indent\par This is clear since (\ref{eqC6.3.15}) is a directed $C_4$.
\\
\\Case 3.2. For each $[i],[j]\in A_2$, $i \neq j$, each $1\le \alpha \le \deg_T([i])-1$, and each $1\le \beta \le \deg_T([j])-1$,
\\(i) $v=(p,[\alpha, i]), w=(q,[\beta,j])$ for $p,q=1,2$, or
\\(ii) $v=(p,[\alpha, i]), w=(q,[j])$ for $p,q=1,2$, or
\\(iii) $v=(p,[i]), w=(q,[\beta,j])$ for $p,q=1,2$.
\indent\par By (\ref{eqC6.3.15})-(\ref{eqC6.3.16}), since $O^\mathtt{c}((p,[i]))=\bar{\lambda}_i \not\subseteq \bar{\lambda}_j=O^\mathtt{c}((q,[j]))$, there exists a vertex $(x, \mathtt{c})\in \bar{\lambda}_i\cap\lambda_j$ such that $(p,[\alpha, i])\rightarrow (p,[i])\rightarrow (x,\mathtt{c}) \rightarrow (3-q,[j]) \rightarrow (q,[\beta,j])$.
\\
\\Case 3.3. $v=(r_1,\mathtt{c})$ and $w=(r_2,\mathtt{c})$ for $r_1\neq r_2$ and $1\le r_1, r_2\le s$.
\indent\par Here, we want to prove a stronger claim, $d_{D_3}((r_1,\mathtt{c}), (r_2,\mathtt{c}))=2$. For $1\le k\le s$, let $x_k=(1,[k])$ and observe from (\ref{eqC6.3.16}) that $\lambda_k\rightarrow x_k\rightarrow \bar{\lambda}_k$. The subgraph induced by $V_1=(\mathbb{N}_s,\mathtt{c})$ and $V_2=\{x_k\mid 1\le k \le s\}$ is a complete bipartite graph $K(V_1,V_2)$. By Lemma \ref{lemC6.2.19}, $d_{D_3}((r_1,\mathtt{c}), (r_2,\mathtt{c}))=2$.
\\
\\Case 3.4. $v\in \{(1,[i]), (2,[i]), (1,[\alpha, i]), (2,[\alpha, i])\}$ for each $[i]\in A_2$ and $1\le \alpha \le \deg_T([i])-1$, and $w=(r,\mathtt{c})$ for $1\le r\le s$.
\indent\par Note that there exists some $1\le k\le s$ such that $d_{D_3}(v,(k,\mathtt{c}))\le 2$ by (\ref{eqC6.3.15})-(\ref{eqC6.3.16}). If $k=r$, we are done. If $k\neq r$, then $d_{D_3}((k,\mathtt{c}),(r,\mathtt{c}))=2$ by Case 3.3. Hence, it follows that $d_{D_3}(v,w)\le d_{D_3}(v,(k,\mathtt{c}))+d_{D_3}((k,\mathtt{c}),w)= 4$.
\\
\\Case 3.5. $v=(r,\mathtt{c})$ for $1\le r\le s$ and $w\in \{(1,[i]), (2,[i]), (1,[\alpha, i]), (2,[\alpha, i])\}$ for each $[i]\in A_2$ and $1\le \alpha \le \deg_T([i])-1$.
\indent\par Note that there exists some $1\le k\le s$ such that $d_{D_3}((k,\mathtt{c}),w)\le 2$ by (\ref{eqC6.3.15})-(\ref{eqC6.3.16}). If $k=r$, we are done. If $k\neq r$, then $d_{D_3}((r,\mathtt{c}),(k,\mathtt{c}))=2$ by Case 3.3. Hence, it follows that $d_{D_3}(v,w)\le d_{D_3}(v,(k,\mathtt{c}))+d_{D_3}((k,\mathtt{c}),w)= 4$.
\\
\\Case 3.6. $v=(p,[i])$ and $w=(q, [j])$, where $1\le p,q\le 2$ and $[i],[j]\in A_2$.
\indent\par This follows from the fact that $|O^\mathtt{c}((p,[i]))|>0$, $|I^\mathtt{c}((q,[j]))|>0$, and $d_{D_3}((r_1,\mathtt{c}), (r_2,\mathtt{c}))$ $=2$ for any $r_1\neq r_2$ and $1\le r_1, r_2\le s$ by Case 3.3.
\\
\\Case 4. $|A_2|<\deg_T(\mathtt{c})$ (i.e., $E\neq\emptyset$) and $s\le |A_2|\le {{s}\choose{\lceil{s/2}\rceil}}-1$. (If $s=3$, this case does not apply and we refer to Case 2 instead.)
\indent\par We define an orientation $D_4$ for $\mathcal{H}$ by making a slight enchancement to $D_3$. Noting that $|A_2|\le {{s}\choose{\lceil{s/2}\rceil}}-1$, we include in $D_4$ these extra arcs:
\begin{align*}
\lambda_{{s}\choose{\lceil{s/2}\rceil}}\rightarrow \{(1,[j]),(2,[j])\}\rightarrow \bar{\lambda}_{{s}\choose{\lceil{s/2}\rceil}}
\end{align*}
for all $[j]\in E$. (See Figure \ref{figC6.3.6} for $D_4$ when $s=5$.)
\\
\\Claim 4: For all $v,w\in V(D_4)$, $d_{D_4}(v,w)\le 4$.
\indent\par In view of the similarity between $D_3$ and $D_4$, it suffices to check the following.
\\
\\Case 4.1. $v\in \{(1,[i]), (2,[i])\}$ for each $[i]\in E$ and $1\le \alpha \le \deg_T([i])-1$, and $w=(r,\mathtt{c})$ for $1\le r\le s$.
\indent\par This follows from the fact that $|O^\mathtt{c}((p,[i]))|>0$ for $p=1,2$, and $d_{D_4}((r_1,\mathtt{c}), (r_2,\mathtt{c}))=2$ for any $r_1\neq r_2$ and $1\le r_1, r_2\le s$ by Case 3.3.
\\
\\Case 4.2. $v=(r,\mathtt{c})$ for $1\le r\le s$ and $w\in \{(1,[i]), (2,[i])\}$ for each $[i]\in E$ and $1\le \alpha \le \deg_T([i])-1$.
\indent\par This follows from the fact that $|I^\mathtt{c}((p,[i]))|>0$ for $p=1,2$, and $d_{D_4}((r_1,\mathtt{c}), (r_2,\mathtt{c}))=2$ for any $r_1\neq r_2$ and $1\le r_1, r_2\le s$ by Case 3.3.
\\
\\Case 4.3. $v=(p,[i])$ and $w=(q, [j])$, where $1\le p,q\le 2$ and $[i],[j]\in E$.
\indent\par This follows from the fact that $|O^\mathtt{c}((p,[i]))|>0$, $|I^\mathtt{c}((q,[j]))|>0$, and $d_{D_4}((r_1,\mathtt{c}), (r_2,\mathtt{c}))$ $=2$ for any $r_1\neq r_2$ and $1\le r_1, r_2\le s$ by Case 3.3.
\\
\indent\par Hence, $d(D_i)=4$ for $i=1,2,3,4$. Since every vertex lies in a directed $C_4$ for $D_i$ and $d(D_i)=4$, $\bar{d}(\mathcal{T})\le \max \{4, d(D_i)\}$ by Lemma \ref{lemC6.1.3}, and thus $\bar{d}(\mathcal{T})=4$ .
\qed

\begin{center}
\begin{tikzpicture}[thick,scale=0.7]%
\draw(-3,3)node[circle, draw, fill=black!100, inner sep=0pt, minimum width=6pt, label={[yshift=0cm] 90:{\small $(1,[1])$}}](1_1){};
\draw(-3,1)node[circle, draw, fill=black!100, inner sep=0pt, minimum width=6pt, label={[yshift=0cm]270:{\small $(2,[1])$}}](2_1){};
\draw(-6,3)node[circle, draw, fill=black!100, inner sep=0pt, minimum width=6pt, label={[] 180:{\small $(1,[1,1])$}}](1_11){};
\draw(-6,1)node[circle, draw, fill=black!100, inner sep=0pt, minimum width=6pt, label={[] 180:{\small $(2,[1,1])$}}](2_11){};

\draw(3,1)node[circle, draw, fill=black!100, inner sep=0pt, minimum width=6pt, label={[yshift=0cm] 90:{\small $(1,[2])$}}](1_2){};
\draw(3,-1)node[circle, draw, fill=black!100, inner sep=0pt, minimum width=6pt, label={[yshift=0cm]270:{\small $(2,[2])$}}](2_2){};
\draw(6,1)node[circle, draw, fill=black!100, inner sep=0pt, minimum width=6pt, label={[] 0:{\small $(1,[1,2])$}}](1_12){};
\draw(6,-1)node[circle, draw, fill=black!100, inner sep=0pt, minimum width=6pt, label={[] 0:{\small $(2,[1,2])$}}](2_12){};

\draw(-6,-1)node[circle, draw, fill=black!100, inner sep=0pt, minimum width=6pt, label={[] 180:{\small $(1,[1,3])$}}](1_13){};
\draw(-6,-3)node[circle, draw, fill=black!100, inner sep=0pt, minimum width=6pt, label={[] 180:{\small $(2,[1,3])$}}](2_13){};
\draw(-6,-5)node[circle, draw, fill=black!100, inner sep=0pt, minimum width=6pt, label={[] 180:{\small $(1,[2,3])$}}](1_23){};
\draw(-6,-7)node[circle, draw, fill=black!100, inner sep=0pt, minimum width=6pt, label={[] 180:{\small $(2,[2,3])$}}](2_23){};

\draw(-3,-1)node[circle, draw, fill=black!100, inner sep=0pt, minimum width=6pt, label={[yshift=-0.1cm, xshift=0.35cm] 270:{\small $(1,[3])$}}](1_3){};
\draw(-3,-3)node[circle, draw, fill=black!100, inner sep=0pt, minimum width=6pt, label={[yshift=0cm, xshift=0.35cm]270:{\small $(2,[3])$}}](2_3){};

\draw(0,2)node[circle, draw, inner sep=0pt, minimum width=3pt](1_u){\scriptsize $(1,\mathtt{c})$};
\draw(0,0)node[circle, draw, inner sep=0pt, minimum width=3pt](2_u){\scriptsize $(2,\mathtt{c})$};
\draw(0,-2)node[circle, draw, inner sep=0pt, minimum width=3pt](3_u){\scriptsize $(3,\mathtt{c})$};
\draw(0,-4)node[circle, draw, inner sep=0pt, minimum width=3pt](4_u){\scriptsize $(4,\mathtt{c})$};
\draw(0,-6)node[circle, draw, inner sep=0pt, minimum width=3pt](5_u){\scriptsize $(5,\mathtt{c})$};

\draw(3,-4)node[circle, draw, fill=black!100, inner sep=0pt, minimum width=6pt, label={[yshift=-0.1cm] 270:{\small $(1,[4])$}}](1_4){};
\draw(3,-6)node[circle, draw, fill=black!100, inner sep=0pt, minimum width=6pt, label= {[yshift=0cm]270:{\small $(2,[4])$}}](2_4){};
\draw(6,-4)node[circle, draw, fill=black!100, inner sep=0pt, minimum width=6pt, label={[] 0:{\small $(1,[1,4])$}}](1_14){};
\draw(6,-6)node[circle, draw, fill=black!100, inner sep=0pt, minimum width=6pt, label={[] 0:{\small $(2,[1,4])$}}](2_14){};

\draw[->, line width=0.3mm, >=latex, shorten <= 0.3cm, shorten >= 0.15cm](1_11)--(1_1);
\draw[->, line width=0.3mm, >=latex, shorten <= 0.3cm, shorten >= 0.15cm](2_11)--(2_1);
\draw[dashed, ->, line width=0.3mm, >=latex, shorten <= 0.3cm, shorten >= 0.15cm](1_1)--(2_11);
\draw[dashed, ->, line width=0.3mm, >=latex, shorten <= 0.3cm, shorten >= 0.15cm](2_1)--(1_11);

\draw[->, line width=0.3mm, >=latex, shorten <= 0.3cm, shorten >= 0.1cm](2_1)--(1_u);
\draw[->, line width=0.3mm, >=latex, shorten <= 0.3cm, shorten >= 0.1cm](1_1)--(1_u);

\draw[->, line width=0.3mm, >=latex, shorten <= 0.3cm, shorten >= 0.15cm](1_12)--(1_2);
\draw[->, line width=0.3mm, >=latex, shorten <= 0.3cm, shorten >= 0.15cm](2_12)--(2_2);
\draw[dashed, ->, line width=0.3mm, >=latex, shorten <= 0.3cm, shorten >= 0.15cm](1_2)--(2_12);
\draw[dashed, ->, line width=0.3mm, >=latex, shorten <= 0.3cm, shorten >= 0.15cm](2_2)--(1_12);

\draw[->, line width=0.3mm, >=latex, shorten <= 0.3cm, shorten >= 0.1cm](2_2)--(2_u);
\draw[->, line width=0.3mm, >=latex, shorten <= 0.3cm, shorten >= 0.1cm](1_2)--(2_u);

\draw[->, line width=0.3mm, >=latex, shorten <= 0.3cm, shorten >= 0.15cm](1_13)--(1_3);
\draw[->, line width=0.3mm, >=latex, shorten <= 0.3cm, shorten >= 0.15cm](2_13)--(2_3);
\draw[dashed, ->, line width=0.3mm, >=latex, shorten <= 0.3cm, shorten >= 0.15cm](1_3)--(2_13);
\draw[dashed, ->, line width=0.3mm, >=latex, shorten <= 0.3cm, shorten >= 0.15cm](2_3)--(1_13);

\draw[->, line width=0.3mm, >=latex, shorten <= 0.3cm, shorten >= 0.15cm](1_23)--(1_3);
\draw[->, line width=0.3mm, >=latex, shorten <= 0.3cm, shorten >= 0.15cm](2_23)--(2_3);
\draw[dashed, ->, line width=0.3mm, >=latex, shorten <= 0.3cm, shorten >= 0.15cm](1_3)--(2_23);
\draw[dashed, ->, line width=0.3mm, >=latex, shorten <= 0.3cm, shorten >= 0.15cm](2_3)--(1_23);

\draw[->, line width=0.3mm, >=latex, shorten <= 0.3cm, shorten >= 0.1cm](2_3)--(3_u);
\draw[->, line width=0.3mm, >=latex, shorten <= 0.3cm, shorten >= 0.1cm](1_3)--(3_u);

\draw[->, line width=0.3mm, >=latex, shorten <= 0.3cm, shorten >= 0.15cm](2_14)--(2_4);
\draw[->, line width=0.3mm, >=latex, shorten <= 0.3cm, shorten >= 0.15cm](1_14)--(1_4);
\draw[dashed, ->, line width=0.3mm, >=latex, shorten <= 0.3cm, shorten >= 0.15cm](1_4)--(2_14);
\draw[dashed, ->, line width=0.3mm, >=latex, shorten <= 0.3cm, shorten >= 0.15cm](2_4)--(1_14);

\draw[->, line width=0.3mm, >=latex, shorten <= 0.3cm, shorten >= 0.1cm](1_4)--(4_u);
\draw[->, line width=0.3mm, >=latex, shorten <= 0.3cm, shorten >= 0.1cm](2_4)--(4_u);
\draw[->, line width=0.3mm, >=latex, shorten <= 0.3cm, shorten >= 0.1cm](1_4)--(4_u);
\draw[->, line width=0.3mm, >=latex, shorten <= 0.3cm, shorten >= 0.1cm](2_4)--(4_u);
\draw[->, line width=0.3mm, >=latex, shorten <= 0.3cm, shorten >= 0.1cm](1_4)--(5_u);
\draw[->, line width=0.3mm, >=latex, shorten <= 0.3cm, shorten >= 0.1cm](2_4)--(5_u);
\end{tikzpicture}
{\captionof{figure}{Orientation $D_1$ for $\mathcal{H}$, Case 1. $s=5$, $\deg_T(\mathtt{c})=4$, $A_2=\{[1],[2],[3],[4]\}$.}\label{figC6.3.3}}
\end{center}

\begin{center}
\begin{tikzpicture}[thick,scale=0.7]%
\draw(-3,3)node[circle, draw, fill=black!100, inner sep=0pt, minimum width=6pt, label={[yshift=0cm] 90:{\small $(1,[1])$}}](1_1){};
\draw(-3,1)node[circle, draw, fill=black!100, inner sep=0pt, minimum width=6pt, label={[yshift=0cm]270:{\small $(2,[1])$}}](2_1){};
\draw(-6,3)node[circle, draw, fill=black!100, inner sep=0pt, minimum width=6pt, label={[] 180:{\small $(1,[1,1])$}}](1_11){};
\draw(-6,1)node[circle, draw, fill=black!100, inner sep=0pt, minimum width=6pt, label={[] 180:{\small $(2,[1,1])$}}](2_11){};

\draw(3,1)node[circle, draw, fill=black!100, inner sep=0pt, minimum width=6pt, label={[yshift=0cm] 90:{\small $(1,[2])$}}](1_2){};
\draw(3,-1)node[circle, draw, fill=black!100, inner sep=0pt, minimum width=6pt, label={[yshift=0cm]270:{\small $(2,[2])$}}](2_2){};
\draw(6,1)node[circle, draw, fill=black!100, inner sep=0pt, minimum width=6pt, label={[] 0:{\small $(1,[1,2])$}}](1_12){};
\draw(6,-1)node[circle, draw, fill=black!100, inner sep=0pt, minimum width=6pt, label={[] 0:{\small $(2,[1,2])$}}](2_12){};

\draw(-6,-1)node[circle, draw, fill=black!100, inner sep=0pt, minimum width=6pt, label={[] 180:{\small $(1,[1,3])$}}](1_13){};
\draw(-6,-3)node[circle, draw, fill=black!100, inner sep=0pt, minimum width=6pt, label={[] 180:{\small $(2,[1,3])$}}](2_13){};
\draw(-6,-5)node[circle, draw, fill=black!100, inner sep=0pt, minimum width=6pt, label={[] 180:{\small $(1,[2,3])$}}](1_23){};
\draw(-6,-7)node[circle, draw, fill=black!100, inner sep=0pt, minimum width=6pt, label={[] 180:{\small $(2,[2,3])$}}](2_23){};

\draw(-3,-1)node[circle, draw, fill=black!100, inner sep=0pt, minimum width=6pt, label={[yshift=-0.1cm, xshift=0.35cm] 270:{\small $(1,[3])$}}](1_3){};
\draw(-3,-3)node[circle, draw, fill=black!100, inner sep=0pt, minimum width=6pt, label={[yshift=0cm, xshift=0.35cm]270:{\small $(2,[3])$}}](2_3){};

\draw(0,2)node[circle, draw, inner sep=0pt, minimum width=3pt](1_u){\scriptsize $(1,\mathtt{c})$};
\draw(0,0)node[circle, draw, inner sep=0pt, minimum width=3pt](2_u){\scriptsize $(2,\mathtt{c})$};
\draw(0,-2)node[circle, draw, inner sep=0pt, minimum width=3pt](3_u){\scriptsize $(3,\mathtt{c})$};
\draw(0,-4)node[circle, draw, inner sep=0pt, minimum width=3pt](4_u){\scriptsize $(4,\mathtt{c})$};
\draw(0,-6)node[circle, draw, inner sep=0pt, minimum width=3pt](5_u){\scriptsize $(5,\mathtt{c})$};

\draw(3,-3)node[circle, draw, fill=black!100, inner sep=0pt, minimum width=6pt, label={[yshift=-0.1cm] 270:{\small $(1,[4])$}}](1_4){};
\draw(3,-5)node[circle, draw, fill=black!100, inner sep=0pt, minimum width=6pt, label= {[yshift=0cm]270:{\small $(2,[4])$}}](2_4){};
\draw(6,-3)node[circle, draw, fill=black!100, inner sep=0pt, minimum width=6pt, label={[] 0:{\small $(1,[1,4])$}}](1_14){};
\draw(6,-5)node[circle, draw, fill=black!100, inner sep=0pt, minimum width=6pt, label={[] 0:{\small $(2,[1,4])$}}](2_14){};

\draw(-3,-7)node[circle, draw, fill=black!100, inner sep=0pt, minimum width=6pt, label={[yshift=0cm] 180:{\small $(1,[5])$}}](1_5){};
\draw(-3,-9)node[circle, draw, fill=black!100, inner sep=0pt, minimum width=6pt, label= {[yshift=0cm]180:{\small $(2,[5])$}}](2_5){};

\draw(3,-7)node[circle, draw, fill=black!100, inner sep=0pt, minimum width=6pt, label={[yshift=0cm] 0:{\small $(1,[6])$}}](1_6){};
\draw(3,-9)node[circle, draw, fill=black!100, inner sep=0pt, minimum width=6pt, label= {[yshift=0cm] 0:{\small $(2,[6])$}}](2_6){};

\draw[->, line width=0.3mm, >=latex, shorten <= 0.3cm, shorten >= 0.15cm](1_11)--(1_1);
\draw[->, line width=0.3mm, >=latex, shorten <= 0.3cm, shorten >= 0.15cm](2_11)--(2_1);
\draw[dashed, ->, line width=0.3mm, >=latex, shorten <= 0.3cm, shorten >= 0.15cm](1_1)--(2_11);
\draw[dashed, ->, line width=0.3mm, >=latex, shorten <= 0.3cm, shorten >= 0.15cm](2_1)--(1_11);

\draw[->, line width=0.3mm, >=latex, shorten <= 0.3cm, shorten >= 0.1cm](2_1)--(1_u);
\draw[->, line width=0.3mm, >=latex, shorten <= 0.3cm, shorten >= 0.1cm](1_1)--(1_u);

\draw[->, line width=0.3mm, >=latex, shorten <= 0.3cm, shorten >= 0.15cm](1_12)--(1_2);
\draw[->, line width=0.3mm, >=latex, shorten <= 0.3cm, shorten >= 0.15cm](2_12)--(2_2);
\draw[dashed, ->, line width=0.3mm, >=latex, shorten <= 0.3cm, shorten >= 0.15cm](1_2)--(2_12);
\draw[dashed, ->, line width=0.3mm, >=latex, shorten <= 0.3cm, shorten >= 0.15cm](2_2)--(1_12);

\draw[->, line width=0.3mm, >=latex, shorten <= 0.3cm, shorten >= 0.1cm](2_2)--(2_u);
\draw[->, line width=0.3mm, >=latex, shorten <= 0.3cm, shorten >= 0.1cm](1_2)--(2_u);

\draw[->, line width=0.3mm, >=latex, shorten <= 0.3cm, shorten >= 0.15cm](1_13)--(1_3);
\draw[->, line width=0.3mm, >=latex, shorten <= 0.3cm, shorten >= 0.15cm](2_13)--(2_3);
\draw[dashed, ->, line width=0.3mm, >=latex, shorten <= 0.3cm, shorten >= 0.15cm](1_3)--(2_13);
\draw[dashed, ->, line width=0.3mm, >=latex, shorten <= 0.3cm, shorten >= 0.15cm](2_3)--(1_13);

\draw[->, line width=0.3mm, >=latex, shorten <= 0.3cm, shorten >= 0.15cm](1_23)--(1_3);
\draw[->, line width=0.3mm, >=latex, shorten <= 0.3cm, shorten >= 0.15cm](2_23)--(2_3);
\draw[dashed, ->, line width=0.3mm, >=latex, shorten <= 0.3cm, shorten >= 0.15cm](1_3)--(2_23);
\draw[dashed, ->, line width=0.3mm, >=latex, shorten <= 0.3cm, shorten >= 0.15cm](2_3)--(1_23);

\draw[->, line width=0.3mm, >=latex, shorten <= 0.3cm, shorten >= 0.1cm](2_3)--(3_u);
\draw[->, line width=0.3mm, >=latex, shorten <= 0.3cm, shorten >= 0.1cm](1_3)--(3_u);

\draw[->, line width=0.3mm, >=latex, shorten <= 0.3cm, shorten >= 0.15cm](2_14)--(2_4);
\draw[->, line width=0.3mm, >=latex, shorten <= 0.3cm, shorten >= 0.15cm](1_14)--(1_4);
\draw[dashed, ->, line width=0.3mm, >=latex, shorten <= 0.3cm, shorten >= 0.15cm](1_4)--(2_14);
\draw[dashed, ->, line width=0.3mm, >=latex, shorten <= 0.3cm, shorten >= 0.15cm](2_4)--(1_14);

\draw[->, line width=0.3mm, >=latex, shorten <= 0.3cm, shorten >= 0.1cm](1_4)--(4_u);
\draw[->, line width=0.3mm, >=latex, shorten <= 0.3cm, shorten >= 0.1cm](2_4)--(4_u);
\draw[->, line width=0.3mm, >=latex, shorten <= 0.3cm, shorten >= 0.1cm](1_4)--(4_u);
\draw[->, line width=0.3mm, >=latex, shorten <= 0.3cm, shorten >= 0.1cm](2_4)--(4_u);

\draw[->, line width=0.3mm, >=latex, shorten <= 0.3cm, shorten >= 0.1cm](1_5)--(5_u);
\draw[->, line width=0.3mm, >=latex, shorten <= 0.3cm, shorten >= 0.1cm](2_5)--(5_u);

\draw[->, line width=0.3mm, >=latex, shorten <= 0.3cm, shorten >= 0.1cm](1_6)--(5_u);
\draw[->, line width=0.3mm, >=latex, shorten <= 0.3cm, shorten >= 0.1cm](2_6)--(5_u);
\end{tikzpicture}
\captionsetup{justification=centering}
{\captionof{figure}{Orientation $D_2$ for $\mathcal{H}$, Case 2. $s=5$, $\deg_T(\mathtt{c})=6$,\\$A_2=\{[1],[2],[3],[4]\}$, $E=\{[5],[6]\}$.}\label{figC6.3.4}}
\end{center}

\begin{center}
\begin{tikzpicture}[thick,scale=0.7]%
\draw(3,2)node[circle, draw, fill=black!100, inner sep=0pt, minimum width=6pt, label={[yshift=0cm, xshift=0.4cm] 90:{\small $(1,[1])$}}](1_1){};
\draw(3,0)node[circle, draw, fill=black!100, inner sep=0pt, minimum width=6pt, label={[yshift=0cm, xshift=0.4cm] 270:{\small $(2,[1])$}}](2_1){};

\draw(6,2)node[circle, draw, fill=black!100, inner sep=0pt, minimum width=6pt, label={[] 0:{\small $(1,[1,1])$}}](1_11){};
\draw(6,0)node[circle, draw, fill=black!100, inner sep=0pt, minimum width=6pt, label={[] 0:{\small $(2,[1,1])$}}](2_11){};

\draw(-6,2)node[circle, draw, fill=black!100, inner sep=0pt, minimum width=6pt, label={[] 180:{\small $(1,[1,2])$}}](1_12){};
\draw(-6,0)node[circle, draw, fill=black!100, inner sep=0pt, minimum width=6pt, label={[] 180:{\small $(2,[1,2])$}}](2_12){};
\draw(-6,-2)node[circle, draw, fill=black!100, inner sep=0pt, minimum width=6pt, label={[] 180:{\small $(1,[2,2])$}}](1_22){};
\draw(-6,-4)node[circle, draw, fill=black!100, inner sep=0pt, minimum width=6pt, label={[] 180:{\small $(2,[2,2])$}}](2_22){};

\draw(-3,0)node[circle, draw, fill=black!100, inner sep=0pt, minimum width=6pt, label={[yshift=0.2cm, xshift=-0.1cm] 90:{\small $(1,[2])$}}](1_2){};
\draw(-3,-2)node[circle, draw, fill=black!100, inner sep=0pt, minimum width=6pt, label={[yshift=-0.2cm, xshift=-0.1cm] 270:{\small $(2,[2])$}}](2_2){};

\draw(0,4)node[circle, draw, inner sep=0pt, minimum width=3pt](1_u){\scriptsize $(1,\mathtt{c})$};
\draw(0,1.5)node[circle, draw, inner sep=0pt, minimum width=3pt](2_u){\scriptsize $(2,\mathtt{c})$};
\draw(0,-1)node[circle, draw, inner sep=0pt, minimum width=3pt](3_u){\scriptsize $(3,\mathtt{c})$};
\draw(0,-3.5)node[circle, draw, inner sep=0pt, minimum width=3pt](4_u){\scriptsize $(4,\mathtt{c})$};
\draw(0,-6)node[circle, draw, inner sep=0pt, minimum width=3pt](5_u){\scriptsize $(5,\mathtt{c})$};

\draw(3,-2)node[circle, draw, fill=black!100, inner sep=0pt, minimum width=6pt, label={[yshift=0cm, xshift=0.4cm] 270:{\small $(1,[3])$}}](1_3){};
\draw(3,-4)node[circle, draw, fill=black!100, inner sep=0pt, minimum width=6pt, label={[yshift=0cm, xshift=0.4cm] 270:{\small $(2,[3])$}}](2_3){};

\draw(6,-2)node[circle, draw, fill=black!100, inner sep=0pt, minimum width=6pt, label={[] 0:{\small $(1,[1,3])$}}](1_13){};
\draw(6,-4)node[circle, draw, fill=black!100, inner sep=0pt, minimum width=6pt, label={[] 0:{\small $(2,[1,3])$}}](2_13){};

\draw[->, line width=0.3mm, >=latex, shorten <= 0.2cm, shorten >= 0.15cm](1_12)--(1_2);
\draw[->, line width=0.3mm, >=latex, shorten <= 0.2cm, shorten >= 0.15cm](2_12)--(2_2);
\draw[dashed, ->, line width=0.3mm, >=latex, shorten <= 0.2cm, shorten >= 0.15cm](1_2)--(2_12);
\draw[dashed, ->, line width=0.3mm, >=latex, shorten <= 0.2cm, shorten >= 0.15cm](2_2)--(1_12);

\draw[->, line width=0.3mm, >=latex, shorten <= 0.2cm, shorten >= 0.15cm](1_22)--(1_2);
\draw[->, line width=0.3mm, >=latex, shorten <= 0.2cm, shorten >= 0.15cm](2_22)--(2_2);
\draw[dashed, ->, line width=0.3mm, >=latex, shorten <= 0.2cm, shorten >= 0.15cm](1_2)--(2_22);
\draw[dashed, ->, line width=0.3mm, >=latex, shorten <= 0.2cm, shorten >= 0.15cm](2_2)--(1_22);

\draw[densely dotted, ->, line width=0.3mm, >=latex, shorten <= 0.2cm, shorten >= 0.1cm](2_2)--(1_u);
\draw[densely dotted, ->, line width=0.3mm, >=latex, shorten <= 0.2cm, shorten >= 0.1cm](1_2)--(1_u);
\draw[densely dotted, ->, line width=0.3mm, >=latex, shorten <= 0.2cm, shorten >= 0.1cm](2_2)--(5_u);
\draw[densely dotted, ->, line width=0.3mm, >=latex, shorten <= 0.2cm, shorten >= 0.1cm](1_2)--(5_u);

\draw[->, line width=0.3mm, >=latex, shorten <= 0.2cm, shorten >= 0.15cm](2_13)--(2_3);
\draw[->, line width=0.3mm, >=latex, shorten <= 0.2cm, shorten >= 0.15cm](1_13)--(1_3);
\draw[dashed, ->, line width=0.3mm, >=latex, shorten <= 0.2cm, shorten >= 0.15cm](1_3) to [out=350, in=115] (2_13);
\draw[dashed, ->, line width=0.3mm, >=latex, shorten <= 0.2cm, shorten >= 0.15cm](2_3)--(1_13);

\draw[densely dotted, ->, line width=0.3mm, >=latex, shorten <= 0.2cm, shorten >= 0.1cm](1_3)--(1_u);
\draw[densely dotted, ->, line width=0.3mm, >=latex, shorten <= 0.2cm, shorten >= 0.1cm](2_3)--(1_u);
\draw[densely dotted, ->, line width=0.3mm, >=latex, shorten <= 0.2cm, shorten >= 0.1cm](1_3)--(2_u);
\draw[densely dotted, ->, line width=0.3mm, >=latex, shorten <= 0.2cm, shorten >= 0.1cm](2_3)--(2_u);

\draw[->, line width=0.3mm, >=latex, shorten <= 0.2cm, shorten >= 0.15cm](2_11)--(2_1);
\draw[->, line width=0.3mm, >=latex, shorten <= 0.2cm, shorten >= 0.15cm](1_11)--(1_1);
\draw[dashed, ->, line width=0.3mm, >=latex, shorten <= 0.2cm, shorten >= 0.15cm](1_1)--(2_11);
\draw[dashed, ->, line width=0.3mm, >=latex, shorten <= 0.2cm, shorten >= 0.15cm](2_1)--(1_11);

\draw[densely dotted, ->, line width=0.3mm, >=latex, shorten <= 0.2cm, shorten >= 0.1cm](1_1)--(4_u);
\draw[densely dotted, ->, line width=0.3mm, >=latex, shorten <= 0.2cm, shorten >= 0.1cm](2_1)--(4_u);
\draw[densely dotted, ->, line width=0.3mm, >=latex, shorten <= 0.2cm, shorten >= 0.1cm](1_1)--(5_u);
\draw[densely dotted, ->, line width=0.3mm, >=latex, shorten <= 0.2cm, shorten >= 0.1cm](2_1)--(5_u);
\end{tikzpicture}
{\captionof{figure}{\protect\footnotemark Orientation $D_3$ for $\mathcal{H}$, Case 3. $s=5$, $\deg_T(\mathtt{c})=5$, $A_2=\{[1],[2],[3],[4],[5]\}$. \\}\label{figC6.3.5}
\caption*{Note: In addition to the simplification noted in Figure \ref{figC6.2.2}, in Figures \ref{figC6.3.5} to \ref{figC6.4.18}, we use \lineexample{densely dotted} (\lineexample{densely dashdotdotted} resp.) lines to elucidate the `out-sets' $B^O_2$ and $B^O_3$ (complements of the `in-sets' $B^I_2$ and $B^I_3$ resp.); and in cases where both coincide, the densely dotted lines take precedent.}}
\end{center}
\footnotetext{For clarity, we only show the vertices $[\alpha,i]$ and $[i]$ for $i=1,2,3$.}

\begin{center}
\begin{tikzpicture}[thick,scale=0.7]%
\draw(3,2)node[circle, draw, fill=black!100, inner sep=0pt, minimum width=6pt, label={[yshift=0cm, xshift=0.4cm] 90:{\small $(1,[1])$}}](1_1){};
\draw(3,0)node[circle, draw, fill=black!100, inner sep=0pt, minimum width=6pt, label={[yshift=0cm, xshift=0.4cm] 270:{\small $(2,[1])$}}](2_1){};

\draw(6,2)node[circle, draw, fill=black!100, inner sep=0pt, minimum width=6pt, label={[] 0:{\small $(1,[1,1])$}}](1_11){};
\draw(6,0)node[circle, draw, fill=black!100, inner sep=0pt, minimum width=6pt, label={[] 0:{\small $(2,[1,1])$}}](2_11){};

\draw(-6,2)node[circle, draw, fill=black!100, inner sep=0pt, minimum width=6pt, label={[] 180:{\small $(1,[1,2])$}}](1_12){};
\draw(-6,0)node[circle, draw, fill=black!100, inner sep=0pt, minimum width=6pt, label={[] 180:{\small $(2,[1,2])$}}](2_12){};
\draw(-6,-2)node[circle, draw, fill=black!100, inner sep=0pt, minimum width=6pt, label={[] 180:{\small $(1,[2,2])$}}](1_22){};
\draw(-6,-4)node[circle, draw, fill=black!100, inner sep=0pt, minimum width=6pt, label={[] 180:{\small $(2,[2,2])$}}](2_22){};

\draw(-3,0)node[circle, draw, fill=black!100, inner sep=0pt, minimum width=6pt, label={[yshift=0.2cm, xshift=-0.25cm] 90:{\small $(1,[2])$}}](1_2){};
\draw(-3,-2)node[circle, draw, fill=black!100, inner sep=0pt, minimum width=6pt, label={[yshift=-0.3cm, xshift=-0.25cm] 270:{\small $(2,[2])$}}](2_2){};

\draw(0,2)node[circle, draw, inner sep=0pt, minimum width=3pt](1_u){\scriptsize $(1,\mathtt{c})$};
\draw(0,-0.5)node[circle, draw, inner sep=0pt, minimum width=3pt] (2_u){\scriptsize $(2,\mathtt{c})$};
\draw(0,-3)node[circle, draw, inner sep=0pt, minimum width=3pt](3_u){\scriptsize $(3,\mathtt{c})$};
\draw(0,-5.5)node[circle, draw, inner sep=0pt, minimum width=3pt](4_u){\scriptsize $(4,\mathtt{c})$};
\draw(0,-8)node[circle, draw, inner sep=0pt, minimum width=3pt](5_u){\scriptsize $(5,\mathtt{c})$};

\draw(3,-2)node[circle, draw, fill=black!100, inner sep=0pt, minimum width=6pt, label={[yshift=0cm, xshift=0.4cm] 270:{\small $(1,[3])$}}](1_3){};
\draw(3,-4)node[circle, draw, fill=black!100, inner sep=0pt, minimum width=6pt, label={[yshift=0cm, xshift=0.4cm] 270:{\small $(2,[3])$}}](2_3){};

\draw(6,-2)node[circle, draw, fill=black!100, inner sep=0pt, minimum width=6pt, label={[] 0:{\small $(1,[1,3])$}}](1_13){};
\draw(6,-4)node[circle, draw, fill=black!100, inner sep=0pt, minimum width=6pt, label={[] 0:{\small $(2,[1,3])$}}](2_13){};

\draw(-3,-6)node[circle, draw, fill=black!100, inner sep=0pt, minimum width=6pt, label={[] 180:{\small $(1,[10])$}}](1_10){};
\draw(-3,-8)node[circle, draw, fill=black!100, inner sep=0pt, minimum width=6pt, label={[] 180:{\small $(2,[10])$}}](2_10){};

\draw(3,-6)node[circle, draw, fill=black!100, inner sep=0pt, minimum width=6pt, label={[] 0:{\small $(1,[11])$}}](1_111){};
\draw(3,-8)node[circle, draw, fill=black!100, inner sep=0pt, minimum width=6pt, label={[] 0:{\small $(2,[11])$}}](2_111){};

\draw[->, line width=0.3mm, >=latex, shorten <= 0.2cm, shorten >= 0.15cm](1_12)--(1_2);
\draw[->, line width=0.3mm, >=latex, shorten <= 0.2cm, shorten >= 0.15cm](2_12)--(2_2);
\draw[dashed, ->, line width=0.3mm, >=latex, shorten <= 0.2cm, shorten >= 0.15cm](1_2)--(2_12);
\draw[dashed, ->, line width=0.3mm, >=latex, shorten <= 0.2cm, shorten >= 0.15cm](2_2)--(1_12);

\draw[->, line width=0.3mm, >=latex, shorten <= 0.2cm, shorten >= 0.15cm](1_22)--(1_2);
\draw[->, line width=0.3mm, >=latex, shorten <= 0.2cm, shorten >= 0.15cm](2_22)--(2_2);
\draw[dashed, ->, line width=0.3mm, >=latex, shorten <= 0.2cm, shorten >= 0.15cm](1_2)--(2_22);
\draw[dashed, ->, line width=0.3mm, >=latex, shorten <= 0.2cm, shorten >= 0.15cm](2_2)--(1_22);

\draw[densely dotted, ->, line width=0.3mm, >=latex, shorten <= 0.2cm, shorten >= 0.1cm](2_2)--(1_u);
\draw[densely dotted, ->, line width=0.3mm, >=latex, shorten <= 0.2cm, shorten >= 0.1cm](1_2)--(1_u);
\draw[densely dotted, ->, line width=0.3mm, >=latex, shorten <= 0.2cm, shorten >= 0.1cm](2_2)--(5_u);
\draw[densely dotted, ->, line width=0.3mm, >=latex, shorten <= 0.2cm, shorten >= 0.1cm](1_2)--(5_u);

\draw[->, line width=0.3mm, >=latex, shorten <= 0.2cm, shorten >= 0.15cm](2_13)--(2_3);
\draw[->, line width=0.3mm, >=latex, shorten <= 0.2cm, shorten >= 0.15cm](1_13)--(1_3);
\draw[dashed, ->, line width=0.3mm, >=latex, shorten <= 0.2cm, shorten >= 0.15cm](1_3) to [out=350, in=115] (2_13);
\draw[dashed, ->, line width=0.3mm, >=latex, shorten <= 0.2cm, shorten >= 0.15cm](2_3)--(1_13);

\draw[densely dotted, ->, line width=0.3mm, >=latex, shorten <= 0.2cm, shorten >= 0.1cm](1_3)--(1_u);
\draw[densely dotted, ->, line width=0.3mm, >=latex, shorten <= 0.2cm, shorten >= 0.1cm](2_3)--(1_u);
\draw[densely dotted, ->, line width=0.3mm, >=latex, shorten <= 0.2cm, shorten >= 0.1cm](1_3)--(2_u);
\draw[densely dotted, ->, line width=0.3mm, >=latex, shorten <= 0.2cm, shorten >= 0.1cm](2_3)--(2_u);

\draw[->, line width=0.3mm, >=latex, shorten <= 0.2cm, shorten >= 0.15cm](2_11)--(2_1);
\draw[->, line width=0.3mm, >=latex, shorten <= 0.2cm, shorten >= 0.15cm](1_11)--(1_1);
\draw[dashed, ->, line width=0.3mm, >=latex, shorten <= 0.2cm, shorten >= 0.15cm](1_1)--(2_11);
\draw[dashed, ->, line width=0.3mm, >=latex, shorten <= 0.2cm, shorten >= 0.15cm](2_1)--(1_11);

\draw[densely dotted, ->, line width=0.3mm, >=latex, shorten <= 0.2cm, shorten >= 0.1cm](1_1)--(4_u);
\draw[densely dotted, ->, line width=0.3mm, >=latex, shorten <= 0.2cm, shorten >= 0.1cm](2_1)--(4_u);
\draw[densely dotted, ->, line width=0.3mm, >=latex, shorten <= 0.2cm, shorten >= 0.1cm](1_1)--(5_u);
\draw[densely dotted, ->, line width=0.3mm, >=latex, shorten <= 0.2cm, shorten >= 0.1cm](2_1)--(5_u);
\draw[->, line width=0.3mm, >=latex, shorten <= 0.2cm, shorten >= 0.1cm](1_10)--(2_u);
\draw[->, line width=0.3mm, >=latex, shorten <= 0.2cm, shorten >= 0.1cm](2_10)--(2_u);
\draw[->, line width=0.3mm, >=latex, shorten <= 0.2cm, shorten >= 0.1cm](1_10)--(4_u);
\draw[->, line width=0.3mm, >=latex, shorten <= 0.2cm, shorten >= 0.1cm](2_10)--(4_u);

\draw[->, line width=0.3mm, >=latex, shorten <= 0.2cm, shorten >= 0.1cm](1_111)--(2_u);
\draw[->, line width=0.3mm, >=latex, shorten <= 0.2cm, shorten >= 0.1cm](2_111)--(2_u);
\draw[->, line width=0.3mm, >=latex, shorten <= 0.2cm, shorten >= 0.1cm](1_111)--(4_u);
\draw[->, line width=0.3mm, >=latex, shorten <= 0.2cm, shorten >= 0.1cm](2_111)--(4_u);
\end{tikzpicture}
\captionsetup{justification=centering}
{
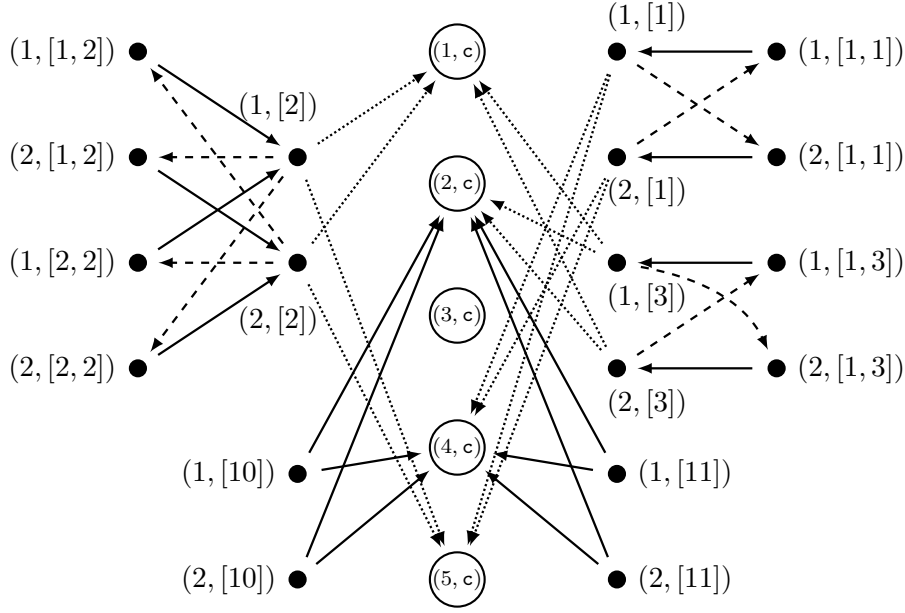
\captionof{figure}{\protect\footnotemark Orientation $D_4$ for $\mathcal{H}$, Case 4. $s=5$, $\deg_T(\mathtt{c})=11$, \\$A_2=\{[1],[2],\ldots,[9]\}$, $E=\{[10],[11]\}$.}\label{figC6.3.6}}
\end{center}
\footnotetext{Assume $ \lambda_{{s}\choose{\lceil{s/2}\rceil}}=\{(1,\mathtt{c}), (3,\mathtt{c}),(5,\mathtt{c})\}$. For clarity, we only show the vertices $[\alpha,i]$ and $[i]$ for $i=1,2,3,10,11$.}
\begin{cor}\label{corC6.3.8}
Suppose $s\ge 3$ for a $\mathcal{T}$. If
\\(i) $|A_{\ge 2}|\le {{s}\choose{\lceil{s/2}\rceil}}-1$, or 
\\(ii) $|A_{\ge 2}|\le {{s}\choose{\lceil{s/2}\rceil}}$ and $|A_{\ge 2}|=\deg_T(\mathtt{c})$,
\\then $\mathcal{T}\in \mathscr{C}_0$.
\end{cor}
\noindent\textit{Proof}: Note in the proof of Proposition \ref{ppnC6.3.5} that every vertex lies in a directed $C_4$ for each orientation $D_i$ and $d(D_i)\le 4$, for $i=1,2,3,4$. By Lemma \ref{lemC6.1.3}, $\bar{d}(\mathcal{T})\le \max \{4, d(D_i)\}$ for $i=1,2,3,4$, and thus $\bar{d}(\mathcal{T})=4$.
\qed

\indent\par For the remaining propositions of this section, we consider even $s\ge 4$. The proof of Proposition \ref{ppnC6.3.9} is centered on a reduction to cross-intersecting antichains and Theorems \ref{thmC6.2.4} and \ref{thmC6.2.6}.
\begin{ppn}\label{ppnC6.3.9}
Suppose $s\ge 4$ is even, $A_2=\emptyset$ and $A_3\neq \emptyset$ for a $\mathcal{T}$. Then, $\mathcal{T}\in \mathscr{C}_0$ if and only if $|A_3|\le {{s}\choose{s/2}}+{{s}\choose{(s/2)+1}}-2$.
\end{ppn}
\noindent\textit{Proof}: ($\Rightarrow$) Since $\mathcal{T}\in \mathscr{C}_0$, there exists an orientation $D$ of $\mathcal{T}$, where $d(D)=4$. As $A_3\neq\emptyset$, we assume (\ref{eqC6.3.3})-(\ref{eqC6.3.5}) here. By Sperner's theorem, $|B^O_3|\le {{s}\choose{s/2}}$ and $|B^I_3|\le {{s}\choose{s/2}}$. If $|B^O_3|=0$ or $|B^I_3|=0$, then $|A_3|=|B^O_3|+|B^I_3|\le {{s}\choose{s/2}}$. Therefore, we assume $|B^O_3|>0$ and $|B^I_3|>0$.
\indent\par Observe also that for each $[i]\in A^O_3$ and each $[j]\in A^I_3$, $d_D((1,[i]), (1,[j]))=4$ implies $X\cap Y\neq \emptyset$ for all $X\in B^O_3$ and all $Y\in B^I_3$. By Theorem \ref{thmC6.2.4}, $|A_3|=|B^O_3|+|B^I_3|\le {{s}\choose{s/2}}+{{s}\choose{(s/2)+1}}$. Suppose $|A_3|>{{s}\choose{s/2}}+{{s}\choose{(s/2)+1}}-2$ for a contradiction. It follows from Theorems \ref{thmC6.2.4} and \ref{thmC6.2.6} that $\{B^O_3, B^I_3\}=\{\mathscr{A}, \mathscr{B}\}$, where
\\(1) $\mathscr{A}={{(\mathbb{N}_s,\mathtt{c})}\choose{s/2}}$, $\mathscr{B}={{(\mathbb{N}_s,\mathtt{c})}\choose{(s/2)+1}}$, or 
\\(2) $\mathscr{A}={{(\mathbb{N}_s,\mathtt{c})}\choose{s/2}}$, $\mathscr{B}\subset{{(\mathbb{N}_s,\mathtt{c})}\choose{(s/2)+1}}$ and $|\mathscr{B}|={{s}\choose{(s/2)+1}}-1$, or
\\(3) $\mathscr{A}\subset{{(\mathbb{N}_s,\mathtt{c})}\choose{s/2}}$, $|\mathscr{A}|={{s}\choose{s/2}}-1$, and $\mathscr{B}={{(\mathbb{N}_s,\mathtt{c})}\choose{(s/2)+1}}$.
\\
\\Case 1. $B^O_3={{(\mathbb{N}_s,\mathtt{c})}\choose{s/2}}$.
\indent\par Let $[i]\in A^I_3$. For all $[j]\in A^O_3$ and $p=1,2$, $d_D((1,[1,j]),(p,[i]))=3$ implies $X\cap I^\mathtt{c}((p,[i]))\neq\emptyset$ for all $X\in B^O_3$. It follows that $|I^\mathtt{c}((p,[i]))|\ge\frac{s}{2}+1$ for all $p=1,2$. As a result, $O^\mathtt{c}((1,[i]))$ and $O^\mathtt{c}((2,[i]))$ are independent. Otherwise, $O^\mathtt{c}((1,[i]))\cup O^\mathtt{c}((2,[i]))\subset X$ for some $X\in B^O_3$, which contradicts Lemma \ref{lemC6.2.17}(b).
\\
\\Subcase 1.1. $B^I_3={{(\mathbb{N}_s,\mathtt{c})}\choose{(s/2)+1}}$. 
\indent\par Let $[i^*]\in A^I_3$. For all $[i]\in A^I_3-\{[i^*]\}$ and $p=1,2$, $d_D((p,[i^*]),(1,[1,i]))=3$ implies $X\cap O^\mathtt{c}((p,[i^*]))\neq \emptyset$ for all $X\in B^I_3-\{I^\mathtt{c}((3,[i^*]))\}$. Consequently, we have either $|O^\mathtt{c}((p,[i^*]))|\ge\frac{s}{2}$, or $O^\mathtt{c}((p,[i^*]))=O^\mathtt{c}((3,[i^*]))$ for each $p=1,2$. Since $|I^\mathtt{c}((p,[i^*]))|\ge\frac{s}{2}+1$, $|O^\mathtt{c}((p,[i^*]))|<\frac{s}{2}$. Hence, $O^\mathtt{c}((1,[i^*]))=O^\mathtt{c}((3,[i^*]))=O^\mathtt{c}((2,[i^*]))$, a contradiction to $O^\mathtt{c}((1,[i^*]))$ and $O^\mathtt{c}((2,[i^*]))$ being independent.
\\
\\Subcase 1.2. $B^I_3\subset{{(\mathbb{N}_s,\mathtt{c})}\choose{(s/2)+1}}$ and $|B^I_3|={{s}\choose{(s/2)+1}}-1$.
\indent\par Let ${{(\mathbb{N}_s,\mathtt{c})}\choose{(s/2)+1}}-B^I_3=\{\psi\}$. If $|O^\mathtt{c}((p,[i]))|\le \frac{s}{2}-2$ for some $[i]\in A^I_3$ and some $p=1,2$, then there are ${{|I^\mathtt{c}((p,[i]))|}\choose{(s/2)+1}}\ge {{(s/2)+2}\choose{(s/2)+1}}=\frac{s}{2}+2\ge 4>2$ $(\frac{s}{2}+1)$-subsets of $I^\mathtt{c}((p,[i]))$, i.e., $X\subset I^\mathtt{c}((p,[i]))$ for some $X\in B^I_3-\{I^\mathtt{c}((3,[i]))\}$, a contradiction to Lemma \ref{lemC6.2.18}(a). So, for each $[i]\in A^I_3$ and each $p=1,2$, we have either $|O^\mathtt{c}((p,[i]))|\ge \frac{s}{2}$, or $O^\mathtt{c}((p,[i]))=O^\mathtt{c}((3,[i]))$, or $O^\mathtt{c}((p,[i]))=\bar{\psi}$.  Since $|I^\mathtt{c}((p,[i]))|\ge\frac{s}{2}+1$ and $O^\mathtt{c}((1,[i]))$ and $O^\mathtt{c}((2,[i]))$ are independent, we may assume without loss of generality that $O^\mathtt{c}((1,[i]))=\bar{\psi}$ and $O^\mathtt{c}((2,[i]))=O^\mathtt{c}((3,[i]))$ for each $[i]\in A^I_3$.
\indent\par Now, we claim that there exists some $[j]\in A^I_3$ such that $|\bar{\psi}\cup O^\mathtt{c}((3,[j]))|=\frac{s}{2}$. Note that $|\bar{\psi}\cup O^\mathtt{c}((3,[j]))|=\frac{s}{2}$ if and only if $|\bar{\psi}\cap O^\mathtt{c}((3,[j]))|=\frac{s}{2}-2$ if and only if $|\psi\cap O^\mathtt{c}((3,[j]))|=1$. Since ${{|\bar{\psi}|}\choose{(s/2)-2}}{{|\psi|}\choose{1}}={{(s/2)-1}\choose{(s/2)-2}}{{(s/2)+1}\choose{1}}=\frac{s^2}{4}-1\ge 3$ and $|B^I_3|={{s}\choose{(s/2)+1}}-1$, the claim follows. Hence, $O^\mathtt{c}((1,[j]))\cup O^\mathtt{c}((2,[j]))=O^\mathtt{c}((1,[j]))\cup O^\mathtt{c}((3,[j]))=\bar{\psi}\cup O^\mathtt{c}((3,[j]))=O^\mathtt{c}((3,[k]))$ for some $[k]\in A^O_3$. This contradicts Lemma \ref{lemC6.2.17}(b).
\\
\\Case 2. $B^I_3={{(\mathbb{N}_s,\mathtt{c})}\choose{s/2}}$.
\indent\par If $B^O_3={{(\mathbb{N}_s,\mathtt{c})}\choose{(s/2)+1}}$ ($B^O_3\subset{{(\mathbb{N}_s,\mathtt{c})}\choose{(s/2)+1}}$ and $|B^O_3|={{s}\choose{(s/2)+1}}-1$ resp.), then the result follows from Subcase 1.1 (Subcase 1.2 resp.) by the Duality Lemma.
\\
\\Case 3. $B^O_3={{(\mathbb{N}_s,\mathtt{c})}\choose{(s/2)+1}}$, $B^I_3\subset{{(\mathbb{N}_s,\mathtt{c})}\choose{s/2}}$ and $|B^I_3|={{s}\choose{s/2}}-1$
\indent\par Let $[i]\in A^I_3$. For all $[j]\in A^O_3$ and $p=1,2$, $d_D((1,[1,j]),(p,[i]))=3$ implies $X\cap I^\mathtt{c}((p,[i]))\neq\emptyset$ for all $X\in B^O_3$. It follows that $|I^\mathtt{c}((p,[i]))|\ge\frac{s}{2}$ for all $p=1,2$. Furthermore, $O^\mathtt{c}((1,[i]))$ and $O^\mathtt{c}((2,[i]))$ are independent. Otherwise, $O^\mathtt{c}((1,[i]))\cup O^\mathtt{c}((2,[i]))\subset X$ for some $X\in B^O_3$, which contradicts Lemma \ref{lemC6.2.17}(b).
\indent\par Let ${{(\mathbb{N}_s,\mathtt{c})}\choose{s/2}}-B^I_3=\{\lambda\}$. If $|O^\mathtt{c}((p,[i]))|\le \frac{s}{2}-1$ for some $p=1,2$, and some $[i]\in A^I_3$, then there are ${{|I^\mathtt{c}((p,[i]))|}\choose{s/2}}\ge {{(s/2)+1}\choose{s/2}}=\frac{s}{2}+1\ge 3$ $\frac{s}{2}$-subsets of $I^\mathtt{c}((p,[i]))$, i.e., $X\subset I^\mathtt{c}((p,[i]))$ for some $X\in B^I_3-\{I^\mathtt{c}((3,[i]))\}$, a contradiction to Lemma \ref{lemC6.2.18}(a). Consequently, we have for each $p=1,2$, $|O^\mathtt{c}((p,[i]))|\ge\frac{s}{2}+1$, or $O^\mathtt{c}((p,[i]))=O^\mathtt{c}((3,[i]))$, or $O^\mathtt{c}((p,[i]))=\bar{\lambda}$. Since $|I^\mathtt{c}((p,[i]))|\ge\frac{s}{2}$ and $O^\mathtt{c}((1,[i]))$ and $O^\mathtt{c}((2,[i]))$ are independent, we may assume without loss of generality that $O^\mathtt{c}((1,[i]))=\bar{\lambda}$ and $O^\mathtt{c}((2,[i]))=O^\mathtt{c}((3,[i]))$ for each $[i]\in A^I_3$.
\indent\par Now, we claim that there exists some $[j]\in A^I_3$ such that $|\bar{\lambda}\cup O^\mathtt{c}((3,[j]))|=\frac{s}{2}+1$. Note that $|\bar{\lambda}\cup O^\mathtt{c}((3,[j]))|=\frac{s}{2}+1$ if and only if $|\bar{\lambda}\cap O^\mathtt{c}((3,[j]))|=\frac{s}{2}-1$ if and only if $|\lambda\cap O^\mathtt{c}((3,[j]))|=1$. Since ${{|\bar{\lambda}|}\choose{(s/2)-1}}{{|\lambda|}\choose{1}}={{s/2}\choose{(s/2)-1}}{{s/2}\choose{1}}=\frac{s^2}{4}\ge 4$ and $|B^I_3|={{s}\choose{(s/2)+1}}-1$, the claim follows. Hence, $O^\mathtt{c}((1,[j]))\cup O^\mathtt{c}((2,[j]))=\bar{\lambda}\cup O^\mathtt{c}((3,[j]))=O^\mathtt{c}((3,[k]))$ for some $[k]\in A^O_3$. This contradicts Lemma \ref{lemC6.2.17}(b).
\\
\\Case 4. $B^I_3={{(\mathbb{N}_s,\mathtt{c})}\choose{(s/2)+1}}$, $B^O_3\subset{{(\mathbb{N}_s,\mathtt{c})}\choose{s/2}}$ and $|B^O_3|={{s}\choose{s/2}}-1$.
\indent\par This follows from Case 3 by the Duality Lemma.
\\
\\($\Leftarrow$) If $|A_{\ge 3}|\le {{s}\choose{s/2}}-1$, then by Corollary \ref{corC6.3.8}(i), $\mathcal{T}\in\mathscr{C}_0$. Hence, we assume $|A_{\ge 3}|\ge{{s}\choose{s/2}}$ hereafter, on top of the hypothesis that $|A_3|\le {{s}\choose{s/2}}+{{s}\choose{(s/2)+1}}-2$. If $|A_3|\ge {{s}\choose{s/2}}-2$, define $A^{\diamond}_3=A_3$. Otherwise, let $A^{\diamond}_3=A_3\cup A^*$, where $A^*$ is an arbitrary subset of $A_{\ge 4}$ such that $|A^{\diamond}_3|={{s}\choose{s/2}}-2$. Then, let $A^{\diamond}_4=A_{\ge 4}-A^{\diamond}_3$. Furthermore, assume without loss of generality that $A^{\diamond}_3=\{[i]\mid i\in\mathbb{N}_{|A^{\diamond}_3|}\}$ and $A^{\diamond}_4=\{[i]\mid i\in\mathbb{N}_{|A^{\diamond}_3|+|A^{\diamond}_4|}-\mathbb{N}_{|A^{\diamond}_3|}\}$.
\indent\par Let $\mathcal{H}=T(t_1,t_2,\ldots, t_n)$ be the subgraph of $\mathcal{T}$, where $t_\mathtt{c}=s$, $t_{[i]}=3$ for all $[i]\in \mathcal{T}(A^{\diamond}_3)$, $t_{[j]}=4$ for all $[j]\in \mathcal{T}(A^{\diamond}_4)$ and $t_v=2$ otherwise. We will use $A_j$ for $\mathcal{H}(A_j)$ for the remainder of this proof. Let ${{(\mathbb{N}_s,\mathtt{c})}\choose{(s/2)+1}}=\{\psi_i\mid i=1,2,\ldots, {{s}\choose{(s/2)+1}}\}$ and recall $\lambda_i$ from Definition \ref{defnC6.3.7}. Define an orientation $D$ of $\mathcal{H}$ as follows.
\begin{align}
& (3,[i])\rightarrow \{(1,[\alpha,i]),(2,[\alpha,i])\}\rightarrow \{(1,[i]),(2,[i])\},\text{ and}\label{eqC6.3.17}\\
& \bar{\lambda}_1=\lambda_{\frac{s}{2}+1} \rightarrow (1,[i])\rightarrow \lambda_1\rightarrow (2,[i])\rightarrow \lambda_{\frac{s}{2}+1}\label{eqC6.3.18}
\end{align}
for all $1\le i\le{{s}\choose{s/2}}-2$ and all $1\le \alpha\le \deg_T([i])-1$.
\begin{align}
&\lambda_{i+1}\rightarrow (3,[i])\rightarrow \bar{\lambda}_{i+1} \label{eqC6.3.19}
\end{align}
for all $1\le i\le \frac{s}{2}-1$.
\begin{align}
&\lambda_{i+2}\rightarrow (3,[i])\rightarrow \bar{\lambda}_{i+2} \label{eqC6.3.20}
\end{align}
for all $\frac{s}{2}\le i\le {{s}\choose{s/2}}-2$, i.e., excluding $\lambda_1$ and $\lambda_{\frac{s}{2}+1}$, the $\frac{s}{2}$-sets $\lambda_i$'s are used as `in-sets' to construct $B^I_3$.
\begin{align}
& \{(1,[j]),(2,[j])\}\rightarrow \{(1,[\beta,j]),(2,[\beta,j])\}\rightarrow (3,[j]),\label{eqC6.3.21}\\
& \lambda_1 \rightarrow (1,[j])\rightarrow \lambda_{\frac{s}{2}+1}\rightarrow (2,[j])\rightarrow \lambda_1,\text{ and}\label{eqC6.3.22}\\
&\bar{\psi}_{j+2-{{s}\choose{s/2}}}\rightarrow (3,[j]) \rightarrow \psi_{j+2-{{s}\choose{s/2}}} \label{eqC6.3.23}
\end{align}
for all ${{s}\choose{s/2}}-1\le j\le |A_3|$ and all $1\le \beta\le \deg_T([j])-1$, i.e., the $(\frac{s}{2}+1)$-sets $\psi_1,\psi_2,\ldots, \psi_{|A_3|+2-{{s}\choose{s/2}}}$ are used as `out-sets' to construct $B^O_3$.
\begin{align}
&(2,[\gamma,k])\rightarrow \{(2,[k]),(4,[k])\}\rightarrow (1,[\gamma,k])\rightarrow \{(1,[k]), (3,[k])\}\rightarrow (2,[\gamma,k]),\label{eqC6.3.24}\\
&\text{and }\lambda_{\frac{s}{2}+1} \rightarrow \{(1,[k]), (4,[k])\}\rightarrow \lambda_1 \rightarrow \{(2,[k]), (3,[k])\}\rightarrow \lambda_{\frac{s}{2}+1}\label{eqC6.3.25}
\end{align}
for all $[k]\in A_4$ and all $1\le  \gamma\le \deg_T([k])-1$.
\begin{align}
&\lambda_1\rightarrow \{(1,[l]), (2,[l])\} \rightarrow \lambda_{\frac{s}{2}+1} \label{eqC6.3.26}
\end{align}
for any $[l]\in E$. (See Figures \ref{figC6.3.7} and \ref{figC6.3.8} when $s=4$.)
\\
\\Claim: $d_{D}(v,w)\le 4$ for all $v,w\in V(D)$.
\\
\\Case 1.1. $v,w \in \{(1,[\alpha,i]),(2,[\alpha,i]),(1,[i]),(2,[i]),(3,[i])\}$ for each $1\le i\le {{s}\choose{s/2}}-2$ and $1\le\alpha\le \deg_T([i])-1$.
\indent\par By (\ref{eqC6.3.18}), $(1,[i])\rightarrow \lambda_1$ and $(2,[i])\rightarrow \lambda_{\frac{s}{2}+1}$. For each $1\le i\le {{s}\choose{s/2}}-2$, since $I^\mathtt{c}((3,[i]))\in {{(\mathbb{N}_s,\mathtt{c})}\choose{s/2}}-\{\lambda_1, \lambda_{\frac{s}{2}+1}\}$ by (\ref{eqC6.3.19})-(\ref{eqC6.3.20}), there exist some $x_i\in\lambda_1\cap I^\mathtt{c}((3,[i]))$ and some $y_i\in\lambda_{\frac{s}{2}+1}\cap I^\mathtt{c}((3,[i]))$ such that $\{x_i,y_i\}\rightarrow (3,[i])$. With (\ref{eqC6.3.17}), this case follows.
\\
\\Case 1.2. $v,w \in \{(1,[\alpha,i]),(2,[\alpha,i]),(1,[i]),(2,[i]),(3,[i])\}$ for each ${{s}\choose{s/2}}-1\le i\le |A_3|$, and $1\le\alpha\le \deg_T([i])-1$.
\indent\par By (\ref{eqC6.3.22}), $\lambda_1\rightarrow (1,[i])$, and $\lambda_{\frac{s}{2}+1}\rightarrow (2,[i])$. For each ${{s}\choose{s/2}}-1\le i \le \deg_T(\mathtt{c})$, since $O^\mathtt{c}((3,[i]))\in  {{(\mathbb{N}_s,\mathtt{c})}\choose{(s/2)+1}}$ and $\min\{|O^\mathtt{c}((3,[i]))|+|\lambda_1|, |O^\mathtt{c}((3,[i]))|+|\lambda_{\frac{s}{2}+1}|\}>s$,  there exist some $x_i\in O^\mathtt{c}((3,[i]))\cap \lambda_1$ and some $y_i\in O^\mathtt{c}((3,[i]))\cap \lambda_{\frac{s}{2}+1}$ such that $(3,[i])\rightarrow \{x_i, y_i\}$. With (\ref{eqC6.3.21}), this case follows.
\\
\\Case 1.3. $v,w \in \{(1,[\alpha,i]),(2,[\alpha,i]),(1,[i]),(2,[i]), (3,[i]), (4,[i])\}$ for each $[i]\in A_4$ and $1\le\alpha\le \deg_T([i])-1$.
\indent\par This is clear since (\ref{eqC6.3.24}) guarantees a directed $C_4$.
\\
\\Case 2.  For each $1\le i,j \le {{s}\choose{s/2}}-2$, $i\neq j$, each $1\le \alpha\le \deg_T([i])-1$, and each $1\le \beta\le \deg_T([j])-1$,
\\(i) $v=(p,[\alpha,i]), w=(q,[j])$ for each $p=1,2$ and $q=1,2,3$.
\indent\par By (\ref{eqC6.3.17})-(\ref{eqC6.3.20}), $\{(1,[\alpha,i]),(2,[\alpha,i])\}\rightarrow \{(1,[i]),(2,[i])\}$, $(1,[i])\rightarrow \lambda_1\rightarrow (2,[j])$, $(2,[i])\rightarrow \lambda_{\frac{s}{2}+1}\rightarrow (1,[j])$, and $|I^\mathtt{c}((q,[j]))|>0$.
\\
\\(ii) $v=(p,[\alpha,i]), w=(q,[\beta,j])$ for each $p,q=1,2$.
\indent\par From (i), $d_D((p,[\alpha,i]),(3,[j]))=3$. Since $(3,[j])\rightarrow \{(1,[\beta,j]),(2,[\beta,j])\}$ by (\ref{eqC6.3.17}), this subcase follows.
\\
\\(iii) $v=(q,[j]), w=(p,[\alpha,i])$ for each $p=1,2$ and $q=1,2,3$.
\indent\par By (\ref{eqC6.3.18}), $(1,[j])\rightarrow \lambda_1$ and $(2,[j])\rightarrow \lambda_{\frac{s}{2}+1}$. Since $I^\mathtt{c}((3,[i]))\in {{(\mathbb{N}_s,\mathtt{c})}\choose{s/2}}-\{\lambda_1, \lambda_{\frac{s}{2}+1}\}$ by (\ref{eqC6.3.19})-(\ref{eqC6.3.20}), there exist some $x_i\in\lambda_1\cap I^\mathtt{c}((3,[i]))$ and $y_i\in\lambda_{\frac{s}{2}+1}\cap I^\mathtt{c}((3,[i]))$ such that $\{x_i, y_i\}\rightarrow (3,[i])\rightarrow \{(1,[\alpha,i]),(2,[\alpha,i])\}$. Also, since $i\neq j$, $O^\mathtt{c}((3,[i]))$ and $I^\mathtt{c}((3,[j]))$ are not complementary sets. So, there exists some $z_{ij}\in O^\mathtt{c}((3,[i]))\cap I^\mathtt{c}((3,[j]))$ such that $(3,[j])\rightarrow z_{ij}\rightarrow (3,[i])\rightarrow \{(1,[\alpha,i]), (2,[\alpha,i])\}$.
\\
\\Case 3.  For each ${{s}\choose{s/2}}-1 \le i,j \le |A_3|$, $i\neq j$, each $1\le\alpha\le \deg_T([i])-1$, and each $1\le\beta\le \deg_T([j])-1$,
\\(i) $v=(p,[\alpha,i]), w=(q,[j])$ for each $p=1,2$ and $q=1,2,3$.
\indent\par By (\ref{eqC6.3.21})-(\ref{eqC6.3.22}), $\{(1,[\alpha,i]),(2,[\alpha,i])\}\rightarrow (3,[i])$, $\lambda_1 \rightarrow (1,[j]), \lambda_{\frac{s}{2}+1}\rightarrow (2,[j])$. By (\ref{eqC6.3.23}) and $\min\{|\psi_{i+2-{{s}\choose{s/2}}}|+|\lambda_1|, |\psi_{i+2-{{s}\choose{s/2}}}|+|\lambda_{\frac{s}{2}+1}|\}>s$, there exist some $x_i\in \psi_{i+2-{{s}\choose{s/2}}}\cap\lambda_1$ and $y_i\in \psi_{i+2-{{s}\choose{s/2}}}\cap\lambda_{\frac{s}{2}+1}$ such that $(3,[i])\rightarrow \{x_i, y_i\}$. Also, since $i\neq j$, $\psi_{i+2-{{s}\choose{s/2}}}$ and $\bar{\psi}_{j+2-{{s}\choose{s/2}}}$ are not complementary sets. So, there exists some $z_{ij}\in \psi_{i+2-{{s}\choose{s/2}}} \cap \bar{\psi}_{j+2-{{s}\choose{s/2}}}$ such that $(3,[i])\rightarrow z_{ij}\rightarrow (3,[j])$.
\\
\\(ii) $v=(p,[\alpha,i]), w=(q,[\beta,j])$ for each $p,q=1,2$.
\indent\par From (i), $d_D((p,[\alpha,i]),(1,[j]))=3$. Since $(1,[j])\rightarrow \{(1, [\beta,j]), (2, [\beta,j])\}$ by (\ref{eqC6.3.21}), this subcase follows.
\\
\\(iii) $v=(q,[j]), w=(p,[\alpha,i])$ for each $p=1,2$ and $q=1,2,3$.
\indent\par By (\ref{eqC6.3.22}), $(1,[j])\rightarrow \lambda_{\frac{s}{2}+1}\rightarrow (2,[i])$, $(2,[j])\rightarrow \lambda_1 \rightarrow (1,[i])$. Furthermore, $|O^\mathtt{c}((3,[j]))|>0$, there exists some $x_j\in\lambda_{\frac{s}{2}+1}\cup\lambda_1$ such that $(3,[j])\rightarrow x_j$. With $\{(1,[i]),(2,[i])\}\rightarrow \{(1,[\alpha,i]),(2,[\alpha,i])\}$, it follows that $d_{D}(v,w)\le 3$.
\\
\\Case 4. For each $[i], [j]\in A_4$, $i \neq j$, each $1\le\alpha\le \deg_T([i])-1$, and each $1\le\beta\le \deg_T([j])-1$,
\\(i) $v=(p,[\alpha,i]), w=(q,[j])$ for $p=1,2$, and $q=1,2,3,4$.
\indent\par By (\ref{eqC6.3.24})-(\ref{eqC6.3.25}), $(2,[\alpha,i])\rightarrow \{(2,[i]),(4,[i])\}$ and $(1,[\alpha,i])\rightarrow \{(1,[i]), (3,[i])\}$, $\{(1,[i]),$ $(4,[i])\}\rightarrow \lambda_1$, $\{(2,[i]), (3,[i])\}\rightarrow \lambda_{\frac{s}{2}+1}$ and $|I^\mathtt{c}((q,[j]))|>0$.
\\
\\(ii) $v=(p,[i]), w=(q,[\beta,j])$ for for $p=1,2,3,4$, and $q=1,2$.
\indent\par By (\ref{eqC6.3.24})-(\ref{eqC6.3.25}), $\{(2,[j]),(4,[j])\}\rightarrow (1,[\beta,j])$, $\{(1,[j]), (3,[j])\}\rightarrow (2,[\beta,j])$, $\lambda_{\frac{s}{2}+1}\rightarrow \{(1,[j]), (4,[j])\}$, $\lambda_1 \rightarrow\{(2,[j]), (3,[j])\}$ and $|O^\mathtt{c}((p,[i]))|>0$.
\\
\\(iii) $v=(p,[\alpha,i]), w=(q,[\beta,j])$ for $p,q=1,2$.
\indent\par From (i), $d_D((p,[\alpha,i]),(r,[j]))=3$ for $r=1,2,3,4$. So, this subcase holds as $\{(2,[j]),$ $(4,[j])\}\rightarrow (1,[\beta,j])$, and $\{(1,[j]), (3,[j])\}\rightarrow (2,[\beta,j])$ by (\ref{eqC6.3.24}).
\\
\\Case 5. For each $1\le i \le {{s}\choose{s/2}}-2$, each ${{s}\choose{s/2}}-1 \le j \le |A_3| $, each $1\le\alpha\le \deg_T([i])-1$, and each $1\le\beta\le \deg_T([j])-1$,
\\(i) $v=(p,[\alpha,i]), w=(q,[j])$ for each $p=1,2$ and $q=1,2,3$.
\indent\par This follows from $\{(1,[\alpha,i]),(2,[\alpha,i])\}\rightarrow \{(1,[i]),(2,[i])\}$, $(1,[i])\rightarrow \lambda_1$, $(2,[i])\rightarrow \lambda_{\frac{s}{2}+1}$ and $|I^\mathtt{c}((q,[j]))|>0$ by (\ref{eqC6.3.17})-(\ref{eqC6.3.18}) and (\ref{eqC6.3.22})-(\ref{eqC6.3.23}).
\\
\\(ii) $v=(q,[j]), w=(p,[\alpha,i])$ for each $p=1,2$ and $q=1,2,3$.
\indent\par By (\ref{eqC6.3.22}), $(1,[j])\rightarrow \lambda_{\frac{s}{2}+1}$, $(2,[j])\rightarrow \lambda_1$. By (\ref{eqC6.3.19})-(\ref{eqC6.3.20}) and (\ref{eqC6.3.17}), $I^\mathtt{c}((3,[i]))\in {{(\mathbb{N}_s,\mathtt{c})}\choose{s/2}}-\{\lambda_1, \lambda_{\frac{s}{2}+1}\}$, there exist some $x_i\in\lambda_1\cap I^\mathtt{c}((3,[i]))$ and some $y_i\in\lambda_{\frac{s}{2}+1}\cap I^\mathtt{c}((3,[i]))$ such that $\{x_i,y_i\}\rightarrow (3,[i])\rightarrow \{(1,[\alpha,i]),(2,[\alpha,i])\}$. Furthermore with (\ref{eqC6.3.23}), since $|O^\mathtt{c}((3,[j]))|+|I^\mathtt{c}((3,[i]))|>s$, there exists some $z_{ij}\in O^\mathtt{c}((3,[j]))\cap I^\mathtt{c}((3,[i]))$ such that $(3,[j])\rightarrow z_{ij} \rightarrow (3,[i])\rightarrow \{(1,[\alpha,i]),(2,[\alpha,i])\}$.
\\
\\(iii) $v=(p,[i]), w=(q,[\beta,j])$ for each $p=1,2,3$ and $q=1,2$.
\indent\par This follows from $(1,[i])\rightarrow \lambda_1\rightarrow (1,[j])$, $(2,[i])\rightarrow \lambda_{\frac{s}{2}+1}\rightarrow (2,[j])$, $|O^\mathtt{c}((3,[i]))|>0$ and $\{(1,[j]),(2,[j])\}\rightarrow \{(1,[\beta,j]),(2,[\beta,j])\}$ by (\ref{eqC6.3.18})-(\ref{eqC6.3.22}).
\\
\\(iv) $v=(q,[\beta,j]), w=(p,[i])$ for each $p=1,2,3$ and $q=1,2$.
\indent\par Note that $\{(1,[\beta,j]),(2,[\beta,j])\}\rightarrow (3,[j])$ by (\ref{eqC6.3.21}), and $\lambda_{\frac{s}{2}+1}\rightarrow (1,[i]) $, $\lambda_1\rightarrow (2,[i])$ by (\ref{eqC6.3.18}). Since $\min\{|\psi_{j+2-{{s}\choose{s/2}}}|+|\lambda_1|, |\psi_{j+2-{{s}\choose{s/2}}}|+|\lambda_{\frac{s}{2}+1}|\}>s$, there exist some $x_i\in\psi_{j+2-{{s}\choose{s/2}}}\cap\lambda_1$ and some $y_i\in\psi_{j+2-{{s}\choose{s/2}}}\cap\lambda_{\frac{s}{2}+1}$ such that $(3,[j])\rightarrow \{x_i,y_i\}$ by (\ref{eqC6.3.23}). Furthermore by (\ref{eqC6.3.19})-(\ref{eqC6.3.20}), since $|O^\mathtt{c}((3,[j]))|+|I^\mathtt{c}((3,[i]))|>s$, there exists some $z_{ij}\in O^\mathtt{c}((3,[j]))\cap I^\mathtt{c}((3,[i]))$ such that $(3,[j])\rightarrow z_{ij} \rightarrow (3,[i])$.
\\
\\(v) $v=(p,[\alpha,i]), w=(q,[\beta,j])$ for each $p,q=1,2$.
\indent\par We have $\{(1,[\alpha,i]),(2,[\alpha,i])\}\rightarrow (1,[i])\rightarrow \lambda_1\rightarrow (1,[j])\rightarrow \{(1,[\beta,j]), (2, [\beta,j])\}$ by (\ref{eqC6.3.17})-(\ref{eqC6.3.18}) and (\ref{eqC6.3.21})-(\ref{eqC6.3.22}).
\\
\\(vi) $v=(q,[\beta,j]), w=(p,[\alpha,i])$ for each $p,q=1,2$.
\indent\par From (iv), $d_D((q,[\beta,j]),(3,[i]))=3$. Since $(3,[i])\rightarrow \{(1,[\alpha,i]),(2,[\alpha,i])\}$ by (\ref{eqC6.3.17}), this subcase follows.
\\
\\Case 6.  For each $1\le i\le {{s}\choose{s/2}}-2$, each $[j]\in A_4$, each $1\le\alpha\le \deg_T([i])-1$, and each $1\le\beta\le \deg_T([j])-1$,
\\(i) $v=(p,[\alpha,i]), w=(q,[j])$ for each $p=1,2$ and $q=1,2,3,4$.
\indent\par By (\ref{eqC6.3.17})-(\ref{eqC6.3.18}) and (\ref{eqC6.3.25}), $\{(1,[\alpha,i]),(2,[\alpha,i])\}\rightarrow \{(1,[i]),(2,[i])\}$, $(1,[i])\rightarrow \lambda_1\rightarrow \{(2,[j]),(3,[j])\}$ and $(2,[i])\rightarrow \lambda_{\frac{s}{2}+1}\rightarrow \{(1,[j]),(4,[j])\}$.
\\
\\(ii) $v=(q,[j]), w=(p,[\alpha,i])$ for each $p=1,2$ and $q=1,2,3,4$.
\indent\par By (\ref{eqC6.3.25}), $\{(1,[j]), (4,[j])\}\rightarrow \lambda_1$ and $\{(2,[j]), (3,[j])\}\rightarrow \lambda_{\frac{s}{2}+1}$. Since $I^\mathtt{c}((3,[i]))\in {{(\mathbb{N}_s,\mathtt{c})}\choose{s/2}}-\{\lambda_1,\lambda_{\frac{s}{2}+1}\}$, there exist some $x_i\in\lambda_1\cap I^\mathtt{c}((3,[i]))$ and $y_i\in\lambda_{\frac{s}{2}+1}\cap I^\mathtt{c}((3,[i]))$ such that $\{x_i, y_i\}\rightarrow (3,[i])\rightarrow \{(1,[\alpha,i]),(2,[\alpha,i])\}$ by (\ref{eqC6.3.17}) and (\ref{eqC6.3.19})-(\ref{eqC6.3.20}).
\\
\\(iii) $v=(p,[i]), w=(q,[\beta,j])$ for each $p=1,2,3$ and $q=1,2$.
\indent\par By (\ref{eqC6.3.18})-(\ref{eqC6.3.20}) and (\ref{eqC6.3.24})-(\ref{eqC6.3.25}), $|O^\mathtt{c}((p,[i]))|>0$, $\lambda_1\rightarrow \{(2,[j]), (3,[j])\}$, $\lambda_{\frac{s}{2}+1}\rightarrow \{(1,[j]), (4,[j])\}$, $\{(2,[j]),(4,[j])\}\rightarrow (1,[\beta,j])$ and $\{(1,[j]), (3,[j])\}\rightarrow (2,[\beta,j])$.
\\
\\(iv) $v=(q,[\beta,j]), w=(p,[i])$ for each $p=1,2,3$ and $q=1,2$.
\indent\par By (\ref{eqC6.3.18})-(\ref{eqC6.3.20}) and (\ref{eqC6.3.24})-(\ref{eqC6.3.25}), $(1,[\beta,j])\rightarrow \{(1,[j]), (3,[j])\}$, $(2,[\beta, j])\rightarrow \{(2,[j]),(4,[j])\}$, $\{(1,[j]), (4,[j])\}\rightarrow \lambda_1$, $\{(2,[j]), (3,[j])\}\rightarrow \lambda_{\frac{s}{2}+1}$ and $|I^\mathtt{c}((p,[i]))|>0$.
\\
\\(v) $v=(p,[\alpha,i]), w=(q,[\beta,j])$ for each $p,q=1,2$.
\indent\par From (i), $d_D((p,[\alpha,i]),(r,[j]))=3$ for $r=3,4$. Since $(3,[j])\rightarrow (2,[\beta,j])$ and $(4,[j])\rightarrow (1,[\beta,j])$ by (\ref{eqC6.3.24}), this subcase follows.
\\
\\(vi) $v=(q,[\beta,j]), w=(p,[\alpha,i])$ for each $p,q=1,2$.
\indent\par From (iv), $d_D((q,[\beta,j]),(3,[i]))=3$. Since $(3,[i])\rightarrow \{(1,[\alpha,i]), (2,[\alpha,i])\}$ by (\ref{eqC6.3.17}), this subcase follows.
\\
\\Case 7.  For each ${{s}\choose{s/2}}-1\le i\le |A_3|$, each $[j]\in A_4$, each $1\le\alpha\le \deg_T([i])-1$, and each $1\le\beta\le \deg_T([j])-1$,
\\(i) $v=(p,[\alpha,i]), w=(q,[j])$ for each $p=1,2$ and $q=1,2,3,4$.
\indent\par By (\ref{eqC6.3.21}), (\ref{eqC6.3.23}) and (\ref{eqC6.3.25}), $\{(1,[\alpha,i]),(2,[\alpha,i])\}\rightarrow (3,[i])\rightarrow \psi_{j+2-{{s}\choose{s/2}}}$, $\lambda_1\rightarrow \{(2,[j]),$ $(3,[j])\}$, and $\lambda_{\frac{s}{2}+1}\rightarrow \{(1,[j]), (4,[j])\}$. Since $\min\{|\psi_{i+2-{{s}\choose{s/2}}}|+|\lambda_1|, |\psi_{i+2-{{s}\choose{s/2}}}|+|\lambda_{\frac{s}{2}+1}|\}>s$, there exist some $x_i\in\psi_{i+2-{{s}\choose{s/2}}}\cap\lambda_1$ and $y_i\in\psi_{i+2-{{s}\choose{s/2}}}\cap\lambda_{\frac{s}{2}+1}$ such that $(3,[i])\rightarrow \{x_i, y_i\}$.
\\
\\(ii) $v=(q,[j]), w=(p,[\alpha,i])$ for each $p=1,2$ and $q=1,2,3,4$.
\indent\par By (\ref{eqC6.3.21})-(\ref{eqC6.3.22}) and (\ref{eqC6.3.25}), $\{(1,[j]), (4,[j])\}\rightarrow \lambda_1\rightarrow (1,[i])$, $\{(2,[j]), (3,[j])\}\rightarrow \lambda_{\frac{s}{2}+1}\rightarrow (2,[i])$, and $\{(1,[i]), (2,[i])\}\rightarrow \{(1,[\alpha,i]),(2,[\alpha,i])\}$.
\\
\\(iii) $v=(p,[i]), w=(q,[\beta,j])$ for each $p=1,2,3$ and $q=1,2$.
\indent\par By (\ref{eqC6.3.22}) and (\ref{eqC6.3.24})-(\ref{eqC6.3.25}), $(2,[i])\rightarrow \lambda_1\rightarrow \{(2,[j]),(3,[j])\}$, $(1,[i])\rightarrow \lambda_{\frac{s}{2}+1}\rightarrow \{(1,[j]), (4,[j])\}$, $\{(2,[j]),(4,[j])\}\rightarrow (1,[\beta,j])$ and $\{(1,[j]), (3,[j])\}\rightarrow (2,[\beta,j])$. Furthermore by (\ref{eqC6.3.23}), since $\min\{|\psi_{i+2-{{s}\choose{s/2}}}|+|\lambda_{1}|, |\psi_{i+2-{{s}\choose{s/2}}}|+|\lambda_{\frac{s}{2}+1}|\}>s$, there exist some $x_i\in\psi_{i+2-{{s}\choose{s/2}}}\cap\lambda_1$ and $y_i\in\psi_{i+2-{{s}\choose{s/2}}}\cap\lambda_{\frac{s}{2}+1}$ such that $(3,[i])\rightarrow \{x_i, y_i\}$.
\\
\\(iv) $v=(q,[\beta,j]), w=(p,[i])$ for each $p=1,2,3$ and $q=1,2$.
\indent\par  By (\ref{eqC6.3.22})-(\ref{eqC6.3.25}), $(1,[\beta,j])\rightarrow \{(1,[j]), (3,[j])\}$, $(2,[\beta,j])\rightarrow \{(2,[j]),(4,[j])\}$, $\{(1,[j]),$ $(4,[j])\}$ $\rightarrow \lambda_1$, $\{(2,[j]), (3,[j])\}\rightarrow \lambda_{\frac{s}{2}+1}$, and $|I^\mathtt{c}((p,[i]))|>0$.
\\
\\(v) $v=(p,[\alpha,i]), w=(q,[\beta,j])$ for each $p,q=1,2$.
\indent\par From (i), $d_D((p,[\alpha,i]),(r,[j]))=3$ for $r=3,4$. Since $(3,[j])\rightarrow (2,[\beta,j])$ and $(4,[j])\rightarrow (1,[\beta,j])$ by (\ref{eqC6.3.24}), this subcase follows.
\\
\\(vi) $v=(q,[\beta,j]), w=(p,[\alpha,i])$ for each $p,q=1,2$.
\indent\par From (iv), $d_D((q,[\beta,j]),(1,[i]))=3$ for $p=1,2$. Since $(1,[i])\rightarrow \{(1,[\alpha,i]), (2,[\alpha,i])\}$ by (\ref{eqC6.3.21}), this subcase follows.
\\
\\Case 8. For each $[i]\in A_3$, each $1\le\alpha\le \deg_T([i])-1$, and each $[j]\in E$,
\\(i) $v=(p,[\alpha,i]), w=(q,[j])$ for each $p=1,2$ and $q=1,2$.
\\(ii) $v=(q,[j]), w=(p,[\alpha,i])$ for each $p=1,2$ and $q=1,2$.
\indent\par Let $[k]\in A_4$. Since $\lambda_1\rightarrow \{(q,[j]),(2,[k])\}\rightarrow \bar{\lambda}_1$ by (\ref{eqC6.3.25})-(\ref{eqC6.3.26}), this case follows from Cases 6(i)-(ii) and Cases 7(i)-(ii).
\\
\\Case 9. For each $[i]\in A_4$, each $1\le\alpha\le \deg_T([i])-1$, and each $[j]\in E$,
\\(i) $v=(p,[\alpha,i]), w=(q,[j])$ for each $p=1,2$ and $q=1,2$.
\\(ii) $v=(q,[j]), w=(p,[\alpha,i])$ for each $p=1,2$ and $q=1,2$.
\indent\par Since $\lambda_1\rightarrow \{(q,[j]),(1,[|A_3|])\}\rightarrow \lambda_{\frac{s}{2}+1}$ by (\ref{eqC6.3.22}) and (\ref{eqC6.3.26}), this case follows from Cases 7(i)-(ii).
\\
\\Case 10. $v=(r_1,\mathtt{c})$ and $w=(r_2,\mathtt{c})$ for $r_1\neq r_2$ and $1\le r_1, r_2\le s$.
\indent\par Here, we want to prove a stronger claim, $d_{D}((r_1,\mathtt{c}), (r_2,\mathtt{c}))=2$. Let $x_1=(2,[1])$, $x_{k+1}=(3,[k])$ for $1\le k\le \frac{s}{2}-1$, $x_{\frac{s}{2}+1}=(1,[1])$, and $x_{k+2}=(3,[k])$ for $ \frac{s}{2}\le k\le s$. Observe from (\ref{eqC6.3.18})-(\ref{eqC6.3.20}) that $\lambda_k\rightarrow x_k\rightarrow \bar{\lambda}_k$ for $1\le k\le s$ and the subgraph induced by $V_1=(\mathbb{N}_s,\mathtt{c})$ and $V_2=\{x_k\mid 1\le k \le s\}$ is a complete bipartite graph $K(V_1,V_2)$. By Lemma \ref{lemC6.2.19}, $d_{D}((r_1,\mathtt{c}), (r_2,\mathtt{c}))=2$.
\\
\\Case 11. $v\in \{(1,[i]), (2,[i]), (3,[i]), (4,[i]), (1,[\alpha,i]), (2,[\alpha,i])\}$ for each $1\le i\le \deg_T(\mathtt{c})$ and $1\le\alpha\le \deg_T([i])-1$, and $w=(r,\mathtt{c})$ for $1\le r\le s$.
\indent\par Note that there exists some $1\le k\le s$ such that $d_{D}(v,(k,\mathtt{c}))\le 2$, and $d_{D}((k,\mathtt{c}),w)\le 2$ by Case 10. Hence, it follows that $d_{D}(v,w)\le d_{D}(v,(k,\mathtt{c}))+d_D((k,\mathtt{c}),w)\le 4$.
\\
\\Case 12. $v=(r,\mathtt{c})$ for $1\le r\le s$ and $w\in \{(1,[i]), (2,[i]), (3,[i]), (4,[i]), (1,[\alpha,i]), (2,[\alpha,i])\}$ for each $1\le i\le \deg_T(\mathtt{c})$ and $1\le\alpha\le \deg_T([i])-1$.
\indent\par Note that there exists some $1\le k\le s$ such that $d_{D}((k,\mathtt{c}), w)\le 2$, and $d_{D}(v,(k,\mathtt{c}))\le 2$ by Case 10. Hence, it follows that $d_{D}(v,w)\le d_{D}(v,(k,\mathtt{c}))+d_{D}((k,\mathtt{c}),w)\le 4$.
\\
\\Case 13. $v=(p,[i])$ and $w=(q, [j])$, where $1\le p,q\le 3$ and $1\le i,j\le \deg_T(\mathtt{c})$.
\noindent\par This follows from the fact that $|O^\mathtt{c}((p,[i]))|>0$, $|I^\mathtt{c}((q,[j]))|>0$, and $d_{D}((r_1,\mathtt{c}), (r_2,\mathtt{c}))$ $=2$ for any $r_1\neq r_2$ and $1\le r_1, r_2\le s$.
\\
\indent\par Therefore, the claim follows. Since every vertex lies in a directed $C_4$ for $D$ and $d(D)=4$, $\bar{d}(\mathcal{T})\le \max \{4, d(D)\}$ by Lemma \ref{lemC6.1.3}, and thus $\bar{d}(\mathcal{T})=4$ .
\qed

\begin{center}
\begin{tikzpicture}[thick,scale=0.7]%
\draw(0,4.5)node[circle, draw, inner sep=0pt, minimum width=3pt](1_u){\scriptsize $(1,\mathtt{c})$};
\draw(0,1.5)node[circle, draw, inner sep=0pt, minimum width=3pt](2_u){\scriptsize $(2,\mathtt{c})$};
\draw(0,-1.5)node[circle, draw, inner sep=0pt, minimum width=3pt](3_u){\scriptsize $(3,\mathtt{c})$};
\draw(0,-4.5)node[circle, draw, inner sep=0pt, minimum width=3pt](4_u){\scriptsize $(4,\mathtt{c})$};

\draw(-6,4)node[circle, draw, fill=black!100, inner sep=0pt, minimum width=6pt, label={[] 180:{\small $(1,[1,1])$}}](1_11){};
\draw(-6,2)node[circle, draw, fill=black!100, inner sep=0pt, minimum width=6pt, label={[] 180:{\small $(2,[1,1])$}}](2_11){};

\draw(-3,5)node[circle, draw, fill=black!100, inner sep=0pt, minimum width=6pt, label={[yshift=0cm, xshift=-0.45cm] 90:{\small $(1,[1])$}}](1_1){};
\draw(-3,3)node[circle, draw, fill=black!100, inner sep=0pt, minimum width=6pt, label={[yshift=-0.1cm, xshift=-0.45cm] 270:{\small $(2,[1])$}}](2_1){};
\draw(-3,1)node[circle, draw, fill=black!100, inner sep=0pt, minimum width=6pt, label={[yshift=0cm, xshift=-0.45cm] 270:{\small $(3,[1])$}}](3_1){};

\draw(6,4)node[circle, draw, fill=black!100, inner sep=0pt, minimum width=6pt, label={[] 0:{\small $(1,[1,2])$}}](1_21){};
\draw(6,2)node[circle, draw, fill=black!100, inner sep=0pt, minimum width=6pt, label={[] 0:{\small $(2,[1,2])$}}](2_21){};

\draw(3,5)node[circle, draw, fill=black!100, inner sep=0pt, minimum width=6pt, label={[yshift=0cm, xshift=0.45cm] 90:{\small $(1,[2])$}}](1_2){};
\draw(3,3)node[circle, draw, fill=black!100, inner sep=0pt, minimum width=6pt, label={[yshift=-0.1cm, xshift=0.45cm] 270:{\small $(2,[2])$}}](2_2){};
\draw(3,1)node[circle, draw, fill=black!100, inner sep=0pt, minimum width=6pt, label={[yshift=0cm, xshift=0.45cm] 270:{\small $(3,[2])$}}](3_2){};

\draw(-6,-2)node[circle, draw, fill=black!100, inner sep=0pt, minimum width=6pt, label={[] 180:{\small $(1,[1,3])$}}](1_31){};
\draw(-6,-4)node[circle, draw, fill=black!100, inner sep=0pt, minimum width=6pt, label={[] 180:{\small $(2,[1,3])$}}](2_31){};

\draw(-3,-1)node[circle, draw, fill=black!100, inner sep=0pt, minimum width=6pt, label={[yshift=-0.15cm, xshift=-0.45cm] 270:{\small $(1,[3])$}}](1_3){};
\draw(-3,-3)node[circle, draw, fill=black!100, inner sep=0pt, minimum width=6pt, label={[yshift=-0.1cm, xshift=-0.45cm] 270:{\small $(2,[3])$}}](2_3){};
\draw(-3,-5)node[circle, draw, fill=black!100, inner sep=0pt, minimum width=6pt, label={[yshift=0cm, xshift=-0.45cm] 270:{\small $(3,[3])$}}](3_3){};

\draw(6,-2)node[circle, draw, fill=black!100, inner sep=0pt, minimum width=6pt, label={[] 0:{\small $(1,[1,4])$}}](1_41){};
\draw(6,-4)node[circle, draw, fill=black!100, inner sep=0pt, minimum width=6pt, label={[] 0:{\small $(2,[1,4])$}}](2_41){};

\draw(3,-1)node[circle, draw, fill=black!100, inner sep=0pt, minimum width=6pt, label={[yshift=-0.15cm, xshift=0.45cm] 270:{\small $(1,[4])$}}](1_4){};
\draw(3,-3)node[circle, draw, fill=black!100, inner sep=0pt, minimum width=6pt, label={[yshift=-0.1cm, xshift=0.45cm] 270:{\small $(2,[4])$}}](2_4){};
\draw(3,-5)node[circle, draw, fill=black!100, inner sep=0pt, minimum width=6pt, label={[yshift=0cm, xshift=0.45cm] 270:{\small $(3,[4])$}}](3_4){};

\draw[dashed, ->, line width=0.3mm, >=latex, shorten <= 0.2cm, shorten >= 0.15cm](3_1) to [out=158, in=275] (1_11);
\draw[dashed, ->, line width=0.3mm, >=latex, shorten <= 0.2cm, shorten >= 0.15cm](3_1)--(2_11);
\draw[->, line width=0.3mm, >=latex, shorten <= 0.2cm, shorten >= 0.15cm](1_11)--(1_1);
\draw[->, line width=0.3mm, >=latex, shorten <= 0.2cm, shorten >= 0.15cm](2_11)--(1_1);
\draw[->, line width=0.3mm, >=latex, shorten <= 0.2cm, shorten >= 0.15cm](1_11)--(2_1);
\draw[->, line width=0.3mm, >=latex, shorten <= 0.2cm, shorten >= 0.15cm](2_11)--(2_1);

\draw[dashed, ->, line width=0.3mm, >=latex, shorten <= 0.2cm, shorten >= 0.15cm](3_2) to [out=22, in=265] (1_21);
\draw[dashed, ->, line width=0.3mm, >=latex, shorten <= 0.2cm, shorten >= 0.15cm](3_2)--(2_21);
\draw[->, line width=0.3mm, >=latex, shorten <= 0.2cm, shorten >= 0.15cm](1_21)--(1_2);
\draw[->, line width=0.3mm, >=latex, shorten <= 0.2cm, shorten >= 0.15cm](2_21)--(1_2);
\draw[->, line width=0.3mm, >=latex, shorten <= 0.2cm, shorten >= 0.15cm](1_21)--(2_2);
\draw[->, line width=0.3mm, >=latex, shorten <= 0.2cm, shorten >= 0.15cm](2_21)--(2_2);

\draw[dashed, ->, line width=0.3mm, >=latex, shorten <= 0.2cm, shorten >= 0.15cm](3_3) to [out=158, in=275] (1_31);
\draw[dashed, ->, line width=0.3mm, >=latex, shorten <= 0.2cm, shorten >= 0.15cm](3_3)--(2_31);
\draw[->, line width=0.3mm, >=latex, shorten <= 0.2cm, shorten >= 0.15cm](1_31)--(1_3);
\draw[->, line width=0.3mm, >=latex, shorten <= 0.2cm, shorten >= 0.15cm](2_31) to [out=85, in=202] (1_3);
\draw[->, line width=0.3mm, >=latex, shorten <= 0.2cm, shorten >= 0.15cm](1_31)--(2_3);
\draw[->, line width=0.3mm, >=latex, shorten <= 0.2cm, shorten >= 0.15cm](2_31)--(2_3);

\draw[dashed, ->, line width=0.3mm, >=latex, shorten <= 0.2cm, shorten >= 0.15cm](3_4) to [out=22, in=265] (1_41);
\draw[dashed, ->, line width=0.3mm, >=latex, shorten <= 0.2cm, shorten >= 0.15cm](3_4)--(2_41);
\draw[->, line width=0.3mm, >=latex, shorten <= 0.2cm, shorten >= 0.15cm](1_41)--(1_4);
\draw[->, line width=0.3mm, >=latex, shorten <= 0.2cm, shorten >= 0.15cm](2_41) to [out=95, in=338] (1_4);
\draw[->, line width=0.3mm, >=latex, shorten <= 0.2cm, shorten >= 0.15cm](1_41)--(2_4);
\draw[->, line width=0.3mm, >=latex, shorten <= 0.2cm, shorten >= 0.15cm](2_41)--(2_4);

\draw[->, line width=0.3mm, >=latex, shorten <= 0.2cm, shorten >= 0.1cm](1_1)--(1_u);
\draw[->, line width=0.3mm, >=latex, shorten <= 0.2cm, shorten >= 0.1cm](1_1)--(2_u);
\draw[->, line width=0.3mm, >=latex, shorten <= 0.2cm, shorten >= 0.1cm](2_1)--(3_u);
\draw[->, line width=0.3mm, >=latex, shorten <= 0.2cm, shorten >= 0.1cm](2_1)--(4_u);
\draw[densely dashdotdotted, ->, line width=0.3mm, >=latex, shorten <= 0.2cm, shorten >= 0.1cm](3_1)--(1_u);
\draw[densely dashdotdotted, ->, line width=0.3mm, >=latex, shorten <= 0.2cm, shorten >= 0.1cm](3_1)--(4_u);

\draw[->, line width=0.3mm, >=latex, shorten <= 0.2cm, shorten >= 0.1cm](1_2)--(1_u);
\draw[->, line width=0.3mm, >=latex, shorten <= 0.2cm, shorten >= 0.1cm](1_2)--(2_u);
\draw[->, line width=0.3mm, >=latex, shorten <= 0.2cm, shorten >= 0.1cm](2_2)--(3_u);
\draw[->, line width=0.3mm, >=latex, shorten <= 0.2cm, shorten >= 0.1cm](2_2)--(4_u);
\draw[densely dashdotdotted, ->, line width=0.3mm, >=latex, shorten <= 0.2cm, shorten >= 0.1cm](3_2)--(2_u);
\draw[densely dashdotdotted, ->, line width=0.3mm, >=latex, shorten <= 0.2cm, shorten >= 0.1cm](3_2)--(3_u);

\draw[->, line width=0.3mm, >=latex, shorten <= 0.2cm, shorten >= 0.1cm](1_3)--(1_u);
\draw[->, line width=0.3mm, >=latex, shorten <= 0.2cm, shorten >= 0.1cm](1_3)--(2_u);
\draw[->, line width=0.3mm, >=latex, shorten <= 0.2cm, shorten >= 0.1cm](2_3)--(3_u);
\draw[->, line width=0.3mm, >=latex, shorten <= 0.2cm, shorten >= 0.1cm](2_3)--(4_u);
\draw[densely dashdotdotted, ->, line width=0.3mm, >=latex, shorten <= 0.2cm, shorten >= 0.1cm](3_3)-- (2_u);
\draw[densely dashdotdotted, ->, line width=0.3mm, >=latex, shorten <= 0.2cm, shorten >= 0.1cm](3_3)--(4_u);

\draw[->, line width=0.3mm, >=latex, shorten <= 0.2cm, shorten >= 0.1cm](1_4)--(1_u);
\draw[->, line width=0.3mm, >=latex, shorten <= 0.2cm, shorten >= 0.1cm](1_4)--(2_u);
\draw[->, line width=0.3mm, >=latex, shorten <= 0.2cm, shorten >= 0.1cm](2_4)--(3_u);
\draw[->, line width=0.3mm, >=latex, shorten <= 0.2cm, shorten >= 0.1cm](2_4)--(4_u);
\draw[densely dashdotdotted, ->, line width=0.3mm, >=latex, shorten <= 0.2cm, shorten >= 0.1cm](3_4)--(1_u);
\draw[densely dashdotdotted, ->, line width=0.3mm, >=latex, shorten <= 0.2cm, shorten >= 0.1cm](3_4)--(3_u);
\end{tikzpicture}
\captionsetup{justification=centering}
{
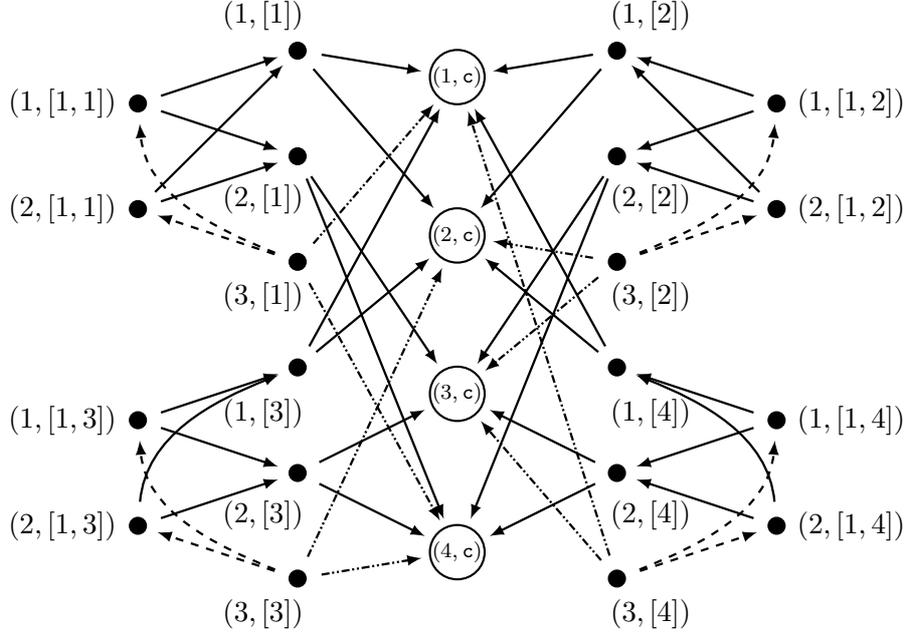
\captionof{figure}{Partial orientation $D$ for $\mathcal{H}$ for $s=4$,\\$A_3=\{[1],[2],\ldots,[6]\}$, $A_4=\{[7],[8]\}$, $E=\{[9],[10]\}$; showing $[i]$ for $1\le i\le {{s}\choose{s/2}}-2$.}\label{figC6.3.7}}
\end{center}

\begin{center}
\begin{tikzpicture}[thick,scale=0.7]%
\draw(0,3.5)node[circle, draw, inner sep=0pt, minimum width=3pt](1_u){\scriptsize $(1,\mathtt{c})$};
\draw(0,-1.5)node[circle, draw, inner sep=0pt, minimum width=3pt](2_u){\scriptsize $(2,\mathtt{c})$};
\draw(0,-5.5)node[circle, draw, inner sep=0pt, minimum width=3pt](3_u){\scriptsize $(3,\mathtt{c})$};
\draw(0,-10.5)node[circle, draw, inner sep=0pt, minimum width=3pt](4_u){\scriptsize $(4,\mathtt{c})$};

\draw(-6,3)node[circle, draw, fill=black!100, inner sep=0pt, minimum width=6pt, label={[] 180:{\small $(1,[1,5])$}}](1_51){};
\draw(-6,1)node[circle, draw, fill=black!100, inner sep=0pt, minimum width=6pt, label={[] 180:{\small $(2,[1,5])$}}](2_51){};

\draw(-3,4)node[circle, draw, fill=black!100, inner sep=0pt, minimum width=6pt, label={[yshift=0cm, xshift=-0.45cm] 90:{\small $(1,[5])$}}](1_5){};
\draw(-3,2)node[circle, draw, fill=black!100, inner sep=0pt, minimum width=6pt, label={[yshift=-0.1cm, xshift=-0.45cm] 270:{\small $(2,[5])$}}](2_5){};
\draw(-3,0)node[circle, draw, fill=black!100, inner sep=0pt, minimum width=6pt, label={[yshift=0cm, xshift=-0.45cm] 270:{\small $(3,[5])$}}](3_5){};

\draw(6,3)node[circle, draw, fill=black!100, inner sep=0pt, minimum width=6pt, label={[] 0:{\small $(1,[1,6])$}}](1_61){};
\draw(6,1)node[circle, draw, fill=black!100, inner sep=0pt, minimum width=6pt, label={[] 0:{\small $(2,[1,6])$}}](2_61){};

\draw(3,4)node[circle, draw, fill=black!100, inner sep=0pt, minimum width=6pt, label={[yshift=0cm, xshift=0.45cm] 90:{\small $(1,[6])$}}](1_6){};
\draw(3,2)node[circle, draw, fill=black!100, inner sep=0pt, minimum width=6pt, label={[yshift=-0.1cm, xshift=0.45cm] 270:{\small $(2,[6])$}}](2_6){};
\draw(3,0)node[circle, draw, fill=black!100, inner sep=0pt, minimum width=6pt, label={[yshift=0cm, xshift=0.45cm] 270:{\small $(3,[6])$}}](3_6){};

\draw(-6,-4)node[circle, draw, fill=black!100, inner sep=0pt, minimum width=6pt, label={[] 180:{\small $(1,[1,7])$}}](1_71){};
\draw(-6,-6)node[circle, draw, fill=black!100, inner sep=0pt, minimum width=6pt, label={[] 180:{\small $(2,[1,7])$}}](2_71){};

\draw(-3,-2)node[circle, draw, fill=black!100, inner sep=0pt, minimum width=6pt, label={[yshift=-0.1cm, xshift=-0.45cm] 270:{\small $(1,[7])$}}](1_7){};
\draw(-3,-4)node[circle, draw, fill=black!100, inner sep=0pt, minimum width=6pt, label={[yshift=0cm, xshift=-0.45cm] 270:{\small $(2,[7])$}}](2_7){};
\draw(-3,-6)node[circle, draw, fill=black!100, inner sep=0pt, minimum width=6pt, label={[yshift=0.05cm, xshift=-0.45cm] 270:{\small $(3,[7])$}}](3_7){};
\draw(-3,-8)node[circle, draw, fill=black!100, inner sep=0pt, minimum width=6pt, label={[yshift=0cm, xshift=-0.45cm] 270:{\small $(4,[7])$}}](4_7){};

\draw(6,-4)node[circle, draw, fill=black!100, inner sep=0pt, minimum width=6pt, label={[] 0:{\small $(1,[1,8])$}}](1_81){};
\draw(6,-6)node[circle, draw, fill=black!100, inner sep=0pt, minimum width=6pt, label={[] 0:{\small $(2,[1,8])$}}](2_81){};

\draw(3,-2)node[circle, draw, fill=black!100, inner sep=0pt, minimum width=6pt, label={[yshift=-0.1cm, xshift=0.45cm] 270:{\small $(1,[8])$}}](1_8){};
\draw(3,-4)node[circle, draw, fill=black!100, inner sep=0pt, minimum width=6pt, label={[yshift=0cm, xshift=0.45cm] 270:{\small $(2,[8])$}}](2_8){};
\draw(3,-6)node[circle, draw, fill=black!100, inner sep=0pt, minimum width=6pt, label={[yshift=0.05cm, xshift=0.45cm] 270:{\small $(3,[8])$}}](3_8){};
\draw(3,-8)node[circle, draw, fill=black!100, inner sep=0pt, minimum width=6pt, label={[yshift=0cm, xshift=0.45cm] 270:{\small $(4,[8])$}}](4_8){};

\draw(-3,-10)node[circle, draw, fill=black!100, inner sep=0pt, minimum width=6pt, label={[] 180:{\small $(1,[9])$}}](1_9){};
\draw(-3,-12)node[circle, draw, fill=black!100, inner sep=0pt, minimum width=6pt, label={[] 180:{\small $(2,[9])$}}](2_9){};

\draw(3,-10)node[circle, draw, fill=black!100, inner sep=0pt, minimum width=6pt, label={[] 0:{\small $(1,[10])$}}](1_10){};
\draw(3,-12)node[circle, draw, fill=black!100, inner sep=0pt, minimum width=6pt, label={[] 0:{\small $(2,[10])$}}](2_10){};

\draw[dashed, ->, line width=0.3mm, >=latex, shorten <= 0.2cm, shorten >= 0.15cm](1_5)--(1_51);
\draw[dashed, ->, line width=0.3mm, >=latex, shorten <= 0.2cm, shorten >= 0.15cm](1_5)--(2_51);
\draw[dashed, ->, line width=0.3mm, >=latex, shorten <= 0.2cm, shorten >= 0.15cm](2_5)--(1_51);
\draw[dashed, ->, line width=0.3mm, >=latex, shorten <= 0.2cm, shorten >= 0.15cm](2_5)--(2_51);
\draw[->, line width=0.3mm, >=latex, shorten <= 0.2cm, shorten >= 0.15cm](1_51) to [out=275, in=158] (3_5);
\draw[->, line width=0.3mm, >=latex, shorten <= 0.2cm, shorten >= 0.15cm](2_51)--(3_5);

\draw[dashed, ->, line width=0.3mm, >=latex, shorten <= 0.2cm, shorten >= 0.15cm](1_6)--(1_61);
\draw[dashed, ->, line width=0.3mm, >=latex, shorten <= 0.2cm, shorten >= 0.15cm](1_6)--(2_61);
\draw[dashed, ->, line width=0.3mm, >=latex, shorten <= 0.2cm, shorten >= 0.15cm](2_6)--(1_61);
\draw[dashed, ->, line width=0.3mm, >=latex, shorten <= 0.2cm, shorten >= 0.15cm](2_6)--(2_61);
\draw[->, line width=0.3mm, >=latex, shorten <= 0.2cm, shorten >= 0.15cm](1_61) to [out=265, in=22] (3_6);
\draw[->, line width=0.3mm, >=latex, shorten <= 0.2cm, shorten >= 0.15cm](2_61)--(3_6);

\draw[dashed, ->, line width=0.3mm, >=latex, shorten <= 0.2cm, shorten >= 0.15cm](2_7)--(1_71);
\draw[dashed, ->, line width=0.3mm, >=latex, shorten <= 0.2cm, shorten >= 0.15cm](1_7) to [out=195, in=80] (2_71);
\draw[dashed, ->, line width=0.3mm, >=latex, shorten <= 0.2cm, shorten >= 0.15cm](4_7) to [out=140, in=290] (1_71);
\draw[dashed, ->, line width=0.3mm, >=latex, shorten <= 0.2cm, shorten >= 0.15cm](3_7)--(2_71);
\draw[->, line width=0.3mm, >=latex, shorten <= 0.2cm, shorten >= 0.15cm](1_71) to [out=65, in=180] (1_7);
\draw[->, line width=0.3mm, >=latex, shorten <= 0.2cm, shorten >= 0.15cm](2_71) to [out=65, in=190] (2_7);
\draw[->, line width=0.3mm, >=latex, shorten <= 0.2cm, shorten >= 0.15cm](1_71)--(3_7);
\draw[->, line width=0.3mm, >=latex, shorten <= 0.2cm, shorten >= 0.15cm](2_71)--(4_7);

\draw[dashed, ->, line width=0.3mm, >=latex, shorten <= 0.2cm, shorten >= 0.15cm](2_8)--(1_81);
\draw[dashed, ->, line width=0.3mm, >=latex, shorten <= 0.2cm, shorten >= 0.15cm](1_8) to [out=345, in=100] (2_81);
\draw[dashed, ->, line width=0.3mm, >=latex, shorten <= 0.2cm, shorten >= 0.15cm](4_8) to [out=40, in=250] (1_81);
\draw[dashed, ->, line width=0.3mm, >=latex, shorten <= 0.2cm, shorten >= 0.15cm](3_8)--(2_81);
\draw[->, line width=0.3mm, >=latex, shorten <= 0.2cm, shorten >= 0.15cm](1_81) to [out=115, in=0] (1_8);
\draw[->, line width=0.3mm, >=latex, shorten <= 0.2cm, shorten >= 0.15cm](2_81) to [out=115, in=350] (2_8);
\draw[->, line width=0.3mm, >=latex, shorten <= 0.2cm, shorten >= 0.15cm](1_81)--(3_8);
\draw[->, line width=0.3mm, >=latex, shorten <= 0.2cm, shorten >= 0.15cm](2_81)--(4_8);

\draw[->, line width=0.3mm, >=latex, shorten <= 0.2cm, shorten >= 0.1cm](1_5)--(3_u);
\draw[->, line width=0.3mm, >=latex, shorten <= 0.2cm, shorten >= 0.1cm](1_5)--(4_u);
\draw[->, line width=0.3mm, >=latex, shorten <= 0.2cm, shorten >= 0.1cm](2_5)--(1_u);
\draw[->, line width=0.3mm, >=latex, shorten <= 0.2cm, shorten >= 0.1cm](2_5)--(2_u);
\draw[densely dotted, ->, line width=0.3mm, >=latex, shorten <= 0.2cm, shorten >= 0.1cm](3_5)--(1_u);
\draw[densely dotted, ->, line width=0.3mm, >=latex, shorten <= 0.2cm, shorten >= 0.1cm](3_5)--(2_u);
\draw[densely dotted, ->, line width=0.3mm, >=latex, shorten <= 0.2cm, shorten >= 0.1cm](3_5)--(3_u);

\draw[->, line width=0.3mm, >=latex, shorten <= 0.2cm, shorten >= 0.1cm](1_6)--(3_u);
\draw[->, line width=0.3mm, >=latex, shorten <= 0.2cm, shorten >= 0.1cm](1_6)--(4_u);
\draw[->, line width=0.3mm, >=latex, shorten <= 0.2cm, shorten >= 0.1cm](2_6)--(1_u);
\draw[->, line width=0.3mm, >=latex, shorten <= 0.2cm, shorten >= 0.1cm](2_6)--(2_u);
\draw[densely dotted, ->, line width=0.3mm, >=latex, shorten <= 0.2cm, shorten >= 0.1cm](3_6)--(2_u);
\draw[densely dotted, ->, line width=0.3mm, >=latex, shorten <= 0.2cm, shorten >= 0.1cm](3_6)--(3_u);
\draw[densely dotted, ->, line width=0.3mm, >=latex, shorten <= 0.2cm, shorten >= 0.1cm](3_6)--(4_u);

\draw[->, line width=0.3mm, >=latex, shorten <= 0.2cm, shorten >= 0.1cm](1_7)--(1_u);
\draw[->, line width=0.3mm, >=latex, shorten <= 0.2cm, shorten >= 0.1cm](4_7)--(1_u);
\draw[->, line width=0.3mm, >=latex, shorten <= 0.2cm, shorten >= 0.1cm](1_7)--(2_u);
\draw[->, line width=0.3mm, >=latex, shorten <= 0.2cm, shorten >= 0.1cm](4_7)--(2_u);

\draw[->, line width=0.3mm, >=latex, shorten <= 0.2cm, shorten >= 0.1cm](1_8)--(1_u);
\draw[->, line width=0.3mm, >=latex, shorten <= 0.2cm, shorten >= 0.1cm](4_8)--(1_u);
\draw[->, line width=0.3mm, >=latex, shorten <= 0.2cm, shorten >= 0.1cm](1_8)--(2_u);
\draw[->, line width=0.3mm, >=latex, shorten <= 0.2cm, shorten >= 0.1cm](4_8)--(2_u);

\draw[->, line width=0.3mm, >=latex, shorten <= 0.2cm, shorten >= 0.1cm](2_7)--(3_u);
\draw[->, line width=0.3mm, >=latex, shorten <= 0.2cm, shorten >= 0.1cm](3_7)--(3_u);
\draw[->, line width=0.3mm, >=latex, shorten <= 0.2cm, shorten >= 0.1cm](2_7)--(4_u);
\draw[->, line width=0.3mm, >=latex, shorten <= 0.2cm, shorten >= 0.1cm](3_7)--(4_u);

\draw[->, line width=0.3mm, >=latex, shorten <= 0.2cm, shorten >= 0.1cm](2_8)--(3_u);
\draw[->, line width=0.3mm, >=latex, shorten <= 0.2cm, shorten >= 0.1cm](3_8)--(3_u);
\draw[->, line width=0.3mm, >=latex, shorten <= 0.2cm, shorten >= 0.1cm](2_8)--(4_u);
\draw[->, line width=0.3mm, >=latex, shorten <= 0.2cm, shorten >= 0.1cm](3_8)--(4_u);

\draw[->, line width=0.3mm, >=latex, shorten <= 0.2cm, shorten >= 0.1cm](1_9)--(3_u);
\draw[->, line width=0.3mm, >=latex, shorten <= 0.2cm, shorten >= 0.1cm](2_9)--(3_u);
\draw[->, line width=0.3mm, >=latex, shorten <= 0.2cm, shorten >= 0.1cm](1_9)--(4_u);
\draw[->, line width=0.3mm, >=latex, shorten <= 0.2cm, shorten >= 0.1cm](2_9)--(4_u);

\draw[->, line width=0.3mm, >=latex, shorten <= 0.2cm, shorten >= 0.1cm](1_10)--(3_u);
\draw[->, line width=0.3mm, >=latex, shorten <= 0.2cm, shorten >= 0.1cm](2_10)--(3_u);
\draw[->, line width=0.3mm, >=latex, shorten <= 0.2cm, shorten >= 0.1cm](1_10)--(4_u);
\draw[->, line width=0.3mm, >=latex, shorten <= 0.2cm, shorten >= 0.1cm](2_10)--(4_u);
\end{tikzpicture}
\captionsetup{justification=centering}
{
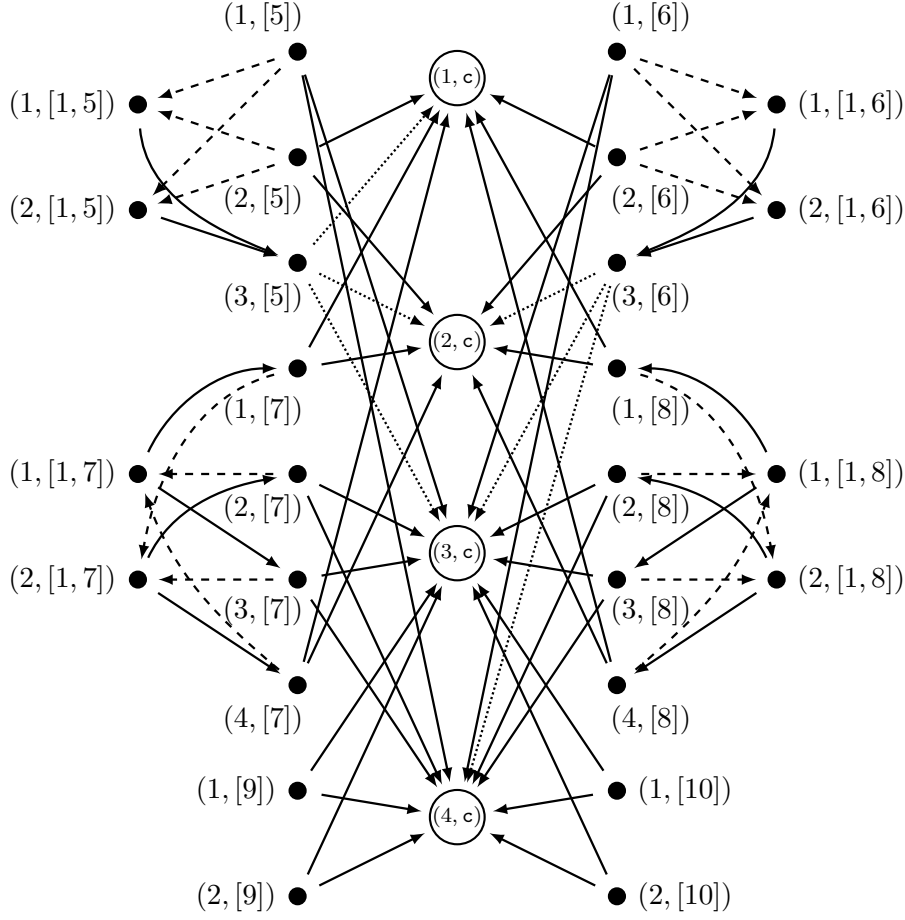
\captionof{figure}{Partial orientation $D$ for $\mathcal{H}$ for $s=4$;\\showing $[i]$ for ${{s}\choose{s/2}}-1\le i\le |A_3|$ or $[i]\in A_4$ or $[i]\in E$.}\label{figC6.3.8}}
\end{center}

\begin{ppn}\label{ppnC6.3.10}
Suppose $s\ge 4$ is even, $A_2\neq \emptyset$, $A_3=\emptyset$ and $A_{\ge 4}\neq \emptyset$ for a $\mathcal{T}$. Then,
\begin{equation}
\mathcal{T}\in \mathscr{C}_0\iff \left\{
  \begin{array}{@{}ll@{}}
    |A_2|\le {{s}\choose{s/2}}-2, & \text{if}\ |A_{\ge 4}|\ge 2 \text{ or }|A_{\ge 2}|<\deg_T(\mathtt{c}), \nonumber\\
    |A_2|\le{{s}\choose{s/2}}-1, & \text{otherwise}. \nonumber
  \end{array}\right.
\end{equation}
\end{ppn}
\noindent\textit{Proof}: $(\Rightarrow)$ Since $\mathcal{T}\in \mathscr{C}_0$, there exists an orientation $D$ of $\mathcal{T}$, where $d(D)=4$. As $A_2\neq\emptyset$, we assume (\ref{eqC6.3.1})-(\ref{eqC6.3.2}) here. Let $[j]\in A_{\ge 4}$. Suppose $|I^\mathtt{c}((p,[j]))|\le \frac{s}{2}$ for some $1\le p\le s_{[j]}$. Then, for any $[i]\in A_2$, $d_D((1,[1,i]),(p,[j]))=3$ implies $O^\mathtt{c}((1,[i]))\cap I^\mathtt{c}((p,[j]))\neq \emptyset$. Hence, by Lih's theorem, $|A_2|=|B^O_2|\le {{s}\choose{s/2}}-{{s-|I^\mathtt{c}((1,[j]))|}\choose{s/2}}\le {{s}\choose{s/2}}-1$. If $|I^\mathtt{c}((p,[j]))|\ge \frac{s}{2}$ for some $1\le p\le s_{[j]}$, then $|A_2|\le {{s}\choose{s/2}}-1$ by the Duality Lemma.
\noindent\par Now, if $|A_{\ge 4}|=1$ and $|A_{\ge 2}|=\deg_T(\mathtt{c})$, then we are done. Consider the case $|A_{\ge 4}|\ge 2$ (or $|A_{\ge 2}|<\deg_T(\mathtt{c})$ resp.). By Griggs' theorem, $|A_2|={{s}\choose{s/2}}-1$ if and only if $\{O^\mathtt{c}((p,[j]))\}={{(\mathbb{N}_s, \mathtt{c})}\choose{s/2}}-B^O_2$ for all $1\le p\le s_{[j]}$ and $[j]\in A_{\ge 4}\cup E$. This implies $O^\mathtt{c}((p,[j_1]))=O^\mathtt{c}((q,[j_2]))$ for any distinct $[j_1], [j_2]\in A_{\ge 4}$ (or $[j_1]\in A_{\ge 4}$, $[j_2]\in E$ resp.), all $1\le p\le s_{[j_1]}$ and all $1\le q\le s_{[j_2]}$. However, it follows that $d_D((1,[\beta,j_1]),(q,[j_2]))>4$, a contradiction.
\\
\\$(\Leftarrow)$ The case of $|A_{\ge 4}|=1$ and $|A_{\ge 2}|=\deg_T(\mathtt{c})$ follows from Corollary \ref{corC6.3.8}(ii). Hence, consider $|A_{\ge 4}|\ge 2$ or $|A_{\ge 2}|<\deg_T(\mathtt{c})$. Furthermore, if $|A_2|+|A_{\ge 4}|\le {{s}\choose{s/2}}-1$, then by Corollary \ref{corC6.3.8}(i), $\mathcal{T}\in\mathscr{C}_0$. Hence, we assume $|A_2|+|A_{\ge 4}|\ge {{s}\choose{s/2}}$ hereafter, on top of the hypothesis that $|A_2|\le{{s}\choose{s/2}}-2$. If $|A_2|\ge s$, define $A^{\diamond}_2=A_2$. Otherwise, $A^{\diamond}_2=A_2\cup A^*$, where $A^*$ is an arbitrary subset of $A_{\ge 4}$ such that $|A^{\diamond}_2|=s$. Then, let $A^{\diamond}_4=A_2\cup A_{\ge 4}-A^{\diamond}_2$. Furthermore, assume without loss of generality that $A^{\diamond}_2=\{[i]\mid i\in \mathbb{N}_{|A^{\diamond}_2|}\}$ and $A^{\diamond}_4=\{[i]\mid i\in\mathbb{N}_{|A^{\diamond}_2|+|A^{\diamond}_4|}-\mathbb{N}_{|A^{\diamond}_2|}\}$.
\indent\par Let $\mathcal{H}=T(t_1,t_2,\ldots, t_n)$ be the subgraph of $\mathcal{T}$, where $t_\mathtt{c}=s$, $t_{[j]}=4$ for all $[j]\in \mathcal{T}(A^{\diamond}_4)$ and $t_v=2$ otherwise. We will use $A_j$ for $\mathcal{H}(A_j)$ for the remainder of this proof. Define an orientation $D$ of $\mathcal{H}$ as follows.
\begin{align}
&(2,[i])\rightarrow (1,[\alpha,i])\rightarrow (1,[i])\rightarrow (2,[\alpha,i])\rightarrow (2,[i]) \label{eqC6.3.27}
\end{align}
for all $[i]\in A_2$, and $1\le \alpha\le \deg_T([i])-1$.
\begin{align}
&\lambda_{i+1}\rightarrow \{(1,[i]),(2,[i])\} \rightarrow \bar{\lambda}_{i+1} \label{eqC6.3.28}
\end{align}
for all $1\le i\le \frac{s}{2}-1$.
\begin{align}
&\lambda_{i+2}\rightarrow \{(1,[i]),(2,[i])\}\rightarrow \bar{\lambda}_{i+2} \label{eqC6.3.29}
\end{align}
for $\frac{s}{2}\le i\le |A_2|$.  i.e., excluding $\lambda_1$ and $\lambda_{\frac{s}{2}+1}$, the $\frac{s}{2}$-sets $\lambda_i$'s ($\bar{\lambda}_i$'s resp.) are used as `in-sets' (`out-sets' resp.) to construct $B^I_2$ ($B^O_2$ resp.).
\begin{align}
&(2,[\beta,j])\rightarrow \{(2,[j]),(4,[j])\}\rightarrow (1,[\beta,j])\rightarrow \{(1,[j]), (3,[j])\}\rightarrow (2,[\beta,j]), \label{eqC6.3.30}\\
&\text{and }\lambda_{\frac{s}{2}+1} \rightarrow \{(1,[j]), (4,[j])\}\rightarrow \lambda_1 \rightarrow \{(2,[j]), (3,[j])\}\rightarrow \lambda_{\frac{s}{2}+1} \label{eqC6.3.31}
\end{align}
for all $[j]\in A_4$, and all $1\le \beta\le \deg_T([j])-1$.
\begin{align}
&\lambda_1\rightarrow \{(1,[l]), (2,[l])\} \rightarrow \lambda_{\frac{s}{2}+1} \label{eqC6.3.32}
\end{align}
for any $[l]\in E$. (See Figure \ref{figC6.3.9} for $D$ when $s=4$.)
\\
\\Claim: $d_{D}(v,w)\le 4$ for all $v,w\in V(D)$.
\\
\\Case 1.1. $v,w \in \{(1,[\alpha,i]),(2,[\alpha,i]),(1,[i]),(2,[i])\}$ for each $[i]\in A_2$ and $1\le\alpha\le \deg_T([i])-1$.
\indent\par Since the orientation defined for $A_2$ (see (\ref{eqC6.3.27})) is similar to that in Proposition \ref{ppnC6.3.5} (see (\ref{eqC6.3.15})), this case follows from Case 3.1 from Proposition \ref{ppnC6.3.5}.
\\
\\Case 1.2. $v,w \in \{(1,[\alpha,i]),(2,[\alpha,i]),(1,[i]),(2,[i]), (3,[i]), (4,[i])\}$ for each $[i]\in A_4$ and $1\le\alpha\le \deg_T([i])-1$.
\indent\par Since the orientation defined for $A_4$ (see (\ref{eqC6.3.30})) is similar to that in Proposition \ref{ppnC6.3.9} (see (\ref{eqC6.3.24})), this case follows from Case 1.3 from Proposition \ref{ppnC6.3.9}.
\\
\\Case 2. For each $[i],[j]\in A_2$, $i \neq j$, each $1\le\alpha\le \deg_T([i])-1$, and each $1\le\beta\le \deg_T([j])-1$,
\\(i) $v=(p,[\alpha,i]), w=(q,[j])$ for $p,q=1,2$.
\indent\par Since $O^\mathtt{c}((p,[i]))$ and $O^\mathtt{c}((q,[j]))$ are independent by (\ref{eqC6.3.28})-(\ref{eqC6.3.29}), there exists a vertex $z_{ij}\in O^\mathtt{c}((p,[i]))-O^\mathtt{c}((q,[j]))$ such that $(p,[\alpha,i])\rightarrow (p,[i])\rightarrow z_{ij} \rightarrow (q,[j])$ for $p,q=1,2$.
\\
\\(ii) $v=(p,[i]), w=(q,[\beta,j])$ for $p,q=1,2$.
\indent\par From (i), $d_D((p,[i]),(q,[j]))=2$. Since $(3-q,[j]) \rightarrow (q,[\beta,j])$ by (\ref{eqC6.3.27}), this subcase follows.
\\
\\(iii) $v=(p,[\alpha,i]), w=(q,[\beta,j])$ for each $p,q=1,2$.
\indent\par From (i), $d_D((p,[\alpha,i]),(q, [j]))=3$. Since $(3-q,[j]) \rightarrow (q,[\beta,j])$ by (\ref{eqC6.3.27}), the subcase follows.
\\
\\Case 3. For each $[i], [j]\in A_4$, $i \neq j$, each $1\le\alpha\le \deg_T([i])-1$, and each $1\le\beta\le \deg_T([j])-1$,
\\(i) $v=(p,[\alpha,i]), w=(q,[j])$ for $p=1,2$, and $q=1,2,3,4$.
\\(ii) $v=(p,[i]), w=(q,[\beta,j])$ for for $p=1,2,3,4$, and $q=1,2$.
\\(iii) $v=(p,[\alpha,i]), w=(q,[\beta,j])$ for $p,q=1,2$.
\indent\par Since the orientation defined for $A_4$ (see (\ref{eqC6.3.30})-(\ref{eqC6.3.31})) is similar to that in Proposition \ref{ppnC6.3.9} (see (\ref{eqC6.3.24})-(\ref{eqC6.3.25})), this case follows from Case 4 from Proposition \ref{ppnC6.3.9}.
\\
\\Case 4.  For each $[i]\in A_2$, each $[j]\in A_4$, each $1\le\alpha\le \deg_T([i])-1$, and each $1\le\beta\le \deg_T([j])-1$,
\\(i) $v=(p,[\alpha,i]), w=(q,[j])$ for each $p=1,2$, and $q=1,2,3,4$.
\indent\par By (\ref{eqC6.3.27}) and (\ref{eqC6.3.31}), $(p,[\alpha, i])\rightarrow (p,[i])$, $\lambda_1\rightarrow \{(2,[j]), (3,[j])\}$, and $\lambda_{\frac{s}{2}+1}\rightarrow \{(1,[j]),$ $(4,[j])\}$. Since $O^\mathtt{c}((p,[i]))\in {{(\mathbb{N}_s,\mathtt{c})}\choose{s/2}}-\{\lambda_1,\lambda_{\frac{s}{2}+1}\}$ by (\ref{eqC6.3.28})-(\ref{eqC6.3.29}), there exist some $x_i\in O^\mathtt{c}((p,[i]))\cap\lambda_1$ and $y_i\in O^\mathtt{c}((p,[i]))\cap\lambda_{\frac{s}{2}+1}$ such that $(p,[i])\rightarrow \{x_i, y_i\}$.
\\
\\(ii) $v=(q,[j]), w=(p,[\alpha,i])$ for each $p=1,2$, and $q=1,2,3,4$.
\indent\par By (\ref{eqC6.3.31}), $\{(1,[j]), (4,[j])\}\rightarrow \lambda_1$ and $\{(2,[j]), (3,[j])\}\rightarrow \lambda_{\frac{s}{2}+1}$. Since $(p,[i])\rightarrow (3-p,[\alpha, i])$ by (\ref{eqC6.3.27}), and $I^\mathtt{c}((p,[i]))\in {{(\mathbb{N}_s,\mathtt{c})}\choose{s/2}}-\{\lambda_1,\lambda_{\frac{s}{2}+1}\}$ by (\ref{eqC6.3.28})-(\ref{eqC6.3.29}), there exist some $x_{i}\in \lambda_1\cap I^\mathtt{c}((p,[i]))$ and $y_{i}\in\lambda_{\frac{s}{2}+1}\cap I^\mathtt{c}((p,[i]))$ such that $\{x_{i}, y_{i}\}\rightarrow (p,[i])$.
\\
\\(iii) $v=(p,[i]), w=(q,[\beta,j])$ for each $p,q=1,2$.
\indent\par By (\ref{eqC6.3.28})-(\ref{eqC6.3.31}), $|O^\mathtt{c}((p,[i]))|>0$, $\lambda_{\frac{s}{2}+1} \rightarrow \{(1,[j]), (4,[j])\}$, $\lambda_1 \rightarrow \{(2,[j]), (3,[j])\}$, \\$\{(2,[j]),(4,[j])\}\rightarrow (1,[\beta,j])$ and $\{(1,[j]), (3,[j])\}\rightarrow (2,[\beta,j])$.
\\
\\(iv) $v=(q,[\beta,j]), w=(p,[i])$ for each $p,q=1,2$.
\indent\par By (\ref{eqC6.3.28})-(\ref{eqC6.3.31}), $(2,[\beta,j])\rightarrow \{(2,[j]),(4,[j])\}$, $(1,[\beta,j])\rightarrow \{(1,[j]), (3,[j])\}$, $\{(1,[j]),$ $(4,[j])\}\rightarrow \lambda_1$, $\{(2,[j]), (3,[j])\}\rightarrow \lambda_{\frac{s}{2}+1}$, and $|I^\mathtt{c}((p,[i]))|>0$ for $p=1,2$.
\\
\\(v) $v=(p,[\alpha,i]), w=(q,[\beta,j])$ for each $p,q=1,2$.
\indent\par From (i), $d_D((p,[\alpha,i]),(q,[j]))=3$. Since $(3-q,[j])\rightarrow (q,[\beta,j])$ by (\ref{eqC6.3.30}), this subcase follows.
\\
\\(vi) $v=(q,[\beta,j]), w=(p,[\alpha,i])$ for each $p,q=1,2$.
\indent\par From (iv), $d_D((q,[\beta,j]),(p,[i]))=3$. Since $(3-p,[i])\rightarrow (p,[\alpha,i])$ by (\ref{eqC6.3.27}), this subcase follows.
\\
\\Case 5. For each $[i]\in A_2$, each $1\le\alpha\le \deg_T([i])-1$, and each $[j]\in E$,
\\(i) $v=(p,[\alpha,i]), w=(q,[j])$ for each $p=1,2$ and $q=1,2$.
\\(ii) $v=(q,[j]), w=(p,[\alpha,i])$ for each $p=1,2$ and $q=1,2$.
\indent\par Let $[k]\in A_4$. By (\ref{eqC6.3.31})-(\ref{eqC6.3.32}), $\lambda_1\rightarrow \{(q,[j]),(2,[k])\}\rightarrow \lambda_{\frac{s}{2}+1}$, this case follows from Cases 4(i)-(ii).
\\
\\Case 6. For each $[i]\in A_4$, each $1\le\alpha\le \deg_T([i])-1$, and each $[j]\in E$,
\\(i) $v=(p,[\alpha,i]), w=(q,[j])$ for each $p=1,2$ and $q=1,2$.
\\(ii) $v=(q,[j]), w=(p,[\alpha,i])$ for each $p=1,2$ and $q=1,2$.
\indent\par (i) follows from $(1,[\alpha,i])\rightarrow \{(1,[i]), (3,[i])\}$, $(2,[\alpha,i])\rightarrow \{(2,[i]),(4,[i])\}$, $\{(1,[i]),$ $(4,[i])\}\rightarrow \lambda_1$, $\{(2,[i]), (3,[i])\}\rightarrow \lambda_{\frac{s}{2}+1}$, and $|I^\mathtt{c}((q,[j]))|>0$ by (\ref{eqC6.3.30})-(\ref{eqC6.3.32}). (ii) follows from $\lambda_{\frac{s}{2}+1} \rightarrow \{(1,[i]), (4,[i])\}$, $\lambda_1 \rightarrow \{(2,[i]), (3,[i])\}$, $\{(2,[i]),(4,[i])\}\rightarrow (1,[\alpha,i])$ and $\{(1,[i]), (3,[i])\}\rightarrow (2,[\alpha,i])$ and $|O^\mathtt{c}((q,[j]))|>0$ by (\ref{eqC6.3.30})-(\ref{eqC6.3.32}).
\\
\\Case 7. $v=(r_1,\mathtt{c})$ and $w=(r_2,\mathtt{c})$ for $r_1\neq r_2$ and $1\le r_1, r_2\le s$.
\indent\par Here, we want to prove a stronger claim, $d_{D}((r_1,\mathtt{c}), (r_2,\mathtt{c}))=2$. Let $[k]\in A_4$, $x_1=(2,[k])$, $x_{i+1}=(1,[i])$ for all $1\le i\le \frac{s}{2}-1$, $x_{\frac{s}{2}+1}=(1,[k])$ and $x_{j+2}=(1,[j])$ for all $\frac{s}{2}\le j\le s$. Observe that $\lambda_i\rightarrow x_i\rightarrow \bar{\lambda}_i$ for $1\le i\le s$ by (\ref{eqC6.3.28})-(\ref{eqC6.3.29}) and (\ref{eqC6.3.31}), the subgraph induced by $V_1=(\mathbb{N}_s,\mathtt{c})$ and $V_2=\{x_i\mid 1\le i \le s\}$ is a complete bipartite graph $K(V_1,V_2)$. By Lemma \ref{lemC6.2.19}, $d_{D}((r_1,\mathtt{c}), (r_2,\mathtt{c}))=2$.
\\
\\Case 8. $v\in \{(1,[i]), (2,[i]), (3,[i]), (4,[i]), (1,[\alpha,i]), (2,[\alpha,i])\}$ for each $1\le i\le \deg_T(\mathtt{c})$ and $1\le\alpha\le \deg_T([i])-1$, and $w=(r,\mathtt{c})$ for $1\le r\le s$.
\indent\par Note that there exists some $1\le k\le s$ such that $d_D(v,(k,\mathtt{c}))\le 2$, and $d_D((k,\mathtt{c}),w)\le 2$ by Case 7. Hence, it follows that $d_D(v,w)\le d_D(v,(k,\mathtt{c}))+d_D((k,\mathtt{c}),w)\le 4$.
\\
\\Case 9. $v=(r,\mathtt{c})$ for $1\le r\le s$ and $w\in \{(1,[i]), (2,[i]), (3,[i]), (4,[i]), (1,[\alpha,i]), (2,[\alpha,i])\}$ for each $1\le i\le \deg_T(\mathtt{c})$ and $1\le\alpha\le \deg_T([i])-1$.
\indent\par Note that there exists some $1\le k\le s$ such that $d_D((k,\mathtt{c}), w)\le 2$, and $d_D(v,(k,\mathtt{c}))\le 2$ by Case 7. Hence, it follows that $d_D(v,w)\le d_D(v,(k,\mathtt{c}))+d_D((k,\mathtt{c}),w)\le 4$.
\\
\\Case 10. $v=(p,[i])$ and $w=(q, [j])$, where $1\le p,q\le 2$ and $1\le i,j\le \deg_T(\mathtt{c})$.
\noindent\par This follows from the fact that $|O^\mathtt{c}((p,[i]))|>0$, $|I^\mathtt{c}((q,[j]))|>0$, and $d_{D}((r_1,\mathtt{c}), (r_2,\mathtt{c}))$ $=2$ for any $r_1\neq r_2$ and $1\le r_1, r_2\le s$.
\\
\indent\par Therefore, the claim follows. Since every vertex lies in a directed $C_4$ for $D$ and $d(D)=4$, $\bar{d}(\mathcal{T})\le \max \{4, d(D)\}$ by Lemma \ref{lemC6.1.3}, and thus $\bar{d}(\mathcal{T})=4$ .
\qed

\begin{center}
\begin{tikzpicture}[thick,scale=0.75]%
\draw(0,5)node[circle, draw, inner sep=0pt, minimum width=3pt](1_u){\scriptsize $(1,\mathtt{c})$};
\draw(0,1)node[circle, draw, inner sep=0pt, minimum width=3pt](2_u){\scriptsize $(2,\mathtt{c})$};
\draw(0,-3)node[circle, draw, inner sep=0pt, minimum width=3pt](3_u){\scriptsize $(3,\mathtt{c})$};
\draw(0,-7)node[circle, draw, inner sep=0pt, minimum width=3pt](4_u){\scriptsize $(4,\mathtt{c})$};

\draw(-6,6)node[circle, draw, fill=black!100, inner sep=0pt, minimum width=6pt, label={[] 180:{\small $(1,[1,1])$}}](1_11){};
\draw(-6,4)node[circle, draw, fill=black!100, inner sep=0pt, minimum width=6pt, label={[] 180:{\small $(2,[1,1])$}}](2_11){};

\draw(-3,6)node[circle, draw, fill=black!100, inner sep=0pt, minimum width=6pt, label={[yshift=0cm, xshift=-0.45cm] 90:{\small $(1,[1])$}}](1_1){};
\draw(-3,4)node[circle, draw, fill=black!100, inner sep=0pt, minimum width=6pt, label={[yshift=0cm, xshift=-0.45cm] 270:{\small $(2,[1])$}}](2_1){};

\draw(6,6)node[circle, draw, fill=black!100, inner sep=0pt, minimum width=6pt, label={[] 0:{\small $(1,[1,2])$}}](1_21){};
\draw(6,4)node[circle, draw, fill=black!100, inner sep=0pt, minimum width=6pt, label={[] 0:{\small $(2,[1,2])$}}](2_21){};

\draw(3,6)node[circle, draw, fill=black!100, inner sep=0pt, minimum width=6pt, label={[yshift=0cm, xshift=0.45cm] 90:{\small $(1,[2])$}}](1_2){};
\draw(3,4)node[circle, draw, fill=black!100, inner sep=0pt, minimum width=6pt, label={[yshift=0cm, xshift=0.45cm] 270:{\small $(2,[2])$}}](2_2){};

\draw(-6,2)node[circle, draw, fill=black!100, inner sep=0pt, minimum width=6pt, label={[] 180:{\small $(1,[1,3])$}}](1_31){};
\draw(-6,0)node[circle, draw, fill=black!100, inner sep=0pt, minimum width=6pt, label={[] 180:{\small $(2,[1,3])$}}](2_31){};

\draw(-3,2)node[circle, draw, fill=black!100, inner sep=0pt, minimum width=6pt, label={[yshift=0cm, xshift=-0.45cm] 270:{\small $(1,[3])$}}](1_3){};
\draw(-3,0)node[circle, draw, fill=black!100, inner sep=0pt, minimum width=6pt, label={[yshift=0cm, xshift=-0.45cm] 270:{\small $(2,[3])$}}](2_3){};

\draw(6,2)node[circle, draw, fill=black!100, inner sep=0pt, minimum width=6pt, label={[] 0:{\small $(1,[1,4])$}}](1_41){};
\draw(6,0)node[circle, draw, fill=black!100, inner sep=0pt, minimum width=6pt, label={[] 0:{\small $(2,[1,4])$}}](2_41){};

\draw(3,2)node[circle, draw, fill=black!100, inner sep=0pt, minimum width=6pt, label={[yshift=0cm, xshift=0.45cm] 270:{\small $(1,[4])$}}](1_4){};
\draw(3,0)node[circle, draw, fill=black!100, inner sep=0pt, minimum width=6pt, label={[yshift=0cm, xshift=0.45cm] 270:{\small $(2,[4])$}}](2_4){};

\draw(-6,-4)node[circle, draw, fill=black!100, inner sep=0pt, minimum width=6pt, label={[] 180:{\small $(1,[1,5])$}}](1_51){};
\draw(-6,-6)node[circle, draw, fill=black!100, inner sep=0pt, minimum width=6pt, label={[] 180:{\small $(2,[1,5])$}}](2_51){};

\draw(-3,-2)node[circle, draw, fill=black!100, inner sep=0pt, minimum width=6pt, label={[yshift=-0.1cm, xshift=-0.45cm] 270:{\small $(1,[5])$}}](1_5){};
\draw(-3,-4)node[circle, draw, fill=black!100, inner sep=0pt, minimum width=6pt, label={[yshift=0cm, xshift=-0.45cm] 270:{\small $(2,[5])$}}](2_5){};
\draw(-3,-6)node[circle, draw, fill=black!100, inner sep=0pt, minimum width=6pt, label={[yshift=0.05cm, xshift=-0.45cm] 270:{\small $(3,[5])$}}](3_5){};
\draw(-3,-8)node[circle, draw, fill=black!100, inner sep=0pt, minimum width=6pt, label={[yshift=0cm, xshift=-0.45cm] 270:{\small $(4,[5])$}}](4_5){};

\draw(6,-4)node[circle, draw, fill=black!100, inner sep=0pt, minimum width=6pt, label={[] 0:{\small $(1,[1,6])$}}](1_61){};
\draw(6,-6)node[circle, draw, fill=black!100, inner sep=0pt, minimum width=6pt, label={[] 0:{\small $(2,[1,6])$}}](2_61){};

\draw(3,-2)node[circle, draw, fill=black!100, inner sep=0pt, minimum width=6pt, label={[yshift=-0.1cm, xshift=0.45cm] 270:{\small $(1,[6])$}}](1_6){};
\draw(3,-4)node[circle, draw, fill=black!100, inner sep=0pt, minimum width=6pt, label={[yshift=0cm, xshift=0.45cm] 270:{\small $(2,[6])$}}](2_6){};
\draw(3,-6)node[circle, draw, fill=black!100, inner sep=0pt, minimum width=6pt, label={[yshift=0.05cm, xshift=0.45cm] 270:{\small $(3,[6])$}}](3_6){};
\draw(3,-8)node[circle, draw, fill=black!100, inner sep=0pt, minimum width=6pt, label={[yshift=0cm, xshift=0.45cm] 270:{\small $(4,[6])$}}](4_6){};

\draw[->, line width=0.3mm, >=latex, shorten <= 0.2cm, shorten >= 0.15cm](1_11)--(1_1);
\draw[->, line width=0.3mm, >=latex, shorten <= 0.2cm, shorten >= 0.15cm](2_11)--(2_1);
\draw[dashed, ->, line width=0.3mm, >=latex, shorten <= 0.2cm, shorten >= 0.15cm](2_1)--(1_11);
\draw[dashed, ->, line width=0.3mm, >=latex, shorten <= 0.2cm, shorten >= 0.15cm](1_1)--(2_11);

\draw[dashed, ->, line width=0.3mm, >=latex, shorten <= 0.2cm, shorten >= 0.15cm](2_2)--(1_21);
\draw[dashed, ->, line width=0.3mm, >=latex, shorten <= 0.2cm, shorten >= 0.15cm](1_2)--(2_21);
\draw[->, line width=0.3mm, >=latex, shorten <= 0.2cm, shorten >= 0.15cm](1_21)--(1_2);
\draw[->, line width=0.3mm, >=latex, shorten <= 0.2cm, shorten >= 0.15cm](2_21)--(2_2);

\draw[dashed, ->, line width=0.3mm, >=latex, shorten <= 0.2cm, shorten >= 0.15cm](2_3)--(1_31);
\draw[dashed, ->, line width=0.3mm, >=latex, shorten <= 0.2cm, shorten >= 0.15cm](1_3) to [out=190, in=65] (2_31);
\draw[->, line width=0.3mm, >=latex, shorten <= 0.2cm, shorten >= 0.15cm](1_31)--(1_3);
\draw[->, line width=0.3mm, >=latex, shorten <= 0.2cm, shorten >= 0.15cm](2_31)--(2_3);

\draw[dashed, ->, line width=0.3mm, >=latex, shorten <= 0.2cm, shorten >= 0.15cm](2_4)--(1_41);
\draw[dashed, ->, line width=0.3mm, >=latex, shorten <= 0.2cm, shorten >= 0.15cm](1_4) to [out=350, in=115] (2_41);
\draw[->, line width=0.3mm, >=latex, shorten <= 0.2cm, shorten >= 0.15cm](1_41)--(1_4);
\draw[->, line width=0.3mm, >=latex, shorten <= 0.2cm, shorten >= 0.15cm](2_41)--(2_4);

\draw[dashed, ->, line width=0.3mm, >=latex, shorten <= 0.2cm, shorten >= 0.15cm](1_6) to [out=345, in=100] (2_61);
\draw[dashed, ->, line width=0.3mm, >=latex, shorten <= 0.2cm, shorten >= 0.15cm](2_6)--(1_61);
\draw[dashed, ->, line width=0.3mm, >=latex, shorten <= 0.2cm, shorten >= 0.15cm](3_6)--(2_61);
\draw[dashed, ->, line width=0.3mm, >=latex, shorten <= 0.2cm, shorten >= 0.15cm](4_6) to [out=40, in=250] (1_61);
\draw[->, line width=0.3mm, >=latex, shorten <= 0.2cm, shorten >= 0.15cm](1_61) to [out=115, in=0] (1_6);
\draw[->, line width=0.3mm, >=latex, shorten <= 0.2cm, shorten >= 0.15cm](2_61) to [out=115, in=350] (2_6);
\draw[->, line width=0.3mm, >=latex, shorten <= 0.2cm, shorten >= 0.15cm](1_61)--(3_6);
\draw[->, line width=0.3mm, >=latex, shorten <= 0.2cm, shorten >= 0.15cm](2_61)--(4_6);

\draw[dashed, ->, line width=0.3mm, >=latex, shorten <= 0.2cm, shorten >= 0.15cm](1_5) to [out=195, in=85] (2_51);
\draw[dashed, ->, line width=0.3mm, >=latex, shorten <= 0.2cm, shorten >= 0.15cm](2_5)--(1_51);
\draw[dashed, ->, line width=0.3mm, >=latex, shorten <= 0.2cm, shorten >= 0.15cm](3_5)--(2_51);
\draw[dashed, ->, line width=0.3mm, >=latex, shorten <= 0.2cm, shorten >= 0.15cm](4_5) to [out=140, in=290] (1_51);
\draw[->, line width=0.3mm, >=latex, shorten <= 0.2cm, shorten >= 0.15cm](1_51) to [out=65, in=180] (1_5);
\draw[->, line width=0.3mm, >=latex, shorten <= 0.2cm, shorten >= 0.15cm](2_51) to [out=65, in=190] (2_5);
\draw[->, line width=0.3mm, >=latex, shorten <= 0.2cm, shorten >= 0.15cm](1_51)--(3_5);
\draw[->, line width=0.3mm, >=latex, shorten <= 0.2cm, shorten >= 0.15cm](2_51)--(4_5);

\draw[densely dotted, ->, line width=0.3mm, >=latex, shorten <= 0.2cm, shorten >= 0.1cm](2_1)--(1_u);
\draw[densely dotted, ->, line width=0.3mm, >=latex, shorten <= 0.2cm, shorten >= 0.1cm](2_1)--(4_u);
\draw[densely dotted, ->, line width=0.3mm, >=latex, shorten <= 0.2cm, shorten >= 0.1cm](1_1)--(1_u);
\draw[densely dotted, ->, line width=0.3mm, >=latex, shorten <= 0.2cm, shorten >= 0.1cm](1_1)--(4_u);

\draw[densely dotted, ->, line width=0.3mm, >=latex, shorten <= 0.2cm, shorten >= 0.1cm](2_2)--(2_u);
\draw[densely dotted, ->, line width=0.3mm, >=latex, shorten <= 0.2cm, shorten >= 0.1cm](2_2)--(3_u);
\draw[densely dotted, ->, line width=0.3mm, >=latex, shorten <= 0.2cm, shorten >= 0.1cm](1_2)--(2_u);
\draw[densely dotted, ->, line width=0.3mm, >=latex, shorten <= 0.2cm, shorten >= 0.1cm](1_2)--(3_u);

\draw[densely dotted, ->, line width=0.3mm, >=latex, shorten <= 0.2cm, shorten >= 0.1cm](2_3)--(1_u);
\draw[densely dotted, ->, line width=0.3mm, >=latex, shorten <= 0.2cm, shorten >= 0.1cm](2_3)--(3_u);
\draw[densely dotted, ->, line width=0.3mm, >=latex, shorten <= 0.2cm, shorten >= 0.1cm](1_3)--(1_u);
\draw[densely dotted, ->, line width=0.3mm, >=latex, shorten <= 0.2cm, shorten >= 0.1cm](1_3)--(3_u);

\draw[densely dotted, ->, line width=0.3mm, >=latex, shorten <= 0.2cm, shorten >= 0.1cm](2_4)--(2_u);
\draw[densely dotted, ->, line width=0.3mm, >=latex, shorten <= 0.2cm, shorten >= 0.1cm](2_4)--(4_u);
\draw[densely dotted, ->, line width=0.3mm, >=latex, shorten <= 0.2cm, shorten >= 0.1cm](1_4)--(2_u);
\draw[densely dotted, ->, line width=0.3mm, >=latex, shorten <= 0.2cm, shorten >= 0.1cm](1_4)--(4_u);

\draw[->, line width=0.3mm, >=latex, shorten <= 0.2cm, shorten >= 0.1cm](2_5)--(3_u);
\draw[->, line width=0.3mm, >=latex, shorten <= 0.2cm, shorten >= 0.1cm](2_5)--(4_u);
\draw[->, line width=0.3mm, >=latex, shorten <= 0.2cm, shorten >= 0.1cm](3_5)--(3_u);
\draw[->, line width=0.3mm, >=latex, shorten <= 0.2cm, shorten >= 0.1cm](3_5)--(4_u);
\draw[->, line width=0.3mm, >=latex, shorten <= 0.2cm, shorten >= 0.1cm](1_5)--(1_u);
\draw[->, line width=0.3mm, >=latex, shorten <= 0.2cm, shorten >= 0.1cm](1_5)--(2_u);
\draw[->, line width=0.3mm, >=latex, shorten <= 0.2cm, shorten >= 0.1cm](4_5)--(1_u);
\draw[->, line width=0.3mm, >=latex, shorten <= 0.2cm, shorten >= 0.1cm](4_5)--(2_u);

\draw[->, line width=0.3mm, >=latex, shorten <= 0.2cm, shorten >= 0.1cm](2_6)--(3_u);
\draw[->, line width=0.3mm, >=latex, shorten <= 0.2cm, shorten >= 0.1cm](3_6)--(3_u);
\draw[->, line width=0.3mm, >=latex, shorten <= 0.2cm, shorten >= 0.1cm](2_6)--(4_u);
\draw[->, line width=0.3mm, >=latex, shorten <= 0.2cm, shorten >= 0.1cm](3_6)--(4_u);
\draw[->, line width=0.3mm, >=latex, shorten <= 0.2cm, shorten >= 0.1cm](1_6)--(1_u);
\draw[->, line width=0.3mm, >=latex, shorten <= 0.2cm, shorten >= 0.1cm](1_6)--(2_u);
\draw[->, line width=0.3mm, >=latex, shorten <= 0.2cm, shorten >= 0.1cm](4_6)--(1_u);
\draw[->, line width=0.3mm, >=latex, shorten <= 0.2cm, shorten >= 0.1cm](4_6)--(2_u);
\end{tikzpicture}
\captionsetup{justification=centering}

{\captionof{figure}{Orientation $D$ for $\mathcal{H}$ for $s=4$, $A_2=\{[1],[2],[3],[4]\}$, $A_4=\{[5],[6]\}$.}\label{figC6.3.9}}
\end{center}

\begin{ppn}\label{ppnC6.3.11}
Suppose $s\ge 4$ is even, $A_2\neq \emptyset$, $A_3=\{[j]\}$, and $A_{\ge 4}=\emptyset$ for a $\mathcal{T}$. Then, 
\begin{align*}
\mathcal{T}\in \mathscr{C}_0\iff \left\{
  \begin{array}{@{}ll@{}}
    |A_2|\le{{s}\choose{s/2}}-2, & \text{if}\ |A_{\ge 2}|< \deg_T(\mathtt{c}),\\
    |A_2|\le{{s}\choose{s/2}}-1, & \text{if}\ |A_{\ge 2}|=\deg_T(\mathtt{c}).
  \end{array}\right.
\end{align*}
\end{ppn}
\noindent\textit{Proof}: $(\Rightarrow)$ Since $\mathcal{T}\in \mathscr{C}_0$, there exists an orientation $D$ of $\mathcal{T}$, where $d(D)=4$. As $A_2\neq\emptyset$, we assume (\ref{eqC6.3.1})-(\ref{eqC6.3.2}) here; (\ref{eqC6.3.3})-(\ref{eqC6.3.5}) is not required though $A_3\neq\emptyset$.
\\
\\Case 1. $|A_2|= \deg_T(\mathtt{c})$.
\indent\par If $|O((1, [1,j]))|=1$, say $O((1, [1,j]))=\{(1, [j])\}$, then by Lemma \ref{lemC6.2.15}, $\{O^\mathtt{c}((1, [i]))\mid [i]$ $\in A_{\ge 2}\}$ is an antichain. So, $|A_{\ge 2}|\le {{s}\choose{s/2}}$ by Sperner's theorem, i.e., $|A_2|\le{{s}\choose{s/2}}-1$. The case of $|I((1,[1,j]))|=1$ follows from the Duality Lemma.
\\
\\Case 2. $|A_2|< \deg_T(\mathtt{c})$.
\indent\par Suppose $|O((1, [1,j]))|=1$, say $O((1, [1,j]))=\{(1, [j])\}$, then by Lemma \ref{lemC6.2.15},  $\{O^\mathtt{c}((1, [i]))\mid [i]\in A_{\ge 2}\}$ is an antichain. Let $[i^*]\in E$.
\indent\par If $|O^\mathtt{c}((1,[i^*]))|\ge \frac{s}{2}$, then $d_D((1,[1,i]),(1,[i^*]))=3$ implies $O^\mathtt{c}((1,[i]))\cap I^\mathtt{c}((1,[i^*]))\neq \emptyset$ for all $[i]\in A_{\ge 2}$. It follows from Lih's theorem that $|A_{\ge 2}|=|\{O^\mathtt{c}((1,[i]))\mid [i]\in A_{\ge 2}\}|\le {{s}\choose{{s/2}}}-{{s-|I^\mathtt{c}((1,[i^*]))|}\choose{{s/2}}}\le {{s}\choose{{s/2}}}-{{{s/2}}\choose{{s/2}}}={{s}\choose{{s/2}}}-1$, i.e., $|A_2|\le {{s}\choose{{s/2}}}-2$. 
\indent\par If $|O^\mathtt{c}((1,[i^*]))|<\frac{s}{2}$, then $d_D((1,[i^*]),(1,[1,i]))=3$ implies $I^\mathtt{c}((2,[i]))\cap O^\mathtt{c}((1,[i^*]))\neq \emptyset$ for all $[i]\in A_2$. By  Lemma \ref{lemC6.2.16}, $\{I^\mathtt{c}((2,[i]))\mid [i]\in A_2\}$ is an antichain. It follows from Lih's theorem that $|A_2|=|\{I^\mathtt{c}((2,[i]))\mid [i]\in A_2\}|\le {{s}\choose{{s/2}}}-{{s-|O^\mathtt{c}((1,[i^*]))|}\choose{{s/2}}}\le {{s}\choose{{s/2}}}-{{{(s/2)+1}}\choose{{s/2}}}={{s}\choose{{s/2}}}-\frac{s}{2}-1\le {{s}\choose{{s/2}}}-3$.
\indent\par The case of $|I((1,[1,j]))|=1$ follows from the Duality Lemma.
\\
\\$(\Leftarrow)$ By Corollary \ref{corC6.3.8}(i), $\mathcal{T}\in\mathscr{C}_0$.
\qed
\begin{ppn}\label{ppnC6.3.12}
Suppose $s\ge 4$ is even, $A_2\neq \emptyset$ and either $|A_3|\ge 2$ or $|A_3|=1$ and $A_{\ge 4}\neq\emptyset$ for a $\mathcal{T}$.
\\(a) If $\mathcal{T}\in \mathscr{C}_0$, then $2|A_2|+|A_3|\le {{s}\choose{s/2}}+{{s}\choose{(s/2)+1}}-\kappa^*_{s,\frac{s}{2}}(k)$ for some $k\le |A_2|+|A_3|$.
\\(b) If there exists some $|A_2|+1\le k\le \min\{|A_2|+|A_3|, {{s}\choose{s/2}}-1\}$ such that $2|A_2|+|A_3|\le {{s}\choose{s/2}}+{{s}\choose{(s/2)+1}}-\kappa_{s,\frac{s}{2}}(k)-3$, then $\mathcal{T}\in \mathscr{C}_0$.
\end{ppn}
\noindent\textit{Proof}: (a) Since $\mathcal{T}\in \mathscr{C}_0$, there exists an orientation $D$ of $\mathcal{T}$, where $d(D)=4$. As $A_2\neq\emptyset$ and $A_3\neq\emptyset$, we assume (\ref{eqC6.3.1})-(\ref{eqC6.3.5}) here. Partition $A^O_3$ ($A^I_3$ resp.) into $A^{O(D)}_3$ and $A^{O(S)}_3$ ($A^{I(D)}_3$ and $A^{I(S)}_3$ resp.), where 
\begin{align}
\left. \begin{array}{@{}ll@{}}
&A^{O(D)}_3=\{[i]\in A^O_3\mid O^\mathtt{c}((1,[i]))\neq O^\mathtt{c}((2,[i]))\},\\
&A^{O(S)}_3=\{[i]\in A^O_3\mid O^\mathtt{c}((1,[i]))=O^\mathtt{c}((2,[i]))\},\\
&A^{I(D)}_3=\{[i]\in A^I_3\mid O^\mathtt{c}((1,[i]))\neq O^\mathtt{c}((2,[i]))\},\\
\text{and} &A^{I(S)}_3=\{[i]\in A^I_3\mid O^\mathtt{c}((1,[i]))=O^\mathtt{c}((2,[i]))\}.
  \end{array}\right\}
\label{eqC6.3.33}
\end{align}
Note that $O^\mathtt{c}((1,[i]))$ and $O^\mathtt{c}((2,[i]))$ are \textit{Different} (\textit{Same} resp.) for $A^{O(D)}_3$ and $A^{I(D)}_3$ ($A^{O(S)}_3$ and $A^{I(S)}_3$ resp.).
\indent\par Both $B^O_2\cup B^O_3\cup \{O^\mathtt{c}((1,[l]))\mid [l]\in A^{I(S)}_3\}$ and $B^I_2\cup B^I_3\cup \{I^\mathtt{c}((2,[l]))\mid [l]\in A^{O(S)}_3\}$ are antichains by Lemmas \ref{lemC6.2.15} and \ref{lemC6.2.16} respectively. Furthermore, $d_D((1,[1,i]),(1,[1,j]))=4$ for each $[i]\in A_2\cup A^O_3\cup A^{I(S)}_3$ and $[j]\in A_2\cup A^I_3\cup A^{O(S)}_3$, $i\neq j$, implies $O^\mathtt{c}((\delta_1,[i]))\cap I^\mathtt{c}((\delta_2,[j]))\neq\emptyset$ where 
\begin{align*}
\delta_1=\left\{
  \begin{array}{@{}ll@{}}
    1, & \text{if}\ [i]\in A_2 \cup A^{I(S)}_3,\\
    3, & \text{if}\ [i]\in A^O_3,\\
  \end{array}\right.
\text{ and }
\delta_2=\left\{
  \begin{array}{@{}ll@{}}
    2, & \text{if}\ [j]\in A_2 \cup A^{O(S)}_3,\\
    3, & \text{if}\ [j]\in A^I_3.\\
  \end{array}\right.
\end{align*}

\noindent Equivalently, for any $O^\mathtt{c}((\delta_1,[i]))\in B^O_2\cup B^O_3\cup \{O^\mathtt{c}((1,[l]))\mid [l]\in A^{I(S)}_3\}$ and $I^\mathtt{c}((\delta_2,[j]))\in B^I_2\cup B^I_3\cup \{I^\mathtt{c}((2,[l]))\mid [l]\in A^{O(S)}_3\}$, $O^\mathtt{c}((\delta_1,[i]))\cap I^\mathtt{c}((\delta_2,[j]))=\emptyset$ only if $[i]=[j]\in A_2\cup A^{O(S)}_3\cup A^{I(S)}_3$. Then,
\begin{align*}
&\ 2|A_2|+|A_3|\\
=&\ 2|A_2|+|A^O_3|+|A^I_3|\\
\le &(|A_2|+|A^O_3|+|A^{I(S)}_3|)+(|A_2|+|A^I_3|+|A^{O(S)}_3|) \\
=&\ |B^O_2\cup B^O_3\cup \{O^\mathtt{c}((1,[l]))\mid [l]\in A^{I(S)}_3\}|+|B^I_2\cup B^I_3\cup \{I^\mathtt{c}((2,[l]))\mid [l]\in A^{O(S)}_3\}|\nonumber\\
\le&\ {{s}\choose{s/2}}+{{s}\choose{(s/2)+1}}-\kappa^*_{s,\frac{s}{2}}(|A_2|+|A^{O(S)}_3|+|A^{I(S)}_3|) \text{ by  Theorem } \ref{thmC6.2.10},
\end{align*} 
where $k=|A_2|+|A^{O(S)}_3|+|A^{I(S)}_3|\le |A_2|+|A_3|$ is as required.
\\
\\(b) Let the $\frac{s}{2}$-subsets of $\mathbb{N}_s$ in squashed order be $X_1, X_2, \ldots, X_{{s}\choose{s/2}}$. Note that $\bar{X}_i=X_{{{s}\choose{s/2}}-i}$ for $i=1,2,\ldots, {{s}\choose{s/2}}$ and $L_{s,\frac{s}{2}}(k)=\{X_{{s}\choose{s/2}}, X_{{{s}\choose{s/2}}-1},\ldots, X_{{{s}\choose{s/2}}-k+1}\}$. Furthermore, let $\mu_i=(X_i,\mathtt{c})$ for $i=1,2,\ldots, {{s}\choose{s/2}}$. We also use the previous notation, ${{(\mathbb{N}_s,\mathtt{c})}\choose{(s/2)+1}}=\{\psi_i\mid i=1,2,\ldots, {{s}\choose{(s/2)+1}}\}$, and further assume $\{\psi_i\mid i=1,2,\ldots,{{s}\choose{(s/2)+1}}-|\nabla L_{s,\frac{s}{2}}(k)|\}=\{(Y,\mathtt{c})\mid Y\in {{\mathbb{N}_s}\choose{(s/2)+1}}-\nabla L_{s,\frac{s}{2}}(k)\}$.
\indent\par If $|A_{\ge 2}|\le {{s}\choose{s/2}}-1$, then by Corollary \ref{corC6.3.8}(i), $\mathcal{T}\in\mathscr{C}_0$. Hence, we assume $|A_{\ge 2}|\ge{{s}\choose{s/2}}$. Let $A^{\diamond}_2=A_2\cup A^*$, where $A^*$ is an arbitrary subset of $A_3$ such that $|A^{\diamond}_2|=k-1$; $A^*=\emptyset$ if $|A_2|=k-1$. Then, let $A^{\diamond}_3=A_3-A^*$. Furthermore, assume without loss of generality that $A^{\diamond}_2=\{[i]\mid i\in\mathbb{N}_{|A^{\diamond}_2|}\}$, $A^{\diamond}_3=\{[i]\mid i\in\mathbb{N}_{|A^{\diamond}_2|+|A^{\diamond}_3|}-\mathbb{N}_{|A^{\diamond}_2|}\}$ and $A^{\diamond}_4=\{[i]\mid i\in\mathbb{N}_{|A^{\diamond}_2|+|A^{\diamond}_3|+|A^{\diamond}_4|}-\mathbb{N}_{|A^{\diamond}_2|+|A^{\diamond}_3|}\}$.
\noindent\par Let $\mathcal{H}=T(t_1,t_2,\ldots, t_n)$ be the subgraph of $\mathcal{T}$, where  $t_\mathtt{c}=s$, $t_{[i]}=3$ for all $[i]\in \mathcal{T}(A^{\diamond}_3)$, $t_{[j]}=4$ for all $[j]\in \mathcal{T}(A^{\diamond}_4)$ and $t_v=2$ otherwise. We will use $A_j$ for $\mathcal{H}(A_j)$ for the remainder of this proof. Define an orientation $D$ of $\mathcal{H}$ as follows.
\begin{align}
&(2,[i])\rightarrow (1,[\alpha,i])\rightarrow (1,[i])\rightarrow (2,[\alpha,i])\rightarrow (2,[i]),\text{ and}\label{eqC6.3.34}\\
&\mu_{i+1}\rightarrow \{(1,[i]),(2,[i])\} \rightarrow \bar{\mu}_{i+1}\label{eqC6.3.35}
\end{align}
for all $1\le i\le |A_2|=k-1$, and $1\le \alpha\le \deg_T([i])-1$, i.e., the $\frac{s}{2}$-sets $\mu_2,\mu_3,\ldots,\mu_k$ ($\bar{\mu}_2,\bar{\mu}_3,\ldots,\bar{\mu}_k$ resp.) are used as `in-sets' (`out-sets' resp.) to construct $B^I_2$ ($B^O_2$ resp.).
\begin{align}
& (3,[j])\rightarrow \{(1,[\beta,j]),(2,[\beta,j])\}\rightarrow \{(1,[j]),(2,[j])\},\label{eqC6.3.36}\\
& \bar{\mu}_1=\mu_{{s}\choose{s/2}} \rightarrow (1,[j])\rightarrow \mu_1\rightarrow (2,[j])\rightarrow \mu_{{s}\choose{s/2}},\text{ and}\label{eqC6.3.37}\\
& \mu_{j+1}\rightarrow (3,[j])\rightarrow \bar{\mu}_{j+1} \label{eqC6.3.38}
\end{align}
for all $|A_2|+1\le j\le {{s}\choose{s/2}}-2$ and all $1\le \alpha\le \deg_T([j])-1$, i.e., the $\frac{s}{2}$-sets $\mu_{k+1},\mu_{k+2},\ldots,\mu_{{{s}\choose{s/2}}-1}$ are used as `in-sets' to construct $B^I_3$.
\begin{align}
& \{(1,[l]),(2,[l])\}\rightarrow \{(1,[\gamma,l]),(2,[\gamma,l])\}\rightarrow (3,[l]),\label{eqC6.3.39}\\
& \mu_1 \rightarrow (1,[l])\rightarrow \mu_{{s}\choose{s/2}}\rightarrow (2,[l])\rightarrow \mu_1,\text{ and}\label{eqC6.3.40}\\
&\bar{\psi}_{l+2-{{s}\choose{s/2}}}\rightarrow (3,[l])\rightarrow \psi_{l+2-{{s}\choose{s/2}}}\label{eqC6.3.41}
\end{align}
for all ${{s}\choose{s/2}}-1\le l\le |A_2|+|A_3|$ and all $1\le \gamma\le \deg_T([l])-1$, i.e.,  the $(\frac{s}{2}+1)$-sets $\psi_1,\psi_2,\ldots,\psi_{{{s}\choose{(s/2)+1}}-|\nabla L_{s,\frac{s}{2}}(k)|}$ are used as `out-sets' to construct $B^O_3$.
\begin{align}
&(2,[\tau,x])\rightarrow \{(2,[x]),(4,[x])\}\rightarrow (1,[\tau,x])\rightarrow \{(1,[x]), (3,[x])\}\rightarrow (2,[\tau,x]),\label{eqC6.3.42}\\
&\text{and }\mu_{{s}\choose{s/2}} \rightarrow \{(1,[x]), (4,[x])\}\rightarrow \mu_1 \rightarrow \{(2,[x]), (3,[x])\}\rightarrow \mu_{{s}\choose{s/2}}\label{eqC6.3.43}
\end{align}
for all $[x]\in A_4$ and all $1\le \tau\le \deg_T([x])-1$.
\begin{align}
&\mu_1\rightarrow \{(1,[y]), (2,[y])\} \rightarrow \mu_{{s}\choose{s/2}}\label{eqC6.3.44}
\end{align}
for any $[y]\in E$. (See Figures \ref{figC6.3.10} and \ref{figC6.3.11} when $s=6$, $k=13$.)
\\
\\Claim: $d_{D}(v,w)\le 4$ for all $v,w\in V(D)$.
\\
\\Case 1. $v \in \{(1,[\alpha,i]),(2,[\alpha,i]),(1,[i]),(2,[i])\}$ and $w\in \{(1,[\beta,j]),(2,[\beta,j]),(1,[j]),$ $(2,[j])\}$ for each $[i], [j]\in A_2$, each $1\le\alpha\le \deg_T([i])-1$, and each $1\le\beta\le \deg_T([j])-1$.
\noindent\par Since the orientation defined for $A_2$ (see (\ref{eqC6.3.34})-(\ref{eqC6.3.35})) is similar to that in Proposition \ref{ppnC6.3.10} (see (\ref{eqC6.3.27})-(\ref{eqC6.3.29})), this case follows from Cases 1.1 and 2 of Proposition \ref{ppnC6.3.10}, except possibly for $v=(p,[i])$ and $w=(q, [j])$, which will be covered in Case 14 later.
\\
\\Case 2. $v \in \{(1,[\alpha,i]),(2,[\alpha,i]),(1,[i]),(2,[i]),(3,[i])\}$ and $w\in \{(1,[\beta,j]),(2,[\beta,j]),$ $(1,[j]),(2,[j]),(3,[j])\}$ for each $[i],[j]\in A_3$, each $1\le\alpha\le \deg_T([i])-1$, and each $1\le\beta\le \deg_T([j])-1$.
\indent\par Since the orientation defined for $A_3$ (see (\ref{eqC6.3.36})-(\ref{eqC6.3.41})) is similar to that in Proposition \ref{ppnC6.3.9} (see (\ref{eqC6.3.17})-(\ref{eqC6.3.23}); note that $\mu_1=\lambda_1$ and $\mu_{{s}\choose{s/2}}=\lambda_{\frac{s}{2}+1}$), this case follows from Cases 1.1-1.2, 2-3 and 5 of Proposition \ref{ppnC6.3.9}, except possibly for $v=(p,[i])$ and $w=(q, [j])$, which will be covered in Case 14 later.
\\
\\Case 3. $v \in \{(1,[\alpha,i]),(2,[\alpha,i]),(1,[i]),(2,[i]),(3,[i]),(4,[i])\}$ and $w\in \{(1,[\beta,j]),(2,[\beta,j]),$ $(1,[j]),(2,[j]),(3,[j]),(4,[j])\}$ for each $[i], [j]\in A_4$, each $1\le\alpha\le \deg_T([i])-1$, and each $1\le\beta\le \deg_T([j])-1$.
\indent\par Since the orientation defined for $A_4$ (see (\ref{eqC6.3.42})-(\ref{eqC6.3.43})) is similar to that in Proposition \ref{ppnC6.3.9} (see (\ref{eqC6.3.24})-(\ref{eqC6.3.25})), this case follows from Cases 1.3 and 4 of Proposition \ref{ppnC6.3.9}, except possibly for $v=(p,[i])$ and $w=(q, [j])$, which will be covered in Case 14 later.
\\
\\Case 4. $v \in \{(1,[\alpha,i]),(2,[\alpha,i]),(1,[i]),(2,[i])\}$ and $w\in \{(1,[\beta,j]),(2,[\beta,j]),(1,[j]),$ $(2,[j]),(3,[j])\}$ for each $1\le i\le |A_2|$, each $|A_2|+1\le j\le {{s}\choose{s/2}}-2$, each $1\le\alpha\le \deg_T([i])-1$, and each $1\le\beta\le \deg_T([j])-1$.
\\(i) $v=(p,[\alpha,i]), w=(q,[j])$ for each $p=1,2$ and $q=1,2,3$.
\indent\par By (\ref{eqC6.3.34}) and (\ref{eqC6.3.37}), $(p,[\alpha,i])\rightarrow (p,[i])$, $\mu_{{s}\choose{s/2}} \rightarrow (1,[j])$, and $\mu_1\rightarrow (2,[j])$. With (\ref{eqC6.3.35}), $O^\mathtt{c}((p,[i]))\in {{(\mathbb{N}_s,\mathtt{c})}\choose{s/2}}-\{\mu_1,\mu_{{s}\choose{s/2}}\}$, there exist some $x_i\in O^\mathtt{c}((p,[i]))\cap\mu_1$ and $y_i\in O^\mathtt{c}((p,[i]))\cap\mu_{{s}\choose{s/2}}$ such that $(p,[i])\rightarrow \{x_i, y_i\}$. Furthermore, with (\ref{eqC6.3.38}), $O^\mathtt{c}((p,[i]))$ and $O^\mathtt{c}((3,[j]))$ are independent. So, there exists some $z_{ij}\in O^\mathtt{c}((p,[i]))\cap I^\mathtt{c}((3,[j]))$ such that $(p,[i])\rightarrow z_{ij} \rightarrow (3,[j])$.
\\
\\(ii) $v=(q,[j]), w=(p,[\alpha,i])$ for each $p=1,2$ and $q=1,2,3$.
\indent\par By (\ref{eqC6.3.34}) and (\ref{eqC6.3.37}), $(3-p,[i])\rightarrow (p,[\alpha, i])$, $(1,[j])\rightarrow \mu_1$ and $(2,[j])\rightarrow \mu_{{s}\choose{s/2}}$. Since $I^\mathtt{c}((p,[i]))\in {{(\mathbb{N}_s,\mathtt{c})}\choose{s/2}}-\{\mu_1,\mu_{{s}\choose{s/2}}\}$, there exist some $x_i\in \mu_1\cap I^\mathtt{c}((p,[i]))$ and $y_i\in \mu_{{s}\choose{s/2}}\cap I^\mathtt{c}((p,[i]))$ such that $\{x_i, y_i\}\rightarrow (p,[i])$. Furthermore, with (\ref{eqC6.3.38}), $O^\mathtt{c}((p,[i]))$ and $O^\mathtt{c}((3,[j]))$ are independent. So, there exists some $z_{ij}\in O^\mathtt{c}((3,[j]))\cap I^\mathtt{c}((p,[i]))$ such that $(3,[j])\rightarrow z_{ij} \rightarrow (p,[i])$.
\\
\\(iii) $v=(p,[i]), w=(q,[\beta,j])$ for each $p=1,2$ and $q=1,2$.
\indent\par From (i), $d_D((p,[i]),(3,[j]))=2$. Since $(3,[j])\rightarrow \{(1,[\beta,j]),(2,[\beta,j])\}$ by (\ref{eqC6.3.36}), this subcase follows.
\\
\\(iv) $v=(q,[\beta,j]), w=(p,[i])$ for each $p=1,2$ and $q=1,2$.
\indent\par By (\ref{eqC6.3.35})-(\ref{eqC6.3.37}), $\{(1,[\beta,j]),(2,[\beta,j])\}\rightarrow \{(1,[j]),(2,[j])\}$, $(1,[j])\rightarrow \mu_1$, $(2,[j])\rightarrow \mu_{{s}\choose{s/2}}$ and $|I^\mathtt{c}((p,[i]))|>0$.
\\
\\(v) $v=(p,[\alpha,i]), w=(q,[\beta,j])$ for each $p,q=1,2$.
\indent\par From (i), $d_D((p,[\alpha,i]),(3,[j]))=3$. Since $(3,[j])\rightarrow \{(1,[\beta,j]),(2,[\beta,j])\}$ by (\ref{eqC6.3.36}), this subcase follows.
\\
\\(vi) $v=(q,[\beta,j]), w=(p,[\alpha,i])$ for each $p,q=1,2$.
\indent\par From (iv), $d_D((q,[\beta,j]),(p,[i]))=3$. Since $(p,[i])\rightarrow (3-p,[\alpha,i])$ by (\ref{eqC6.3.34}), this subcase follows.
\\
\\Case 5. $v \in \{(1,[\alpha,i]),(2,[\alpha,i]),(1,[i]),(2,[i])\}$ and $w\in \{(1,[\beta,j]),(2,[\beta,j]),(1,[j]),$ $(2,[j]),(3,[j])\}$ for each $1\le i\le |A_2|$, each ${{s}\choose{s/2}}-1\le j\le |A_3|$, each $1\le\alpha\le \deg_T([i])-1$, and each $1\le\beta\le \deg_T([j])-1$.
\\(i) $v=(p,[\alpha,i]), w=(q,[j])$ for each $p=1,2$ and $q=1,2,3$.
\indent\par By (\ref{eqC6.3.34}) and (\ref{eqC6.3.40}), $(p,[\alpha,i])\rightarrow (p,[i])$, $\mu_1 \rightarrow (1,[j])$ and $\mu_{{s}\choose{s/2}}\rightarrow (2,[j])$. With (\ref{eqC6.3.35}), $O^\mathtt{c}((p,[i]))\in {{(\mathbb{N}_s,\mathtt{c})}\choose{s/2}}-\{\mu_1,\mu_{{s}\choose{s/2}}\}$, there exist some $x_i\in O^\mathtt{c}((p,[i]))\cap\mu_1$ and $y_i\in O^\mathtt{c}((p,[i]))\cap\mu_{{s}\choose{s/2}}$ such that $(p,[i])\rightarrow \{x_i, y_i\}$. Furthermore, with (\ref{eqC6.3.41}), $O^\mathtt{c}((p,[i]))=\bar{\mu}_{i+1}=(\bar{X}_{i+1},\mathtt{c})$, where $\bar{X}_{i+1}\in L_{n,\frac{n}{2}}(k)$ and $O^\mathtt{c}((3,[j]))=\psi_{j+2-{{s}\choose{s/2}}}=(Y,\mathtt{c})$ for some $Y\in{{\mathbb{N}_s}\choose{(s/2)+1}}-\nabla L_{n,\frac{n}{2}}(k)$ so that  $\bar{\mu}_{i+1}\not\subseteq \psi_{j+2-{{s}\choose{s/2}}}$. That is, there exists some $z_{ij}\in O^\mathtt{c}((p,[i])) \cap I^\mathtt{c}((3,[j]))$ such that $(p,[i])\rightarrow z_{ij}\rightarrow (3,[j])$.
\\
\\(ii) $v=(q,[j]), w=(p,[\alpha,i])$ for each $p=1,2$ and $q=1,2,3$.
\indent\par By (\ref{eqC6.3.34}) and  (\ref{eqC6.3.40}),  $(p,[i])\rightarrow (3-p,[\alpha,i])$, $(1,[j])\rightarrow \mu_{{s}\choose{s/2}}$, and $(2,[j])\rightarrow \mu_1$. With (\ref{eqC6.3.35}), $I^\mathtt{c}((p,[i]))\in {{(\mathbb{N}_s,\mathtt{c})}\choose{s/2}}-\{\mu_1,\mu_{{s}\choose{s/2}}\}$ and there exist some $x_i\in \mu_1\cap I^\mathtt{c}((p,[i]))$ and $y_i\in \mu_{{s}\choose{s/2}}\cap I^\mathtt{c}((p,[i]))$ such that $\{x_i, y_i\}\rightarrow (p,[i])$. Furthermore, with (\ref{eqC6.3.41}), $|O^\mathtt{c}((3,[j]))|+|I^\mathtt{c}((p,[i]))|=(\frac{s}{2}+1)+\frac{s}{2}>s$. That is, there exists some $z_{ij}\in O^\mathtt{c}((3,[j])) \cap I^\mathtt{c}((p,[i]))$ such that $(3,[j])\rightarrow z_{ij} \rightarrow (p,[i])$.
\\
\\(iii) $v=(p,[i]), w=(q,[\beta,j])$ for each $p=1,2$ and $q=1,2$.
\indent\par From (i), $d_D((p,[i]),(1,[j]))=2$. Since $(1,[j])\rightarrow \{(1,[\beta,j]), (2,[\beta,j])\}$ by (\ref{eqC6.3.39}), this subcase follows.
\\
\\(iv) $v=(q,[\beta,j]), w=(p,[i])$ for each $p=1,2$ and $q=1,2$.
\indent\par From (ii), $d_D((3,[j]),(p,[i]))=2$. Since $\{(1,[\beta,j]), (2,[\beta,j])\} \rightarrow (3,[j])$ by (\ref{eqC6.3.39}), this subcase follows.
\\
\\(v) $v=(p,[\alpha,i]), w=(q,[\beta,j])$ for each $p,q=1,2$.
\indent\par From (i), $d_D((p,[\alpha,i]),(1,[j]))=3$. Since $(1,[j])\rightarrow \{(1,[\beta,j]), (2,[\beta,j])\}$ by (\ref{eqC6.3.39}), this subcase follows.
\\
\\(vi) $v=(q,[\beta,j]), w=(p,[\alpha,i])$ for each $p,q=1,2$.
\indent\par From (iv), $d_D((q,[\beta,j]),(p,[i]))=3$. Since $(3-p,[i])\rightarrow (p,[\alpha,i])$ by (\ref{eqC6.3.34}), this subcase follows.
\\
\\Case 6. $v \in \{(1,[\alpha,i]),(2,[\alpha,i]),(1,[i]),(2,[i])\}$ and $w\in \{(1,[\beta,j]),(2,[\beta,j]),(1,[j]),$ $(2,[j]),(3,[j]),(4,[j])\}$ for each $[i]\in A_2$, each $[j]\in A_4$, each $1\le\alpha\le \deg_T([i])-1$, and each $1\le\beta\le \deg_T([j])-1$.
\indent\par Since the orientation defined for $A_2$ and $A_4$ (see (\ref{eqC6.3.34})-(\ref{eqC6.3.35}), (\ref{eqC6.3.42})-(\ref{eqC6.3.43})) are similar to those in Proposition \ref{ppnC6.3.10} (see (\ref{eqC6.3.27})-(\ref{eqC6.3.31})), this case follows from Case 4 of Proposition \ref{ppnC6.3.10}, except possibly for $v=(p,[i])$ and $w=(q, [j])$, which will be covered in Case 14 later.
\\
\\Case 7. $v \in \{(1,[\alpha,i]),(2,[\alpha,i]),(1,[i]),(2,[i]),(3,[i])\}$ and $w\in \{(1,[\beta,j]),(2,[\beta,j]),(1,[j]),$ $(2,[j]),(3,[j]),(4,[j])\}$ for each $[i]\in A_3$, each $[j]\in A_4$, each $1\le\alpha\le \deg_T([i])-1$, and each $1\le\beta\le \deg_T([j])-1$.
\indent\par Since the orientation defined for $A_3$ and $A_4$ (see (\ref{eqC6.3.36})-(\ref{eqC6.3.43})) are similar to those in Proposition \ref{ppnC6.3.9} (see (\ref{eqC6.3.17})-(\ref{eqC6.3.25})), this case follows from Cases 6 and 7 of Proposition \ref{ppnC6.3.9}, except possibly for $v=(p,[i])$ and $w=(q, [j])$, which will be covered in Case 14 later.
\\
\\Case 8. For each $1\le i\le |A_2|$, each $1\le\alpha\le \deg_T([i])-1$, and each $[j]\in E$,
\\(i) $v=(p,[\alpha,i]), w=(q,[j])$ for each $p=1,2$ and $q=1,2$.
\\(ii) $v=(q,[j]), w=(p,[\alpha,i])$ for each $p=1,2$ and $q=1,2$.
\indent\par Let $[x]\in A_4$. By (\ref{eqC6.3.43})-(\ref{eqC6.3.44}), $\mu_1\rightarrow \{(q,[j]),(2,[x])\}\rightarrow  \mu_{{s}\choose{s/2}}$, this case follows from Case 6.
\\
\\Case 9. For each $|A_2|+1\le i\le |A_2|+|A_3|$, each $1\le\alpha\le \deg_T([i])-1$, and each $[j]\in E$,
\\(i) $v=(p,[\alpha,i]), w=(q,[j])$ for each $p=1,2$ and $q=1,2$.
\\(ii) $v=(q,[j]), w=(p,[\alpha,i])$ for each $p=1,2$ and $q=1,2$.
\indent\par Let $[x]\in A_4$. By (\ref{eqC6.3.43})-(\ref{eqC6.3.44}), $\mu_1\rightarrow \{(q,[j]),(2,[x])\}\rightarrow  \mu_{{s}\choose{s/2}}$, this case follows from  Case 7.
\\
\\Case 10. For each $[i]\in A_4$, each $1\le\alpha\le \deg_T([i])-1$, and each $[j]\in E$,
\\(i) $v=(p,[\alpha,i]), w=(q,[j])$ for each $p=1,2$ and $q=1,2$.
\\(ii) $v=(q,[j]), w=(p,[\alpha,i])$ for each $p=1,2$ and $q=1,2$.
\indent\par Let $x=|A_2|+1$. By (\ref{eqC6.3.37}) and (\ref{eqC6.3.44}), $\mu_1\rightarrow \{(q,[j]),(2,[x])\}\rightarrow  \mu_{{s}\choose{s/2}}$, this case follows from Case 7.
\\
\\Case 11. $v=(r_1,\mathtt{c})$ and $w=(r_2,\mathtt{c})$ for $r_1\neq r_2$ and $1\le r_1, r_2\le s$.
\indent\par Here, we want to prove a stronger claim, $d_{D}((r_1,\mathtt{c}), (r_2,\mathtt{c}))=2$. Let $V_2=\{(1,[|A_2|+1]), (2,[|A_2|+1])\}\cup\{(1,[i])\mid 1\le i\le |A_2|\}\cup \{(3,[j])\mid |A_2|+1\le j\le {{s}\choose{s/2}}-2\}$. Then, $\{I^\mathtt{c}(v)\mid v\in V_2\}={{\mathbb{N}_s}\choose{s/2}}$ and the subgraph induced by $V_1=(\mathbb{N}_s,\mathtt{c})$ and $V_2$ is a complete bipartite graph $K(V_1,V_2)$. By Lemma \ref{lemC6.2.19}, $d_{D}((r_1,\mathtt{c}), (r_2,\mathtt{c}))=2$.
\\
\\Case 12. $v\in \{(1,[i]), (2,[i]), (3,[i]), (4,[i]), (1,[\alpha,i]), (2,[\alpha,i])\}$ for each $1\le i\le \deg_T(\mathtt{c})$ and $1\le\alpha\le \deg_T([i])-1$, and $w=(r,\mathtt{c})$ for $1\le r\le s$.
\indent\par Note that there exists some $1\le l\le s$ such that $d_D(v,(l,\mathtt{c}))\le 2$, and $d_D((l,\mathtt{c}),w)\le 2$ by Case 11. Hence, it follows that $d_D(v,w)\le d_D(v,(l,\mathtt{c}))+d_D((l,\mathtt{c}),w)\le 4$.
\\
\\Case 13. $v=(r,\mathtt{c})$ for $1\le r\le s$ and $w\in \{(1,[i]), (2,[i]), (3,[i]), (4,[i]), (1,[\alpha,i]), (2,[\alpha,i])\}$ for each $1\le i\le \deg_T(\mathtt{c})$ and $1\le\alpha\le \deg_T([i])-1$.
\indent\par Note that there exists some $1\le l\le s$ such that $d_D((l,\mathtt{c}), w)\le 2$, and $d_D(v,(l,\mathtt{c}))\le 2$ by Case 11. Hence, it follows that $d_D(v,w)\le d_D(v,(l,\mathtt{c}))+d_D((l,\mathtt{c}),w)\le 4$.
\\
\\Case 14. $v=(p,[i])$ and $w=(q, [j])$, where $1\le p,q\le 4$ and $1\le i,j\le \deg_T(\mathtt{c})$.
\noindent\par This follows from the fact that $|O^\mathtt{c}((p,[i]))|>0$, $|I^\mathtt{c}((q,[j]))|>0$, and $d_{D}((r_1,\mathtt{c}), (r_2,\mathtt{c}))$ $=2$ for any $r_1\neq r_2$ and $1\le r_1, r_2\le s$.
\\
\indent\par Therefore, the claim follows. Since every vertex lies in a directed $C_4$ for $D$ and $d(D)=4$, $\bar{d}(\mathcal{T})\le \max \{4, d(D)\}$ by Lemma \ref{lemC6.1.3} and thus $\bar{d}(\mathcal{T})=4$ .
\qed
\begin{center}
\begin{tikzpicture}[thick,scale=0.8]%
\draw(0,6)node[circle, draw, inner sep=0pt, minimum width=3pt](1_u){\scriptsize $(1,\mathtt{c})$};
\draw(0,3.6)node[circle, draw, inner sep=0pt, minimum width=3pt](2_u){\scriptsize $(2,\mathtt{c})$};
\draw(0,1.2)node[circle, draw, inner sep=0pt, minimum width=3pt](3_u){\scriptsize $(3,\mathtt{c})$};
\draw(0,-1.2)node[circle, draw, inner sep=0pt, minimum width=3pt](4_u){\scriptsize $(4,\mathtt{c})$};
\draw(0,-3.6)node[circle, draw, inner sep=0pt, minimum width=3pt](5_u){\scriptsize $(5,\mathtt{c})$};
\draw(0,-6)node[circle, draw, inner sep=0pt, minimum width=3pt](6_u){\scriptsize $(6,\mathtt{c})$};

\draw(-6,6)node[circle, draw, fill=black!100, inner sep=0pt, minimum width=6pt, label={[] 180:{\small $(1,[1,1])$}}](1_11){};
\draw(-6,4)node[circle, draw, fill=black!100, inner sep=0pt, minimum width=6pt, label={[] 180:{\small $(2,[1,1])$}}](2_11){};

\draw(-3,6)node[circle, draw, fill=black!100, inner sep=0pt, minimum width=6pt, label={[yshift=0cm, xshift=-0.45cm] 90:{\small $(1,[1])$}}](1_1){};
\draw(-3,4)node[circle, draw, fill=black!100, inner sep=0pt, minimum width=6pt, label={[yshift=0cm, xshift=-0.45cm] 270:{\small $(2,[1])$}}](2_1){};

\draw(6,6)node[circle, draw, fill=black!100, inner sep=0pt, minimum width=6pt, label={[] 0:{\small $(1,[1,2])$}}](1_21){};
\draw(6,4)node[circle, draw, fill=black!100, inner sep=0pt, minimum width=6pt, label={[] 0:{\small $(2,[1,2])$}}](2_21){};

\draw(3,6)node[circle, draw, fill=black!100, inner sep=0pt, minimum width=6pt, label={[yshift=0cm, xshift=0.45cm] 90:{\small $(1,[2])$}}](1_2){};
\draw(3,4)node[circle, draw, fill=black!100, inner sep=0pt, minimum width=6pt, label={[yshift=0cm, xshift=0.45cm] 270:{\small $(2,[2])$}}](2_2){};

\draw(-6,2)node[circle, draw, fill=black!100, inner sep=0pt, minimum width=6pt, label={[] 180:{\small $(1,[1,11])$}}](1_31){};
\draw(-6,0)node[circle, draw, fill=black!100, inner sep=0pt, minimum width=6pt, label={[] 180:{\small $(2,[1,11])$}}](2_31){};

\draw(-3,2)node[circle, draw, fill=black!100, inner sep=0pt, minimum width=6pt, label={[yshift=0cm, xshift=-0.45cm] 270:{\small $(1,[11])$}}](1_3){};
\draw(-3,0)node[circle, draw, fill=black!100, inner sep=0pt, minimum width=6pt, label={[yshift=0cm, xshift=-0.45cm] 270:{\small $(2,[11])$}}](2_3){};

\draw(6,2)node[circle, draw, fill=black!100, inner sep=0pt, minimum width=6pt, label={[] 0:{\small $(1,[1,12])$}}](1_41){};
\draw(6,0)node[circle, draw, fill=black!100, inner sep=0pt, minimum width=6pt, label={[] 0:{\small $(2,[1,12])$}}](2_41){};

\draw(3,2)node[circle, draw, fill=black!100, inner sep=0pt, minimum width=6pt, label={[yshift=0cm, xshift=0.45cm] 270:{\small $(1,[12])$}}](1_4){};
\draw(3,0)node[circle, draw, fill=black!100, inner sep=0pt, minimum width=6pt, label={[yshift=0cm, xshift=0.45cm] 270:{\small $(2,[12])$}}](2_4){};

\draw(-6,-3)node[circle, draw, fill=black!100, inner sep=0pt, minimum width=6pt, label={[] 180:{\small $(1,[1,13])$}}](1_51){};
\draw(-6,-5)node[circle, draw, fill=black!100, inner sep=0pt, minimum width=6pt, label={[] 180:{\small $(2,[1,13])$}}](2_51){};

\draw(-3,-2)node[circle, draw, fill=black!100, inner sep=0pt, minimum width=6pt, label={[yshift=-0.2cm, xshift=-0.45cm] 270:{\small $(1,[13])$}}](1_5){};
\draw(-3,-4)node[circle, draw, fill=black!100, inner sep=0pt, minimum width=6pt, label={[yshift=-0.1cm, xshift=-0.45cm] 270:{\small $(2,[13])$}}](2_5){};
\draw(-3,-6)node[circle, draw, fill=black!100, inner sep=0pt, minimum width=6pt, label={[yshift=0cm, xshift=-0.45cm] 270:{\small $(3,[13])$}}](3_5){};

\draw(6,-3)node[circle, draw, fill=black!100, inner sep=0pt, minimum width=6pt, label={[] 0:{\small $(1,[1,14])$}}](1_61){};
\draw(6,-5)node[circle, draw, fill=black!100, inner sep=0pt, minimum width=6pt, label={[] 0:{\small $(2,[1,14])$}}](2_61){};

\draw(3,-2)node[circle, draw, fill=black!100, inner sep=0pt, minimum width=6pt, label={[yshift=-0.2cm, xshift=0.45cm] 270:{\small $(1,[14])$}}](1_6){};
\draw(3,-4)node[circle, draw, fill=black!100, inner sep=0pt, minimum width=6pt, label={[yshift=-0.1cm, xshift=0.45cm] 270:{\small $(2,[14])$}}](2_6){};
\draw(3,-6)node[circle, draw, fill=black!100, inner sep=0pt, minimum width=6pt, label={[yshift=0cm, xshift=0.45cm] 270:{\small $(3,[14])$}}](3_6){};

\draw[->, line width=0.3mm, >=latex, shorten <= 0.2cm, shorten >= 0.15cm](1_11)--(1_1);
\draw[->, line width=0.3mm, >=latex, shorten <= 0.2cm, shorten >= 0.15cm](2_11)--(2_1);
\draw[dashed, ->, line width=0.3mm, >=latex, shorten <= 0.2cm, shorten >= 0.15cm](2_1)--(1_11);
\draw[dashed, ->, line width=0.3mm, >=latex, shorten <= 0.2cm, shorten >= 0.15cm](1_1)--(2_11);

\draw[dashed, ->, line width=0.3mm, >=latex, shorten <= 0.2cm, shorten >= 0.15cm](2_2)--(1_21);
\draw[dashed, ->, line width=0.3mm, >=latex, shorten <= 0.2cm, shorten >= 0.15cm](1_2)--(2_21);
\draw[->, line width=0.3mm, >=latex, shorten <= 0.2cm, shorten >= 0.15cm](1_21)--(1_2);
\draw[->, line width=0.3mm, >=latex, shorten <= 0.2cm, shorten >= 0.15cm](2_21)--(2_2);

\draw[dashed, ->, line width=0.3mm, >=latex, shorten <= 0.2cm, shorten >= 0.15cm](2_3)--(1_31);
\draw[dashed, ->, line width=0.3mm, >=latex, shorten <= 0.2cm, shorten >= 0.15cm](1_3) to [out=190, in=65] (2_31);
\draw[->, line width=0.3mm, >=latex, shorten <= 0.2cm, shorten >= 0.15cm](1_31)--(1_3);
\draw[->, line width=0.3mm, >=latex, shorten <= 0.2cm, shorten >= 0.15cm](2_31)--(2_3);

\draw[dashed, ->, line width=0.3mm, >=latex, shorten <= 0.2cm, shorten >= 0.15cm](2_4)--(1_41);
\draw[dashed, ->, line width=0.3mm, >=latex, shorten <= 0.2cm, shorten >= 0.15cm](1_4) to [out=350, in=115] (2_41);
\draw[->, line width=0.3mm, >=latex, shorten <= 0.2cm, shorten >= 0.15cm](1_41)--(1_4);
\draw[->, line width=0.3mm, >=latex, shorten <= 0.2cm, shorten >= 0.15cm](2_41)--(2_4);

\draw[dashed, ->, line width=0.3mm, >=latex, shorten <= 0.2cm, shorten >= 0.15cm](3_6) to [out=22, in=265] (1_61);
\draw[dashed, ->, line width=0.3mm, >=latex, shorten <= 0.2cm, shorten >= 0.15cm](3_6)--(2_61);
\draw[->, line width=0.3mm, >=latex, shorten <= 0.2cm, shorten >= 0.15cm](1_61)--(1_6);
\draw[->, line width=0.3mm, >=latex, shorten <= 0.2cm, shorten >= 0.15cm](2_61) to [out=95, in=338] (1_6);
\draw[->, line width=0.3mm, >=latex, shorten <= 0.2cm, shorten >= 0.15cm](1_61)--(2_6);
\draw[->, line width=0.3mm, >=latex, shorten <= 0.2cm, shorten >= 0.15cm](2_61)--(2_6);

\draw[dashed, ->, line width=0.3mm, >=latex, shorten <= 0.2cm, shorten >= 0.15cm](3_5) to [out=158, in=275] (1_51);
\draw[dashed, ->, line width=0.3mm, >=latex, shorten <= 0.2cm, shorten >= 0.15cm](3_5)--(2_51);
\draw[->, line width=0.3mm, >=latex, shorten <= 0.2cm, shorten >= 0.15cm](1_51)--(1_5);
\draw[->, line width=0.3mm, >=latex, shorten <= 0.2cm, shorten >= 0.15cm](2_51) to [out=85, in=202] (1_5);
\draw[->, line width=0.3mm, >=latex, shorten <= 0.2cm, shorten >= 0.15cm](1_51)--(2_5);
\draw[->, line width=0.3mm, >=latex, shorten <= 0.2cm, shorten >= 0.15cm](2_51)--(2_5);

\draw[densely dotted, ->, line width=0.3mm, >=latex, shorten <= 0.2cm, shorten >= 0.1cm](2_1)--(3_u);
\draw[densely dotted, ->, line width=0.3mm, >=latex, shorten <= 0.2cm, shorten >= 0.1cm](2_1)--(5_u);
\draw[densely dotted, ->, line width=0.3mm, >=latex, shorten <= 0.2cm, shorten >= 0.1cm](2_1)--(6_u);
\draw[densely dotted, ->, line width=0.3mm, >=latex, shorten <= 0.2cm, shorten >= 0.1cm](1_1)--(3_u);
\draw[densely dotted, ->, line width=0.3mm, >=latex, shorten <= 0.2cm, shorten >= 0.1cm](1_1)--(5_u);
\draw[densely dotted, ->, line width=0.3mm, >=latex, shorten <= 0.2cm, shorten >= 0.1cm](1_1)--(6_u);

\draw[densely dotted, ->, line width=0.3mm, >=latex, shorten <= 0.2cm, shorten >= 0.1cm](2_2)--(2_u);
\draw[densely dotted, ->, line width=0.3mm, >=latex, shorten <= 0.2cm, shorten >= 0.1cm](2_2)--(5_u);
\draw[densely dotted, ->, line width=0.3mm, >=latex, shorten <= 0.2cm, shorten >= 0.1cm](2_2)--(6_u);
\draw[densely dotted, ->, line width=0.3mm, >=latex, shorten <= 0.2cm, shorten >= 0.1cm](1_2)--(2_u);
\draw[densely dotted, ->, line width=0.3mm, >=latex, shorten <= 0.2cm, shorten >= 0.1cm](1_2)--(5_u);
\draw[densely dotted, ->, line width=0.3mm, >=latex, shorten <= 0.2cm, shorten >= 0.1cm](1_2)--(6_u);

\draw[densely dotted, ->, line width=0.3mm, >=latex, shorten <= 0.2cm, shorten >= 0.1cm](2_3)--(2_u);
\draw[densely dotted, ->, line width=0.3mm, >=latex, shorten <= 0.2cm, shorten >= 0.1cm](2_3)--(4_u);
\draw[densely dotted, ->, line width=0.3mm, >=latex, shorten <= 0.2cm, shorten >= 0.1cm](2_3)--(5_u);
\draw[densely dotted, ->, line width=0.3mm, >=latex, shorten <= 0.2cm, shorten >= 0.1cm](1_3)--(2_u);
\draw[densely dotted, ->, line width=0.3mm, >=latex, shorten <= 0.2cm, shorten >= 0.1cm](1_3)--(4_u);
\draw[densely dotted, ->, line width=0.3mm, >=latex, shorten <= 0.2cm, shorten >= 0.1cm](1_3)--(5_u);

\draw[densely dotted, ->, line width=0.3mm, >=latex, shorten <= 0.2cm, shorten >= 0.1cm](2_4)--(1_u);
\draw[densely dotted, ->, line width=0.3mm, >=latex, shorten <= 0.2cm, shorten >= 0.1cm](2_4)--(4_u);
\draw[densely dotted, ->, line width=0.3mm, >=latex, shorten <= 0.2cm, shorten >= 0.1cm](2_4)--(5_u);
\draw[densely dotted, ->, line width=0.3mm, >=latex, shorten <= 0.2cm, shorten >= 0.1cm](1_4)--(1_u);
\draw[densely dotted, ->, line width=0.3mm, >=latex, shorten <= 0.2cm, shorten >= 0.1cm](1_4)--(4_u);
\draw[densely dotted, ->, line width=0.3mm, >=latex, shorten <= 0.2cm, shorten >= 0.1cm](1_4)--(5_u);

\draw[->, line width=0.3mm, >=latex, shorten <= 0.2cm, shorten >= 0.1cm](1_5)--(1_u);
\draw[->, line width=0.3mm, >=latex, shorten <= 0.2cm, shorten >= 0.1cm](1_5)--(2_u);
\draw[->, line width=0.3mm, >=latex, shorten <= 0.2cm, shorten >= 0.1cm](1_5)--(3_u);
\draw[->, line width=0.3mm, >=latex, shorten <= 0.2cm, shorten >= 0.1cm](2_5)--(4_u);
\draw[->, line width=0.3mm, >=latex, shorten <= 0.2cm, shorten >= 0.1cm](2_5)--(5_u);
\draw[->, line width=0.3mm, >=latex, shorten <= 0.2cm, shorten >= 0.1cm](2_5)--(6_u);
\draw[densely dashdotdotted, ->, line width=0.3mm, >=latex, shorten <= 0.2cm, shorten >= 0.1cm](3_5)--(2_u);
\draw[densely dashdotdotted, ->, line width=0.3mm, >=latex, shorten <= 0.2cm, shorten >= 0.1cm](3_5)--(3_u);
\draw[densely dashdotdotted, ->, line width=0.3mm, >=latex, shorten <= 0.2cm, shorten >= 0.1cm](3_5)--(5_u);

\draw[->, line width=0.3mm, >=latex, shorten <= 0.2cm, shorten >= 0.1cm](1_6)--(1_u);
\draw[->, line width=0.3mm, >=latex, shorten <= 0.2cm, shorten >= 0.1cm](1_6)--(2_u);
\draw[->, line width=0.3mm, >=latex, shorten <= 0.2cm, shorten >= 0.1cm](1_6)--(3_u);
\draw[->, line width=0.3mm, >=latex, shorten <= 0.2cm, shorten >= 0.1cm](2_6)--(4_u);
\draw[->, line width=0.3mm, >=latex, shorten <= 0.2cm, shorten >= 0.1cm](2_6)--(5_u);
\draw[->, line width=0.3mm, >=latex, shorten <= 0.2cm, shorten >= 0.1cm](2_6)--(6_u);
\draw[densely dashdotdotted, ->, line width=0.3mm, >=latex, shorten <= 0.2cm, shorten >= 0.1cm](3_6)--(1_u);
\draw[densely dashdotdotted, ->, line width=0.3mm, >=latex, shorten <= 0.2cm, shorten >= 0.1cm](3_6)--(3_u);
\draw[densely dashdotdotted, ->, line width=0.3mm, >=latex, shorten <= 0.2cm, shorten >= 0.1cm](3_6)--(5_u);
\end{tikzpicture}
\captionsetup{justification=centering}
\captionof{figure}{Partial orientation $D$ for $\mathcal{H}$ for $s=6$, $k=13$;\\
 $A_2=\{[1],[2],\ldots,[12] \}$, $A_3=\{[13],[14],\ldots, [20]\}$,  $A_4=\{[21],[22]\}$, $E=\{[23],[24]\}$.\label{figC6.3.10}}
\captionsetup{justification=justified}
\caption*{Note that $\bm{123} <_s \bm{124} <_s \bm{134} <_s \bm{234} <_s \bm{125} <_s \bm{135} <_s \bm{235} <_s \bm{145} <_s \bm{245} <_s \bm{345}<_s \bm{126} <_s \bm{136} <_s \bm{236} <_s \bm{146} <_s \bm{246} <_s \bm{346} <_s \bm{156} <_s \bm{256} <_s \bm{356} <_s \bm{456}$. So, $L_{6,3}(13)=\{\bm{145}, \bm{245},\ldots, \bm{456}\}$, ${{\mathbb{N}_6}\choose{4}}-\nabla L_{6,3}(13)=\{\bm{1234},\bm{1235}\}$ and $\kappa_{6,3}(13)=0$.}
\end{center}

\begin{center}
\begin{tikzpicture}[thick,scale=0.8]%
\draw(0,6)node[circle, draw, inner sep=0pt, minimum width=3pt](1_u){\scriptsize $(1,\mathtt{c})$};
\draw(0,2.8)node[circle, draw, inner sep=0pt, minimum width=3pt](2_u){\scriptsize $(2,\mathtt{c})$};
\draw(0,-0.4)node[circle, draw, inner sep=0pt, minimum width=3pt](3_u){\scriptsize $(3,\mathtt{c})$};
\draw(0,-3.6)node[circle, draw, inner sep=0pt, minimum width=3pt](4_u){\scriptsize $(4,\mathtt{c})$};
\draw(0,-6.8)node[circle, draw, inner sep=0pt, minimum width=3pt](5_u){\scriptsize $(5,\mathtt{c})$};
\draw(0,-10)node[circle, draw, inner sep=0pt, minimum width=3pt](6_u){\scriptsize $(6,\mathtt{c})$};

\draw(-6,5)node[circle, draw, fill=black!100, inner sep=0pt, minimum width=6pt, label={[] 180:{\small $(1,[1,19])$}}](1_119){};
\draw(-6,3)node[circle, draw, fill=black!100, inner sep=0pt, minimum width=6pt, label={[] 180:{\small $(2,[1,19])$}}](2_119){};

\draw(-3,6)node[circle, draw, fill=black!100, inner sep=0pt, minimum width=6pt, label={[yshift=0cm, xshift=-0.5cm] 90:{\small $(1,[19])$}}](1_19){};
\draw(-3,4)node[circle, draw, fill=black!100, inner sep=0pt, minimum width=6pt, label={[yshift=-0.15cm, xshift=-0.5cm] 270:{\small $(2,[19])$}}](2_19){};
\draw(-3,2)node[circle, draw, fill=black!100, inner sep=0pt, minimum width=6pt, label={[yshift=0cm, xshift=-0.5cm] 270:{\small $(3,[19])$}}](3_19){};

\draw(6,5)node[circle, draw, fill=black!100, inner sep=0pt, minimum width=6pt, label={[] 0:{\small $(1,[1,20])$}}](1_120){};
\draw(6,3)node[circle, draw, fill=black!100, inner sep=0pt, minimum width=6pt, label={[] 0:{\small $(2,[1,20])$}}](2_120){};

\draw(3,6)node[circle, draw, fill=black!100, inner sep=0pt, minimum width=6pt, label={[yshift=0cm, xshift=0.5cm] 90:{\small $(1,[20])$}}](1_20){};
\draw(3,4)node[circle, draw, fill=black!100, inner sep=0pt, minimum width=6pt, label={[yshift=-0.15cm, xshift=0.5cm] 270:{\small $(2,[20])$}}](2_20){};
\draw(3,2)node[circle, draw, fill=black!100, inner sep=0pt, minimum width=6pt, label={[yshift=0cm, xshift=0.5cm] 270:{\small $(3,[20])$}}](3_20){};

\draw(-6,-2)node[circle, draw, fill=black!100, inner sep=0pt, minimum width=6pt, label={[] 180:{\small $(1,[1,21])$}}](1_121){};
\draw(-6,-4)node[circle, draw, fill=black!100, inner sep=0pt, minimum width=6pt, label={[] 180:{\small $(2,[1,21])$}}](2_121){};

\draw(-3,0)node[circle, draw, fill=black!100, inner sep=0pt, minimum width=6pt, label={[yshift=-0.15cm, xshift=-0.5cm] 270:{\small $(1,[21])$}}](1_21){};
\draw(-3,-2)node[circle, draw, fill=black!100, inner sep=0pt, minimum width=6pt, label={[yshift=-0.05cm, xshift=-0.5cm] 270:{\small $(2,[21])$}}](2_21){};
\draw(-3,-4)node[circle, draw, fill=black!100, inner sep=0pt, minimum width=6pt, label={[yshift=0.05cm, xshift=-0.5cm] 270:{\small $(3,[21])$}}](3_21){};
\draw(-3,-6)node[circle, draw, fill=black!100, inner sep=0pt, minimum width=6pt, label={[yshift=0cm, xshift=-0.5cm] 270:{\small $(4,[21])$}}](4_21){};

\draw(6,-2)node[circle, draw, fill=black!100, inner sep=0pt, minimum width=6pt, label={[] 0:{\small $(1,[1,22])$}}](1_122){};
\draw(6,-4)node[circle, draw, fill=black!100, inner sep=0pt, minimum width=6pt, label={[] 0:{\small $(2,[1,22])$}}](2_122){};

\draw(3,0)node[circle, draw, fill=black!100, inner sep=0pt, minimum width=6pt, label={[yshift=-0.15cm, xshift=0.5cm] 270:{\small $(1,[22])$}}](1_22){};
\draw(3,-2)node[circle, draw, fill=black!100, inner sep=0pt, minimum width=6pt, label={[yshift=-0.05cm, xshift=0.5cm] 270:{\small $(2,[22])$}}](2_22){};
\draw(3,-4)node[circle, draw, fill=black!100, inner sep=0pt, minimum width=6pt, label={[yshift=0.05cm, xshift=0.5cm] 270:{\small $(3,[22])$}}](3_22){};
\draw(3,-6)node[circle, draw, fill=black!100, inner sep=0pt, minimum width=6pt, label={[yshift=0cm, xshift=0.5cm] 270:{\small $(4,[22])$}}](4_22){};

\draw(-3,-8)node[circle, draw, fill=black!100, inner sep=0pt, minimum width=6pt, label={[] 180:{\small $(1,[23])$}}](1_23){};
\draw(-3,-10)node[circle, draw, fill=black!100, inner sep=0pt, minimum width=6pt, label={[] 180:{\small $(2,[23])$}}](2_23){};

\draw(3,-8)node[circle, draw, fill=black!100, inner sep=0pt, minimum width=6pt, label={[] 0:{\small $(1,[24])$}}](1_24){};
\draw(3,-10)node[circle, draw, fill=black!100, inner sep=0pt, minimum width=6pt, label={[] 0:{\small $(2,[24])$}}](2_24){};

\draw[dashed, ->, line width=0.3mm, >=latex, shorten <= 0.2cm, shorten >= 0.15cm](1_20)--(1_120);
\draw[dashed, ->, line width=0.3mm, >=latex, shorten <= 0.2cm, shorten >= 0.15cm](1_20)--(2_120);
\draw[dashed, ->, line width=0.3mm, >=latex, shorten <= 0.2cm, shorten >= 0.15cm](2_20)--(1_120);
\draw[dashed, ->, line width=0.3mm, >=latex, shorten <= 0.2cm, shorten >= 0.15cm](2_20)--(2_120);
\draw[->, line width=0.3mm, >=latex, shorten <= 0.2cm, shorten >= 0.15cm](1_120) to [out=265, in=22] (3_20);
\draw[->, line width=0.3mm, >=latex, shorten <= 0.2cm, shorten >= 0.15cm](2_120)--(3_20);

\draw[dashed, ->, line width=0.3mm, >=latex, shorten <= 0.2cm, shorten >= 0.15cm](1_19)--(1_119);
\draw[dashed, ->, line width=0.3mm, >=latex, shorten <= 0.2cm, shorten >= 0.15cm](1_19)--(2_119);
\draw[dashed, ->, line width=0.3mm, >=latex, shorten <= 0.2cm, shorten >= 0.15cm](2_19)--(1_119);
\draw[dashed, ->, line width=0.3mm, >=latex, shorten <= 0.2cm, shorten >= 0.15cm](2_19)--(2_119);
\draw[->, line width=0.3mm, >=latex, shorten <= 0.2cm, shorten >= 0.15cm](1_119) to [out=275, in=158] (3_19);
\draw[->, line width=0.3mm, >=latex, shorten <= 0.2cm, shorten >= 0.15cm](2_119)--(3_19);

\draw[dashed, ->, line width=0.3mm, >=latex, shorten <= 0.2cm, shorten >= 0.15cm](1_22) to [out=345, in=100] (2_122);
\draw[dashed, ->, line width=0.3mm, >=latex, shorten <= 0.2cm, shorten >= 0.15cm](2_22)--(1_122);
\draw[dashed, ->, line width=0.3mm, >=latex, shorten <= 0.2cm, shorten >= 0.15cm](3_22)--(2_122);
\draw[dashed, ->, line width=0.3mm, >=latex, shorten <= 0.2cm, shorten >= 0.15cm](4_22) to [out=40, in=250] (1_122);
\draw[->, line width=0.3mm, >=latex, shorten <= 0.2cm, shorten >= 0.15cm](1_122) to [out=115, in=0] (1_22);
\draw[->, line width=0.3mm, >=latex, shorten <= 0.2cm, shorten >= 0.15cm](2_122) to [out=115, in=350] (2_22);
\draw[->, line width=0.3mm, >=latex, shorten <= 0.2cm, shorten >= 0.15cm](1_122)--(3_22);
\draw[->, line width=0.3mm, >=latex, shorten <= 0.2cm, shorten >= 0.15cm](2_122)--(4_22);

\draw[dashed, ->, line width=0.3mm, >=latex, shorten <= 0.2cm, shorten >= 0.15cm](1_21) to [out=195, in=85] (2_121);
\draw[dashed, ->, line width=0.3mm, >=latex, shorten <= 0.2cm, shorten >= 0.15cm](2_21)--(1_121);
\draw[dashed, ->, line width=0.3mm, >=latex, shorten <= 0.2cm, shorten >= 0.15cm](3_21)--(2_121);
\draw[dashed, ->, line width=0.3mm, >=latex, shorten <= 0.2cm, shorten >= 0.15cm](4_21) to [out=140, in=290] (1_121);
\draw[->, line width=0.3mm, >=latex, shorten <= 0.2cm, shorten >= 0.15cm](1_121) to [out=65, in=180] (1_21);
\draw[->, line width=0.3mm, >=latex, shorten <= 0.2cm, shorten >= 0.15cm](2_121) to [out=65, in=190] (2_21);
\draw[->, line width=0.3mm, >=latex, shorten <= 0.2cm, shorten >= 0.15cm](1_121)--(3_21);
\draw[->, line width=0.3mm, >=latex, shorten <= 0.2cm, shorten >= 0.15cm](2_121)--(4_21);
\draw[->, line width=0.3mm, >=latex, shorten <= 0.2cm, shorten >= 0.1cm](2_19)--(1_u);
\draw[->, line width=0.3mm, >=latex, shorten <= 0.2cm, shorten >= 0.1cm](2_19)--(2_u);
\draw[->, line width=0.3mm, >=latex, shorten <= 0.2cm, shorten >= 0.1cm](2_19)--(3_u);
\draw[->, line width=0.3mm, >=latex, shorten <= 0.2cm, shorten >= 0.1cm](1_19)--(4_u);
\draw[->, line width=0.3mm, >=latex, shorten <= 0.2cm, shorten >= 0.1cm](1_19)--(5_u);
\draw[->, line width=0.3mm, >=latex, shorten <= 0.2cm, shorten >= 0.1cm](1_19)--(6_u);
\draw[densely dotted, ->, line width=0.3mm, >=latex, shorten <= 0.2cm, shorten >= 0.1cm](3_19)--(1_u);
\draw[densely dotted, ->, line width=0.3mm, >=latex, shorten <= 0.2cm, shorten >= 0.1cm](3_19)--(2_u);
\draw[densely dotted, ->, line width=0.3mm, >=latex, shorten <= 0.2cm, shorten >= 0.1cm](3_19)--(3_u);
\draw[densely dotted, ->, line width=0.3mm, >=latex, shorten <= 0.2cm, shorten >= 0.1cm](3_19)--(5_u);

\draw[->, line width=0.3mm, >=latex, shorten <= 0.2cm, shorten >= 0.1cm](2_20)--(1_u);
\draw[->, line width=0.3mm, >=latex, shorten <= 0.2cm, shorten >= 0.1cm](2_20)--(2_u);
\draw[->, line width=0.3mm, >=latex, shorten <= 0.2cm, shorten >= 0.1cm](2_20)--(3_u);
\draw[->, line width=0.3mm, >=latex, shorten <= 0.2cm, shorten >= 0.1cm](1_20)--(4_u);
\draw[->, line width=0.3mm, >=latex, shorten <= 0.2cm, shorten >= 0.1cm](1_20)--(5_u);
\draw[->, line width=0.3mm, >=latex, shorten <= 0.2cm, shorten >= 0.1cm](1_20)--(6_u);
\draw[densely dotted, ->, line width=0.3mm, >=latex, shorten <= 0.2cm, shorten >= 0.1cm](3_20)--(1_u);
\draw[densely dotted, ->, line width=0.3mm, >=latex, shorten <= 0.2cm, shorten >= 0.1cm](3_20)--(2_u);
\draw[densely dotted, ->, line width=0.3mm, >=latex, shorten <= 0.2cm, shorten >= 0.1cm](3_20)--(3_u);
\draw[densely dotted, ->, line width=0.3mm, >=latex, shorten <= 0.2cm, shorten >= 0.1cm](3_20)--(4_u);

\draw[->, line width=0.3mm, >=latex, shorten <= 0.2cm, shorten >= 0.1cm](2_21)--(4_u);
\draw[->, line width=0.3mm, >=latex, shorten <= 0.2cm, shorten >= 0.1cm](2_21)--(5_u);
\draw[->, line width=0.3mm, >=latex, shorten <= 0.2cm, shorten >= 0.1cm](2_21)--(6_u);
\draw[->, line width=0.3mm, >=latex, shorten <= 0.2cm, shorten >= 0.1cm](3_21)--(4_u);
\draw[->, line width=0.3mm, >=latex, shorten <= 0.2cm, shorten >= 0.1cm](3_21)--(5_u);
\draw[->, line width=0.3mm, >=latex, shorten <= 0.2cm, shorten >= 0.1cm](3_21)--(6_u);
\draw[->, line width=0.3mm, >=latex, shorten <= 0.2cm, shorten >= 0.1cm](1_21)--(1_u);
\draw[->, line width=0.3mm, >=latex, shorten <= 0.2cm, shorten >= 0.1cm](1_21)--(2_u);
\draw[->, line width=0.3mm, >=latex, shorten <= 0.2cm, shorten >= 0.1cm](1_21)--(3_u);
\draw[->, line width=0.3mm, >=latex, shorten <= 0.2cm, shorten >= 0.1cm](4_21)--(1_u);
\draw[->, line width=0.3mm, >=latex, shorten <= 0.2cm, shorten >= 0.1cm](4_21)--(2_u);
\draw[->, line width=0.3mm, >=latex, shorten <= 0.2cm, shorten >= 0.1cm](4_21)--(3_u);

\draw[->, line width=0.3mm, >=latex, shorten <= 0.2cm, shorten >= 0.1cm](2_22)--(4_u);
\draw[->, line width=0.3mm, >=latex, shorten <= 0.2cm, shorten >= 0.1cm](2_22)--(5_u);
\draw[->, line width=0.3mm, >=latex, shorten <= 0.2cm, shorten >= 0.1cm](2_22)--(6_u);
\draw[->, line width=0.3mm, >=latex, shorten <= 0.2cm, shorten >= 0.1cm](3_22)--(4_u);
\draw[->, line width=0.3mm, >=latex, shorten <= 0.2cm, shorten >= 0.1cm](3_22)--(5_u);
\draw[->, line width=0.3mm, >=latex, shorten <= 0.2cm, shorten >= 0.1cm](3_22)--(6_u);
\draw[->, line width=0.3mm, >=latex, shorten <= 0.2cm, shorten >= 0.1cm](1_22)--(1_u);
\draw[->, line width=0.3mm, >=latex, shorten <= 0.2cm, shorten >= 0.1cm](1_22)--(2_u);
\draw[->, line width=0.3mm, >=latex, shorten <= 0.2cm, shorten >= 0.1cm](1_22)--(3_u);
\draw[->, line width=0.3mm, >=latex, shorten <= 0.2cm, shorten >= 0.1cm](4_22)--(1_u);
\draw[->, line width=0.3mm, >=latex, shorten <= 0.2cm, shorten >= 0.1cm](4_22)--(2_u);
\draw[->, line width=0.3mm, >=latex, shorten <= 0.2cm, shorten >= 0.1cm](4_22)--(3_u);

\draw[->, line width=0.3mm, >=latex, shorten <= 0.2cm, shorten >= 0.1cm](1_23)--(4_u);
\draw[->, line width=0.3mm, >=latex, shorten <= 0.2cm, shorten >= 0.1cm](1_23)--(5_u);
\draw[->, line width=0.3mm, >=latex, shorten <= 0.2cm, shorten >= 0.1cm](1_23)--(6_u);
\draw[->, line width=0.3mm, >=latex, shorten <= 0.2cm, shorten >= 0.1cm](2_23)--(4_u);
\draw[->, line width=0.3mm, >=latex, shorten <= 0.2cm, shorten >= 0.1cm](2_23)--(5_u);
\draw[->, line width=0.3mm, >=latex, shorten <= 0.2cm, shorten >= 0.1cm](2_23)--(6_u);

\draw[->, line width=0.3mm, >=latex, shorten <= 0.2cm, shorten >= 0.1cm](1_24)--(4_u);
\draw[->, line width=0.3mm, >=latex, shorten <= 0.2cm, shorten >= 0.1cm](1_24)--(5_u);
\draw[->, line width=0.3mm, >=latex, shorten <= 0.2cm, shorten >= 0.1cm](1_24)--(6_u);
\draw[->, line width=0.3mm, >=latex, shorten <= 0.2cm, shorten >= 0.1cm](2_24)--(4_u);
\draw[->, line width=0.3mm, >=latex, shorten <= 0.2cm, shorten >= 0.1cm](2_24)--(5_u);
\draw[->, line width=0.3mm, >=latex, shorten <= 0.2cm, shorten >= 0.1cm](2_24)--(6_u);
\end{tikzpicture}
\captionsetup{justification=centering}
{
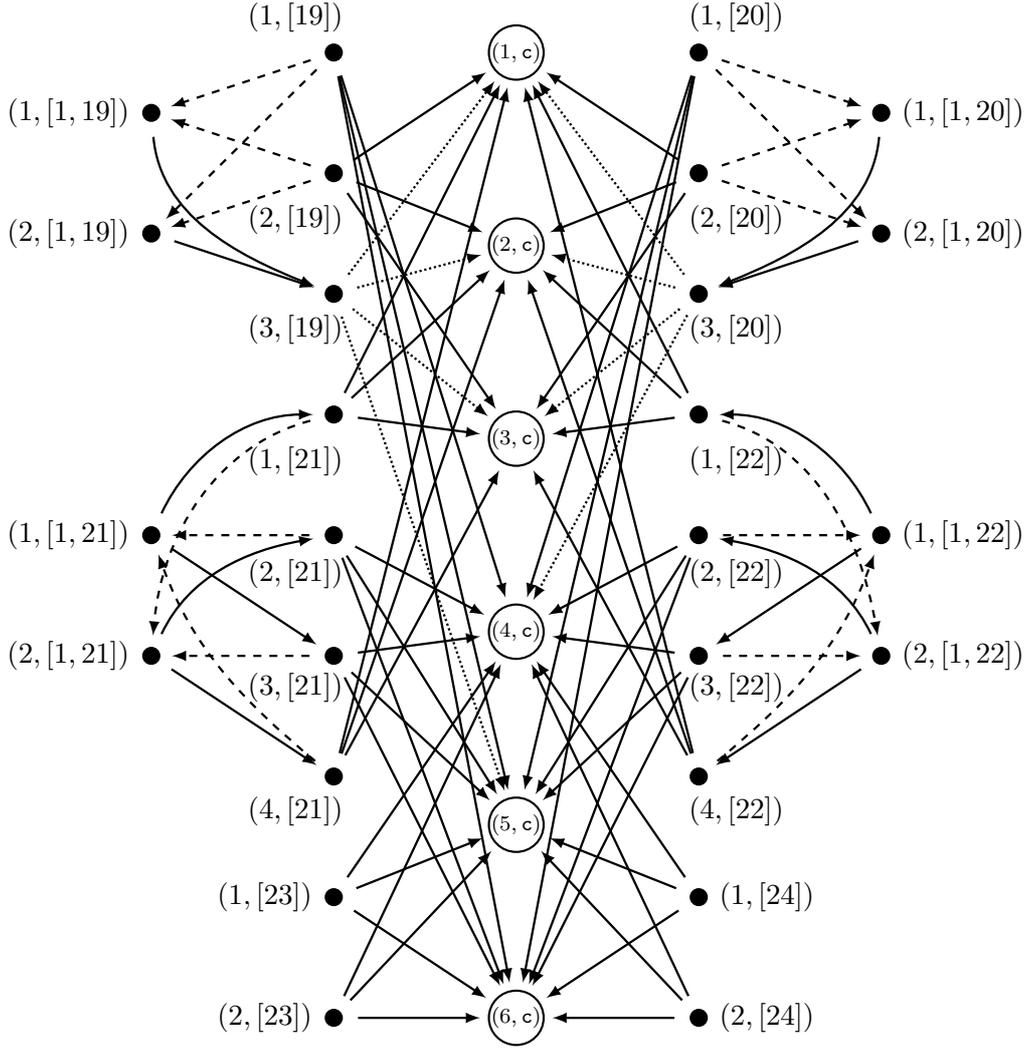
\captionof{figure}{Partial orientation $D$ for $\mathcal{H}$ for $s=6$, $k=13$;\\
 $A_2=\{[1],[2],\ldots,[12] \}$, $A_3=\{[13],[14],\ldots, [20]\}$, $A_4=\{[21],[22]\}$, $E=\{[23],[24]\}$.}\label{figC6.3.11}}
\end{center}
\indent\par This concludes the proof of Theorem \ref{thmC6.1.7}. Unfortunately, we were not able to give a complete characterisation for Proposition \ref{ppnC6.3.12}. Our core idea is the consideration of cross-intersecting antichains with at most $k$ disjoint pairs, thus, invoking Theorem \ref{thmC6.2.10}. We believe the gap between the necessary and sufficient conditions ($\kappa_{s,\frac{s}{2}}^*(\cdot)$ and $\kappa_{s,\frac{s}{2}}(\cdot)$ resp.) may be further tightened if there is an analogue of Theorem \ref{thmC6.2.10} on \textit{exactly} $k$ disjoint pairs; more discussion may be found in the last section of \cite{WHW TEG 5}. 

\section{Proof of Theorem \ref{thmC6.1.8}}
Similar to the previous section, we prove Theorem \ref{thmC6.1.8} by collating several propositions.
\begin{ppn}\label{ppnC6.4.1}
Suppose $s\ge 3$ is odd, $A_2=\emptyset$ and $A_3\neq \emptyset$ for a $\mathcal{T}$. Then, $\mathcal{T}\in \mathscr{C}_0$ if and only if $|A_3|\le 2{{s}\choose{\lceil{s/2}\rceil}}-2$.
\end{ppn}
\noindent\textit{Proof}: $(\Rightarrow)$ Since $\mathcal{T}\in \mathscr{C}_0$, there exists an orientation $D$ of $\mathcal{T}$, where $d(D)=4$. As $A_3\neq\emptyset$, we assume (\ref{eqC6.3.3})-(\ref{eqC6.3.5}) here. By Sperner's theorem, $|B^O_3|\le {{s}\choose{\lfloor{s/2}\rfloor}}$ and $|B^I_3|\le {{s}\choose{\lfloor{s/2}\rfloor}}$. So, if $|B^O_3|=0$ or $|B^I_3|=0$, then $|A_3|=|B^O_3|+|B^I_3|\le {{s}\choose{\lfloor{s/2}\rfloor}}\le 2{{s}\choose{\lceil{s/2}\rceil}}-2$. Therefore, we assume $|B^O_3|>0$ and $|B^I_3|>0$.
\indent\par In what follows, we first show that for any $[i]\in A_3$ and any $p=1,2,3$, if $|O^\mathtt{c}((p,[i]))|$ is too big ($>\lceil\frac{s}{2}\rceil$) or too small ($<\lfloor\frac{s}{2}\rfloor$), then $|A_3|\le 2{{s}\choose{\lceil{s/2}\rceil}}-2$.
\\
\\Case 1. There exists some $[i]\in A_3$ such that  $|O^\mathtt{c}((p,[i]))|>\lceil\frac{s}{2}\rceil$ for some $p=1,2,3$.
\indent\par For any $[j]\in A^O_3-\{[i]\}$, $d_D((1,[1,j]),(p,[i]))=3$ implies $O^\mathtt{c}((3,[j]))\cap I^\mathtt{c}((p,[i]))\neq \emptyset$. By Lih's theorem, $|B^O_3|-1\le |B^O_3-\{O^\mathtt{c}((3,[i]))\}|\le {{s}\choose{\lceil{s/2}\rceil}}-{{s-|I^\mathtt{c}((p,[i]))|}\choose{\lceil{s/2}\rceil}}\le {{s}\choose{\lceil{s/2}\rceil}}-{{\lceil{s/2}\rceil+1}\choose{\lceil{s/2}\rceil}}={{s}\choose{\lceil{s/2}\rceil}}-(\lceil{\frac{s}{2}}\rceil+1)\le {{s}\choose{\lceil{s/2}\rceil}}-3$. It follows that $|B^O_3|\le {{s}\choose{\lceil{s/2}\rceil}}-2$ and $|A_3|=|B^O_3|+|B^I_3|\le [{{s}\choose{\lceil{s/2}\rceil}}-2]+{{s}\choose{\lfloor{s/2}\rfloor}}=2{{s}\choose{\lceil{s/2}\rceil}}-2$.
\\
\\Case 2. There exists some $[i]\in A_3$ such that $|O^\mathtt{c}((p,[i]))|<\lfloor\frac{s}{2}\rfloor$ for some $p=1,2,3$.
\indent\par In other words, $|I^\mathtt{c}((p,[i]))|>\lceil\frac{s}{2}\rceil$. Hence, this case follows from Case 1 by the Duality Lemma.
\\
\\\noindent Case 3. For all $[i]\in A_3$ and all $p=1,2,3$, $\lfloor\frac{s}{2}\rfloor\le|O^\mathtt{c}((p,[i]))|\le \lceil\frac{s}{2}\rceil$.
\indent\par Note that for all $[i]\in A^O_3$ and $[j]\in A^I_3$, $d_D((1,[1,i]),(1,[1,j]))=4$ implies $X\cap Y\neq\emptyset$ for all $X\in B^O_3$ and $Y\in B^I_3$. Now, it suffices to consider the case where $B^O_3\cup B^I_3 \subseteq {{(\mathbb{N}_s,\mathtt{c})}\choose{\lceil{s/2}\rceil}}$. For otherwise, $|A_3|=|B^O_3|+|B^I_3|\le 2{{s}\choose{\lceil{s/2}\rceil}}-2$ by Theorems \ref{thmC6.2.4} and \ref{thmC6.2.5}. Partition $A^O_3$ ($A^I_3$ resp.) into $A^{O(D)}_3$ and $A^{O(S)}_3$ ($A^{I(D)}_3$ and $A^{I(S)}_3$ resp.) as in (\ref{eqC6.3.33}). 
\begin{rmk}\label{rmkC6.4.2}
Now, we shall make a series of assumptions on the structure of $D$, on which we will derive $|A_3|\le 2{{s}\choose{\lceil{s/2}\rceil}}-2$ if any one fails to hold. We will then show that under all these assumptions, we still arrive at the same required conclusion.
\end{rmk}
\noindent Assumption 1: $|A^{O(D)}_3|\ge 2$ and $|A^{I(D)}_3|\ge 2$.
\indent\par Suppose $|A^{O(D)}_3|\le 1$. By Lemma \ref{lemC6.2.16}, $\{I^\mathtt{c}((1,[j]))\mid [j]\in A^{O(S)}_3\}\cup B^I_3$ is an antichain. By Sperner's theorem, $|A_3|=|B^O_3|+|B^I_3|=|A^{O(D)}_3|+|\{I^\mathtt{c}((1,[j]))\mid [j]\in A^{O(S)}_3\}|+|B^I_3|\le 1+{{s}\choose{\lfloor{s/2}\rfloor}}\le 2{{s}\choose{\lceil{s/2}\rceil}}-2$. A similar argument follows if $|A^{I(D)}_3|\le 1$.
\\
\\Subcase 3.1. $|O^\mathtt{c}((1,[i]))|=|O^\mathtt{c}((2,[i]))|=\lfloor\frac{s}{2}\rfloor$ for some $[i]\in A^{O(D)}_3$.
\indent\par For any $[j]\in A^I_3$ and $p=1,2$, $d_D((p,[i]),(1,[1,j]))=3$ implies $O^\mathtt{c}((p,[i]))\cap I^\mathtt{c}((3,[j]))\neq \emptyset$, i.e., $I^\mathtt{c}((3,[j]))\neq I^\mathtt{c}((p,[i]))$. It follows that $|B^I_3|\le {{s}\choose{\lceil{s/2}\rceil}}-2$. Hence, $|A_3|=|B^O_3|+|B^I_3|\le {{s}\choose{\lceil{s/2}\rceil}}+[{{s}\choose{\lceil{s/2}\rceil}}-2]=2{{s}\choose{\lceil{s/2}\rceil}}-2$.
\\
\\Subcase 3.2. $|O^\mathtt{c}((1,[i]))|=|O^\mathtt{c}((2,[i]))|=\lceil\frac{s}{2}\rceil$ for some $[i]\in A^{I(D)}_3$.
\indent\par This follows from Subcase 3.1 by the Duality Lemma.
\\
\\Subcase 3.3. $|O^\mathtt{c}((1,[i^*]))|=\lfloor\frac{s}{2}\rfloor$ and $|O^\mathtt{c}((2,[i^*]))|=\lceil\frac{s}{2}\rceil$ for some $[i^*]\in A^{O(D)}_3$.
\indent\par $(\star)$ For any $[j]\in A^I_3$ and $p=1,2$, $d_D((1,[i^*]),(1,[1,j]))=3$ implies $I^\mathtt{c}((3,[j]))\neq I^\mathtt{c}((1,[i^*]))$. It follows that $|B^I_3|\le {{s}\choose{\lceil{s/2}\rceil}}-1$. Now, we are going to establish some assumptions regarding $A^{O(D)}_3$ and $A^{I(D)}_3$, and provide justifications accordingly.
\\
\\Assumption 2A: $|O^\mathtt{c}((1,[i]))|=\lfloor\frac{s}{2}\rfloor$ and $|O^\mathtt{c}((2,[i]))|=\lceil\frac{s}{2}\rceil$ for all $[i]\in A^{O(D)}_3$.
\indent\par Suppose there exists some $[i]\in A^{O(D)}_3-\{[i^*]\}$ such that $|O^\mathtt{c}((1,[i]))|=|O^\mathtt{c}((2,[i]))|=\lceil\frac{s}{2}\rceil$. Note by definition of $A^{O(D)}_3$ that $O^\mathtt{c}((3,[i]))$ is equal to at most one of $O^\mathtt{c}((1,[i]))$ and $O^\mathtt{c}((2,[i]))$, say $O^\mathtt{c}((3,[i]))\neq O^\mathtt{c}((1,[i]))$. Also, for any $[j]\in A^O_3-\{[i]\}$, $d_D((1,[1,j]),$ $(1,[i]))=3$ implies $O^\mathtt{c}((3,[j]))\neq O^\mathtt{c}((1,[i]))$. It follows that $|B^O_3|\le {{s}\choose{\lceil{s/2}\rceil}}-1$. Hence, $|A_3|=|B^O_3|+|B^I_3|\le 2[{{s}\choose{\lceil{s/2}\rceil}}-1]=2{{s}\choose{\lceil{s/2}\rceil}}-2$. Therefore, and in view of Subcase 3.1, we may now assume $|O^\mathtt{c}((1,[i]))|=\lfloor\frac{s}{2}\rfloor$ and $|O^\mathtt{c}((2,[i]))|=\lceil\frac{s}{2}\rceil$ for all $[i]\in A^{O(D)}_3$.
\\
\\Assumption 3A: $O^\mathtt{c}((1,[i]))=O^\mathtt{c}((1,[i^*]))$ for all $[i]\in A^{O(D)}_3$.
\indent\par Suppose there exists some $[i]\in A^{O(D)}_3-\{[i^*]\}$ such that $O^\mathtt{c}((1,[i]))\neq O^\mathtt{c}((1,[i^*]))$. For any $[j]\in A^I_3$ and $x=i,i^*$, $d_D((1,[x]),(1,[1,j]))=3$ implies $I^\mathtt{c}((3,[j]))\neq I^\mathtt{c}((1,[x]))$. It follows that $|B^I_3|\le {{s}\choose{\lceil{s/2}\rceil}}-2$. Hence, $|A_3|=|B^O_3|+|B^I_3|\le {{s}\choose{\lceil{s/2}\rceil}}+[{{s}\choose{\lceil{s/2}\rceil}}-2]=2{{s}\choose{\lceil{s/2}\rceil}}-2$. Thus, we may assume $O^\mathtt{c}((1,[i]))=O^\mathtt{c}((1,[i^*]))$ for all $[i]\in A^{O(D)}_3$.
\\
\\Assumption 4A: $O^\mathtt{c}((2,[i]))=O^\mathtt{c}((3,[i]))$ for all $[i]\in A^{O(D)}_3$.
\indent\par Suppose there exists some $[i]\in A^{O(D)}_3$ such that $O^\mathtt{c}((2,[i]))\neq O^\mathtt{c}((3,[i]))$. Also, for any $[j]\in A^O_3-\{[i]\}$, $d_D((1,[1,j]),(2,[i]))=3$ implies $O^\mathtt{c}((3,[j]))\neq O^\mathtt{c}((2,[i]))$. It follows that $|B^O_3|\le {{s}\choose{\lceil{s/2}\rceil}}-1$, and $|A_3|=|B^O_3|+|B^I_3|\le 2[{{s}\choose{\lceil{s/2}\rceil}}-1]=2{{s}\choose{\lceil{s/2}\rceil}}-2$. Therefore, the assumption follows.
\\
\\Assumption 5A: $|O^\mathtt{c}((1,[i]))|=|O^\mathtt{c}((2,[i]))|=\lfloor\frac{s}{2}\rfloor$ for all $[i]\in A^{I(D)}_3$.
\indent\par Suppose there exists some $[i]\in A^{I(D)}_3$ such that $|O^\mathtt{c}((1,[i]))|=\lfloor\frac{s}{2}\rfloor$ and $|O^\mathtt{c}((2,[i]))|=\lceil\frac{s}{2}\rceil$. For any $[j]\in A^O_3$, $d_D((1,[1,j]),(2,[i]))=3$ implies $O^\mathtt{c}((3,[j]))\neq O^\mathtt{c}((2,[i]))$. So, $|B^O_3|\le {{s}\choose{\lceil{s/2}\rceil}}-1$. Hence, $|A_3|=|B^O_3|+|B^I_3|\le 2[{{s}\choose{\lceil{s/2}\rceil}}-1]=2{{s}\choose{\lceil{s/2}\rceil}}-2$. Therefore, and in view of Subcase 3.2, we may now assume $|O^\mathtt{c}((1,[i]))|=|O^\mathtt{c}((2,[i]))|=\lfloor\frac{s}{2}\rfloor$ for all $[i]\in A^{I(D)}_3$.
\\
\\Assumption 6A: $O^\mathtt{c}((1,[i]))= O^\mathtt{c}((1,[i^*]))$ and $O^\mathtt{c}((2,[i]))=O^\mathtt{c}((3,[i]))$ for all $[i]\in A^{I(D)}_3$.
\noindent\par Suppose there exists some $[i]\in A^{I(D)}_3$ and some $p=1,2$, such that $O^\mathtt{c}((p,[i]))\neq O^\mathtt{c}((1,[i^*]))$ and $O^\mathtt{c}((p,[i]))\neq O^\mathtt{c}((3,[i]))$. Also, for any $[j]\in A^I_3-\{[i]\}$, $d_D((p,[i]),(1,[1,j]))$ $=3$ implies $I^\mathtt{c}((3,[j]))\neq I^\mathtt{c}((p,[i]))$. Therefore, for all $X\in B^I_3$, $X\neq I^\mathtt{c}((p,[i]))$ and recall from ($\star$) that $X\neq I^\mathtt{c}((1,[i^*]))$. It follows that $|B^I_3|\le {{s}\choose{\lceil{s/2}\rceil}}-2$ and $|A_3|=|B^O_3|+|B^I_3|\le {{s}\choose{\lceil{s/2}\rceil}}+[{{s}\choose{\lceil{s/2}\rceil}}-2]=2{{s}\choose{\lceil{s/2}\rceil}}-2$.
\indent\par Therefore, for each $[i]\in A^{I(D)}_3$ and each $p=1,2$, either $O^\mathtt{c}((p,[i]))=O^\mathtt{c}((3,[i]))$ or $O^\mathtt{c}((p,[i]))= O^\mathtt{c}((1,[i^*]))$. By the definition of $A^{I(D)}_3$, we may assume without loss of generality that $O^\mathtt{c}((1,[i]))= O^\mathtt{c}((1,[i^*]))$ and $O^\mathtt{c}((2,[i]))=O^\mathtt{c}((3,[i]))$ for all $[i]\in A^{I(D)}_3$.

\indent\par Now, with Assumptions 1, 2A-6A in place, consider $[j]\in A^{I(D)}_3$. For any $[k]\in A^{O(D)}_3$, $O^\mathtt{c}((1,[k]))=O^\mathtt{c}((1,[i^*]))=O^\mathtt{c}((1,[j]))$ and $d_D((1,[j]),(1,[1,k]))=3$ imply $O^\mathtt{c}((1,[j]))\cap I^\mathtt{c}((2,[k]))\neq \emptyset$. Equivalently, $O^\mathtt{c}((1,[i^*]))=O^\mathtt{c}((1,[j]))\not\subseteq O^\mathtt{c}((2,[k]))$. Note also that there are $\lceil\frac{s}{2}\rceil$ number of $\lceil\frac{s}{2}\rceil$-supersets of $O^\mathtt{c}((1,[i^*]))$. Recall that $O^\mathtt{c}((2,[k]))=O^\mathtt{c}((3,[k]))$, so that $\{O^\mathtt{c}((2,[k]))\mid [k]\in A^{O(D)}_3\}=\{O^\mathtt{c}((3,[k]))\mid [k]\in A^{O(D)}_3\}\subseteq {{(\mathbb{N}_s,\mathtt{c})}\choose{\lceil{s/2}\rceil}}$. So, $|A^{O(D)}_3|\le {{s}\choose{\lceil{s/2}\rceil}}-\lceil\frac{s}{2}\rceil\le {{s}\choose{\lceil{s/2}\rceil}}-2$. Since $\{I^\mathtt{c}((1,[i]))\mid [i]\in A^{O(S)}_3\}\cup B^I_3$ is an antichain by Lemma \ref{lemC6.2.16}, $|A^{O(S)}_3|+|A^I_3|=|\{I^\mathtt{c}((1,[i]))\mid [i]\in A^{O(S)}_3\}\cup B^I_3|\le {{s}\choose{\lfloor{s/2}\rfloor}}$ by Sperner's theorem. Hence, $|A_3|=|A^O_3|+|A^I_3|=|A^{O(D)}_3|+(|A^{O(S)}_3|+|A^I_3|)\le [{{s}\choose{\lceil{s/2}\rceil}}-2]+{{s}\choose{\lfloor{s/2}\rfloor}}=2{{s}\choose{\lceil{s/2}\rceil}}-2$.
\\
\\Subcase 3.4. $|O^\mathtt{c}((1,[i^*]))|=\lceil\frac{s}{2}\rceil$ and $|O^\mathtt{c}((2,[i^*]))|=\lfloor\frac{s}{2}\rfloor$ for some $[i^*]\in A^{I(D)}_3$.
\indent\par This follows from Subcase 3.3 by the Duality Lemma. For clarity, we state the analogous versions of Assumptions 2A-6A below.
\\
\\Assumption 2B: $|O^\mathtt{c}((1,[i]))|=\lceil\frac{s}{2}\rceil$ and $|O^\mathtt{c}((2,[i]))|=\lfloor\frac{s}{2}\rfloor$ for all $[i]\in A^{I(D)}_3$.
\\Assumption 3B: $O^\mathtt{c}((1,[i]))=O^\mathtt{c}((1,[i^*]))$ for all $[i]\in A^{I(D)}_3$.
\\Assumption 4B: $O^\mathtt{c}((2,[i]))=O^\mathtt{c}((3,[i]))$ for all $[i]\in A^{I(D)}_3$.
\\Assumption 5B: $|O^\mathtt{c}((1,[j]))|=|O^\mathtt{c}((2,[j]))|=\lceil\frac{s}{2}\rceil$ for all $[j]\in A^{O(D)}_3$.
\\Assumption 6B: $O^\mathtt{c}((1,[j]))= O^\mathtt{c}((1,[i^*]))$ and $O^\mathtt{c}((2,[j]))=O^\mathtt{c}((3,[j]))$ for all $[j]\in A^{O(D)}_3$. 
\\
\\Subcase 3.5. $|O^\mathtt{c}((1,[i]))|=|O^\mathtt{c}((2,[i]))|=\lceil\frac{s}{2}\rceil$ and $|O^\mathtt{c}((1,[j]))|=|O^\mathtt{c}((2,[j]))|=\lfloor\frac{s}{2}\rfloor$ for some $[i]\in A^{O(D)}_3$ and $[j]\in A^{I(D)}_3$.
\indent\par Note by definition of $A^{O(D)}_3$ that $O^\mathtt{c}((3,[i]))$ is equal to at most one of $O^\mathtt{c}((1,[i]))$ and $O^\mathtt{c}((2,[i]))$, say $O^\mathtt{c}((3,[i]))\neq O^\mathtt{c}((1,[i]))$. Also, for any $[k]\in A^O_3-\{[i]\}$, $d_D((1,[1,k]),(1,[i]))$ $=3$ implies $O^\mathtt{c}((3,[k]))\neq O^\mathtt{c}((1,[i]))$. It follows that $|B^O_3|\le {{s}\choose{\lceil{s/2}\rceil}}-1$.
\indent\par Similarly, by definition of $A^{I(D)}_3$, $I^\mathtt{c}((3,[j]))$ is equal to at most one of $I^\mathtt{c}((1,[j]))$ and $I^\mathtt{c}((2,[j]))$, say $I^\mathtt{c}((3,[j]))\neq I^\mathtt{c}((1,[j]))$. Also, for any $[k]\in A^I_3-\{[j]\}$, $d_D((1,[j]),(1,[1,k]))$ $=3$ implies $I^\mathtt{c}((3,[k]))\neq I^\mathtt{c}((1,[j]))$. It follows that $|B^I_3 |\le {{s}\choose{\lceil{s/2}\rceil}}-1$. Hence, $|A_3|=|B^O_3|+|B^I_3|\le 2[{{s}\choose{\lceil{s/2}\rceil}}-1]=2{{s}\choose{\lceil{s/2}\rceil}}-2$.

\indent\par In view of the above, it is intuitive to construct an optimal orientation $D$ of $\mathcal{T}$ with 
$B^O_3\cup B^I_3\subseteq {{(\mathbb{N}_s,\mathtt{c})}\choose{\lceil{s/2}\rceil}}$. To this end, we recall Definition \ref{defnC6.3.7}.
\\
\\$(\Leftarrow)$ If $|A_{\ge 3}|\le {{s}\choose{\lceil{s/2}\rceil}}-1$, then by Corollary \ref{corC6.3.8}(i), $\mathcal{T}\in\mathscr{C}_0$. Hence, we assume $|A_{\ge 3}|\ge {{s}\choose{\lceil{s/2}\rceil}}$ hereafter, on top of the hypothesis that $|A_3|\le 2{{s}\choose{\lceil{s/2}\rceil}}-2$. If $|A_3|\ge {{s}\choose{\lceil{s/2}\rceil}}$, define $A^{\diamond}_3=A_3$. Otherwise, let $A^{\diamond}_3=A_3\cup A^*$, where $A^*$ is an arbitrary subset of $A_{\ge 4}$ such that $|A^{\diamond}_3|={{s}\choose{\lceil{s/2}\rceil}}$. Then, let $A^{\diamond}_4=A_3\cup A_{\ge 4}-A^{\diamond}_3$. Furthermore, assume without loss of generality that $A^{\diamond}_3=\{[i]\mid i\in\mathbb{N}_{|A^{\diamond}_3|}\}$ and $A^{\diamond}_4=\{[i]\mid i\in\mathbb{N}_{|A^{\diamond}_3|+|A^{\diamond}_4|}-\mathbb{N}_{|A^{\diamond}_3|}\}$.
\indent\par Let $\mathcal{H}=T(t_1,t_2,\ldots, t_n)$ be the subgraph of $\mathcal{T}$, where $t_\mathtt{c}=s$, $t_{[i]}=3$ for all $[i]\in \mathcal{T}(A^{\diamond}_3)$, $t_{[j]}=4$ for all $[j]\in \mathcal{T}(A^{\diamond}_4)$ and $t_v=2$ otherwise. We will use $A_j$ for $\mathcal{H}(A_j)$ for the remainder of this proof. Define an orientation $D$ of $\mathcal{H}$ as follows.
\begin{align}
& (3,[i])\rightarrow \{(1,[\alpha,i]),(2,[\alpha,i])\}\rightarrow \{(1,[i]),(2,[i])\},\label{eqC6.4.1}\\
& \bar{\lambda}_1 \rightarrow (1,[i])\rightarrow \lambda_1\rightarrow (2,[i])\rightarrow \bar{\lambda}_1,\text{ and}\label{eqC6.4.2}\\
&\lambda_{i+1}\rightarrow (3,[i]) \rightarrow \bar{\lambda}_{i+1}\label{eqC6.4.3}
\end{align}
for all $1\le i\le {{s}\choose{\lceil{s/2}\rceil}}-1$ and all $1\le \alpha\le \deg_T([i])-1$, i.e., excluding $\lambda_1$, the $\lceil\frac{s}{2}\rceil$-sets $\lambda_i$'s are used as `in-sets' to construct $B^I_3$.
\begin{align}
& \{(1,[j]),(2,[j])\}\rightarrow \{(1,[\beta,j]),(2,[\beta,j])\}\rightarrow (3,[j]),\label{eqC6.4.4}\\
& \lambda_1 \rightarrow (1,[j])\rightarrow \bar{\lambda}_1\rightarrow (2,[j])\rightarrow \lambda_1,\text{ and}\label{eqC6.4.5}\\
&\bar{\lambda}_{j+2-{{s}\choose{\lceil{s/2}\rceil}}}\rightarrow (3,[j]) \rightarrow \lambda_{j+2-{{s}\choose{\lceil{s/2}\rceil}}}\label{eqC6.4.6}
\end{align}
for all ${{s}\choose{\lceil{s/2}\rceil}}\le j\le |A_3|$ and all $1\le \beta\le \deg_T([j])-1$, i.e., the $\lceil\frac{s}{2}\rceil$-sets $\lambda_2,\lambda_3,\ldots,\lambda_{|A_3|+2-{{s}\choose{\lceil{s/2}\rceil}}}$ are used as `out-sets' to construct $B^O_3$.
\begin{align}
&(2,[\gamma,k])\rightarrow \{(2,[k]),(4,[k])\}\rightarrow (1,[\gamma,k])\rightarrow \{(1,[k]), (3,[k])\}\rightarrow (2,[\gamma,k]),\label{eqC6.4.7}\\
&\text{and }\bar{\lambda}_1 \rightarrow \{(1,[k]), (4,[k])\}\rightarrow \lambda_1 \rightarrow\{(2,[k]), (3,[k])\}\rightarrow \bar{\lambda}_1\label{eqC6.4.8}
\end{align}
for all $[k]\in A_4$ and all $1\le\gamma\le \deg_T([k])-1$.
\begin{align}
\lambda_1 \rightarrow \{(1,[l]), (2,[l])\} \rightarrow \bar{\lambda}_1 \label{eqC6.4.9}
\end{align}
for all $[l]\in E$. (See Figure \ref{figC6.4.12} for $D$ when $s=3$.)
\\
\\Claim: $d_{D}(v,w)\le 4$ for all $v,w\in V(D)$.
\\
\\Case 1.1. $v,w \in \{(1,[\alpha,i]),(2,[\alpha,i]),(1,[i]),(2,[i]),(3,[i])\}$ for each $1\le i\le {{s}\choose{\lceil{s/2}\rceil}}-1$ and $1\le\alpha\le \deg_T([i])-1$.
\indent\par By (\ref{eqC6.4.2}), $(1,[i])\rightarrow \lambda_1$, and $(2,[i])\rightarrow \bar{\lambda}_1$. Since $\lambda_{i+1}\in {{(\mathbb{N}_s,\mathtt{c})}\choose{\lceil s/2\rceil}}-\{\lambda_1\}$ by (\ref{eqC6.4.3}), there exist some $x_i\in\lambda_1\cap \lambda_{i+1}$ and some $y_i\in\bar{\lambda}_1\cap \lambda_{i+1}$ such that $\{x_i,y_i\}\rightarrow (3,[i])$. With (\ref{eqC6.4.1}), this case follows.
\\
\\Case 1.2. $v,w \in \{(1,[\alpha,i]),(2,[\alpha,i]),(1,[i]),(2,[i]),(3,[i])\}$ for each ${{s}\choose{\lceil{s/2}\rceil}}\le i\le |A_3|$ and $1\le\alpha\le \deg_T([i])-1$.
\indent\par By (\ref{eqC6.4.5}), $\lambda_1\rightarrow (1,[i])$, and $\bar{\lambda}_1\rightarrow (2,[i])$. Since $\lambda_{i+2-{{s}\choose{\lceil{s/2}\rceil}}}\in {{(\mathbb{N}_s,\mathtt{c})}\choose{\lceil s/2\rceil}}-\{\lambda_1\}$ by (\ref{eqC6.4.6}), there exist some $x_i\in\lambda_1\cap\lambda_{i+2-{{s}\choose{\lceil{s/2}\rceil}}}$ and some $y_i\in\bar{\lambda}_1\cap\lambda_{i+2-{{s}\choose{\lceil{s/2}\rceil}}}$ such that $(3,[i])\rightarrow \{x_i, y_i\}$. With (\ref{eqC6.4.4}), this case follows.
\\
\\Case 1.3. $v,w \in \{(1,[\alpha,i]),(2,[\alpha,i]),(1,[i]),(2,[i]), (3,[i]), (4,[i])\}$ for each $[i]\in A_4$ and $1\le\alpha\le \deg_T([i])-1$.
\indent\par This is clear since (\ref{eqC6.4.7}) guarantees a directed $C_4$.
\\
\\Case 2.  For each $1\le i,j \le {{s}\choose{\lceil{s/2}\rceil}}-1$, $i\neq j$, each $1\le \alpha\le \deg_T([i])-1$, and each $1\le \beta\le \deg_T([j])-1$,
\\(i) $v=(p,[\alpha,i]), w=(q,[j])$ for each $p=1,2$ and $q=1,2,3$.
\indent\par By (\ref{eqC6.4.1})-(\ref{eqC6.4.2}), $\{(1,[\alpha,i]),(2,[\alpha,i])\}\rightarrow \{(1,[i]),(2,[i])\}$, $(1,[i])\rightarrow \lambda_1\rightarrow$ $ (2,[j])$ and $(2,[i])\rightarrow \bar{\lambda}_1\rightarrow (1,[j])$. Since $|\lambda_1|=|\lambda_{j+1}|=\lceil{\frac{s}{2}}\rceil$, there exists some $x_j\in \lambda_1\cap\lambda_{j+1}$ such that $x_j\rightarrow (3,[j])$ by (\ref{eqC6.4.3}).
\\
\\(ii) $v=(p,[\alpha,i]), w=(q,[\beta,j])$ for each $p,q=1,2$.
\indent\par From (i), $d_D((p,[\alpha,i]),(3,[j]))=3$ by (i). Since $(3,[j])\rightarrow \{(1,[\beta,j]),(2,[\beta,j])\}$ by (\ref{eqC6.4.1}), this subcase follows.
\\
\\(iii) $v=(q,[j]), w=(p,[\alpha,i])$ for each $p=1,2$ and $q=1,2,3$.
\indent\par By (\ref{eqC6.4.1})-(\ref{eqC6.4.2}), $(1,[j])\rightarrow \lambda_1$, $(2,[j])\rightarrow \bar{\lambda}_1$ and $(3,[i])\rightarrow \{(1,[\alpha,i]),(2,[\alpha,i])\}$. Since $\lambda_{i+1}\in {{(\mathbb{N}_s,\mathtt{c})}\choose{\lceil s/2\rceil}}-\{\lambda_1\}$, there exist some $x_i\in\lambda_1\cap\lambda_{i+1}$ and $y_i\in\bar{\lambda}_1\cap\lambda_{i+1}$ such that $\{x_i, y_i\}\rightarrow (3,[i])$ by (\ref{eqC6.4.3}). Also, since $i\neq j$, there exists $z_{ij}\in \bar{\lambda}_{j+1}\cap \lambda_{i+1}$ such that $(3,[j])\rightarrow z_{ij}\rightarrow (3,[i])$.
\\
\\Case 3.  For each ${{s}\choose{\lceil{s/2}\rceil}} \le i,j \le |A_3|$, $i\neq j$, each $1\le\alpha\le \deg_T([i])-1$, and each $1\le\beta\le \deg_T([j])-1$,
\\(i) $v=(p,[\alpha,i]), w=(q,[j])$ for each $p=1,2$ and $q=1,2,3$.
\indent\par By (\ref{eqC6.4.4})-(\ref{eqC6.4.5}), $\{(1,[\alpha,i]),(2,[\alpha,i])\}\rightarrow (3,[i])$, $\lambda_1 \rightarrow (1,[j]), \bar{\lambda}_1\rightarrow (2,[j])$. Since $\lambda_{i+2-{{s}\choose{\lceil{s/2}\rceil}}}$ $\in {{(\mathbb{N}_s,\mathtt{c})}\choose{\lceil s/2\rceil}}-\{\lambda_1\}$, there exist some $x_i\in\lambda_1\cap\lambda_{i+2-{{s}\choose{\lceil{s/2}\rceil}}}$ and $y_i\in\bar{\lambda}_1\cap\lambda_{i+2-{{s}\choose{\lceil{s/2}\rceil}}}$ such that $(3,[i])\rightarrow \{x_i, y_i\}$ by (\ref{eqC6.4.6}). Also, since $i\neq j$, there exists $z_{ij}\in \lambda_{i+2-{{s}\choose{\lceil{s/2}\rceil}}} \cap \bar{\lambda}_{j+2-{{s}\choose{\lceil{s/2}\rceil}}}$ such that $(3,[i])\rightarrow z_{ij}\rightarrow (3,[j])$.
\\
\\(ii) $v=(p,[\alpha,i]), w=(q,[\beta,j])$ for each $p,q=1,2$.
\indent\par From (i), $d_D((p,[\alpha,i]),(1,[j]))=3$. Since $(1,[j])\rightarrow \{(1, [\beta,j]), (2, [\beta,j])\}$ by (\ref{eqC6.4.4}), this subcase follows.
\\
\\(iii) $v=(q,[j]), w=(p,[\alpha,i])$ for each $p=1,2$ and $q=1,2,3$.
\indent\par By (\ref{eqC6.4.5}), $(1,[j])\rightarrow \bar{\lambda}_1\rightarrow (2,[i])$, $(2,[j])\rightarrow \lambda_1 \rightarrow (1,[i])$. With (\ref{eqC6.4.6}), since $|\lambda_{j+2-{{s}\choose{\lceil{s/2}\rceil}}}|=|\lambda_1|=\lceil{\frac{s}{2}}\rceil$, there exists some $x_j\in\lambda_1\cap \lambda_{j+2-{{s}\choose{\lceil{s/2}\rceil}}}$ such that $(3,[j])\rightarrow x_j$. With $\{(1,[i]),(2,[i])\}\rightarrow \{(1,[\alpha,i]),(2,[\alpha,i])\}$ by (\ref{eqC6.4.4}), this subcase follows.
\\
\\Case 4. For each $[i],[j]\in A_4$, $i \neq j$, each $1\le\alpha\le \deg_T([i])-1$, and each $1\le\beta\le \deg_T([j])-1$,
\\(i) $v=(p,[\alpha,i]), w=(q,[j])$ for $p=1,2$, and $q=1,2,3,4$.
\indent\par By (\ref{eqC6.4.7})-(\ref{eqC6.4.8}), $(2,[\alpha,i])\rightarrow \{(2,[i]),(4,[i])\}$ and $(1,[\alpha,i])\rightarrow \{(1,[i]), (3,[i])\}$, $\{(1,[i]),$ $(4,[i])\}\rightarrow \lambda_1$, $\{(2,[i]), (3,[i])\}\rightarrow \bar{\lambda}_1$ and $|I^\mathtt{c}((q,[j]))|>0$.
\\
\\(ii) $v=(p,[i]), w=(q,[\beta,j])$ for for $p=1,2,3,4$, and $q=1,2$.
\indent\par By (\ref{eqC6.4.7})-(\ref{eqC6.4.8}), $\{(2,[j]),(4,[j])\}\rightarrow (1,[\beta,j])$, $\{(1,[j]), (3,[j])\}\rightarrow (2,[\beta,j])$, $\bar{\lambda}_1 \rightarrow \{(1,[j]), (4,[j])\}$, $\lambda_1 \rightarrow\{(2,[j]), (3,[j])\}$ and $|O^\mathtt{c}((p,[i]))|>0$.
\\
\\(iii) $v=(p,[\alpha,i]), w=(q,[\beta,j])$ for $p,q=1,2$.
\indent\par From (i), $d_D((p,[\alpha,i]),(r,[j]))=3$ for $r=1,2,3,4$. So, this subcase holds as $\{(2,[j]),(4,[j])\}\rightarrow (1,[\beta,j])$, and $\{(1,[j]), (3,[j])\}\rightarrow (2,[\beta,j])$ by (\ref{eqC6.4.7}).
\\
\\Case 5. For each $1\le i \le {{s}\choose{\lceil{s/2}\rceil}}-1$, each ${{s}\choose{\lceil{s/2}\rceil}} \le j \le |A_3|$, each $1\le\alpha\le \deg_T([i])-1$, and each $1\le\beta\le \deg_T([j])-1$,
\\(i) $v=(p,[\alpha,i]), w=(q,[j])$ for each $p=1,2$ and $q=1,2,3$.
\indent\par By (\ref{eqC6.4.1})-(\ref{eqC6.4.2}) and (\ref{eqC6.4.5})-(\ref{eqC6.4.6}), $\{(1,[\alpha,i]),(2,[\alpha,i])\}\rightarrow \{(1,[i]),(2,[i])\}$, $(1,[i])\rightarrow \lambda_1\rightarrow (1,[j])$, $(2,[i])\rightarrow \bar{\lambda}_1\rightarrow (2,[j])$, and $|I^\mathtt{c}((3,[j]))|>0$.
\\
\\(ii) $v=(q,[j]), w=(p,[\alpha,i])$ for each $p=1,2$ and $q=1,2,3$.
\indent\par By (\ref{eqC6.4.5}), $(1,[j])\rightarrow \bar{\lambda}_1$ and $(2,[j])\rightarrow \lambda_1$. Since $\lambda_{i+1}\in {{(\mathbb{N}_s,\mathtt{c})}\choose{\lceil s/2\rceil}}-\{\lambda_1\}$, there exist some $x_i\in\lambda_1\cap\lambda_{i+1}$ and some $y_i\in\bar{\lambda}_1\cap\lambda_{i+1}$ such that $\{x_i,y_i\}\rightarrow (3,[i])$ by (\ref{eqC6.4.3}). With (\ref{eqC6.4.6}), since $|\lambda_{j+2-{{s}\choose{\lceil{s/2}\rceil}}}|=|\lambda_{i+1}|=\lceil\frac{s}{2}\rceil$, there exists some $z_{ij}\in\lambda_{j+2-{{s}\choose{\lceil{s/2}\rceil}}}\cap \lambda_{i+1}$ such that $(3,[j])\rightarrow z_{ij} \rightarrow (3,[i])$. With $(3,[i])\rightarrow \{(1,[\alpha,i]),(2,[\alpha,i])\}$ by ({\ref{eqC6.4.1}), this subcase follows.
\\
\\(iii) $v=(p,[i]), w=(q,[\beta,j])$ for each $p=1,2,3$ and $q=1,2$.
\indent\par By (\ref{eqC6.4.2})-(\ref{eqC6.4.5}), $(1,[i])\rightarrow \lambda_1\rightarrow (1,[j])$, $(2,[i])\rightarrow \bar{\lambda}_1\rightarrow (2,[j])$, $|O^\mathtt{c}((3,[i]))|>0$ and $\{(1,[j]),(2,[j])\}\rightarrow \{(1,[\beta,j]),(2,[\beta,j])\}$.
\\
\\(iv) $v=(q,[\beta,j]), w=(p,[i])$ for each $p=1,2,3$ and $q=1,2$.
\indent\par By (\ref{eqC6.4.4}) and (\ref{eqC6.4.2}), $\{(1,[\beta,j]),(2,[\beta,j])\}\rightarrow (3,[j])$, $\bar{\lambda}_1\rightarrow (1,[i]) $, $\lambda_1\rightarrow (2,[i])$. Since $\lambda_{j+2-{{s}\choose{\lceil{s/2}\rceil}}}$ $\in {{(\mathbb{N}_s,\mathtt{c})}\choose{\lceil s/2\rceil}}-\{\lambda_1\}$, there exists some $x_j\in\lambda_{j+2-{{s}\choose{\lceil{s/2}\rceil}}}\cap \lambda_1$ and some $y_j\in\lambda_{j+2-{{s}\choose{\lceil{s/2}\rceil}}}\cap \bar{\lambda}_1$ such that $(3,[j])\rightarrow \{x_j,y_j\}$ by (\ref{eqC6.4.6}). With (\ref{eqC6.4.3}), since $|\lambda_{j+2-{{s}\choose{\lceil{s/2}\rceil}}}|=|\lambda_{i+1}|=\lceil\frac{s}{2}\rceil$, there exists some $z_{ij}\in\lambda_{j+2-{{s}\choose{\lceil{s/2}\rceil}}}\cap \lambda_{i+1}$ such that $(3,[j])\rightarrow z_{ij}\rightarrow (3,[i])$.
\\
\\(v) $v=(p,[\alpha,i]), w=(q,[\beta,j])$ for each $p,q=1,2$.
\indent\par By (\ref{eqC6.4.1})-(\ref{eqC6.4.2}) and (\ref{eqC6.4.4})-(\ref{eqC6.4.5}), $\{(1,[\alpha,i]),(2,[\alpha,i])\}\rightarrow (1,[i])\rightarrow \lambda_1\rightarrow (1,[j])\rightarrow \{(1,[\beta,j]),$ $(2, [\beta,j])\}$.
\\
\\(vi) $v=(q,[\beta,j]), w=(p,[\alpha,i])$ for each $p,q=1,2$.
\indent\par From (iv), $d_D((q,[\beta,j]),(3,[i]))=3$. Since $(3,[i])\rightarrow \{(1, [\alpha,i]), (2, [\alpha,i])\}$ by (\ref{eqC6.4.1}), this subcase follows.
\\
\\Case 6.  For each $1\le i\le {{s}\choose{\lceil{s/2}\rceil}}-1$, each $[j]\in A_4$, each $1\le\alpha\le \deg_T([i])-1$, and each $1\le\beta\le \deg_T([j])-1$,
\\(i) $v=(p,[\alpha,i]), w=(q,[j])$ for each $p=1,2$ and $q=1,2,3,4$.
\indent\par By (\ref{eqC6.4.1})-(\ref{eqC6.4.2}) and (\ref{eqC6.4.8}), $\{(1,[\alpha,i]),(2,[\alpha,i])\}\rightarrow \{(1,[i]),(2,[i])\}$, $(1,[i])\rightarrow \lambda_1\rightarrow \{(2,[j]),$ $(3,[j])\}$ and $(2,[i])\rightarrow \bar{\lambda}_1\rightarrow \{(1,[j]),(4,[j])\}$.
\\
\\(ii) $v=(q,[j]), w=(p,[\alpha,i])$ for each $p=1,2$ and $q=1,2,3,4$.
\indent\par By (\ref{eqC6.4.8}), $\{(1,[j]), (4,[j])\}\rightarrow \lambda_1$ and $\{(2,[j]), (3,[j])\}\rightarrow \bar{\lambda}_1$. Since $\lambda_{i+1}\in {{(\mathbb{N}_s,\mathtt{c})}\choose{\lceil s/2\rceil}}-\{\lambda_1\}$, there exist some $x_i\in\lambda_{i+1}\cap\lambda_1$ and $y_i\in\lambda_{i+1}\cap\bar{\lambda}_1$ such that $\{x_i, y_i\}\rightarrow (3,[i])\rightarrow \{(1,[\alpha,i]),(2,[\alpha,i])\}$ by (\ref{eqC6.4.1}) and (\ref{eqC6.4.3}). 
\\
\\(iii) $v=(p,[i]), w=(q,[\beta,j])$ for each $p=1,2,3$ and $q=1,2$.
\indent\par By (\ref{eqC6.4.7})-(\ref{eqC6.4.8}) and (\ref{eqC6.4.2})-(\ref{eqC6.4.3}), $\lambda_1\rightarrow \{(2,[j]), (3,[j])\}$, $\bar{\lambda}_1\rightarrow \{(1,[j]), (4,[j])\}$, $\{(2,[j]),$ $(4,[j])\}\rightarrow (1,[\beta,j])$, $\{(1,[j]), (3,[j])\}\rightarrow (2,[\beta,j]),$ and $|O^\mathtt{c}((p,[i]))|>0$ for $p=1,2,3$.
\\
\\(iv) $v=(q,[\beta,j]), w=(p,[i])$ for each $p=1,2,3$ and $q=1,2$.
\indent\par By (\ref{eqC6.4.7})-(\ref{eqC6.4.8}) and (\ref{eqC6.4.2})-(\ref{eqC6.4.3}), $(2,[\beta,j])\rightarrow \{(2,[j]),(4,[j])\}$ and $(1,[\beta,j])\rightarrow \{(1,[j]),$ $(3,[j])\}$, $\{(1,[j]), (4,[j])\}$ $\rightarrow \lambda_1$, $\{(2,[j]), (3,[j])\}\rightarrow \bar{\lambda}_1$, and $|I^\mathtt{c}((p,[i]))|>0$ for $p=1,2,3$.
\\
\\(v) $v=(p,[\alpha,i]), w=(q,[\beta,j])$ for each $p,q=1,2$.
\indent\par From (i), $d_D((p,[\alpha,i]),(r,[j]))=3$ for $r=2,3$. Since  $(2,[j])\rightarrow (1,[\beta,j])$ and $(3,[j])\rightarrow (2,[\beta,j])$ by (\ref{eqC6.4.7}), this subcase follows.
\\
\\(vi) $v=(q,[\beta,j]), w=(p,[\alpha,i])$ for each $p,q=1,2$.
\indent\par From (iv), $d_D((q,[\beta,j]),(3,[i]))=3$. Since $(3,[i])\rightarrow \{(1,[\alpha,i]), (2,[\alpha,i])\}$ by (\ref{eqC6.4.1}), this subcase follows.
\\
\\Case 7.  For each ${{s}\choose{\lceil{s/2}\rceil}}\le i\le |A_3|$, each $[j]\in A_4$, each $1\le\alpha\le \deg_T([i])-1$, and each $1\le\beta\le \deg_T([j])-1$,
\\(i) $v=(p,[\alpha,i]), w=(q,[j])$ for each $p=1,2$ and $q=1,2,3,4$.
\indent\par By (\ref{eqC6.4.4}) and (\ref{eqC6.4.8}), $\{(1,[\alpha,i]),(2,[\alpha,i])\}\rightarrow (3,[i])$, $\lambda_1\rightarrow \{(2,[j]), (3,[j])\}$, and $\bar{\lambda}_1\rightarrow \{(1,[j]), (4,[j])\}$. Since $\lambda_{i+2-{{s}\choose{\lceil{s/2}\rceil}}}\in {{(\mathbb{N}_s,\mathtt{c})}\choose{\lceil s/2\rceil}}-\{\lambda_1\}$, there exist some $x_i\in\lambda_{i+2-{{s}\choose{\lceil{s/2}\rceil}}}\cap\lambda_1$ and $y_i\in\lambda_{i+2-{{s}\choose{\lceil{s/2}\rceil}}}\cap\bar{\lambda}_1$ such that $(3,[i])\rightarrow \{x_i, y_i\}$ by (\ref{eqC6.4.6}).
\\
\\(ii) $v=(q,[j]), w=(p,[\alpha,i])$ for each $p=1,2$ and $q=1,2,3,4$.
\indent\par By (\ref{eqC6.4.4})-(\ref{eqC6.4.5}) and (\ref{eqC6.4.8}), $\{(1,[j]), (4,[j])\}\rightarrow \lambda_1\rightarrow (1,[i])$, $\{(2,[j]), (3,[j])\}\rightarrow \bar{\lambda}_1\rightarrow (2,[i])$, and $\{(1,[i]), (2,[i])\}\rightarrow \{(1,[\alpha,i]),(2,[\alpha,i])\}$.
\\
\\(iii) $v=(p,[i]), w=(q,[\beta,j])$ for each $p=1,2,3$ and $q=1,2$.
\indent\par By (\ref{eqC6.4.5}) and (\ref{eqC6.4.8}), $(2,[i])\rightarrow \lambda_1\rightarrow \{(2,[j]),(3,[j])\}$ and $(1,[i])\rightarrow \bar{\lambda}_1\rightarrow \{(1,[j]),(4,[j])\}$. With (\ref{eqC6.4.6}), since $|\lambda_{i+2-{{s}\choose{\lceil{s/2}\rceil}}}|=|\lambda_{1}|=\lceil{\frac{s}{2}}\rceil$, there exists some $x_i\in \lambda_{i+2-{{s}\choose{\lceil{s/2}\rceil}}}\cap\lambda_{1}$, so that $(3,[i])\rightarrow x_i\rightarrow \{(2,[j]), (3,[j])\}$. With $\{(2,[j]),(4,[j])\}\rightarrow (1,[\beta,j])$ and $\{(1,[j]), (3,[j])\}$ $\rightarrow (2,[\beta,j])$ by (\ref{eqC6.4.7}), this subcase follows.
\\
\\(iv) $v=(q,[\beta,j]), w=(p,[i])$ for each $p=1,2,3$ and $q=1,2$.
\indent\par By (\ref{eqC6.4.5})-(\ref{eqC6.4.8}), $(2,[\beta,j])\rightarrow \{(2,[j]),(4,[j])\}$, $(1,[\beta,j])\rightarrow \{(1,[j]), (3,[j])\}$, $\{(1,[j]), (4,[j])\}$ $\rightarrow \lambda_1$, $\{(2,[j]), (3,[j])\}\rightarrow \bar{\lambda}_1$ and $|I^\mathtt{c}((p,[i]))|>0$.
\\
\\(v) $v=(p,[\alpha,i]), w=(q,[\beta,j])$ for each $p,q=1,2$.
\indent\par From (i), $d_D((p,[\alpha,i]),(r,[j]))=3$ for $r=2,3$. Since $(2,[j])\rightarrow (1,[\beta,j])$ and $(3,[j])\rightarrow (2,[\beta,j])$ by (\ref{eqC6.4.7}), this subcase follows.
\\
\\(vi) $v=(q,[\beta,j]), w=(p,[\alpha,i])$ for each $p,q=1,2$.
\indent\par From (iv), $d_D((q,[\beta,j]),(1,[i]))=3$. Since $(1,[i])\rightarrow \{(1,[\alpha,i]), (2,[\alpha,i])\}$ by (\ref{eqC6.4.4}), this subcase follows.
\\
\\Case 8.  For each $[i]\in A_3$, each $1\le\alpha\le \deg_T([i])-1$, and each $[j]\in E$,
\\(i) $v=(p,[\alpha,i]), w=(q,[j])$ for each $p,q=1,2$.
\\(ii) $v=(q,[j]), w=(p,[\alpha,i])$ for each $p,q=1,2$.
\indent\par Let $[k]\in A_4$. Since $\lambda_1\rightarrow \{(q,[j]),(2,[k])\}\rightarrow \bar{\lambda}_1$ by (\ref{eqC6.4.8})-(\ref{eqC6.4.9}), this case follows from Cases 6(i)-(ii) and 7(i)-(ii).
\\
\\Case 9. For each $[i]\in A_4$, each $1\le \alpha\le \deg_T([i])-1$, and each $[j]\in E$,
\\(i) $v=(p,[\alpha,i]), w=(q,[j])$ for each $p=1,2$ and $q=1,2$.
\\(ii) $v=(q,[j]), w=(p,[\alpha,i])$ for each $p=1,2$ and $q=1,2$.
\indent\par Since $\lambda_1\rightarrow \{(q,[j]),(2,[1])\}\rightarrow \bar{\lambda}_1$ by (\ref{eqC6.4.2}) and (\ref{eqC6.4.9}), this case follows from Cases 6(iii)-(iv).
\\
\\Case 10. $v=(r_1,\mathtt{c})$ and $w=(r_2,\mathtt{c})$ for $r_1\neq r_2$ and $1\le r_1, r_2\le s$.
\indent\par Here, we want to prove a stronger claim, $d_{D}((r_1,\mathtt{c}), (r_2,\mathtt{c}))=2$. Let $x_1=(2,[1])$, $x_{k+1}=(3,[k])$ for all $1\le k\le s-1$. Observe from (\ref{eqC6.4.2})-(\ref{eqC6.4.3}) that $\lambda_k\rightarrow x_k\rightarrow \bar{\lambda}_k$ for $1\le k\le s$ and the subgraph induced by $V_1=(\mathbb{N}_s,\mathtt{c})$ and $V_2=\{x_i\mid 1\le i \le s\}$ is a complete bipartite graph $K(V_1,V_2)$. By Lemma \ref{lemC6.2.19}, $d_{D}((r_1,\mathtt{c}), (r_2,\mathtt{c}))=2$.
\\
\\Case 11. $v\in \{(1,[i]), (2,[i]), (3,[i]), (4,[i]), (1,[\alpha,i]), (2,[\alpha,i])\}$ for each $1\le i\le \deg_T(\mathtt{c})$ and $1\le\alpha\le \deg_T([i])-1$, and $w=(r,\mathtt{c})$ for $1\le r\le s$.
\indent\par Note that there exists some $1\le k\le s$ such that $d_{D}(v,(k,\mathtt{c}))\le 2$, and $d_{D}((k,\mathtt{c}),w)\le 2$ by Case 10. Hence, it follows that $d_{D}(v,w)\le d_{D}(v,(k,\mathtt{c}))+d_D((k,\mathtt{c}),w)\le 4$.
\\
\\Case 12. $v=(r,\mathtt{c})$ for $r=1,2,\ldots, s$ and $w\in \{(1,[i]), (2,[i]), (3,[i]), (4,[i]), (1,[\alpha,i]), (2,[\alpha,i])\}$ for each $1\le i\le \deg_T(\mathtt{c})$ and $1\le\alpha\le \deg_T([i])-1$.
\indent\par Note that there exists some $1\le k\le s$ such that $d_{D}((k,\mathtt{c}), w)\le 2$, and $d_{D}(v,(k,\mathtt{c}))\le 2$ by Case 10. Hence, it follows that $d_{D}(v,w)\le d_{D}(v,(k,\mathtt{c}))+d_{D}((k,\mathtt{c}),w)\le 4$.
\\
\\Case 13. $v=(p,[i])$ and $w=(q, [j])$, where $1\le p,q\le 4$ and $1\le i,j\le \deg_T(\mathtt{c})$.
\noindent\par This follows from the fact that $|O^\mathtt{c}((p,[i]))|>0$, $|I^\mathtt{c}((q,[j]))|>0$, and $d_{D}((r_1,\mathtt{c}), (r_2,\mathtt{c}))$ $=2$ for any $r_1\neq r_2$ and $1\le r_1, r_2\le s$.
\\
\indent\par Therefore, the claim follows. Since every vertex lies in a directed $C_4$ for $D$ and $d(D)=4$, $\bar{d}(\mathcal{T})\le \max \{4, d(D)\}$ by Lemma \ref{lemC6.1.3}, and thus $\bar{d}(\mathcal{T})=4$.
\qed

\begin{center}
\begin{tikzpicture}[thick,scale=0.7]%
\draw(0,2)node[circle, draw, inner sep=0pt, minimum width=3pt](1_u){\scriptsize $(1,\mathtt{c})$};
\draw(0,-4)node[circle, draw, inner sep=0pt, minimum width=3pt](2_u){\scriptsize $(2,\mathtt{c})$};
\draw(0,-10)node[circle, draw, inner sep=0pt, minimum width=3pt](3_u){\scriptsize $(3,\mathtt{c})$};

\draw(-6,5)node[circle, draw, fill=black!100, inner sep=0pt, minimum width=6pt, label={[] 180:{\small $(1,[1,1])$}}](1_11){};
\draw(-6,3)node[circle, draw, fill=black!100, inner sep=0pt, minimum width=6pt, label={[] 180:{\small $(2,[1,1])$}}](2_11){};

\draw(-3,6)node[circle, draw, fill=black!100, inner sep=0pt, minimum width=6pt, label={[yshift=0cm, xshift=-0.45cm] 90:{\small $(1,[1])$}}](1_1){};
\draw(-3,4)node[circle, draw, fill=black!100, inner sep=0pt, minimum width=6pt, label={[yshift=-0.1cm, xshift=-0.45cm] 270:{\small $(2,[1])$}}](2_1){};
\draw(-3,2)node[circle, draw, fill=black!100, inner sep=0pt, minimum width=6pt, label={[yshift=0cm, xshift=-0.45cm] 270:{\small $(3,[1])$}}](3_1){};

\draw(6,5)node[circle, draw, fill=black!100, inner sep=0pt, minimum width=6pt, label={[] 0:{\small $(1,[1,2])$}}](1_21){};
\draw(6,3)node[circle, draw, fill=black!100, inner sep=0pt, minimum width=6pt, label={[] 0:{\small $(2,[1,2])$}}](2_21){};

\draw(3,6)node[circle, draw, fill=black!100, inner sep=0pt, minimum width=6pt, label={[yshift=0cm, xshift=0.45cm] 90:{\small $(1,[2])$}}](1_2){};
\draw(3,4)node[circle, draw, fill=black!100, inner sep=0pt, minimum width=6pt, label={[yshift=-0.1cm, xshift=0.45cm] 270:{\small $(2,[2])$}}](2_2){};
\draw(3,2)node[circle, draw, fill=black!100, inner sep=0pt, minimum width=6pt, label={[yshift=0cm, xshift=0.45cm] 270:{\small $(3,[2])$}}](3_2){};

\draw(-6,-1)node[circle, draw, fill=black!100, inner sep=0pt, minimum width=6pt, label={[] 180:{\small $(1,[1,3])$}}](1_31){};
\draw(-6,-3)node[circle, draw, fill=black!100, inner sep=0pt, minimum width=6pt, label={[] 180:{\small $(2,[1,3])$}}](2_31){};

\draw(-3,-0)node[circle, draw, fill=black!100, inner sep=0pt, minimum width=6pt, label={[yshift=-0.15cm, xshift=-0.45cm] 270:{\small $(1,[3])$}}](1_3){};
\draw(-3,-2)node[circle, draw, fill=black!100, inner sep=0pt, minimum width=6pt, label={[yshift=-0.1cm, xshift=-0.45cm] 270:{\small $(2,[3])$}}](2_3){};
\draw(-3,-4)node[circle, draw, fill=black!100, inner sep=0pt, minimum width=6pt, label={[yshift=0cm, xshift=-0.45cm] 270:{\small $(3,[3])$}}](3_3){};

\draw(6,-1)node[circle, draw, fill=black!100, inner sep=0pt, minimum width=6pt, label={[] 0:{\small $(1,[1,4])$}}](1_41){};
\draw(6,-3)node[circle, draw, fill=black!100, inner sep=0pt, minimum width=6pt, label={[] 0:{\small $(2,[1,4])$}}](2_41){};

\draw(3,-0)node[circle, draw, fill=black!100, inner sep=0pt, minimum width=6pt, label={[yshift=-0.15cm, xshift=0.45cm] 270:{\small $(1,[4])$}}](1_4){};
\draw(3,-2)node[circle, draw, fill=black!100, inner sep=0pt, minimum width=6pt, label={[yshift=-0.1cm, xshift=0.45cm] 270:{\small $(2,[4])$}}](2_4){};
\draw(3,-4)node[circle, draw, fill=black!100, inner sep=0pt, minimum width=6pt, label={[yshift=0cm, xshift=0.45cm] 270:{\small $(3,[4])$}}](3_4){};

\draw(-6,-8)node[circle, draw, fill=black!100, inner sep=0pt, minimum width=6pt, label={[] 180:{\small $(1,[1,5])$}}](1_51){};
\draw(-6,-10)node[circle, draw, fill=black!100, inner sep=0pt, minimum width=6pt, label={[] 180:{\small $(2,[1,5])$}}](2_51){};

\draw(-3,-6)node[circle, draw, fill=black!100, inner sep=0pt, minimum width=6pt, label={[yshift=-0.1cm, xshift=-0.45cm] 270:{\small $(1,[5])$}}](1_5){};
\draw(-3,-8)node[circle, draw, fill=black!100, inner sep=0pt, minimum width=6pt, label={[yshift=0cm, xshift=-0.45cm] 270:{\small $(2,[5])$}}](2_5){};
\draw(-3,-10)node[circle, draw, fill=black!100, inner sep=0pt, minimum width=6pt, label={[yshift=0.05cm, xshift=-0.45cm] 270:{\small $(3,[5])$}}](3_5){};
\draw(-3,-12)node[circle, draw, fill=black!100, inner sep=0pt, minimum width=6pt, label={[yshift=0cm, xshift=-0.45cm] 270:{\small $(4,[5])$}}](4_5){};

\draw(6,-8)node[circle, draw, fill=black!100, inner sep=0pt, minimum width=6pt, label={[] 0:{\small $(1,[1,6])$}}](1_61){};
\draw(6,-10)node[circle, draw, fill=black!100, inner sep=0pt, minimum width=6pt, label={[] 0:{\small $(2,[1,6])$}}](2_61){};

\draw(3,-6)node[circle, draw, fill=black!100, inner sep=0pt, minimum width=6pt, label={[yshift=-0.1cm, xshift=0.45cm] 270:{\small $(1,[6])$}}](1_6){};
\draw(3,-8)node[circle, draw, fill=black!100, inner sep=0pt, minimum width=6pt, label={[yshift=0cm, xshift=0.45cm] 270:{\small $(2,[6])$}}](2_6){};
\draw(3,-10)node[circle, draw, fill=black!100, inner sep=0pt, minimum width=6pt, label={[yshift=0.05cm, xshift=0.45cm] 270:{\small $(3,[6])$}}](3_6){};
\draw(3,-12)node[circle, draw, fill=black!100, inner sep=0pt, minimum width=6pt, label={[yshift=0cm, xshift=0.45cm] 270:{\small $(4,[6])$}}](4_6){};

\draw(-3,-14)node[circle, draw, fill=black!100, inner sep=0pt, minimum width=6pt, label={[] 180:{\small $(1,[7])$}}](1_7){};
\draw(-3,-16)node[circle, draw, fill=black!100, inner sep=0pt, minimum width=6pt, label={[] 180:{\small $(2,[7])$}}](2_7){};

\draw(3,-14)node[circle, draw, fill=black!100, inner sep=0pt, minimum width=6pt, label={[] 0:{\small $(1,[8])$}}](1_8){};
\draw(3,-16)node[circle, draw, fill=black!100, inner sep=0pt, minimum width=6pt, label={[] 0:{\small $(2,[8])$}}](2_8){};

\draw[->, line width=0.3mm, >=latex, shorten <= 0.2cm, shorten >= 0.15cm](1_11)--(1_1);
\draw[->, line width=0.3mm, >=latex, shorten <= 0.2cm, shorten >= 0.15cm](1_11)--(2_1);
\draw[dashed,->, line width=0.3mm, >=latex, shorten <= 0.2cm, shorten >= 0.15cm](3_1) to [out=158, in=275] (1_11);
\draw[->, line width=0.3mm, >=latex, shorten <= 0.2cm, shorten >= 0.15cm](2_11)--(1_1);
\draw[->, line width=0.3mm, >=latex, shorten <= 0.2cm, shorten >= 0.15cm](2_11)--(2_1);
\draw[dashed,->, line width=0.3mm, >=latex, shorten <= 0.2cm, shorten >= 0.15cm](3_1)--(2_11);

\draw[->, line width=0.3mm, >=latex, shorten <= 0.2cm, shorten >= 0.15cm](1_21)--(1_2);
\draw[->, line width=0.3mm, >=latex, shorten <= 0.2cm, shorten >= 0.15cm](1_21)--(2_2);
\draw[dashed,->, line width=0.3mm, >=latex, shorten <= 0.2cm, shorten >= 0.15cm](3_2) to [out=22, in=265] (1_21);
\draw[->, line width=0.3mm, >=latex, shorten <= 0.2cm, shorten >= 0.15cm](2_21)--(1_2);
\draw[->, line width=0.3mm, >=latex, shorten <= 0.2cm, shorten >= 0.15cm](2_21)--(2_2);
\draw[dashed,->, line width=0.3mm, >=latex, shorten <= 0.2cm, shorten >= 0.15cm](3_2)--(2_21);

\draw[->, line width=0.3mm, >=latex, shorten <= 0.2cm, shorten >= 0.15cm](1_31) to [out=275, in=158] (3_3);
\draw[->, line width=0.3mm, >=latex, shorten <= 0.2cm, shorten >= 0.15cm](2_31)--(3_3);
\draw[dashed,->, line width=0.3mm, >=latex, shorten <= 0.2cm, shorten >= 0.15cm](1_3)--(1_31);
\draw[dashed,->, line width=0.3mm, >=latex, shorten <= 0.2cm, shorten >= 0.15cm](1_3) to [out=202, in=85] (2_31);
\draw[dashed,->, line width=0.3mm, >=latex, shorten <= 0.2cm, shorten >= 0.15cm](2_3)--(1_31);
\draw[dashed,->, line width=0.3mm, >=latex, shorten <= 0.2cm, shorten >= 0.15cm](2_3)--(2_31);

\draw[->, line width=0.3mm, >=latex, shorten <= 0.2cm, shorten >= 0.15cm](1_41) to [out=265, in=22] (3_4);
\draw[->, line width=0.3mm, >=latex, shorten <= 0.2cm, shorten >= 0.15cm](2_41)--(3_4);
\draw[dashed,->, line width=0.3mm, >=latex, shorten <= 0.2cm, shorten >= 0.15cm](1_4)--(1_41);
\draw[dashed,->, line width=0.3mm, >=latex, shorten <= 0.2cm, shorten >= 0.15cm](1_4) to [out=338, in=95] (2_41);
\draw[dashed,->, line width=0.3mm, >=latex, shorten <= 0.2cm, shorten >= 0.15cm](2_4)--(1_41);
\draw[dashed,->, line width=0.3mm, >=latex, shorten <= 0.2cm, shorten >= 0.15cm](2_4)--(2_41);

\draw[dashed,->, line width=0.3mm, >=latex, shorten <= 0.2cm, shorten >= 0.15cm](2_5)--(1_51);
\draw[dashed,->, line width=0.3mm, >=latex, shorten <= 0.2cm, shorten >= 0.15cm](1_5) to [out=195, in=80] (2_51);
\draw[dashed,->, line width=0.3mm, >=latex, shorten <= 0.2cm, shorten >= 0.15cm](4_5) to [out=140, in=290] (1_51);
\draw[dashed,->, line width=0.3mm, >=latex, shorten <= 0.2cm, shorten >= 0.15cm](3_5)--(2_51);
\draw[->, line width=0.3mm, >=latex, shorten <= 0.2cm, shorten >= 0.15cm](1_51) to [out=65, in=180] (1_5);
\draw[->, line width=0.3mm, >=latex, shorten <= 0.2cm, shorten >= 0.15cm](2_51) to [out=65, in=190] (2_5);
\draw[->, line width=0.3mm, >=latex, shorten <= 0.2cm, shorten >= 0.15cm](1_51)--(3_5);
\draw[->, line width=0.3mm, >=latex, shorten <= 0.2cm, shorten >= 0.15cm](2_51)--(4_5);

\draw[dashed,->, line width=0.3mm, >=latex, shorten <= 0.2cm, shorten >= 0.15cm](2_6)--(1_61);
\draw[dashed,->, line width=0.3mm, >=latex, shorten <= 0.2cm, shorten >= 0.15cm](1_6) to [out=345, in=100] (2_61);
\draw[dashed,->, line width=0.3mm, >=latex, shorten <= 0.2cm, shorten >= 0.15cm](4_6) to [out=40, in=250] (1_61);
\draw[dashed,->, line width=0.3mm, >=latex, shorten <= 0.2cm, shorten >= 0.15cm](3_6)--(2_61);
\draw[->, line width=0.3mm, >=latex, shorten <= 0.2cm, shorten >= 0.15cm](1_61) to [out=115, in=0] (1_6);
\draw[->, line width=0.3mm, >=latex, shorten <= 0.2cm, shorten >= 0.15cm](2_61) to [out=115, in=350] (2_6);
\draw[->, line width=0.3mm, >=latex, shorten <= 0.2cm, shorten >= 0.15cm](1_61)--(3_6);
\draw[->, line width=0.3mm, >=latex, shorten <= 0.2cm, shorten >= 0.15cm](2_61)--(4_6);

\draw[->, line width=0.3mm, >=latex, shorten <= 0.2cm, shorten >= 0.1cm](1_1)--(1_u);
\draw[->, line width=0.3mm, >=latex, shorten <= 0.2cm, shorten >= 0.1cm](1_1)--(2_u);
\draw[->, line width=0.3mm, >=latex, shorten <= 0.2cm, shorten >= 0.1cm](2_1)--(3_u);
\draw[densely dashdotdotted, ->, line width=0.3mm, >=latex, shorten <= 0.2cm, shorten >= 0.1cm](3_1)--(1_u);

\draw[->, line width=0.3mm, >=latex, shorten <= 0.2cm, shorten >= 0.1cm](1_2)--(1_u);
\draw[->, line width=0.3mm, >=latex, shorten <= 0.2cm, shorten >= 0.1cm](1_2)--(2_u);
\draw[->, line width=0.3mm, >=latex, shorten <= 0.2cm, shorten >= 0.1cm](2_2)--(3_u);
\draw[densely dashdotdotted, ->, line width=0.3mm, >=latex, shorten <= 0.2cm, shorten >= 0.1cm](3_2)--(2_u);

\draw[->, line width=0.3mm, >=latex, shorten <= 0.2cm, shorten >= 0.1cm](1_3)--(3_u);
\draw[->, line width=0.3mm, >=latex, shorten <= 0.2cm, shorten >= 0.1cm](2_3)--(1_u);
\draw[->, line width=0.3mm, >=latex, shorten <= 0.2cm, shorten >= 0.1cm](2_3)--(2_u);
\draw[densely dotted, ->, line width=0.3mm, >=latex, shorten <= 0.2cm, shorten >= 0.1cm](3_3)--(2_u);
\draw[densely dotted, ->, line width=0.3mm, >=latex, shorten <= 0.2cm, shorten >= 0.1cm](3_3)--(3_u);

\draw[->, line width=0.3mm, >=latex, shorten <= 0.2cm, shorten >= 0.1cm](1_4)--(3_u);
\draw[->, line width=0.3mm, >=latex, shorten <= 0.2cm, shorten >= 0.1cm](2_4)--(1_u);
\draw[->, line width=0.3mm, >=latex, shorten <= 0.2cm, shorten >= 0.1cm](2_4)--(2_u);
\draw[densely dotted, ->, line width=0.3mm, >=latex, shorten <= 0.2cm, shorten >= 0.1cm](3_4)--(1_u);
\draw[densely dotted, ->, line width=0.3mm, >=latex, shorten <= 0.2cm, shorten >= 0.1cm](3_4)--(3_u);

\draw[->, line width=0.3mm, >=latex, shorten <= 0.2cm, shorten >= 0.1cm](2_5)--(3_u);
\draw[->, line width=0.3mm, >=latex, shorten <= 0.2cm, shorten >= 0.1cm](3_5)--(3_u);
\draw[->, line width=0.3mm, >=latex, shorten <= 0.2cm, shorten >= 0.1cm](1_5)--(1_u);
\draw[->, line width=0.3mm, >=latex, shorten <= 0.2cm, shorten >= 0.1cm](1_5)--(2_u);
\draw[->, line width=0.3mm, >=latex, shorten <= 0.2cm, shorten >= 0.1cm](4_5)--(1_u);
\draw[->, line width=0.3mm, >=latex, shorten <= 0.2cm, shorten >= 0.1cm](4_5)--(2_u);

\draw[->, line width=0.3mm, >=latex, shorten <= 0.2cm, shorten >= 0.1cm](2_6)--(3_u);
\draw[->, line width=0.3mm, >=latex, shorten <= 0.2cm, shorten >= 0.1cm](3_6)--(3_u);
\draw[->, line width=0.3mm, >=latex, shorten <= 0.2cm, shorten >= 0.1cm](1_6)--(1_u);
\draw[->, line width=0.3mm, >=latex, shorten <= 0.2cm, shorten >= 0.1cm](1_6)--(2_u);
\draw[->, line width=0.3mm, >=latex, shorten <= 0.2cm, shorten >= 0.1cm](4_6)--(1_u);
\draw[->, line width=0.3mm, >=latex, shorten <= 0.2cm, shorten >= 0.1cm](4_6)--(2_u);
\draw[->, line width=0.3mm, >=latex, shorten <= 0.2cm, shorten >= 0.1cm](1_7)--(3_u);
\draw[->, line width=0.3mm, >=latex, shorten <= 0.2cm, shorten >= 0.1cm](2_7)--(3_u);
\draw[->, line width=0.3mm, >=latex, shorten <= 0.2cm, shorten >= 0.1cm](1_8)--(3_u);
\draw[->, line width=0.3mm, >=latex, shorten <= 0.2cm, shorten >= 0.1cm](2_8)--(3_u);
\end{tikzpicture}
\captionsetup{justification=centering}
{\captionof{figure}{Orientation $D$ for $\mathcal{H}$ for $s=3$,
\\$A_3=\{[1],[2],[3],[4]\}$, $A_4=\{[5],[6]\}$, E=\{[7],[8]\}.}\label{figC6.4.12}}
\end{center}

\begin{ppn}\label{ppnC6.4.3}
Suppose $s\ge 3$ is odd, $A_2\neq \emptyset$, $A_3\neq \emptyset$, and $A_{\ge 4}=\emptyset$ for a $\mathcal{T}$. Then, $\mathcal{T}\in \mathscr{C}_0$ if and only if
\begin{equation}
\left\{
  \begin{array}{@{}ll@{}}
    |A_2|\le{{s}\choose{\lceil{s/2}\rceil}}-1, & \text{if}\ |A_3|=1,\nonumber\\
    (i)\ 2|A_2|+|A_3|\le 2{{s}\choose{\lceil{s/2}\rceil}}-2, \text{ or } \nonumber\\
    (ii)\ 2|A_2|+|A_3|= 2{{s}\choose{\lceil{s/2}\rceil}}-1, |A_2|\ge \lceil\frac{s}{2}\rceil \lfloor\frac{s}{2}\rfloor \text{ and } s\ge 5, & \text{if}\ |A_3|\ge 2.\nonumber\\
  \end{array}\right.
\end{equation}
\end{ppn}
\noindent\textit{Proof}: $(\Rightarrow)$ Since $\mathcal{T}\in \mathscr{C}_0$, there exists an orientation $D$ of $\mathcal{T}$, where $d(D)=4$. As $A_2\neq\emptyset$ and $A_3\neq\emptyset$, we assume (\ref{eqC6.3.1})-(\ref{eqC6.3.5}) here, unless stated otherwise.
\\
\\Case 1: $|A_3|=1$.
\indent\par Let $A_3=\{[j]\}$. By Lemma \ref{lemC6.2.13}(b), either $|O((1,[j]))|=1$ or $|I((1,[j]))|=1$. If $|O((1, [1,j]))|=1$, say $O((1, [1,j]))=\{(1, [j])\}$ (instead of (\ref{eqC6.3.4})), then by Lemma \ref{lemC6.2.15}, $\{O^\mathtt{c}((1,[i]))\mid [i]\in A_2\cup A_3\}$ is an antichain. So, $|A_2|+|A_3|\le {{s}\choose{\lfloor{s/2}\rfloor}}$ by Sperner's theorem, i.e., $|A_2|\le{{s}\choose{\lfloor{s/2}\rfloor}}-1$. If $|I((1, [j]))|=1$, then this case follows by the Duality Lemma.
\begin{rmk}
The outline of Case 2 is largely similar to Proposition \ref{ppnC6.4.1}. In particular, Remark \ref{rmkC6.3.3} applies here. Hence, we will make reference to Proposition \ref{ppnC6.4.1} as we proceed. We first show that for any $[i]\in A_2\cup A_3$ and any $p=1,2,3$ (wherever applicable), if $|O^\mathtt{c}((p,[i]))|$ is too big ($>\lceil\frac{s}{2}\rceil$) or too small ($<\lfloor\frac{s}{2}\rfloor$), then $2|A_2|+|A_3|\le 2{{s}\choose{\lceil{s/2}\rceil}}-2$. Here, observe that $B^O_2\cup B^O_3$ ($B^I_2\cup B^I_3$ resp.) plays an analogous role of $B^O_3$ ($B^I_3$ resp.) in Proposition \ref{ppnC6.4.1}.
\end{rmk}
\noindent Case 2: $|A_3|\ge 2$.
\indent\par Note that $B^O_2\cup B^O_3$ and $B^I_2\cup B^I_3$ are antichains by Lemmas \ref{lemC6.2.15} and \ref{lemC6.2.16} respectively. Hence, $|A_2|+|A^O_3|=|B^O_2|+|B^O_3|\le {{s}\choose{\lfloor{s/2}\rfloor}}$ and $|A_2|+|A^I_3|=|B^I_2|+|B^I_3|\le {{s}\choose{\lfloor{s/2}\rfloor}}$ by Sperner's theorem. So, if $|A^O_3|=0$ ($|A^I_3|=0$ resp.), then $|A_2|+|A^I_3|\le {{s}\choose{\lfloor{s/2}\rfloor}}$ ($|A_2|+|A^O_3|\le {{s}\choose{\lfloor{s/2}\rfloor}}$ resp.) implies that $|A_2|\le {{s}\choose{\lfloor{s/2}\rfloor}}-2$. It follows that  $2|A_2|+|A_3|=(|A_2|+|A^O_3|)+(|A_2|+|A^I_3|)\le [{{s}\choose{\lfloor{s/2}\rfloor}}-2]+{{s}\choose{\lfloor{s/2}\rfloor}}=2{{s}\choose{\lceil{s/2}\rceil}}-2$. Therefore, we assume $|A^O_3|=|B^O_3|>0$ and $|A^I_3|=|B^I_3|>0$.
\\
\\Subcase 2.1. There exists some $1\le i\le \deg_T(\mathtt{c})$ such that  $|O^\mathtt{c}((p,[i]))|>\lceil\frac{s}{2}\rceil$ for some $p=1,2,3$.
\indent\par If $[i]\in A_2\cup A^O_3$, then $d_D((1,[1,j]),(p,[i]))=3$ for any $[j]\in A_2\cup A^O_3-\{[i]\}$ implies $X\cap I^\mathtt{c}((p,[i]))\neq \emptyset$ for each $X\in B^O_2\cup B^O_3-\{O^\mathtt{c}((\delta,[i]))\}$, where
\begin{align}
\delta=\left\{
  \begin{array}{@{}ll@{}}
    1, & \text{if}\ [i]\in A_2,\\
    3, & \text{if}\ [i]\in A^O_3.\\
  \end{array}\right.
\label{eqC6.4.10}
\end{align}
By Lih's theorem, $|A_2|+|A^O_3|-1=|B^O_2\cup B^O_3-\{O^\mathtt{c}((\delta,[i]))\}|\le {{s}\choose{\lceil{s/2}\rceil}}-{{s-|I^\mathtt{c}((p,[i]))|}\choose{\lceil{s/2}\rceil}}\le {{s}\choose{\lceil{s/2}\rceil}}-{{\lceil{s/2}\rceil+1}\choose{\lceil{s/2}\rceil}}={{s}\choose{\lceil{s/2}\rceil}}-(\lceil{\frac{s}{2}}\rceil+1)\le {{s}\choose{\lceil{s/2}\rceil}}-3$. It follows that $|A_2|+|A^O_3|\le {{s}\choose{\lceil{s/2}\rceil}}-2$.
\noindent\par If $[i]\not\in A_2\cup A^O_3$,  then $d_D((1,[1,j]),(p,[i]))=3$ for each $[j]\in A_2\cup A^O_3$,  implies $X\cap I^\mathtt{c}((p,[i]))\neq \emptyset$ for each $X\in B^O_2\cup B^O_3$. Hence, by Lih's theorem, $|A_2|+|A^O_3|=|B^O_2\cup B^O_3|\le {{s}\choose{\lceil{s/2}\rceil}}-{{s-|I^\mathtt{c}((p,[i]))|}\choose{\lceil{s/2}\rceil}}\le {{s}\choose{\lceil{s/2}\rceil}}-{{\lceil{s/2}\rceil+1}\choose{\lceil{s/2}\rceil}}={{s}\choose{\lceil{s/2}\rceil}}-(\lceil{\frac{s}{2}}\rceil+1)\le {{s}\choose{\lceil{s/2}\rceil}}-3$. 
\noindent\par In both subcases, $2|A_2|+|A_3|=(|A_2|+|A^O_3|)+(|A_2|+|A^I_3|)\le {{s}\choose{\lceil{s/2}\rceil}}-2+{{s}\choose{\lfloor{s/2}\rfloor}}=2{{s}\choose{\lceil{s/2}\rceil}}-2$.
\\
\\Subcase 2.2. There exists some $1\le i\le \deg_T(\mathtt{c})$ such that $|O^\mathtt{c}((p,[i]))|<\lfloor\frac{s}{2}\rfloor$ for some $p=1,2,3$.
\indent\par In other words, $|I^\mathtt{c}((p,[i]))|>\lfloor\frac{s}{2}\rfloor$. The result follows from Subcase 2.1 by the Duality Lemma.

\begin{rmk}
At this stage of Proposition \ref{ppnC6.4.1}, we invoked Theorems \ref{thmC6.2.4} and \ref{thmC6.2.5} to conclude $B^O_3\cup B^I_3\subseteq {{(\mathbb{N}_s,\mathtt{c})}\choose{\lceil{s/2}\rceil}}$. However, the two theorems cannot apply here because for $[i]\in A_2$, and $d_D((1,[i]),(2,[i]))\le 4$, it is not necessary that $O^\mathtt{c}((1,[i]))\cap I^\mathtt{c}((2,[i]))\neq\emptyset$. Consequently, $B^O_2\cup B^O_3$ and $B^I_2\cup B^I_3$ may not be cross-intersecting.
\indent\par So, we need more work (i.e., Subcases 2.3.1-2.3.5 below) before arriving at the desired last case of $B^O_2\cup B^O_3\cup B^I_2\cup B^I_3 \subseteq {{(\mathbb{N}_s,\mathtt{c})}\choose{\lceil{s/2}\rceil}}$. Retrospectively, we could have used such a similar, yet more lengthy, argument in Proposition \ref{ppnC6.4.1}.
\end{rmk}

\noindent Subcase 2.3. For all $1\le i\le \deg_T(\mathtt{c})$ and all $p=1,2,3$, $\lfloor\frac{s}{2}\rfloor\le|O^\mathtt{c}((p,[i]))|\le \lceil\frac{s}{2}\rceil$.
Partition $A^O_3$ ($A^I_3$ resp.) into $A^{O(D)}_3$ and $A^{O(S)}_3$ ($A^{I(D)}_3$ and $A^{I(S)}_3$ resp.) as in  (\ref{eqC6.3.33}).
\\
\\Assumption 1: $|A^{O(D)}_3|\ge 1$ and $|A^{I(D)}_3|\ge 1$.
\indent\par Suppose $A^{O(D)}_3=\emptyset$, i.e., $O^\mathtt{c}((1,[i]))=O^\mathtt{c}((2,[i]))$ for all $[i]\in A^O_3$. By Lemma \ref{lemC6.2.16}, $\{I^\mathtt{c}((1,[i]))\mid [i]\in A^O_3\}\cup B^I_2\cup B^I_3$ is an antichain. By Sperner's theorem, $|A_2|+|A_3|=|B^O_3|+|B^I_2|+|B^I_3|=|\{I^\mathtt{c}((1,[i]))\mid [i]\in A^O_3\}\cup B^I_2\cup B^I_3|\le {{s}\choose{\lfloor{s/2}\rfloor}}$. Since $|A_3|\ge 2$, this implies $|A_2|\le {{s}\choose{\lfloor{s/2}\rfloor}}-2$. Therefore, $2|A_2|+|A_3|\le {{s}\choose{\lfloor{s/2}\rfloor}}+[{{s}\choose{\lfloor{s/2}\rfloor}}-2]=2{{s}\choose{\lceil{s/2}\rceil}}-2$. A similar argument follows if $A^{I(D)}_3=\emptyset$.
\\
\\Subcase 2.3.1. There exist $X^O, Y^O\in B^O_2\cup B^O_3$ and $X^I, Y^I\in B^I_2\cup B^I_3$, such that $|X^O|=|X^I|=\lfloor\frac{s}{2}\rfloor$ and $|Y^O|=|Y^I|=\lceil\frac{s}{2}\rceil$.
\indent\par By Sperner's theorem, $|A_2|+|A^O_3|=|B^O_2\cup B^O_3|\le {{s}\choose{\lfloor{s/2}\rfloor}}-1$ and $|A_2|+|A^I_3|=|B^I_2\cup B^I_3|\le {{s}\choose{\lfloor{s/2}\rfloor}}-1$. Hence, $2|A_2|+|A_3|=(|A_2|+|A^O_3|)+(|A_2|+|A^I_3|)\le 2{{s}\choose{\lceil{s/2}\rceil}}-2$.
\\
\\Subcase 2.3.2. $|X^*|=\lfloor\frac{s}{2}\rfloor$, $|Y^*|=\lceil\frac{s}{2}\rceil$ for some $X^*, Y^*\in B^O_2\cup B^O_3$ and $B^I_2\cup B^I_3\subseteq {{(\mathbb{N}_s,\mathtt{c})}\choose{\lfloor{s/2}\rfloor}}$, or for some $X^*, Y^*\in B^I_2\cup B^I_3$ and $B^O_2\cup B^O_3\subseteq {{(\mathbb{N}_s,\mathtt{c})}\choose{\lfloor{s/2}\rfloor}}$.
\indent\par Consider the former case; a similar argument follows for the latter. By Sperner's theorem, $|A_2|+|A^O_3|=|B^O_2\cup B^O_3|\le {{s}\choose{\lfloor{s/2}\rfloor}}-1$. Let $[i]\in A^{O(D)}_3$. Suppose $|O^\mathtt{c}((p,[i]))|=\lfloor\frac{s}{2}\rfloor$ for some $p=1,2,3$. Note that $d_D((p,[i]),(1,[1,j]))=3$ for each $[j]\in A_2\cup A^I_3$ implies $X\cap O^\mathtt{c}((p,[i]))$ for all $X\in B^I_2\cup B^I_3$. Since there are ${{\lceil{s/2}\rceil}\choose{\lfloor{s/2}\rfloor}}=\lceil{\frac{s}{2}}\rceil$ number of $\lfloor{\frac{s}{2}}\rfloor$-subsets of $I^\mathtt{c}((p,[i]))$, $|A_2|+|A^I_3|=|B^I_2\cup B^I_3|\le {{s}\choose{\lfloor{s/2}\rfloor}}-\lceil{\frac{s}{2}}\rceil\le {{s}\choose{\lfloor{s/2}\rfloor}}-2$. Therefore, $2|A_2|+|A_3|=(|A_2|+|A^O_3|)+(|A_2|+|A^I_3|)\le [{{s}\choose{\lfloor{s/2}\rfloor}}-1]+[{{s}\choose{\lfloor{s/2}\rfloor}}-2]=2{{s}\choose{\lceil{s/2}\rceil}}-3$. 
\indent\par Suppose $|O^\mathtt{c}((p,[i]))|=\lceil\frac{s}{2}\rceil$ for all $p=1,2,3$. Similarly, for $p=1,2$, $d_D((p,[i]),$ $(1,[1,j]))=3$ for each $[i]\in A_2\cup A^I_3$ implies $X\neq I^\mathtt{c}((p,[i]))$ for all $X\in B^I_2\cup B^I_3$ and $p=1,2$. So, $|A_2|+|A^I_3|=|B^I_2\cup B^I_3|\le {{s}\choose{\lfloor{s/2}\rfloor}}-2$, which implies $2|A_2|+|A_3|=(|A_2|+|A^O_3|)+(|A_2|+|A^I_3|)\le [{{s}\choose{\lfloor{s/2}\rfloor}}-1]+[{{s}\choose{\lfloor{s/2}\rfloor}}-2]=2{{s}\choose{\lceil{s/2}\rceil}}-3$.
\\
\\Subcase 2.3.3. $|X^*|=\lfloor\frac{s}{2}\rfloor$, $|Y^*|=\lceil\frac{s}{2}\rceil$ for some $X^*, Y^*\in B^O_2\cup B^O_3$ and $B^I_2\cup B^I_3\subseteq {{(\mathbb{N}_s,\mathtt{c})}\choose{\lceil{s/2}\rceil}}$, or for some $X^*, Y^*\in B^I_2\cup B^I_3$ and $B^O_2\cup B^O_3\subseteq {{(\mathbb{N}_s,\mathtt{c})}\choose{\lceil{s/2}\rceil}}$.
\indent\par Consider the former case; a similar argument follows for the latter. By Sperner's theorem, $|A_2|+|A^O_3|=|B^O_2\cup B^O_3|\le {{s}\choose{\lfloor{s/2}\rfloor}}-1$. Let $[i]\in A^{O(D)}_3$. Suppose $|O^\mathtt{c}((p,[i]))|=\lfloor\frac{s}{2}\rfloor$ for some $p=1,2,3$. Then, $d_D((p,[i]),(1,[1,j]))=3$ for each $[j]\in A_2\cup A^I_3$ implies $X\neq I^\mathtt{c}((p,[i]))$ for all $X\in B^I_2\cup B^I_3$. So, $|A_2|+|A^I_3|=|B^I_2\cup B^I_3|\le {{s}\choose{\lceil{s/2}\rceil}}-1$, which implies $2|A_2|+|A_3|=(|A_2|+|A^O_3|)+(|A_2|+|A^I_3|)\le [{{s}\choose{\lfloor{s/2}\rfloor}}-1]+[{{s}\choose{\lceil{s/2}\rceil}}-1]=2{{s}\choose{\lceil{s/2}\rceil}}-2$.
\nolinebreak\indent\par Suppose $|O^\mathtt{c}((p,[i]))|=\lceil\frac{s}{2}\rceil$ for all $p=1,2,3$. Note by definition of $A^{O(D)}_3$ that $O^\mathtt{c}((3,[i]))$ is equal to at most one of $O^\mathtt{c}((1,[i]))$ and $O^\mathtt{c}((2,[i]))$, say $O^\mathtt{c}((3,[i]))\neq O^\mathtt{c}((1,[i]))$. Then, $d_D((1,[1,j]),(1,[i]))=3$ for each $[j]\in A_2\cup A^O_3-\{[i]\}$ implies $X\cap I^\mathtt{c}((1,[i])$ for all $X\in B^O_2\cup B^O_3-\{O^\mathtt{c}((3,[i]))\}$. So, $X\not\subseteq O^\mathtt{c}((1,[i]))$ for all $X\in B^O_2\cup B^O_3$. It follows that $R=B^O_2\cup B^O_3\cup\{O^\mathtt{c}((1,[i]))\}$ is an antichain. Since there exist $X,Y\in R$ such that $|X|=\lfloor\frac{s}{2}\rfloor$, and $|Y|=\lceil\frac{s}{2}\rceil$, it follows from Sperner's theorem that $|B^O_2|+|B^O_3|+1=|R|\le {{s}\choose{\lfloor{s/2}\rfloor}}-1$. Hence, $|A_2|+|A^O_3|=|B^O_2|+|B^O_3|\le {{s}\choose{\lfloor{s/2}\rfloor}}-2$. This implies $2|A_2|+|A_3|=(|A_2|+|A^O_3|)+(|A_2|+|A^I_3|)\le [{{s}\choose{\lfloor{s/2}\rfloor}}-2]+{{s}\choose{\lfloor{s/2}\rfloor}}=2{{s}\choose{\lceil{s/2}\rceil}}-2$.
\\
\\Subcase 2.3.4. $B^O_2\cup B^O_3\cup B^I_2\cup B^I_3\subseteq {{(\mathbb{N}_s,\mathtt{c})}\choose{\lfloor{s/2}\rfloor}}$.
\indent\par Let $[i]\in A^O_3$. For any $[j]\in A_2\cup A^I_3$, $d_D((3,[i]),(1,[1,j]))=3$ implies $O^\mathtt{c}((3,[i]))\cap X\neq \emptyset$ for all $X\in B^I_2\cup B^I_3$. In particular, $X\not\subseteq I^\mathtt{c}((3,[i]))$ for all $X\in B^I_2\cup B^I_3$. Since there are ${{\lceil{s/2}\rceil}\choose{\lfloor{s/2}\rfloor}}=\lceil\frac{s}{2}\rceil$ number of $\lfloor \frac{s}{2}\rfloor$-subsets of $I^\mathtt{c}((3,[i]))$, it follows that $|A_2|+|A^I_3|=|B^I_2\cup B^I_3|\le{{s}\choose{\lfloor{s/2}\rfloor}}-\lceil\frac{s}{2}\rceil\le {{s}\choose{\lfloor{s/2}\rfloor}}-2$. Hence, $2|A_2|+|A_3|=(|A_2|+|A^O_3|)+(|A_2|+|A^I_3|)\le {{s}\choose{\lfloor{s/2}\rfloor}}+[{{s}\choose{\lfloor{s/2}\rfloor}}-2]=2{{s}\choose{\lceil{s/2}\rceil}}-2$.
\\
\\Subcase 2.3.5. $B^O_2\cup B^O_3\subseteq {{(\mathbb{N}_s,\mathtt{c})}\choose{\lfloor{s/2}\rfloor}}$ and $B^I_2\cup B^I_3\subseteq {{(\mathbb{N}_s,\mathtt{c})}\choose{\lceil{s/2}\rceil}}$, or $B^O_2\cup B^O_3\subseteq {{(\mathbb{N}_s,\mathtt{c})}\choose{\lceil{s/2}\rceil}}$ and $B^I_2\cup B^I_3\subseteq {{(\mathbb{N}_s,\mathtt{c})}\choose{\lfloor{s/2}\rfloor}}$.
\indent\par Consider the former case; a similar argument follows for the latter. Let $[i]\in A^O_3$, and $[j]\in A^I_3$. For any $[k]\in A_2\cup A^I_3$, $d_D((3,[i]),$ $(1,[1,k]))=3$ implies $X\neq I^\mathtt{c}((3,[i]))$ for all $X\in B^I_2\cup B^I_3$. This implies $|A_2|+|A^I_3|=|B^I_2\cup B^I_3|\le {{s}\choose{\lceil{s/2}\rceil}}-1$. Similarly, for any $[k]\in A_2\cup A^O_3$, $d_D((1,[1,k]),(3,[j]))=3$ implies $X\neq O^\mathtt{c}((3,[j]))$ for all $X\in B^O_2\cup B^O_3$. Thus, $|A_2|+|A^O_3|=|B^O_2\cup B^O_3|\le {{s}\choose{\lfloor{s/2}\rfloor}}-1$. Hence, $2|A_2|+|A_3|=(|A_2|+|A^O_3|)+(|A_2|+|A^I_3|)\le [{{s}\choose{\lfloor{s/2}\rfloor}}-1]+[{{s}\choose{\lceil{s/2}\rceil}}-1]=2{{s}\choose{\lceil{s/2}\rceil}}-2$.
\begin{rmk}
An approach similar to the one mentioned in Remark \ref{rmkC6.4.2} is employed here. We shall establish a series of claims on the structure of $D$, from which we will derive $2|A_2|+|A_3|\le 2{{s}\choose{\lceil{s/2}\rceil}}-2$ if any one fails to hold. In other words, $2|A_2|+|A_3|=2{{s}\choose{\lceil{s/2}\rceil}}-1$ is only possible in the last scenario where all these claims hold.
\end{rmk}
\noindent Subcase 2.3.6. $B^O_2\cup B^O_3\cup B^I_2\cup B^I_3\subseteq {{(\mathbb{N}_s,\mathtt{c})}\choose{\lceil{s/2}\rceil}}$.
\\
\\Subcase 2.3.6.1. $|O^\mathtt{c}((1,[i]))|=|O^\mathtt{c}((2,[i]))|=\lfloor\frac{s}{2}\rfloor$ for some $[i]\in A^{O(D)}_3$.
\indent\par For any $[j]\in A_2\cup A^I_3$ and $p=1,2$, $d_D((p,[i]),(1,[1,j]))=3$ implies $X\neq I^\mathtt{c}((p,[i]))$ for all $X\in B^I_2\cup B^I_3$. It follows that $|A_2|+|A^I_3|=|B^I_2\cup B^I_3|\le {{s}\choose{\lceil{s/2}\rceil}}-2$. Hence, $2|A_2|+|A_3|=(|A_2|+|A^O_3|)+(|A_2|+|A^I_3|)\le {{s}\choose{\lceil{s/2}\rceil}}+[{{s}\choose{\lceil{s/2}\rceil}}-2]=2{{s}\choose{\lceil{s/2}\rceil}}-2$.
\\
\\Subcase 2.3.6.2. $|O^\mathtt{c}((1,[i]))|=|O^\mathtt{c}((2,[i]))|=\lceil\frac{s}{2}\rceil$ for some $[i]\in A^{I(D)}_3$.
\indent\par This subcase follows from the Duality Lemma and Subcase 2.3.6.1.
\\
\\Subcase 2.3.6.3. $|O^\mathtt{c}((1,[i^*]))|=\lfloor\frac{s}{2}\rfloor$ and $|O^\mathtt{c}((2,[i^*]))|=\lceil\frac{s}{2}\rceil$ for some $[i^*]\in A^{O(D)}_3$.
\noindent\par $(\star\star)$ For any $[j]\in A_2\cup A^I_3$, $d_D((1,[i^*]),(1,[1,j]))=3$ implies $X\neq I^\mathtt{c}((1,[i^*]))$ for all $X\in B^I_2\cup B^I_3$. It follows that $|A_2|+|A^I_3|=|B^I_2\cup B^I_3|\le {{s}\choose{\lceil{s/2}\rceil}}-1$.  Therefore, $2|A_2|+|A_3|=(|B^O_2|+|B^O_3|)+(|B^I_2|+|B^I_3|)\le {{s}\choose{\lceil{s/2}\rceil}}+[{{s}\choose{\lceil{s/2}\rceil}}-1]=2{{s}\choose{\lceil{s/2}\rceil}}-1$. 

\indent\par Subcase 2.3.6.3 is done for $2|A_2|+|A_3|\le 2{{s}\choose{\lceil{s/2}\rceil}}-2$. Now, our aim is to prove $|A_2|\ge \lceil\frac{s}{2}\rceil \lfloor\frac{s}{2}\rfloor$ and $s\ge 5$ in the event that $2|A_2|+|A_3|= 2{{s}\choose{\lceil{s/2}\rceil}}-1$. The following claims will help us to achieve the said aim.
\noindent\par Suppose $2|A_2|+|A_3|= 2{{s}\choose{\lceil{s/2}\rceil}}-1$. Equivalently, 
\begin{align}
|B^O_2|+|B^O_3|= {{s}\choose{\lceil{s/2}\rceil}}\text{ and }|B^I_2|+|B^I_3|= {{s}\choose{\lceil{s/2}\rceil}}-1.\label{eqC6.4.11}
\end{align}
Then, Claims 2A-6A hold.
\\
\\Claim 2A: $|O^\mathtt{c}((1,[i]))|=\lfloor\frac{s}{2}\rfloor$ and $|O^\mathtt{c}((2,[i]))|=\lceil\frac{s}{2}\rceil$ for all $[i]\in A^{O(D)}_3$.
\indent\par Suppose there exists some $[i]\in A^{O(D)}_3-\{[i^*]\}$ such that $|O^\mathtt{c}((1,[i]))|=|O^\mathtt{c}((2,[i]))|=\lceil\frac{s}{2}\rceil$. Note by definition of $A^{O(D)}_3$ that $O^\mathtt{c}((3,[i]))$ is equal to at most one of $O^\mathtt{c}((1,[i]))$ and $O^\mathtt{c}((2,[i]))$, say $O^\mathtt{c}((3,[i]))\neq O^\mathtt{c}((1,[i]))$. Also, for any $[j]\in A_2\cup A^O_3-\{[i]\}$, $d_D((1,[1,j]),$ $(1,[i]))=3$ implies $X\neq O^\mathtt{c}((1,[i]))$  for all $X\in B^O_2\cup B^O_3-\{O^\mathtt{c}((3,[i]))\}$. It follows that $|B^O_2|+|B^O_3|\le {{s}\choose{\lceil{s/2}\rceil}}-1$, a contradiction to (\ref{eqC6.4.11}). In view of Subcase 2.3.6.1, the claim follows.
\\
\\Claim 3A: $O^\mathtt{c}((1,[i]))=O^\mathtt{c}((1,[i^*]))$ for all $[i]\in A^{O(D)}_3$.
\indent\par Suppose there exists some $[i]\in A^{O(D)}_3-\{[i^*]\}$ such that $O^\mathtt{c}((1,[i]))\neq O^\mathtt{c}((1,[i^*]))$. For any $[j]\in A_2\cup A^I_3$ and $x=i,i^*$, $d_D((1,[x]),(1,[1,j]))=3$ implies $X\neq I^\mathtt{c}((1,[x]))$ for all $X\in B^I_2\cup B^I_3$. It follows that $|B^I_2|+|B^I_3|\le {{s}\choose{\lceil{s/2}\rceil}}-2$, a contradiction to (\ref{eqC6.4.11}).
\\
\\Claim 4A: $O^\mathtt{c}((2,[i]))=O^\mathtt{c}((3,[i]))$ for all $[i]\in A^{O(D)}_3$.
\indent\par Suppose there exists some $[i]\in A^{O(D)}_3$ such that $O^\mathtt{c}((2,[i]))\neq O^\mathtt{c}((3,[i]))$. Then, for any $[j]\in A_2\cup A^O_3-\{[i]\}$, $d_D((1,[1,j]),(2,[i]))=3$ implies $X\neq O^\mathtt{c}((2,[i]))$ for all $X\in B^O_2\cup B^O_3-\{O^\mathtt{c}((3,[i]))\}$. It follows that $|B^O_2|+|B^O_3|\le {{s}\choose{\lceil{s/2}\rceil}}-1$, a contradiction to (\ref{eqC6.4.11}).
\\
\\Claim 5A: $|O^\mathtt{c}((1,[i]))|=|O^\mathtt{c}((2,[i]))|=\lfloor\frac{s}{2}\rfloor$ for all $[i]\in A^{I(D)}_3$.
\indent\par Suppose there exists some $[i]\in A^{I(D)}_3$ such that $|O^\mathtt{c}((1,[i]))|=\lfloor\frac{s}{2}\rfloor$ and $|O^\mathtt{c}((2,[i]))|=\lceil\frac{s}{2}\rceil$. For any $[j]\in A_2\cup A^O_3$, $d_D((1,[1,j]),(2,[i]))=3$ implies $X\neq O^\mathtt{c}((2,[i]))$ for all $X\in B^O_2\cup B^O_3$. So, $|B^O_2|+|B^O_3|\le {{s}\choose{\lceil{s/2}\rceil}}-1$, a contradiction to (\ref{eqC6.4.11}). In view of Subcase 2.3.6.2, the claim follows.
\\
\\Claim 6A: $O^\mathtt{c}((1,[i]))= O^\mathtt{c}((1,[i^*]))$ and $O^\mathtt{c}((2,[i]))=O^\mathtt{c}((3,[i]))$ for all $[i]\in A^{I(D)}_3$. 
\noindent\par Suppose there exist some $[i]\in A^{I(D)}_3$ and some $p=1,2$, such that $O^\mathtt{c}((p,[i]))\neq O^\mathtt{c}((1,[i^*]))$ and $O^\mathtt{c}((p,[i]))\neq O^\mathtt{c}((3,[i]))$. Also, for any $[j]\in A_2\cup A^I_3-\{[i]\}$, $d_D((p,[i]),$ $(1,[1,j]))=3$ implies $X\cap O^\mathtt{c}((p,[i]))\neq \emptyset$ for all $X\in B^I_2\cup B^I_3-\{I^\mathtt{c}((3,[i]))\}$. Therefore, for all $X\in B^I_2\cup B^I_3$, $X\neq I^\mathtt{c}((p,[i]))$ and recall from $(\star\star)$ that $X\neq I^\mathtt{c}((1,[i^*]))$. It follows that $|B^I_2|+|B^I_3|\le {{s}\choose{\lceil{s/2}\rceil}}-2$, a contradiction to (\ref{eqC6.4.11}).
\indent\par Therefore, it must follow that for each $[i]\in A^{I(D)}_3$ and each $p=1,2$, either $O^\mathtt{c}((p,[i]))=O^\mathtt{c}((3,[i]))$ or $O^\mathtt{c}((p,[i]))= O^\mathtt{c}((1,[i^*]))$. By the definition of $A^{I(D)}_3$, we have $O^\mathtt{c}((1,[i]))\neq O^\mathtt{c}((2,[i]))$ for each $[i]\in A^{I(D)}_3$. Hence, we may assume without loss of generality that $O^\mathtt{c}((1,[i]))= O^\mathtt{c}((1,[i^*]))$ and $O^\mathtt{c}((2,[i]))=O^\mathtt{c}((3,[i]))$ for all $[i]\in A^{I(D)}_3$.

\indent\par Now, we are ready to prove $|A_2|\ge \lceil\frac{s}{2}\rceil \lfloor\frac{s}{2}\rfloor$ and $s\ge 5$ in the event that $2|A_2|+|A_3|= 2{{s}\choose{\lceil{s/2}\rceil}}-1$.
\noindent\par Suppose $s=3$. By Claims 5A-6A, $|O^\mathtt{c}((1,[j]))\cup O^\mathtt{c}((2,[j]))|\le\lceil\frac{s}{2}\rceil$ for any $[j]\in A^{I(D)}_3$. Since $B^O_2\cup B^O_3= {{(\mathbb{N}_s,\mathtt{c})}\choose{\lceil{s/2}\rceil}}$, we have $O^\mathtt{c}((1,[j]))\cup O^\mathtt{c}((2,[j]))\subseteq X$ for some $X\in B^O_2\cup B^O_3$, a contradiction to Lemma \ref{lemC6.2.17}(b). Hence, $s\ge 5$.
\indent\par Let $[j]\in A^{I(D)}_3$, $[i]\in A_2\cup A^O_3$ and recall $\delta$ from (\ref{eqC6.4.10}). Then, Claim 6A and $d_D((1,[1,j]),(\delta,[i]))=3$ imply $O^\mathtt{c}((1,[i^*]))\cup O^\mathtt{c}((2,[j]))=O^\mathtt{c}((1,[j]))\cup O^\mathtt{c}((2,[j]))\neq O^\mathtt{c}((\delta,[i]))$, i.e., $O^\mathtt{c}((2,[j]))$ cannot be a $\lfloor\frac{s}{2}\rfloor$-set whose union with $O^\mathtt{c}((1,[i^*]))$ forms a $\lceil\frac{s}{2}\rceil$-set. So, $\{O^\mathtt{c}((3,[j]))\mid [j]\in A^{I(D)}_3\}\cap R=\{O^\mathtt{c}((2,[j]))\mid [j]\in A^{I(D)}_3\}\cap R=\emptyset$, where $R=\{X\in {{(\mathbb{N}_s,\mathtt{c})}\choose{\lfloor{s/2}\rfloor}}\mid |X\cup O^\mathtt{c}((1,[i^*]))|=\lceil\frac{s}{2}\rceil\}$. Note that there are $\lceil\frac{s}{2}\rceil$ number of $\lceil\frac{s}{2}\rceil$-supersets of $O^\mathtt{c}((1,[i^*]))$ for there are $\lceil\frac{s}{2}\rceil$ elements not in $O^\mathtt{c}((1,[i^*]))$. For each of these supersets, there are exactly $\lfloor\frac{s}{2}\rfloor$ number of $\lfloor\frac{s}{2}\rfloor$-sets whose union with $O^\mathtt{c}((1,[i^*]))$ forms the superset, as each $\lfloor\frac{s}{2}\rfloor$-set differs from $O^\mathtt{c}((1,[i^*]))$ by exactly one element and there are $\lfloor\frac{s}{2}\rfloor$ choices for such an element in $O^\mathtt{c}((1,[i^*]))$. So, $|R|=\lceil\frac{s}{2}\rceil \lfloor\frac{s}{2}\rfloor$. 
\indent\par Furthermore, $d_D((1,[i^*]),(1,[1,j]))=3$ and the definition of $R$ imply $O^\mathtt{c}((1,[i^*]))\not\in \{O^\mathtt{c}((3,[j]))\mid [j]\in A^{I(D)}_3\}\cup R$. It follows that $\{O^\mathtt{c}((3,[j]))\mid [j]\in A^{I(D)}_3\}\cup R\subseteq {{(\mathbb{N}_s,\mathtt{c})}\choose{\lfloor{s/2}\rfloor}}-\{O^\mathtt{c}((1,[i^*]))\}$. So, $|A^{I(D)}_3|\le {{s}\choose{\lceil{s/2}\rceil}}-1-|R|={{s}\choose{\lceil{s/2}\rceil}}-1-\lceil{\frac{s}{2}}\rceil \lfloor{\frac{s}{2}}\rfloor$.
Since $\{O^\mathtt{c}((1,[i]))\mid [i]\in A^{I(S)}_3\}\cup B^O_2\cup B^O_3$ is an antichain by Lemma \ref{lemC6.2.15}, $|A^{I(S)}_3|+|B^O_2|+|B^O_3|\le {{s}\choose{\lceil{s/2}\rceil}}$. So, $|A_2|+|A_3|=|A^{I(D)}_3|+(|A^{I(S)}_3|+|B^O_3|+|B^O_2|)\le 2{{s}\choose{\lceil{s/2}\rceil}}-1-\lceil\frac{s}{2}\rceil \lfloor\frac{s}{2}\rfloor$. Using $2|A_2|+|A_3|=2{{s}\choose{\lfloor{s/2}\rfloor}}-1$, we derive $|A_2|\ge \lceil\frac{s}{2}\rceil\lfloor\frac{s}{2}\rfloor$.
\begin{rmk}
We are done for Subcase 2.3.6.3 but we want to obtain more information to construct the optimal orientations for the case $2|A_2|+|A_3|=2{{s}\choose{\lfloor{s/2}\rfloor}}-1$ when proving the converse later. (It is also neater to present the argument here.) In fact, such an optimal orientation has to satisfy the above set of Claims 2A-6A or the analogous set of Claims 2B-6B in Subcase 2.3.6.4. (see Corollary \ref{corC6.4.10}).
\end{rmk}
\noindent Claim 7A: $A^O_3=A^{O(D)}_3$ and $A^I_3=A^{I(D)}_3$. 
\indent\par Note that $|B^O_2|+|B^O_3|={{s}\choose{\lceil{s/2}\rceil}}$, and $|A^{I(S)}_3|+|B^O_2|+|B^O_3|\le {{s}\choose{\lceil{s/2}\rceil}}$ imply $|A^{I(S)}_3|=0$ and $A^I_3=A^{I(D)}_3$. Also, $\{I^\mathtt{c}((1,[i]))\mid [i]\in A^{O(S)}_3\}\cup B^I_2\cup B^I_3$ is an antichain by Lemma \ref{lemC6.2.16}. For all $[j]\in A^{O(S)}_3\cup A_2\cup A^I_3$, $d_D((1,[i^*]),(1,[1,j]))=3$ implies $O^\mathtt{c}((1,[i^*]))\cap X\neq\emptyset$ for all $X\in\{I^\mathtt{c}((1,[i]))\mid [i]\in A^{O(S)}_3\}\cup B^I_2\cup B^I_3$. By Lih's theorem, $|A^{O(S)}_3|+|B^I_2|+|B^I_3|\le {{s}\choose{\lceil{s/2}\rceil}}-{{s-|O^\mathtt{c}((1,[i^*]))|}\choose{\lceil{s/2}\rceil}}={{s}\choose{\lceil{s/2}\rceil}}-1$. Since $|B^I_2|+|B^I_3|={{s}\choose{\lceil{s/2}\rceil}}-1$, it follows that $|A^{O(S)}_3|=0$ and $A^O_3=A^{O(D)}_3$.
\\
\\Subcase 2.3.6.4. $|O^\mathtt{c}((1,[i^*]))|=\lceil\frac{s}{2}\rceil$ and $|O^\mathtt{c}((2,[i^*]))|=\lfloor\frac{s}{2}\rfloor$ for some $[i^*]\in A^{I(D)}_3$.
\noindent\par This follows from Subcase 3.3.6.3 by the Duality Lemma. For clarity, we state the analogous versions of Claims 2A-6A below that hold if $2|A_2|+|A_3|=2{{s}\choose{\lceil{s/2}\rceil}}-1$.
\\
\\Claim 2B: $|O^\mathtt{c}((1,[i]))|=\lceil\frac{s}{2}\rceil$ and $|O^\mathtt{c}((2,[i]))|=\lfloor\frac{s}{2}\rfloor$ for all $[i]\in A^{I(D)}_3$.
\\Claim 3B: $O^\mathtt{c}((1,[i]))=O^\mathtt{c}((1,[i^*]))$ for all $[i]\in A^{I(D)}_3$.
\\Claim 4B: $O^\mathtt{c}((2,[i]))=O^\mathtt{c}((3,[i]))$ for all $[i]\in A^{I(D)}_3$.
\\Claim 5B: $|O^\mathtt{c}((1,[j]))|=|O^\mathtt{c}((2,[j]))|=\lceil\frac{s}{2}\rceil$ for all $[j]\in A^{O(D)}_3$.
\\Claim 6B: $O^\mathtt{c}((1,[j]))= O^\mathtt{c}((1,[i^*]))$ and $O^\mathtt{c}((2,[j]))=O^\mathtt{c}((3,[j]))$ for all $[j]\in A^{O(D)}_3$.
\\Claim 7B: $A^O_3=A^{O(D)}_3$ and $A^I_3=A^{I(D)}_3$.
\\
\\Subcase 2.3.6.5. $|O^\mathtt{c}((1,[i]))|=|O^\mathtt{c}((2,[i]))|=\lceil\frac{s}{2}\rceil$ and $|O^\mathtt{c}((1,[j]))|=|O^\mathtt{c}((2,[j]))|=\lfloor\frac{s}{2}\rfloor$ for some $[i]\in A^{O(D)}_3$ and $[j]\in A^{I(D)}_3$.
\indent\par By definition of $A^{O(D)}_3$, $O^\mathtt{c}((3,[i]))$ is equal to at most one of $O^\mathtt{c}((1,[i]))$ and $O^\mathtt{c}((2,[i]))$, say $O^\mathtt{c}((3,[i]))\neq O^\mathtt{c}((1,[i]))$. Also, for any $[k]\in A_2\cup A^O_3-\{[i]\}$, $d_D((1,[1,k]),$ $(1,[i]))=3$ implies $X\neq O^\mathtt{c}((1,[i]))$ for all $X\in B^O_2\cup B^O_3-\{O^\mathtt{c}((3,[i]))\}$. Thus, $|B^O_2|+|B^O_3|\le {{s}\choose{\lceil{s/2}\rceil}}-1$. 
\indent\par Similarly, by definition of $A^{I(D)}_3$, $I^\mathtt{c}((3,[j]))$ is equal to at most one of $I^\mathtt{c}((1,[j]))$ and $I^\mathtt{c}((2,[j]))$, say $I^\mathtt{c}((3,[j]))\neq I^\mathtt{c}((1,[j]))$. Also, for any $[k]\in A_2\cup A^I_3-\{[j]\}$, $d_D((1,[j]),$ $(1,[1,k]))=3$ implies $X\neq I^\mathtt{c}((1,[j]))$ for all $X\in B^I_2\cup B^I_3-\{I^\mathtt{c}((3,[i]))\}$. Thus, $|B^I_2|+|B^I_3 |\le {{s}\choose{\lceil{s/2}\rceil}}-1$. Hence, $2|A_2|+|A_3|=(|B^O_2|+|B^O_3|)+(|B^I_2|+|B^I_3|)\le 2[{{s}\choose{\lceil{s/2}\rceil}}-1]=2{{s}\choose{\lceil{s/2}\rceil}}-2$.
\\
\noindent\\$(\Leftarrow)$ By Corollary \ref{corC6.3.8}, $\mathcal{T}\in \mathscr{C}_0$ if $|A_2|\le{{s}\choose{\lceil{s/2}\rceil}}-2$ and $|A_3|=1$, or $|A_2|+|A_3|\le {{s}\choose{\lceil{s/2}\rceil}}-1$ and $|A_3|\ge 2$. Hence, it suffices to consider the following three cases. Let $\mathcal{H}=T(t_1,t_2,\ldots, t_n)$ be the subgraph of $\mathcal{T}$, where $t_\mathtt{c}=s$, $t_{[i]}=3$ for all $[i]\in \mathcal{T}(A_3)$ and $t_v=2$ otherwise. We will use $A_j$ for $\mathcal{H}(A_j)$ for the remainder of this proof.
\\
\\Case 1. $|A_2|={{s}\choose{\lceil{s/2}\rceil}}-1$ and $|A_3|=1$.
\indent\par Assume without loss of generality that $A_3=\{[1]\}$ and $A_2=\{[i]\mid i\in\mathbb{N}_{{{s}\choose{\lceil{s/2}\rceil}}}\}-A_3$. Define an orientation $D_1$ of $\mathcal{H}$ as follows.
\begin{align}
& (2,[i])\rightarrow \{(1,[\alpha,i]),(2,[\alpha,i])\}\rightarrow (1,[i]),\text{ and}\label{eqC6.4.12}\\
& \bar{\lambda}_{i}\rightarrow (1,[i])\rightarrow \lambda_{i}\rightarrow (2,[i]) \rightarrow \bar{\lambda}_{i}\label{eqC6.4.13}
\end{align}
for all $[i]\in A_2$ and all $1\le \alpha\le \deg_T([i])-1$, i.e., excluding $\lambda_1$, the $\lceil \frac{s}{2}\rceil$-sets $\lambda_i$'s are used as `in-sets' (`out-sets' resp.) to construct $B^I_2$ ($B^O_2$ resp.).
\begin{align}
& \{(1,[1]),(2,[1])\}\rightarrow (1,[\beta,1])\rightarrow (3,[1]),\label{eqC6.4.14}\\
& \{(1,[1]),(3,[1])\}\rightarrow (2,[\beta,1])\rightarrow (2,[1]),\text{ and}\label{eqC6.4.15}\\
& \bar{\lambda}_1\rightarrow \{(2,[1]),(3,[1])\} \rightarrow\lambda_1\rightarrow (1,[1])\rightarrow \bar{\lambda}_1\label{eqC6.4.16}
\end{align}
for all $1\le \beta\le \deg_T([1])-1$. Furthermore,
\begin{align}
\lambda_1\rightarrow \{(1,[j]),(2,[j])\}\rightarrow \bar{\lambda}_1\label{eqC6.4.17}
\end{align}
for all $[j]\in E$. (See Figure \ref{figC6.4.13} for $D_1$ when $s=3$.)
\\
\\Claim: $d_{D_1}(v,w)\le 4$ for all $v,w\in V(D_1)$.
\\
\\Case 1.1.1. $v,w \in \{(1,[\alpha,i]),(2,[\alpha,i]),(1,[i]),(2,[i])\}$ for each $[i]\in A_2$ and $1\le\alpha\le \deg_T([i])-1$.
\indent\par By (\ref{eqC6.4.12})-(\ref{eqC6.4.13}), $(2,[i])\rightarrow \{(1,[\alpha,i]),(2,[\alpha,i])\}\rightarrow (1,[i])\rightarrow \lambda_{i}\rightarrow (2,[i])$ guarantees a directed $C_4$.
\\
\\Case 1.1.2. $v,w \in \{(1,[\alpha,1]),(2,[\alpha,1]),(1,[1]),(2,[1]),(3,[1])\}$ for each $1\le\alpha\le \deg_T([1])-1$.
\indent\par By (\ref{eqC6.4.14})-(\ref{eqC6.4.16}), $(1,[1])\rightarrow \{(1,[\alpha,1]),(2,[\alpha,1])\}$, $(1,[\alpha,1])\rightarrow (3,[1])\rightarrow \lambda_1\rightarrow (1,[1])$ and $(2,[\alpha,1])\rightarrow (2,[1])\rightarrow \lambda_1\rightarrow (1,[1])$ guarantee the directed $C_4$'s required.
\\
\\Case 1.2. For each $[i],[j]\in A_2$, $i\neq j$, each $1\le \alpha\le \deg_T([i])-1$, and each $1\le \beta\le \deg_T([j])-1$,
\\(i) $v=(p,[\alpha,i]), w=(q,[j])$ for any $p,q=1,2$.
\indent\par By (\ref{eqC6.4.12})-(\ref{eqC6.4.13}), $\{(1,[\alpha,i]),(2,[\alpha,i])\}\rightarrow (1,[i])\rightarrow \lambda_{i}$, and $\lambda_{j}\rightarrow (2,[j])$. Since $|\lambda_{i}|=|\lambda_{j}|=\lceil\frac{s}{2}\rceil$, there exists some $x_{ij}\in\lambda_{i}\cap\lambda_{j} $ such that $(1,[i])\rightarrow x_{ij}\rightarrow (2,[j])$. Also, $\lambda_{i}\cap\bar{\lambda}_{j}\neq\emptyset$ for $i\neq j$. There exists some $y_{ij}\in\lambda_{i}\cap\bar{\lambda}_{j} $ such that $(1,[i])\rightarrow y_{ij}\rightarrow (1,[j])$.
\\
\\(ii) $v=(q,[j]), w=(p,[\alpha,i])$ for any $p,q=1,2$.
\indent\par By (\ref{eqC6.4.12})-(\ref{eqC6.4.13}), $(2,[i])\rightarrow \{(1,[\alpha,i]), (2,[\alpha,i])\}$. Since $|\lambda_j|=|\lambda_{i}|=\lceil\frac{s}{2}\rceil$, there exists some $x_{ij}\in\lambda_j\cap\lambda_{i} $ such that $(1,[j])\rightarrow x_{ij}\rightarrow (2,[i])$. Also, since $i\neq j$, there exists some $y_{ij}\in\bar{\lambda}_j\cap\lambda_{i} $ such that $(2,[j])\rightarrow y_{ij}\rightarrow (2,[i])$.
\\
\\(iii) $v=(p,[\alpha,i]), w=(q,[\beta,j])$ for each $p,q=1,2$.
\indent\par From (i), $d_{D_1}((p,[\alpha,i]),(2, [j]))=3$. Since $(2,[j])\rightarrow \{(1,[\beta,j]), (2,[\beta,j])\}$ by (\ref{eqC6.4.12}), this subcase follows.
\\
\\Case 1.3. For each $[i]\in A_2$, each $1\le\alpha\le \deg_T([i])-1$, and each $1\le\beta\le \deg_T([1])-1$,
\\(i) $v=(p,[\alpha,i]), w=(q,[1])$ for each $p=1,2$ and $q=1,2,3$.
\indent\par By (\ref{eqC6.4.12}), $\{(1,[\alpha,i]),(2,[\alpha,i])\}\rightarrow (1,[i])$. Since $|\lambda_i|=|\lambda_1|=\lceil\frac{s}{2}\rceil$, there exists some $x_i\in \lambda_1\cap\lambda_{i}$ such that $(1,[i])\rightarrow x_i\rightarrow (1,[1])$ by (\ref{eqC6.4.13}) and (\ref{eqC6.4.16}). Also, since $i\neq 1$, there exists some $y_i\in \bar{\lambda}_1\cap\lambda_{i}$ such that $(1,[i])\rightarrow y_i\rightarrow \{(2,[1]), (3,[1])\}$.
\\
\\(ii) $v=(q,[1]), w=(p,[\alpha,i])$ for each $p=1,2$ and $q=1,2,3$.
\indent\par By (\ref{eqC6.4.12}), $(2,[i])\rightarrow \{(1,[\alpha,i]),(2,[\alpha,i])\}$. Since $i\neq 1$, there exists some $x_i\in \bar{\lambda}_1\cap\lambda_{i}$ such that $(1,[1])\rightarrow x_i\rightarrow (2,[i])$ by (\ref{eqC6.4.13}) and (\ref{eqC6.4.16}). Also, since $|\lambda_1|=|\lambda_{i}|=\lceil\frac{s}{2}\rceil$, there exists some $y_i\in \lambda_1\cap\lambda_{i}$ such that $\{(2,[1]),(3,[1])\}\rightarrow y_i\rightarrow (2,[i])$.
\\
\\(iii) $v=(p,[i]), w=(q,[\beta,1])$ for each $p,q=1,2$.
\indent\par From (i), $d_{D_1}((1,[i]), (1,[1]))=2$ and by (\ref{eqC6.4.14})-(\ref{eqC6.4.15}), $(1,[1])\rightarrow\{(1,[\beta,1]),$ $(2,[\beta,1])\}$. And, since $i\neq 1$, there exists some $x_i\in \bar{\lambda}_{i}\cap \lambda_1$ such that $(2,[i])\rightarrow x_i\rightarrow (1,[1])$ by (\ref{eqC6.4.13}) and (\ref{eqC6.4.16}).
\\
\\(iv) $v=(q,[\beta,1]), w=(p,[i])$ for each $p,q=1,2$.
\indent\par By (\ref{eqC6.4.14})-(\ref{eqC6.4.16}), $(1,[\beta,1])\rightarrow (3,[1])\rightarrow \lambda_1$ and $(2,[\beta,1])\rightarrow (2,[1])\rightarrow \lambda_1$. Since $|\lambda_1|=|\lambda_{i}|=\lceil\frac{s}{2}\rceil$, there exists some $x_i\in \lambda_1\cap\lambda_{i}$ such that $x_i \rightarrow (2,[i])$ by (\ref{eqC6.4.13}). And, since $i\neq 1$, there exists some $y_i\in \lambda_1\cap\bar{\lambda}_{i}$ such that $y_i \rightarrow (1,[i])$.
\\
\\(v) $v=(p,[\alpha,i]), w=(q,[\beta,1])$ for each $p,q=1,2$.
\indent\par From (i), $d_{D_1}((p,[\alpha,i]), (1,[1]))=3$. Since $(1,[1])\rightarrow\{(1,[\beta,1]),(2,[\beta,1])\}$ by (\ref{eqC6.4.14})-(\ref{eqC6.4.15}), this subcase follows.
\\
\\(vi) $v=(q,[\beta,1]), w=(p,[\alpha,i])$ for each $p,q=1,2$.
\indent\par From (iv), $d_{D_1}((q,[\beta,1]),(2,[i]))=3$. Since $(2,[i])\rightarrow \{(1,[\alpha,i]), (2,[\alpha,i])\}$ by (\ref{eqC6.4.12}), this subcase follows.
\\
\\Case 1.4. For each $[i]\in A_2$, each $1\le\alpha\le \deg_T([i])-1$, and each $[j]\in E$,
\\(i) $v=(p,[\alpha,i]), w=(q,[j])$ for each $p,q=1,2$.
\\(ii) $v=(q,[j]), w=(p,[\alpha,i])$ for each $p,q=1,2$.
\indent\par Since $\lambda_1\rightarrow \{(q,[j]),(1,[1])\}\rightarrow \bar{\lambda}_1$ by (\ref{eqC6.4.16})-(\ref{eqC6.4.17}), this case follows from Cases 1.3(i)-(ii).
\\
\\Case 1.5. For each $1\le\alpha\le deg_T([1])-1$ and each $[j]\in E$,
\\(i) $v=(p,[\alpha,1]), w=(q,[j])$ for each $p,q=1,2$.
\noindent\par By (\ref{eqC6.4.14})-(\ref{eqC6.4.17}), $(1,[\alpha,1])\rightarrow (3,[1])$, $(2,[\alpha,1])\rightarrow (2,[1])$, and $\{(2,[1]), (3,[1])\}\rightarrow \lambda_1\rightarrow \{(1,[j]), (2,[j])\}$.
\\
\\(ii) $v=(q,[j]), w=(p,[\alpha,1])$ for each $p,q=1,2$.
\noindent\par By (\ref{eqC6.4.14})-(\ref{eqC6.4.17}), $\{(1,[j]), (2,[j])\}\rightarrow \bar{\lambda}_1\rightarrow \{(2,[1]), (3,[1])\}$, $(2,[1])\rightarrow (1,[\alpha,1])$ and $(3,[1])\rightarrow (2,[\alpha,1])$.
\\
\\Case 1.6. $v=(r_1,\mathtt{c})$ and $w=(r_2,\mathtt{c})$ for $r_1\neq r_2$ and $1\le r_1, r_2\le s$.
\indent\par Here, we want to prove a stronger claim, $d_{D_1}((r_1,\mathtt{c}), (r_2,\mathtt{c}))=2$. Let $x_1=(1,[1])$, $x_k=(2,[k])$ for $2\le k\le {{s}\choose{\lceil{s/2}\rceil}}$. Observe that $\lambda_k\rightarrow x_k\rightarrow \bar{\lambda}_k$ for all $1\le k\le s$ and the subgraph induced by $V_1=(\mathbb{N}_s,\mathtt{c})$ and $V_2=\{x_i\mid 1\le i \le s\}$ is a complete bipartite graph $K(V_1,V_2)$. By Lemma \ref{lemC6.2.19}, $d_{D_1}((r_1,\mathtt{c}), (r_2,\mathtt{c}))=2$.
\\
\\Case 1.7. $v\in \{(1,[i]), (2,[i]), (3,[i]), (1,[\alpha,i]), (2,[\alpha,i])\}$ for each $1\le i\le deg_T(\mathtt{c})$ and $1\le\alpha\le deg_T([i])-1$, and $w=(r,\mathtt{c})$ for $1\le r\le s$.
\indent\par Note that there exists some $1\le k\le s$ such that $d_{D_1}(v,(k,\mathtt{c}))\le 2$, and $d_{D_1}((k,\mathtt{c}),w)\le 2$ by Case 1.6. Hence, it follows that $d_{D_1}(v,w)\le d_{D_1}(v,(k,\mathtt{c}))+d_{D_1}((k,\mathtt{c}),w)\le 4$.
\\
\\Case 1.8. $v=(r,\mathtt{c})$ for $1\le r\le s$ and $w\in \{(1,[i]), (2,[i]), (3,[i]), (1,[\alpha,i]), (2,[\alpha,i])\}$ for each $1\le i\le deg_T(\mathtt{c})$ and $1\le\alpha\le deg_T([i])-1$.
\indent\par Note that there exists some $1\le k\le s$ such that $d_{D_1}((k,\mathtt{c}), w)\le 2$, and $d_{D_1}(v,(k,\mathtt{c}))\le 2$ by Case 1.6. Hence, it follows that $d_{D_1}(v,w)\le d_{D_1}(v,(k,\mathtt{c}))+d_{D_1}((k,\mathtt{c}),w)\le 4$.
\\
\\Case 1.9. $v=(p,[i])$ and $w=(q, [j])$, where $1\le p,q\le 3$ and $1\le i,j\le deg_T(\mathtt{c})$.
\noindent\par This follows from the fact that $|O^\mathtt{c}((p,[i]))|>0$, $|I^\mathtt{c}((q,[j]))|>0$, and $d_{D_1}((r_1,\mathtt{c}), (r_2,\mathtt{c}))$ $=2$ for any $r_1\neq r_2$ and $1\le r_1, r_2\le s$.
\\
\\Case 2. $|A_3|\ge 2$ and $2|A_2|+|A_3|\le 2{{s}\choose{\lceil{s/2}\rceil}}-2$.
\indent\par By  Corollary \ref{corC6.3.8}, we may assume $|A_2|+|A_3|\ge {{s}\choose{\lceil{s/2}\rceil}}$. Furthermore, assume without loss of generality that $A_2=\{[i]\mid i\in\mathbb{N}_{|A_2|}\}$, and $A_3=\{[i]\mid i\in\mathbb{N}_{|A_2|+|A_3|}-\mathbb{N}_{|A_2|}\}$. Define an orientation $D_2$ of $\mathcal{H}$ as follows.
\begin{align}
& (2,[i])\rightarrow \{(1,[\alpha,i]),(2,[\alpha,i])\}\rightarrow (1,[i]),\text{ and}\label{eqC6.4.18}\\
& \bar{\lambda}_{i+1}\rightarrow (1,[i])\rightarrow \lambda_{i+1}\rightarrow (2,[i]) \rightarrow \bar{\lambda}_{i+1},\label{eqC6.4.19}
\end{align}
for all $1\le i\le |A_2|$ and all $1\le \alpha\le deg_T([i])-1$, i.e., the $\lceil\frac{s}{2}\rceil$-sets $\lambda_2,\lambda_3,\ldots,\lambda_{|A_2|+1}$ are used as `in-sets' (`out-sets' resp.) to construct $B^I_2$ ($B^O_2$ resp.).
\begin{align}
& (3,[j])\rightarrow \{(1,[\beta,j]),(2,[\beta,j])\}\rightarrow \{(1,[j]),(2,[j])\},\label{eqC6.4.20}\\
& \bar{\lambda}_1 \rightarrow (1,[j])\rightarrow \lambda_1\rightarrow (2,[j])\rightarrow \bar{\lambda}_1,\text{ and}\label{eqC6.4.21}\\
&\lambda_{j+1}\rightarrow (3,[j]) \rightarrow \bar{\lambda}_{j+1}\label{eqC6.4.22}
\end{align}
for all $|A_2|+1\le j\le {{s}\choose{\lceil{s/2}\rceil}}-1$ and all $1\le \beta \le deg_T([j])-1$, i.e., the $\lceil\frac{s}{2}\rceil$-sets $\lambda_{|A_2|+2},\lambda_{|A_2|+3},\ldots,\lambda_{{{s}\choose{\lceil{s/2}\rceil}}}$ are used as `in-sets' to construct $B^I_3$.
\begin{align}
& \{(1,[k]),(2,[k])\}\rightarrow \{(1,[\gamma,k]),(2,[\gamma,k])\}\rightarrow (3,[k]),\label{eqC6.4.23}\\
& \lambda_1 \rightarrow (1,[k])\rightarrow \bar{\lambda}_1\rightarrow (2,[k])\rightarrow \lambda_1,\text{ and}\label{eqC6.4.24}\\
&\bar{\lambda}_{k-{{s}\choose{\lceil{s/2}\rceil}}+|A_2|+2}\rightarrow (3,[k]) \rightarrow \lambda_{{k-{{s}\choose{\lceil{s/2}\rceil}}+|A_2|+2}}\label{eqC6.4.25}
\end{align}
for all ${{s}\choose{\lceil{s/2}\rceil}}\le k\le|A_2|+|A_3|$ and all $1\le \gamma\le deg_T([k])-1$, i.e., the $\lceil\frac{s}{2}\rceil$-sets $\lambda_{|A_2|+2},\lambda_{|A_2|+3},\ldots,\lambda_{2|A_2|+|A_3|+2-{{s}\choose{\lceil{s/2}\rceil}}}$ are used as `out-sets' to construct $B^O_3$.
\begin{align}
\lambda_1 \rightarrow \{(1,[l]), (2,[l])\} \rightarrow \bar{\lambda}_1\label{eqC6.4.26}
\end{align}
for all $[l]\in E$. (See Figure \ref{figC6.4.14} for $D_2$ when $s=3$.)
\\
\\Claim: $d_{D_2}(v,w)\le 4$ for all $v,w\in V(D_2)$.
\\
\\Case 2.1.1. $v,w \in \{(1,[\alpha,i]),(2,[\alpha,i]),(1,[i]),(2,[i])\}$ for each $1\le i\le |A_2|$ and $1\le\alpha\le deg_T([i])-1$.
\indent\par In view of the similarity in the definitions of $D_2$ (see (\ref{eqC6.4.18})-(\ref{eqC6.4.19})) and $D_1$ (see (\ref{eqC6.4.12})-(\ref{eqC6.4.13})), this case is similar to Case 1.1.1.
\\
\\Case 2.1.2. $v,w \in \{(1,[\alpha,i]),(2,[\alpha,i]),(1,[i]),(2,[i]), (3,[i]), (4,[i])\}$ for each $[i]\in A_3$ and $1\le\alpha\le deg_T([i])-1$.
\indent\par Since the orientation defined for $A_3$ (see (\ref{eqC6.4.20})-(\ref{eqC6.4.25})) is similar to that in Proposition \ref{ppnC6.4.1} (see (\ref{eqC6.4.1})-(\ref{eqC6.4.6})), this case follows from Cases 1.1-1.2 of Proposition \ref{ppnC6.4.1}.
\\
\\Case 2.2.  For each $[i],[j]\in A_3$, $i \neq j$, each $1\le\alpha\le deg_T([i])-1$, and each $1\le\beta\le deg_T([j])-1$,
\\(i) $v=(p,[\alpha,i]), w=(q,[j])$ for $p=1,2$, and $q=1,2,3$.
\\(ii) $v=(p,[i]), w=(q,[\beta,j])$ for for $p=1,2,3$, and $q=1,2$.
\\(iii) $v=(p,[\alpha,i]), w=(q,[\beta,j])$ for $p,q=1,2$.
\indent\par Since the orientation defined for $A_3$ (see (\ref{eqC6.4.20})-(\ref{eqC6.4.25})) is similar to that in Proposition \ref{ppnC6.4.1} (see (\ref{eqC6.4.1})-(\ref{eqC6.4.6})), this case follows from Cases 2-3 and 5 of Proposition \ref{ppnC6.4.1}.
\\
\\Case 2.3. For each $1\le i,j \le  |A_2|$, $i\neq j$, each $1\le \alpha\le deg_T([i])-1$, and each $1\le \beta\le deg_T([j])-1$,
\\(i) $v=(p,[\alpha,i]), w=(q,[j])$ for any $p,q=1,2$.
\\(ii) $v=(q,[j]), w=(p,[\alpha,i])$ for any $p,q=1,2$.
\\(iii) $v=(p,[\alpha,i]), w=(q,[\beta,j])$ for each $p,q=1,2$.
\indent\par In view of the similarity in the definitions of $D_2$ (see (\ref{eqC6.4.18})-(\ref{eqC6.4.19})) and $D_1$ (see (\ref{eqC6.4.12})-(\ref{eqC6.4.13})), this case is similar to Case 1.2.
\\
\\Case 2.4. For each $1\le i\le |A_2|$, each $|A_2|+1\le j\le {{s}\choose{\lceil{s/2}\rceil}}-1$, each $1\le\alpha\le deg_T([i])-1$, and each $1\le\beta\le deg_T([j])-1$,
\\(i) $v=(p,[\alpha,i]), w=(q,[j])$ for each $p=1,2$ and $q=1,2,3$.
\indent\par Since $|\lambda_1|=|\lambda_{i+1}|=|\lambda_{j+1}|=\lceil\frac{s}{2}\rceil$, there exist some $x_i\in \lambda_{i+1}\cap\lambda_1$ and $y_{ij}\in \lambda_{i+1}\cap\lambda_{j+1}$. So, $(1,[i])\rightarrow x_i\rightarrow (2,[j])$ and $(1,[i])\rightarrow y_{ij}\rightarrow (3,[j])$ by (\ref{eqC6.4.19}), (\ref{eqC6.4.21})-(\ref{eqC6.4.22}). Also, there exists some $z_i\in\lambda_{i+1}\cap\bar{\lambda}_1$ as $i+1\neq 1$ such that $(1,[i])\rightarrow z_i\rightarrow (1,[j])$. Since $\{(1,[\alpha,i]),(2,[\alpha,i])\}\rightarrow (1,[i])$ by (\ref{eqC6.4.18}), this subcase follows.
\\
\\(ii) $v=(q,[j]), w=(p,[\alpha,i])$ for each $p=1,2$ and $q=1,2,3$.
\indent\par Since $|\lambda_1|=|\lambda_{i+1}|=\lceil\frac{s}{2}\rceil$, there exists some $x_i\in \lambda_1\cap\lambda_{i+1}$ such that $(1,[j])\rightarrow x_i\rightarrow (2,[i])$ by (\ref{eqC6.4.19}) and (\ref{eqC6.4.21}). Also, there exist some $y_i\in\bar{\lambda}_1\cap\lambda_{i+1}$ and some $z_{ij}\in\bar{\lambda}_{j+1}\cap\lambda_{i+1}$ as $i+1\neq 1$ and $i\neq j$ respectively. So, $(2,[j])\rightarrow y_i\rightarrow (2,[i])$ and $(3,[j])\rightarrow z_{ij}\rightarrow (2,[i])$ with (\ref{eqC6.4.22}). Since $(2,[i])\rightarrow \{(1,[\alpha,i]),(2,[\alpha,i])\}$ by (\ref{eqC6.4.18}), this subcase follows.
\\
\\(iii) $v=(p,[i]), w=(q,[\beta,j])$ for each $p,q=1,2$.
\indent\par From (i), $d_{D_2}((1,[i]), (3,[j]))=2$. Since $i\neq j$, there exists some $x_{ij}\in \bar{\lambda}_{i+1}\cap \lambda_{j+1}$ such that $(2,[i])\rightarrow x_{ij}\rightarrow (3,[j])$ by (\ref{eqC6.4.19}) and (\ref{eqC6.4.22}). Since $(3,[j])\rightarrow\{(1,[\beta,j]),(2,[\beta,j])\}$ by (\ref{eqC6.4.20}), this subcase follows.
\\
\\(iv) $v=(q,[\beta,j]), w=(p,[i])$ for each $p,q=1,2$.
\indent\par By (\ref{eqC6.4.19})-(\ref{eqC6.4.21}), $\{(1, [\beta,j]), (2,[\beta,j])\}\rightarrow \{(1,[j]), (2,[j])\}$, $(1,[j])\rightarrow \lambda_1$, and $(2,[j])\rightarrow \bar{\lambda}_1$ and $|I^\mathtt{c}((p,[i]))|>0$.
\\
\\(v) $v=(p,[\alpha,i]), w=(q,[\beta,j])$ for each $p,q=1,2$.
\indent\par From (i), $d_{D_2}((p,[\alpha,i]),(3,[j]))=3$. Since $(3,[j])\rightarrow \{(1,[\beta,j]), (2,[\beta,j])\}$ by (\ref{eqC6.4.20}), this subcase follows.
\\
\\(vi) $v=(q,[\beta,j]), w=(p,[\alpha,i])$ for each $p,q=1,2$.
\indent\par From (iv), $d_{D_2}((q,[\beta,j]),(2,[i]))=3$. Since $(2,[i])\rightarrow \{(1,[\alpha,i]), (2,[\alpha,i])\}$ by (\ref{eqC6.4.18}), this subcase follows.
\\
\\Case 2.5. For each $1\le i\le |A_2|$, each ${{s}\choose{\lceil{s/2}\rceil}}\le j \le |A_2|+|A_3|$, each $1\le\alpha\le deg_T([i])-1$, and each $1\le\beta\le deg_T([j])-1$,
\\(i) $v=(p,[\alpha,i]), w=(q,[j])$ for each $p=1,2$ and $q=1,2,3$.
\indent\par By (\ref{eqC6.4.18}), $\{(1,[\alpha,i]),(2,[\alpha,i])\}\rightarrow (1,[i])$. Since $|\lambda_1|=|\lambda_{i+1}|=\lceil\frac{s}{2}\rceil$, there exists some $x_i\in\lambda_{i+1}\cap\lambda_1$ such that $(1,[i])\rightarrow x_i \rightarrow (1,[j])$ by (\ref{eqC6.4.19}) and (\ref{eqC6.4.24}). Also, there exist some $y_i\in\lambda_{i+1}\cap\bar{\lambda}_1$ and some $z_{ij}\in\lambda_{i+1}\cap\bar{\lambda}_{j-{{s}\choose{\lceil{s/2}\rceil}}+|A_2|+2}$ as $i+1\neq 1$ and $i+1\neq j-{{s}\choose{\lceil{s/2}\rceil}}+|A_2|+2$ respectively. So, $(1,[i])\rightarrow y_i\rightarrow (2,[j])$ and $(1,[i])\rightarrow z_{ij}\rightarrow (3,[j])$ with (\ref{eqC6.4.25}).
\\
\\(ii) $v=(q,[j]), w=(p,[\alpha,i])$ for each $p=1,2$ and $q=1,2,3$.
\indent\par Since $|\lambda_1|=|\lambda_{j-{{s}\choose{\lceil{s/2}\rceil}}+|A_2|+2}|=|\lambda_{i+1}|=\lceil\frac{s}{2}\rceil$, there exist some $x_i\in \lambda_1\cap\lambda_{i+1}$ and some $y_{ij}\in\lambda_{j-{{s}\choose{\lceil{s/2}\rceil}}+|A_2|+2}\cap\lambda_{i+1}$. So, $(2,[j])\rightarrow x_i\rightarrow (2,[i])$ and $(3,[j])\rightarrow y_{ij}\rightarrow (2,[i])$ by (\ref{eqC6.4.19}) and (\ref{eqC6.4.24})-(\ref{eqC6.4.25}). Also, $i+1\neq 1$ implies that there exists some $z_i\in \bar{\lambda}_1\cap\lambda_{i+1}$. So, $(1,[j])\rightarrow z_i\rightarrow (2,[i])$. Since $(2,[i])\rightarrow \{(1,[\alpha,i]),(2,[\alpha,i])\}$ by (\ref{eqC6.4.18}), this subcase follows.
\\
\\(iii) $v=(p,[i]), w=(q,[\beta,j])$ for each $p,q=1,2$.
\indent\par By (\ref{eqC6.4.19}) and (\ref{eqC6.4.23})-(\ref{eqC6.4.24}), $\lambda_1\rightarrow (1,[j])$, $\bar{\lambda}_1\rightarrow (2,[j])$, $\{(1,[j]),(2,[j])\} \rightarrow \{(1,[\beta,j]),$ $(2,[\beta,j])\}$ and $|O^\mathtt{c}((p,[i]))|>0$.
\\
\\(iv) $v=(q,[\beta,j]), w=(p,[i])$ for each $p,q=1,2$.
\indent\par By (\ref{eqC6.4.23}), $\{(1, [\beta,j]), (2,[\beta,j])\}\rightarrow (3,[j])$. Since $i+1\neq{{j-{{s}\choose{\lceil{s/2}\rceil}}+|A_2|+2}}$, there exists some $x_{ij}\in\lambda_{j-{{s}\choose{\lceil{s/2}\rceil}}+|A_2|+2}\cap \bar{\lambda}_{i+1}$ such that $(3,[j])\rightarrow x_{ij} \rightarrow (1, [i])$ by (\ref{eqC6.4.19}) and (\ref{eqC6.4.25}). Also, since $|\lambda_{{j-{{s}\choose{\lceil{s/2}\rceil}}+|A_2|+2}}|=|\lambda_{i+1}|=\lceil\frac{s}{2}\rceil$, there exists some $y_{ij}\in \lambda_{i+1}\cap \lambda_{{j-{{s}\choose{\lceil{s/2}\rceil}}+|A_2|+2}}$. So, $(3,[j])\rightarrow y_{ij}\rightarrow (2,[i])$.
\\
\\(v) $v=(p,[\alpha,i]), w=(q,[\beta,j])$ for each $p,q=1,2$.
\indent\par From (i), $d_{D_2}((p,[\alpha,i]),(1,[j]))=3$. Since $(1,[j])\rightarrow \{(1,[\beta,j]), (2,[\beta,j])\}$ by (\ref{eqC6.4.23}), this subcase follows.
\\
\\(vi) $v=(q,[\beta,j]), w=(p,[\alpha,i])$ for each $p,q=1,2$.
\indent\par From (iv), $d_{D_2}((q,[\beta,j]),(2,[i]))=3$. Since $(2,[i])\rightarrow \{(1,[\alpha,i]), (2,[\alpha,i])\}$ by (\ref{eqC6.4.18}), this subcase follows.
\\
\\Case 2.6. For each $1\le i\le |A_2|$, each $1\le\alpha\le deg_T([i])-1$, and each $[j]\in E$,
\\(i) $v=(p,[\alpha,i]), w=(q,[j])$ for each $p,q=1,2$.
\\(ii) $v=(q,[j]), w=(p,[\alpha,i])$ for each $p,q=1,2$.
\indent\par Since $\lambda_1\rightarrow \{(q,[j]),(2,[|A_2|+1])\}\rightarrow \bar{\lambda}_1$ by (\ref{eqC6.4.21}) and (\ref{eqC6.4.26}), this case follows from Cases 2.4(i)-(ii).
\\
\\Case 2.7. For each $|A_2|+1\le i\le{{s}\choose{\lceil{s/2}\rceil}}-1$, each $1\le\alpha\le deg_T([i])-1$, and each $[j]\in E$,
\\(i) $v=(p,[\alpha,i]), w=(q,[j])$ for each $p,q=1,2$.
\\(ii) $v=(q,[j]), w=(p,[\alpha,i])$ for each $p,q=1,2$.
\indent\par Since $\lambda_1\rightarrow \{(q,[j]),(1,[{{s}\choose{\lceil{s/2}\rceil}}])\}\rightarrow \bar{\lambda}_1$ by (\ref{eqC6.4.24}) and (\ref{eqC6.4.26}), this case follows from Case 2.2.
\\
\\Case 2.8. For each ${{s}\choose{\lceil{s/2}\rceil}}\le i\le|A_2|+|A_3|$, each $1\le\alpha\le deg_T([i])-1$, and each $[j]\in E$,
\\(i) $v=(p,[\alpha,i]), w=(q,[j])$ for each $p,q=1,2$.
\\(ii) $v=(q,[j]), w=(p,[\alpha,i])$ for each $p,q=1,2$.
\indent\par Since $\lambda_1\rightarrow \{(q,[j]),(2,[|A_2|+1])\}\rightarrow \bar{\lambda}_1$ by (\ref{eqC6.4.21}) and (\ref{eqC6.4.26}), this case follows from Case 2.2.
\\
\\Case 2.9. $v=(r_1,\mathtt{c})$ and $w=(r_2,\mathtt{c})$ for $r_1\neq r_2$ and $1\le r_1, r_2\le s$.
\indent\par Here, we want to prove a stronger claim, $d_{D_2}((r_1,\mathtt{c}), (r_2,\mathtt{c}))=2$. Let $x_1=(2,[|A_2|+1])$, $x_{k+1}=(2,[k])$ for $1\le k\le |A_2|$, and $x_{k+1}=(3,[k])$ for $|A_2|+1\le k\le {{s}\choose{\lceil{s/2}\rceil}}-1$. Observe that $\lambda_k\rightarrow x_k\rightarrow \bar{\lambda}_k$ for all $1\le k\le s$ and the subgraph induced by $V_1=(\mathbb{N}_s,\mathtt{c})$ and $V_2=\{x_i\mid 1\le i \le s\}$ is a complete bipartite graph $K(V_1,V_2)$. By Lemma \ref{lemC6.2.19}, $d_{D_2}((r_1,\mathtt{c}), (r_2,\mathtt{c}))=2$.
\\
\\Case 2.10. $v\in \{(1,[i]), (2,[i]), (3,[i]), (1,[\alpha,i]), (2,[\alpha,i])\}$ for each $1\le i\le deg_T(\mathtt{c})$ and $1\le\alpha\le deg_T([i])-1$, and $w=(r,\mathtt{c})$ for $1\le r\le s$.
\indent\par Note that there exists some $1\le k\le s$ such that $d_{D_2}(v,(k,\mathtt{c}))\le 2$, and $d_{D_2}((k,\mathtt{c}),w)\le 2$ by Case 2.9. Hence, it follows that $d_{D_2}(v,w)\le d_{D_2}(v,(k,\mathtt{c}))+d_{D_2}((k,\mathtt{c}),w)\le 4$.
\\
\\Case 2.11. $v=(r,\mathtt{c})$ for $r=1,2,\ldots, s$ and $w\in \{(1,[i]), (2,[i]), (3,[i]), (1,[\alpha,i]), (2,[\alpha,i])\}$ for each $1\le i\le deg_T(\mathtt{c})$ and $1\le\alpha\le deg_T([i])-1$.
\indent\par Note that there exists some $1\le k\le s$ such that $d_{D_2}((k,\mathtt{c}), w)\le 2$, and $d_{D_2}(v,(k,\mathtt{c}))\le 2$ by Case 2.9. Hence, it follows that $d_{D_2}(v,w)\le d_{D_2}(v,(k,\mathtt{c}))+d_{D_2}((k,\mathtt{c}),w)\le 4$.
\\
\\Case 2.12. $v=(p,[i])$ and $w=(q, [j])$, where $1\le p,q\le 3$ and $1\le i,j\le deg_T(\mathtt{c})$.
\noindent\par This follows from the fact that $|O^\mathtt{c}((p,[i]))|>0$, $|I^\mathtt{c}((q,[j]))|>0$, and $d_{D_2}((r_1,\mathtt{c}), (r_2,\mathtt{c}))$ $=2$ for any $r_1\neq r_2$ and $1\le r_1, r_2\le s$.
\\
\\Case 3. $|A_3|\ge 2$, $|A_2|\ge \lceil\frac{s}{2}\rceil \lfloor\frac{s}{2}\rfloor$, $2|A_2|+|A_3|= 2{{s}\choose{\lceil{s/2}\rceil}}-1$, and $s\ge 5$.
\indent\par Let $\psi=(\mathbb{N}_{\lfloor\frac{s}{2}\rfloor},\mathtt{c})$ and $I_{\psi}=\{\lambda\in {{(\mathbb{N}_s,\mathtt{c})}\choose{\lceil s/2\rceil}}\mid |\lambda\cap \psi|=1\}$. Since there are $\lfloor\frac{s}{2}\rfloor$ choices in $\psi$ as the common element of $\lambda$ and $\psi$, and for each such choice, there are $\lceil\frac{s}{2}\rceil$ choices (elements not in $\psi$) for the remaining $\lfloor\frac{s}{2}\rfloor$ elements of $\lambda$, $|I_{\psi}|=\lceil\frac{s}{2}\rceil\lfloor\frac{s}{2}\rfloor$. Particularly, note that $\bar{\psi}\not\in I_{\psi}$. Let $O_{\psi}=\{\lambda\in{{(\mathbb{N}_s,\mathtt{c})}\choose{\lceil s/2\rceil}}\mid \psi\subset\lambda\}$ and observe that $|O_{\psi}|=\lceil\frac{s}{2}\rceil$. 
\indent\par Our aim is to design an orientation in which the elements of $I_{\psi}$ and $O_{\psi}$ are used as $I^\mathtt{c}((2,[i]))$ and $O^\mathtt{c}((1,[i]))$ respectively, where $[i]\in A_2$, i.e., $I_\psi\subseteq B^I_2$ and $O_\psi\subseteq B^O_2$. To achieve this, we introduce two new listings of the elements in ${{(\mathbb{N}_s,\mathtt{c})}\choose{\lceil s/2\rceil}}$. Let ${{(\mathbb{N}_s,\mathtt{c})}\choose{\lceil s/2\rceil}}=\{\gamma_1, \gamma_2, \ldots, \gamma_{{s}\choose{\lceil{s/2}\rceil}}\}=\{\mu_1, \mu_2, \ldots, \mu_{{s}\choose{\lceil{s/2}\rceil}}\}$ such that $\bar{\psi}=\gamma_{{s}\choose{\lceil{s/2}\rceil}}$, $I_{\psi}=\{\gamma_1, \gamma_2, \ldots, \gamma_{\lceil\frac{s}{2}\rceil\lfloor\frac{s}{2}\rfloor}\}$ and $O_{\psi}=\{\mu_1, \mu_2, \ldots, \mu_{\lceil\frac{s}{2}\rceil}\}$. The denotation of the remaining $\gamma_i$'s and $\mu_j$'s can be arbitrary. Assume further that $A_2=\{[i]\mid i\in\mathbb{N}_{|A_2|}\}$ and $A_3=\{[i]\mid i\in\mathbb{N}_{|A_2|+|A_3|}-\mathbb{N}_{|A_2|}\}$. Define an orientation $D_3$ of $\mathcal{H}$ as follows.
\begin{align}
& (2,[i])\rightarrow \{(1,[\alpha,i]),(2,[\alpha,i])\}\rightarrow (1,[i]),\label{eqC6.4.27}\\
& \bar{\mu}_i\rightarrow (1,[i]) \rightarrow \mu_i,\text{ and } \gamma_i\rightarrow (2,[i]) \rightarrow \bar{\gamma}_i,\label{eqC6.4.28}
\end{align}
for all $1\le i\le |A_2|$ and all $1\le \alpha\le deg_T([i])-1$, i.e., the $\lceil\frac{s}{2}\rceil$-sets $\gamma_1,\gamma_2,\ldots,\gamma_{|A_2|}$ ($\mu_1,\mu_2,\ldots,\mu_{|A_2|}$ resp.) are used as `in-sets' (`out-sets' resp.) to construct $B^I_2$ ($B^O_2$ resp.).
\begin{align}
& \{(1,[j]),(2,[j])\}\rightarrow (1,[\beta,j])\rightarrow (3,[j]),\label{eqC6.4.29}\\
& \{(1,[j]),(3,[j])\}\rightarrow (2,[\beta,j])\rightarrow (2,[j]),\label{eqC6.4.30}\\
& \bar{\psi}\rightarrow (1,[j])\rightarrow \psi, \label{eqC6.4.31}\text{ and}\\
& \bar{\mu}_j\rightarrow \{(2,[j]), (3,[j])\}\rightarrow \mu_j\label{eqC6.4.32}
\end{align}
for all $|A_2|+1\le j\le{{s}\choose{\lceil{s/2}\rceil}}$ and all $1\le \beta\le deg_T([j])-1$, i.e., the $\lceil\frac{s}{2}\rceil$-sets $\mu_{|A_2|+1},\mu_{|A_2|+2},$ $\ldots,\mu_{{{s}\choose{\lceil{s/2}\rceil}}}$ are used as `out-sets' to construct $B^O_3$.
\begin{align}
& (3,[k])\rightarrow (1,[\theta,k])\rightarrow \{(1,[k]),(2,[k])\},\label{eqC6.4.33}\\
& (2,[k])\rightarrow (2,[\theta,k])\rightarrow \{(1,[k]),(3,[k])\},\label{eqC6.4.34}\\
& \bar{\psi} \rightarrow (1,[k])\rightarrow \psi,\label{eqC6.4.35}\text{ and}\\
& \gamma_{k-{{s}\choose{\lceil{s/2}\rceil}}+|A_2|}\rightarrow \{(2,[k]), (3,[k])\}\rightarrow \bar{\gamma}_{k-{{s}\choose{\lceil{s/2}\rceil}}+|A_2|}\label{eqC6.4.36}
\end{align}
for all ${{s}\choose{\lceil{s/2}\rceil}}+1\le k\le |A_2|+|A_3|$ and all $1\le \theta\le deg_T([k])-1$, i.e., the $\lceil\frac{s}{2}\rceil$-sets $\gamma_{|A_2|+1},\gamma_{|A_2|+2},$ $\ldots,\gamma_{{{s}\choose{\lceil{s/2}\rceil}}-1}$ are used as `in-sets' to construct $B^I_3$.
\begin{align}
\bar{\psi} \rightarrow \{(1,[l]), (2,[l])\} \rightarrow \psi \label{eqC6.4.37}
\end{align}
for all $[l]\in E$. (See Figures \ref{figC6.4.15} and \ref{figC6.4.16} for $D_3$ when $s=5$.)
\\
\\Claim: $d_{D_3}(v,w)\le 4$ for all $v,w\in V(D_3)$.
\\
\\Case 3.1.1. $v,w \in \{(1,[\alpha,i]),(2,[\alpha,i]),(1,[i]),(2,[i])\}$ for each $1\le i\le |A_2|$ and $1\le\alpha\le deg_T([i])-1$.
\indent\par For each $1\le i\le |A_2|$, since $|\gamma_i|=|\mu_i|=\lceil\frac{s}{2}\rceil$, there exists some $x_i\in \gamma_i\cap\mu_i$. Then, by (\ref{eqC6.4.27})-(\ref{eqC6.4.28}), $(2,[i])\rightarrow \{(1,[\alpha,i]),(2,[\alpha,i])\}\rightarrow (1,[i])\rightarrow x_i\rightarrow (2,[i])$ guarantees a directed $C_4$.
\\
\\Case 3.1.2. $v,w \in \{(1,[\alpha,i]),(2,[\alpha,i]),(1,[i]),(2,[i]),(3,[i])\}$ for each $|A_2|+1\le i\le {{s}\choose{\lceil{s/2}\rceil}}$, and $1\le\alpha\le deg_T([i])-1$.
\indent\par By (\ref{eqC6.4.29})-(\ref{eqC6.4.32}), since $|\mu_i|=|\bar{\psi}|=\lceil\frac{s}{2}\rceil$, there exists some $x_i\in \mu_i\cap\bar{\psi}$ so that $\{(2,[i]),(3,[i])\}\rightarrow x_i\rightarrow (1,[i])\rightarrow \{(1,[\alpha,i]),(2,[\alpha,i])\}$, $(1,[\alpha,i])\rightarrow (3,[i])$ and $(2,[\alpha,i])\rightarrow (2,[i])$.
\\
\\Case 3.1.3. $v,w \in \{(1,[\alpha,i]),(2,[\alpha,i]),(1,[i]),(2,[i]),(3,[i])\}$ for each ${{s}\choose{\lceil{s/2}\rceil}}+1\le i\le |A_2|+|A_3|$, and $1\le\alpha\le deg_T([i])-1$.
\indent\par By (\ref{eqC6.4.33})-(\ref{eqC6.4.36}), since $2|A_2|+|A_3|= 2{{s}\choose{\lceil{s/2}\rceil}}-1$ implies $\bar{\psi}=\gamma_{{s}\choose{\lceil{s/2}\rceil}}\neq\gamma_{i-{{s}\choose{\lceil{s/2}\rceil}}+|A_2|}$, there exists some $x_i\in \gamma_{i-{{s}\choose{\lceil{s/2}\rceil}}+|A_2|}\cap\psi$ so that $\{(1,[\alpha,i]),(2,[\alpha,i])\}\rightarrow (1,[i])\rightarrow x_i\rightarrow \{(2,[i]), (3,[i])\}$, $(3,[i])\rightarrow (1,[\alpha,i])$ and $(2,[i])\rightarrow (2,[\alpha,i])$.
\\
\\Case 3.2. For each $1\le i,j \le  |A_2|$, $i\neq j$, each $1\le \alpha\le deg_T([i])-1$, and each $1\le \beta\le deg_T([j])-1$,
\\(i) $v=(p,[\alpha,i]), w=(q,[j])$ for any $p,q=1,2$.
\indent\par Since $i\neq j$, there exists $x_{ij}\in\mu_i\cap\bar{\mu}_j$ such that $\{(1,[\alpha,i]),(2,[\alpha,i])\}\rightarrow (1,[i])\rightarrow x_{ij}\rightarrow (1,[j])$ by (\ref{eqC6.4.27})-(\ref{eqC6.4.28}). Also, $|\mu_i|=|\gamma_j|=\lceil\frac{s}{2}\rceil$, there exists some $y_{ij}\in \mu_i\cap\gamma_j$ so that $\{(1,[\alpha,i]),(2,[\alpha,i])\}\rightarrow (1,[i])\rightarrow y_{ij}\rightarrow (2,[j])$.
\\
\\(ii) $v=(q,[j]), w=(p,[\alpha,i])$ for any $p,q=1,2$.
\indent\par By (\ref{eqC6.4.27})-(\ref{eqC6.4.28}), $|\mu_j|=|\gamma_i|=\lceil\frac{s}{2}\rceil$, there exists some $x_{ij}\in \mu_j\cap\gamma_i$ so that $(1,[j]) \rightarrow x_{ij}\rightarrow (2,[i])\rightarrow \{(1,[\alpha,i]),(2,[\alpha,i])\}$. Since $i\neq j$, there exists some $y_{ij}\in \bar{\gamma}_j\cap\gamma_i$ so that $(2,[j]) \rightarrow y_{ij}\rightarrow (2,[i])\rightarrow \{(1,[\alpha,i]),(2,[\alpha,i])\}$.
\\
\\(iii) $v=(p,[\alpha,i]), w=(q,[\beta,j])$ for each $p,q=1,2$.
\indent\par From (i), $d_{D_3}((p,[\alpha,i]),(2, [j]))=3$. Since $(2,[j])\rightarrow \{(1,[\beta,j]), (2,[\beta,j])\}$ by (\ref{eqC6.4.27}), this subcase follows.
\\
\\Case 3.3. For each $|A_2|+1 \le i,j \le {{s}\choose{\lceil{s/2}\rceil}}$, $i\neq j$, each $1\le\alpha\le deg_T([i])-1$, and each $1\le\beta\le deg_T([j])-1$,
\\(i) $v=(p,[\alpha,i]), w=(q,[j])$ for each $p=1,2$ and $q=1,2,3$.
\indent\par Since $i\neq j$, there exists some $x_{ij}\in \mu_i\cap\bar{\mu}_j$. Also, $|\bar{\psi}|=|\mu_i|=\lceil\frac{s}{2}\rceil$ implies there exists some $y_{ij}\in\mu_i\cap\bar{\psi}$. Then, by (\ref{eqC6.4.29})-(\ref{eqC6.4.32}), $(1,[\alpha,i])\rightarrow (3,[i])$, $(2,[\alpha,i])\rightarrow (2,[i])$, $\{(2,[i]),(3,[i])\}\rightarrow x_{ij}\rightarrow\{(2,[j]),(3,[j])\}$, and $\{(2,[i]),(3,[i])\}\rightarrow y_{ij}\rightarrow (1,[j])$.
\\
\\(ii) $v=(q,[j]), w=(p,[\alpha,i])$ for each $p=1,2$ and $q=1,2,3$.
\indent\par By definition of $O_{\psi}$ and $i\ge |A_2|+1>\lceil\frac{s}{2}\rceil \lfloor\frac{s}{2}\rfloor$, $\mu_i\not\in O_{\psi}$. So, there exists some $x_j\in\psi\cap\bar{\mu}_i$ so that $(1,[j])\rightarrow x_j\rightarrow\{(2,[i]), (3,[i])\}$ by (\ref{eqC6.4.31})-(\ref{eqC6.4.32}). Also, since $i\neq j$, there exists some $y_{ij}\in \bar{\mu}_i\cap\mu_j$ so that $\{(2,[j]), (3,[j])\}\rightarrow y_{ij}\rightarrow \{(2,[i]), (3,[i])\}$. With $(2,[i])\rightarrow (1,[\alpha,i])$ and $(3,[i])\rightarrow (2,[\alpha,i])$ by (\ref{eqC6.4.29})-(\ref{eqC6.4.30}), this subcase follows.
\\
\\(iii) $v=(p,[\alpha,i]), w=(q,[\beta,j])$ for each $p,q=1,2$.
\indent\par From (i), $d_{D_3}((p,[\alpha,i]),(1, [j]))=3$. Since $(1,[j])\rightarrow \{(1,[\beta,j]), (2,[\beta,j])\}$ by (\ref{eqC6.4.29})-(\ref{eqC6.4.30}), this subcase follows.

\begin{rmk} In Subcase 3.3(ii), we invoked the condition $|A_2|\ge \lceil\frac{s}{2}\rceil \lfloor\frac{s}{2}\rfloor$. This is the reason for the need of two different constructions $D_2$ and $D_3$ in proving sufficiency. In Case 2, without $|A_2|\ge \lceil\frac{s}{2}\rceil \lfloor\frac{s}{2}\rfloor$, the construction $D_3$ may not satisfy $d(D_3)=4$.
\end{rmk}
\noindent Case 3.4. For each ${{s}\choose{\lceil{s/2}\rceil}}+1 \le i,j \le |A_2|+|A_3|$, $i\neq j$, each $1\le\alpha\le deg_T([i])-1$, and each $1\le\beta\le deg_T([j])-1$,
\\(i) $v=(p,[\alpha,i]), w=(q,[j])$ for each $p=1,2$ and $q=1,2,3$.
\indent\par Since $j\le |A_2|+|A_3|$ and $2|A_2|+|A_3|= 2{{s}\choose{\lceil{s/2}\rceil}}-1$, we have $\bar{\psi}=\gamma_{{s}\choose{\lceil{s/2}\rceil}}\neq\gamma_{j-{{s}\choose{\lceil{s/2}\rceil}}+|A_2|}$. Note also that $|\psi|=\lfloor\frac{s}{2}\rfloor$ and $| \gamma_{j-{{s}\choose{\lceil{s/2}\rceil}}+|A_2|}|=\lceil \frac{s}{2}\rceil$. So, there exists some $x_j\in \psi\cap \gamma_{j-{{s}\choose{\lceil{s/2}\rceil}}+|A_2|}$ such that $\{(1,[\alpha,i]),(2,[\alpha,i])\}\rightarrow(1,[i])\rightarrow x_j \rightarrow \{(2,[j]),(3,[j])\}$ by (\ref{eqC6.4.33})-(\ref{eqC6.4.36}).  Similarly, there exists some $y_j\in \bar{\gamma}_{j-{{s}\choose{\lceil{s/2}\rceil}}+|A_2|}\cap\bar{\psi}$ such that $\{(2,[i]), (3,[i])\}\rightarrow y_j\rightarrow(1,[j])$. Since $(1,[\alpha,i])\rightarrow (2,[i])$ and $(2,[\alpha,i])\rightarrow (3,[i])$, this subcase follows.
\\
\\(ii) $v=(q,[j]), w=(p,[\alpha,i])$ for each $p=1,2$ and $q=1,2,3$.
\indent\par Note that $\bar{\psi}=\gamma_{{s}\choose{\lceil{s/2}\rceil}}\neq\gamma_{j-{{s}\choose{\lceil{s/2}\rceil}}+|A_2|}$, which implies there exists some $x_i\in \psi\cap \gamma_{i-{{s}\choose{\lceil{s/2}\rceil}}+|A_2|}$. So, $(1,[j])\rightarrow x_i\rightarrow\{(2,[i]),(3,[i])\}$ by (\ref{eqC6.4.35})-(\ref{eqC6.4.36}). Since $i\neq j$, there exists some $y_{ij}\in \bar{\gamma}_{j-{{s}\choose{\lceil{s/2}\rceil}}+|A_2|}\cap \gamma_{i-{{s}\choose{\lceil{s/2}\rceil}}+|A_2|}$ so that $\{(2,[j]),(3,[j])\}\rightarrow y_{ij}\rightarrow \{(2,[i]),(3,[i])\}$. Since $(2,[i])\rightarrow (2,[\alpha,i])$ and $(3,[i])\rightarrow (1,[\alpha,i])$ by (\ref{eqC6.4.33})-(\ref{eqC6.4.34}), this subcase follows.
\\
\\(iii) $v=(p,[\alpha,i]), w=(q,[\beta,j])$ for each $p,q=1,2$.
\indent\par From (i), $d_{D_3}((p,[\alpha,i]),(r, [j]))=3$ for $r=2,3$. Since $(2,[j])\rightarrow (2,[\beta,j])$ and $(3,[j])\rightarrow (1,[\beta,j])$ by (\ref{eqC6.4.33})-(\ref{eqC6.4.34}), this subcase follows.
\\
\\Case 3.5. For each $1\le i\le |A_2|$, each $|A_2|+1\le j\le {{s}\choose{\lceil{s/2}\rceil}}$, each $1\le\alpha\le deg_T([i])-1$, and each $1\le\beta\le deg_T([j])-1$,
\\(i) $v=(p,[\alpha,i]), w=(q,[j])$ for each $p=1,2$ and $q=1,2,3$.
\indent\par Since $i\neq j$, there exists some $x_{ij}\in \mu_i\cap\bar{\mu}_j$. Also, $|\bar{\psi}|=|\mu_i|=\lceil\frac{s}{2}\rceil$ implies there exists some $y_i\in \mu_i\cap\bar{\psi}$. Then by (\ref{eqC6.4.28}) and (\ref{eqC6.4.31})-(\ref{eqC6.4.32}), $\{(1,[\alpha,i]),(2,[\alpha,i])\}\rightarrow (1,[i])\rightarrow x_{ij}\rightarrow \{(2,[j]),(3,[j])\}$, and $\{(1,[\alpha,i]),(2,[\alpha,i])\}\rightarrow(1,[i])\rightarrow y_i\rightarrow (1,[j])$.
\\
\\(ii) $v=(q,[j]), w=(p,[\alpha,i])$ for each $p=1,2$ and $q=1,2,3$.
\indent\par Since $\bar{\psi}=\gamma_{{s}\choose{\lceil{s/2}\rceil}}\neq \gamma_i$, there exists some $x_i\in \psi\cap\gamma_i$. Also, $|\gamma_i|=|\mu_j|=\lceil\frac{s}{2}\rceil$ implies there exists some $y_{ij}\in \mu_j\cap\gamma_i$. So by (\ref{eqC6.4.27})-(\ref{eqC6.4.28}) and (\ref{eqC6.4.31})-(\ref{eqC6.4.32}), $(1,[j])\rightarrow x_i\rightarrow (2,[i])\rightarrow\{(1,[\alpha,i]),(2,[\alpha,i])\}$ and $\{(2,[j]),(3,[j])\}\rightarrow y_{ij}\rightarrow (2,[i])\rightarrow\{(1,[\alpha,i]),(2,[\alpha,i])\}$.
\\
\\(iii) $v=(p,[i]), w=(q,[\beta,j])$ for each $p,q=1,2$.
\indent\par Since $i\neq j$, there exists some $x_{ij}\in \mu_i\cap\bar{\mu}_j$ such that $(1,[i])\rightarrow x_{ij}\rightarrow \{(2,[j]),(3,[j])\}$ by (\ref{eqC6.4.28}) and (\ref{eqC6.4.32}). Furthermore, $(2,[j])\rightarrow (1,[\beta,j])$, and $(3,[j])\rightarrow (2,[\beta,j])$ by (\ref{eqC6.4.29})-(\ref{eqC6.4.30}). Also with (\ref{eqC6.4.31}), since $\bar{\psi}=\gamma_{{s}\choose{\lceil{s/2}\rceil}}\neq \gamma_i$, there exists some $y_i\in \bar{\gamma}_i\cap\bar{\psi}$. So, $(2,[i])\rightarrow y_i\rightarrow (1,[j])\rightarrow \{(1,[\beta,j]),(2,[\beta,j])\}$.
\\
\\(iv) $v=(q,[\beta,j]), w=(p,[i])$ for each $p,q=1,2$.
\indent\par Since $|\mu_j|=|\gamma_i|=\lceil\frac{s}{2}\rceil$ implies there exists some $x_{ij}\in \mu_j\cap\gamma_i$. Then, $(1,[\beta,j])\rightarrow (3,[j])$, $(2,[\beta,j])\rightarrow (2,[j])$, and $\{(2,[j]),(3,[j])\}\rightarrow x_{ij}\rightarrow (2,[i])$ by (\ref{eqC6.4.28})-(\ref{eqC6.4.30}) and (\ref{eqC6.4.32}). Since $i\neq j$, there exists some $y_{ij}\in \mu_j\cap\bar{\mu}_i$ such that $\{(2,[j]),(3,[j])\}\rightarrow y_{ij}\rightarrow (1,[i])$.
\\
\\(v) $v=(p,[\alpha,i]), w=(q,[\beta,j])$ for each $p,q=1,2$.
\indent\par From (i), $d_{D_3}((p,[\alpha,i]),(1,[j]))=3$. Since $(1,[j])\rightarrow \{(1,[\beta,j]), (2,[\beta,j])\}$ by (\ref{eqC6.4.29})-(\ref{eqC6.4.30}), this subcase follows.
\\
\\(vi) $v=(q,[\beta,j]), w=(p,[\alpha,i])$ for each $p,q=1,2$.
\indent\par From (iv), $d_{D_3}((q,[\beta,j]),(2,[i]))=3$. Since $(2,[i])\rightarrow \{(1,[\alpha,i]), (2,[\alpha,i])\}$ by (\ref{eqC6.4.27}), this subcase follows.
\\
\\Case 3.6. For each $1\le i\le |A_2|$, each ${{s}\choose{\lceil{s/2}\rceil}}+1 \le j \le |A_2|+|A_3|$, each $1\le\alpha\le deg_T([i])-1$, and each $1\le\beta\le deg_T([j])-1$,
\\(i) $v=(p,[\alpha,i]), w=(q,[j])$ for each $p=1,2$ and $q=1,2,3$.
\indent\par Since $|\bar{\psi}|=|\gamma_{j-{{s}\choose{\lceil{s/2}\rceil}}+|A_2|}|=|\mu_i|=\lceil\frac{s}{2}\rceil$, there exist some $x_i\in \mu_i\cap\bar{\psi}$ and $y_{ij}\in \mu_i\cap\gamma_{j-{{s}\choose{\lceil{s/2}\rceil}}+|A_2|}$. So, $\{(1,[\alpha,i]),(2,[\alpha,i])\}\rightarrow (1,[i])\rightarrow \{x_i,y_{ij}\}$, $x_i\rightarrow (1,[j])$ and $y_{ij}\rightarrow \{(2,[j]),(3,[j])\}$ by (\ref{eqC6.4.27})-(\ref{eqC6.4.28}) and (\ref{eqC6.4.35})-(\ref{eqC6.4.36}).
\\
\\(ii) $v=(q,[j]), w=(p,[\alpha,i])$ for each $p=1,2$ and $q=1,2,3$.
\indent\par Since $\bar{\psi}=\gamma_{{s}\choose{\lceil{s/2}\rceil}}\neq \gamma_i$, there exists some $x_i\in \psi\cap\gamma_i$. So, $(1,[j])\rightarrow x_i\rightarrow (2,[i])\rightarrow \{(1,[\alpha,i]),(2,[\alpha,i])\}$ by (\ref{eqC6.4.27})-(\ref{eqC6.4.28}) and (\ref{eqC6.4.35}). Also with (\ref{eqC6.4.36}), since $i\neq j-{{s}\choose{\lceil{s/2}\rceil}}+|A_2|$, there exists some $y_{ij}\in \bar{\gamma}_{j-{{s}\choose{\lceil{s/2}\rceil}}+|A_2|}\cap \gamma_i$ so that $\{(2,[j]),(3,[j])\}\rightarrow y_{ij}\rightarrow(2,[i])\rightarrow \{(1,[\alpha,i]),(2,[\alpha,i])\}$.
\\
\\(iii) $v=(p,[i]), w=(q,[\beta,j])$ for each $p,q=1,2$.
\indent\par Since $|\gamma_{j-{{s}\choose{\lceil{s/2}\rceil}}+|A_2|}|=|\mu_i|=\lceil\frac{s}{2}\rceil$, there exists some $x_{ij}\in \mu_i\cap\gamma_{j-{{s}\choose{\lceil{s/2}\rceil}}+|A_2|}$. Then, $(1,[i])\rightarrow x_{ij}\rightarrow \{(2,[j]),(3,[j])\}$ by (\ref{eqC6.4.28}) and (\ref{eqC6.4.36}). Since $i\neq j-{{s}\choose{\lceil{s/2}\rceil}}+|A_2|$, there exists some $y_{ij}\in \bar{\gamma}_i\cap\gamma_{j-{{s}\choose{\lceil{s/2}\rceil}}+|A_2|}$ so that $(2,[i]) \rightarrow y_{ij}\rightarrow \{(2,[j]),(3,[j])\}$. With $(3,[j])\rightarrow (1,[\beta,j])$ and $(2,[j])\rightarrow (2,[\beta,j])$ by (\ref{eqC6.4.33})-(\ref{eqC6.4.34}), this subcase follows.
\\
\\(iv) $v=(q,[\beta,j]), w=(p,[i])$ for each $p,q=1,2$.
\indent\par Since $\bar{\psi}=\gamma_{{s}\choose{\lceil{s/2}\rceil}}\neq \gamma_i$, there exists some $x_i\in \psi\cap\gamma_i$. So, $\{(1,[\beta,j]),(2,[\beta,j])\}\rightarrow (1,[j])\rightarrow x_i\rightarrow (2,[i])$ by (\ref{eqC6.4.28}) and (\ref{eqC6.4.33})-(\ref{eqC6.4.35}).
\indent\par Now, observe that $|\gamma_k\cap \psi|=0$ if and only if $\gamma_k=\bar{\psi}$ if and only if $k={{s}\choose{\lceil{s/2}\rceil}}$. And, by definition of $I_{\psi}$, $\gamma_{j-{{s}\choose{\lceil{s/2}\rceil}}+|A_2|}\not\in I_{\psi}$ so that $|\gamma_{j-{{s}\choose{\lceil{s/2}\rceil}}+|A_2|}\cap\psi|\ge 2$ for all ${{s}\choose{\lceil{s/2}\rceil}}+1 \le j \le |A_2|+|A_3|$. So, $|\gamma_{j-{{s}\choose{\lceil{s/2}\rceil}}+|A_2|} \cup \psi|=|\gamma_{j-{{s}\choose{\lceil{s/2}\rceil}}+|A_2|}|+|\psi|-|\gamma_{j-{{s}\choose{\lceil{s/2}\rceil}}+|A_2|}\cap\psi|\le \lceil\frac{s}{2}\rceil+\lfloor\frac{s}{2}\rfloor-2=s-2$; there exist at least two elements neither in $\gamma_{j-{{s}\choose{\lceil{s/2}\rceil}}+|A_2|}$, nor in $\psi$. So, $|\bar{\gamma}_{j-{{s}\choose{\lceil{s/2}\rceil}}+|A_2|}\cup\psi|\ge \lfloor\frac{s}{2}\rfloor+2$. Since $|\bar{\gamma}_{j-{{s}\choose{\lceil{s/2}\rceil}}+|A_2|}\cup\psi|+|\bar{\mu}_i|\ge \lfloor\frac{s}{2}\rfloor+2+\lfloor\frac{s}{2}\rfloor>s$, it follows that there exists some $y_{ij}\in (\bar{\gamma}_{j-{{s}\choose{\lceil{s/2}\rceil}}+|A_2|}\cup\psi)\cap\bar{\mu}_i$. Hence, $(1,[\beta,j])\rightarrow \{(1,[j]),(2,[j])\}$, $(2,[\beta,j])\rightarrow \{(1,[j]),(3,[j])\}$, and either $\{(2,[j]),(3,[j])\}\rightarrow y_{ij}\rightarrow (1,[i])$ or $(1,[j])\rightarrow y_{ij}\rightarrow (1,[i])$ by (\ref{eqC6.4.28}) and (\ref{eqC6.4.33})-(\ref{eqC6.4.36}).
\begin{rmk}
Cases 3.6(iv) and 3.7(iv) justify the way $\psi$ and $I_{\psi}$ are defined. 
\end{rmk}
\noindent (v) $v=(p,[\alpha,i]), w=(q,[\beta,j])$ for each $p,q=1,2$.
\indent\par From (i), $d_{D_3}((p,[\alpha,i]),(r,[j]))=3$ for $r=2,3$. Since $(2,[j])\rightarrow (2,[\beta,j])$ and $(3,[j])\rightarrow (1,[\beta,j])$ by (\ref{eqC6.4.33})-(\ref{eqC6.4.34}), this subcase follows.
\\
\\(vi) $v=(q,[\beta,j]), w=(p,[\alpha,i])$ for each $p,q=1,2$.
\indent\par From (iv), $d_{D_3}((q,[\beta,j]),(2,[i]))=3$. Since $(2,[i])\rightarrow \{(1,[\alpha,i]), (2,[\alpha,i])\}$ by (\ref{eqC6.4.27}), this subcase follows.
\\
\\Case 3.7. For each $|A_2|+1\le i\le {{s}\choose{\lceil{s/2}\rceil}}$, each ${{s}\choose{\lceil{s/2}\rceil}}+1 \le j \le |A_2|+|A_3|$, each $1\le\alpha\le deg_T([i])-1$, and each $1\le\beta\le deg_T([j])-1$,
\\(i) $v=(p,[\alpha,i]), w=(q,[j])$ for each $p=1,2$ and $q=1,2,3$.
\indent\par Since $|\bar{\psi}|=|\gamma_{j-{{s}\choose{\lceil{s/2}\rceil}}+|A_2|}|=|\mu_i|=\lceil\frac{s}{2}\rceil$, there exist some $x_i\in \mu_i\cap\bar{\psi}$ and $y_{ij}\in \mu_i\cap\gamma_{j-{{s}\choose{\lceil{s/2}\rceil}}+|A_2|}$. So, $(1,[\alpha,i])\rightarrow (3,[i])$, $(2,[\alpha,i])\rightarrow (2,[i])$, $\{(2,[i]),(3,[i])\}\rightarrow x_i\rightarrow (1,[j])$ and $\{(2,[i]),(3,[i])\}\rightarrow y_{ij}\rightarrow \{(2,[j]),(3,[j])\}$ by (\ref{eqC6.4.29})-(\ref{eqC6.4.32}) and (\ref{eqC6.4.35})-(\ref{eqC6.4.36}).
\\
\\(ii) $v=(q,[j]), w=(p,[\alpha,i])$ for each $p=1,2$ and $q=1,2,3$.
\indent\par By definition of $O_{\psi}$, $\mu_i\not\in O_{\psi}$. So, $\psi \not\subset \mu_i$. Equivalently, there exists some $x_i\in\psi\cap \bar{\mu}_i$. So, $(1,[j])\rightarrow x_i\rightarrow \{(2,[i]), (3,[i])\}$, $(2,[i])\rightarrow (1,[\alpha,i])$ and $(3,[i])\rightarrow (2,[\alpha,i])$ by (\ref{eqC6.4.29})-(\ref{eqC6.4.30}), (\ref{eqC6.4.32}) and (\ref{eqC6.4.35}). Also, since $\bar{\psi}=\gamma_{{s}\choose{\lceil{s/2}\rceil}}\neq \gamma_{j-{{s}\choose{\lceil{s/2}\rceil}}+|A_2|}$, there exists some $y_j\in \bar{\gamma}_{j-{{s}\choose{\lceil{s/2}\rceil}}+|A_2|}\cap\bar{\psi}$. So, $\{(2,[j]),(3,[j])\}\rightarrow y_j\rightarrow (1,[i])\rightarrow \{(1,[\alpha,i]),(2,[\alpha,i])\}$ with (\ref{eqC6.4.31}) and (\ref{eqC6.4.36}).
\\
\\(iii) $v=(p,[i]), w=(q,[\beta,j])$ for each $p=1,2,3$ and $q=1,2$.
\indent\par Since $\bar{\psi}=\gamma_{{s}\choose{\lceil{s/2}\rceil}}\neq \gamma_{j-{{s}\choose{\lceil{s/2}\rceil}}+|A_2|}$, there exists some $x_j \in \psi\cap\gamma_{j-{{s}\choose{\lceil{s/2}\rceil}}+|A_2|}$ such that $(1,[i])\rightarrow x_j\rightarrow \{(2,[j]),(3,[j])\}$ by (\ref{eqC6.4.31}) and (\ref{eqC6.4.36}). Also, since $|\gamma_{j-{{s}\choose{\lceil{s/2}\rceil}}+|A_2|}|=|\mu_i|=\lceil\frac{s}{2}\rceil$, there exists some $y_{ij}\in \mu_i\cap\gamma_{j-{{s}\choose{\lceil{s/2}\rceil}}+|A_2|}$. So with (\ref{eqC6.4.32}), $\{(2,[i]),(3,[i])\}\rightarrow y_{ij}\rightarrow\{(2,[j]),(3,[j])\}$. With $(3,[j])\rightarrow (1,[\beta,j])$ and $(2,[j])\rightarrow (2,[\beta,j])$ by (\ref{eqC6.4.33})-(\ref{eqC6.4.34}), this subcase follows.
\\
\\(iv) $v=(q,[\beta,j]), w=(p,[i])$ for each $p=1,2$ and $q=1,2,3$.
\indent\par Since $\bar{\psi}=\gamma_{{s}\choose{\lceil{s/2}\rceil}}\neq \gamma_{j-{{s}\choose{\lceil{s/2}\rceil}}+|A_2|}$, there exists some $x_j \in \bar{\gamma}_{j-{{s}\choose{\lceil{s/2}\rceil}}+|A_2|}\cap \bar{\psi}$. So, $(1,[\beta,j])\rightarrow (2,[j])$, $(2,[\beta,j])\rightarrow (3,[j])$ and $\{(2,[j]),(3,[j])\}\rightarrow x_j\rightarrow (1,[i])$ by (\ref{eqC6.4.31}), (\ref{eqC6.4.33})-(\ref{eqC6.4.34}) and (\ref{eqC6.4.36}).
\indent\par Now, observe that $|\gamma_k\cap \psi|=0$ if and only if $\gamma_k=\bar{\psi}$ if and only if $k={{s}\choose{\lceil{s/2}\rceil}}$. And, by definition of $I_{\psi}$, $\gamma_{j-{{s}\choose{\lceil{s/2}\rceil}}+|A_2|}\not\in I_{\psi}$ so that $|\gamma_{j-{{s}\choose{\lceil{s/2}\rceil}}+|A_2|}\cap\psi|\ge 2$ for all ${{s}\choose{\lceil{s/2}\rceil}}+1 \le j \le |A_2|+|A_3|$. As in Case 3.6(iv), we deduce that there exists some $y_{ij}\in (\bar{\gamma}_{j-{{s}\choose{\lceil{s/2}\rceil}}+|A_2|}\cup\psi)\cap\bar{\mu}_i$. Hence, $(1,[\beta,j])\rightarrow \{(1,[j]),(2,[j])\}$, $(2,[\beta,j])\rightarrow \{(1,[j]),(3,[j])\}$, and either $\{(2,[j]),(3,[j])\}\rightarrow y_{ij}\rightarrow \{(2,[i]),(3,[i])\}$ or $(1,[j])\rightarrow y_{ij}\rightarrow \{(2,[i]),(3,[i])\}$ by (\ref{eqC6.4.32})-(\ref{eqC6.4.36}).
\\
\\(v) $v=(p,[\alpha,i]), w=(q,[\beta,j])$ for each $p,q=1,2$.
\indent\par From (i), $d_{D_3}((p,[\alpha,i]),(r,[j]))=3$ for $r=2,3$. Since  $(3,[j])\rightarrow (1,[\beta,j])$ and $(2,[j])\rightarrow (2,[\beta,j])$ by (\ref{eqC6.4.33})-(\ref{eqC6.4.34}), this subcase follows.
\\
\\(vi) $v=(q,[\beta,j]), w=(p,[\alpha,i])$ for each $p,q=1,2$.
\indent\par From (iv), $d_{D_3}((q,[\beta,j]),(1,[i]))=3$. Since $(1,[i])\rightarrow \{(1,[\alpha,i]), (2,[\alpha,i])\}$ by (\ref{eqC6.4.29})-(\ref{eqC6.4.30}), this subcase follows.
\\
\\Case 3.8. For each $1\le i\le |A_2|$, each $1\le\alpha\le deg_T([i])-1$, and each $[j]\in E$,
\\(i) $v=(p,[\alpha,i]), w=(q,[j])$ for each $p,q=1,2$.
\\(ii) $v=(q,[j]), w=(p,[\alpha,i])$ for each $p,q=1,2$.
\indent\par Since $\bar{\psi}\rightarrow \{(q,[j]),(1,[|A_2|+1])\}\rightarrow \psi$ by (\ref{eqC6.4.31}) and (\ref{eqC6.4.37}), this case follows from Cases 3.5(i)-(ii).
\\
\\Case 3.9. For each $|A_2|+1\le i\le {{s}\choose{\lceil{s/2}\rceil}}$, each $1\le\alpha\le deg_T([i])-1$, and each $[j]\in E$,
\\(i) $v=(p,[\alpha,i]), w=(q,[j])$ for each $p,q=1,2$.
\\(ii) $v=(q,[j]), w=(p,[\alpha,i])$ for each $p,q=1,2$.
\indent\par Since $\bar{\psi}\rightarrow \{(q,[j]),(1,[{{s}\choose{\lceil{s/2}\rceil}}+1])\}\rightarrow \psi$ by (\ref{eqC6.4.35}) and (\ref{eqC6.4.37}), this case follows from Cases 3.7(i)-(ii).
\\
\\Case 3.10. For each ${{s}\choose{\lceil{s/2}\rceil}}+1\le i\le|A_2|+|A_3|$, each $1\le\alpha\le deg_T([i])-1$, and each $[j]\in E$,
\\(i) $v=(p,[\alpha,i]), w=(q,[j])$ for each $p,q=1,2$.
\\(ii) $v=(q,[j]), w=(p,[\alpha,i])$ for each $p,q=1,2$.
\indent\par Since $\bar{\psi}\rightarrow \{(q,[j]),(1,[|A_2|+1])\}\rightarrow \psi$ by (\ref{eqC6.4.31}) and (\ref{eqC6.4.37}), this case follows from Cases 3.7(iii)-(iv).
\\
\\Case 3.11. $v=(r_1,\mathtt{c})$ and $w=(r_2,\mathtt{c})$ for $r_1\neq r_2$ and $1\le r_1, r_2\le s$.
\indent\par As before in $D_2$, we want to prove $d_{D_3}((r_1,\mathtt{c}), (r_2,\mathtt{c}))=2$. Let $x_i=(1,[i])$ for $1\le i\le |A_2|$ and $x_i=(2,[i])$ for $i=|A_2|+1\le i\le {{s}\choose{\lceil{s/2}\rceil}}$. Observe that $\bar{\mu}_i\rightarrow x_i\rightarrow \mu_i$ for all $1\le i\le{{s}\choose{\lceil{s/2}\rceil}}$. Since ${{(\mathbb{N}_s,\mathtt{c})}\choose{\lceil{s/2}\rceil}}=\{\mu_i\mid 1\le i\le{{s}\choose{\lceil{s/2}\rceil}}\}$, for every $r_1\neq r_2$, there exists some $1\le i\le{{s}\choose{\lceil{s/2}\rceil}}$ such that $(r_1,\mathtt{c})\not\in\mu_i$, $(r_2,\mathtt{c})\in\mu_i$ and $(r_1,\mathtt{c})\rightarrow x_i\rightarrow (r_2,\mathtt{c})$.
\\
\\Case 3.12. $v\in \{(1,[i]), (2,[i]), (3,[i]), (1,[\alpha,i]), (2,[\alpha,i])\}$ for each $1\le i\le deg_T(\mathtt{c})$ and $1\le\alpha\le deg_T([i])-1$, and $w=(r,\mathtt{c})$ for $1\le r\le s$.
\indent\par Note that there exists some $1\le k\le s$ such that $d_{D_3}(v,(k,\mathtt{c}))\le 2$, and $d_{D_3}((k,\mathtt{c}),w)\le 2$ by Case 3.11. Hence, it follows that $d_{D_3}(v,w)\le d_{D_3}(v,(k,\mathtt{c}))+d_{D_3}((k,\mathtt{c}),w)\le 4$.
\\
\\Case 3.13. $v=(r,\mathtt{c})$ for $r=1,2,\ldots, s$ and $w\in \{(1,[i]), (2,[i]), (3,[i]), (1,[\alpha,i]), (2,[\alpha,i])\}$ for each $1\le i\le deg_T(\mathtt{c})$ and $1\le\alpha\le deg_T([i])-1$.
\indent\par Note that there exists some $1\le k\le s$ such that $d_{D_3}((k,\mathtt{c}), w)\le 2$, and $d_{D_3}(v,(k,\mathtt{c}))\le 2$ by Case 3.11. Hence, it follows that $d_{D_3}(v,w)\le d_{D_3}(v,(k,\mathtt{c}))+d_{D_3}((k,\mathtt{c}),w)\le 4$.
\\
\\Case 3.14. $v=(p,[i])$ and $w=(q, [j])$, where $1\le p,q\le 3$ and $1\le i,j\le deg_T(\mathtt{c})$.
\noindent\par This follows from the fact that $|O^\mathtt{c}((p,[i]))|>0$, $|I^\mathtt{c}((q,[j]))|>0$, and $d_{D_3}((r_1,\mathtt{c}), (r_2,\mathtt{c}))$ $=2$ for any $r_1\neq r_2$ and $1\le r_1, r_2\le s$.

\indent\par Hence, $d(D_i)=4$ for $i=1,2,3$. Since every vertex in $D_i$ lies in a directed $C_4$, i.e., $\bar{d}(\mathcal{T})\le \max \{4, d(D_i)\}$ by Lemma \ref{lemC6.1.3}, and thus $\bar{d}(\mathcal{T})=4$.
\qed
\begin{center}
\begin{tikzpicture}[thick,scale=0.75]%
\draw(0,2)node[circle, draw, inner sep=0pt, minimum width=3pt](1_u){\scriptsize $(1,\mathtt{c})$};
\draw(0,-1)node[circle, draw, inner sep=0pt, minimum width=3pt](2_u){\scriptsize $(2,\mathtt{c})$};
\draw(0,-4)node[circle, draw, inner sep=0pt, minimum width=3pt](3_u){\scriptsize $(3,\mathtt{c})$};

\draw(6,2)node[circle, draw, fill=black!100, inner sep=0pt, minimum width=6pt, label={[] 0:{\small $(1,[1,1])$}}](1_11){};
\draw(6,0)node[circle, draw, fill=black!100, inner sep=0pt, minimum width=6pt, label={[] 0:{\small $(2,[1,1])$}}](2_11){};

\draw(3,3)node[circle, draw, fill=black!100, inner sep=0pt, minimum width=6pt, label={[yshift=0cm, xshift=0.45cm] 90:{\small $(1,[1])$}}](1_1){};
\draw(3,1)node[circle, draw, fill=black!100, inner sep=0pt, minimum width=6pt, label={[yshift=-0.1cm, xshift=0.45cm] 270:{\small $(2,[1])$}}](2_1){};
\draw(3,-1)node[circle, draw, fill=black!100, inner sep=0pt, minimum width=6pt, label={[yshift=0cm, xshift=0.45cm] 270:{\small $(3,[1])$}}](3_1){};

\draw(-6,1)node[circle, draw, fill=black!100, inner sep=0pt, minimum width=6pt, label={[] 180:{\small $(1,[1,2])$}}](1_12){};
\draw(-6,-1)node[circle, draw, fill=black!100, inner sep=0pt, minimum width=6pt, label={[] 180:{\small $(2,[1,2])$}}](2_12){};

\draw(-3,1)node[circle, draw, fill=black!100, inner sep=0pt, minimum width=6pt, label={[yshift=0cm, xshift=-0.45cm] 90:{\small $(1,[2])$}}](1_2){};
\draw(-3,-1)node[circle, draw, fill=black!100, inner sep=0pt, minimum width=6pt, label={[yshift=0cm, xshift=-0.45cm] 270:{\small $(2,[2])$}}](2_2){};

\draw(-3,-3)node[circle, draw, fill=black!100, inner sep=0pt, minimum width=6pt, label={[yshift=0cm] 180:{\small $(1,[3])$}}](1_3){};
\draw(-3,-5)node[circle, draw, fill=black!100, inner sep=0pt, minimum width=6pt, label={[yshift=0cm] 180:{\small $(2,[3])$}}](2_3){};

\draw(3,-3)node[circle, draw, fill=black!100, inner sep=0pt, minimum width=6pt, label={[yshift=0cm] 0:{\small $(1,[4])$}}](1_4){};
\draw(3,-5)node[circle, draw, fill=black!100, inner sep=0pt, minimum width=6pt, label={[yshift=0cm] 0:{\small $(2,[4])$}}](2_4){};

\draw[->, line width=0.3mm, >=latex, shorten <= 0.2cm, shorten >= 0.15cm](1_12)--(1_2);
\draw[->, line width=0.3mm, >=latex, shorten <= 0.2cm, shorten >= 0.15cm](2_12)--(1_2);
\draw[dashed,->, line width=0.3mm, >=latex, shorten <= 0.2cm, shorten >= 0.15cm](2_2)--(1_12);
\draw[dashed,->, line width=0.3mm, >=latex, shorten <= 0.2cm, shorten >= 0.15cm](2_2)--(2_12);

\draw[dashed,->, line width=0.3mm, >=latex, shorten <= 0.2cm, shorten >= 0.15cm](1_1)--(1_11);
\draw[dashed,->, line width=0.3mm, >=latex, shorten <= 0.2cm, shorten >= 0.15cm](1_1)--(2_11);
\draw[dashed,->, line width=0.3mm, >=latex, shorten <= 0.2cm, shorten >= 0.15cm](2_1)--(1_11);
\draw[dashed,->, line width=0.3mm, >=latex, shorten <= 0.2cm, shorten >= 0.15cm](3_1)--(2_11);
\draw[->, line width=0.3mm, >=latex, shorten <= 0.2cm, shorten >= 0.15cm](1_11) to [out=265, in=22] (3_1);
\draw[->, line width=0.3mm, >=latex, shorten <= 0.2cm, shorten >= 0.15cm](2_11)--(2_1);

\draw[densely dotted, ->, line width=0.3mm, >=latex, shorten <= 0.2cm, shorten >= 0.1cm](1_2)--(2_u);
\draw[densely dotted, ->, line width=0.3mm, >=latex, shorten <= 0.2cm, shorten >= 0.1cm](1_2)--(3_u);
\draw[densely dashdotdotted, ->, line width=0.3mm, >=latex, shorten <= 0.2cm, shorten >= 0.1cm](2_2)--(1_u);

\draw[->, line width=0.3mm, >=latex, shorten <= 0.2cm, shorten >= 0.1cm](1_1)--(3_u);
\draw[->, line width=0.3mm, >=latex, shorten <= 0.2cm, shorten >= 0.1cm](2_1)--(1_u);
\draw[->, line width=0.3mm, >=latex, shorten <= 0.2cm, shorten >= 0.1cm](2_1)--(2_u);
\draw[->, line width=0.3mm, >=latex, shorten <= 0.2cm, shorten >= 0.1cm](3_1)--(1_u);
\draw[->, line width=0.3mm, >=latex, shorten <= 0.2cm, shorten >= 0.1cm](3_1)--(2_u);
\draw[->, line width=0.3mm, >=latex, shorten <= 0.2cm, shorten >= 0.1cm](1_3)--(3_u);
\draw[->, line width=0.3mm, >=latex, shorten <= 0.2cm, shorten >= 0.1cm](2_3)--(3_u);
\draw[->, line width=0.3mm, >=latex, shorten <= 0.2cm, shorten >= 0.1cm](1_4)--(3_u);
\draw[->, line width=0.3mm, >=latex, shorten <= 0.2cm, shorten >= 0.1cm](2_4)--(3_u);
\end{tikzpicture}
\captionsetup{justification=centering}
{\captionof{figure}{Orientation $D_1$ for $\mathcal{H}$ for $s=3$, $A_2=\{[2]\}$, $A_3=\{[1]\}$, $E=\{[3],[4]\}$.}\label{figC6.4.13}}
\end{center}
\begin{center}
\begin{tikzpicture}[thick,scale=0.75]%
\draw(0,2)node[circle, draw, inner sep=0pt, minimum width=3pt](1_u){\scriptsize $(1,\mathtt{c})$};
\draw(0,-2)node[circle, draw, inner sep=0pt, minimum width=3pt](2_u){\scriptsize $(2,\mathtt{c})$};
\draw(0,-6)node[circle, draw, inner sep=0pt, minimum width=3pt](3_u){\scriptsize $(3,\mathtt{c})$};

\draw(6,3)node[circle, draw, fill=black!100, inner sep=0pt, minimum width=6pt, label={[] 0:{\small $(1,[1,1])$}}](1_11){};
\draw(6,1)node[circle, draw, fill=black!100, inner sep=0pt, minimum width=6pt, label={[] 0:{\small $(2,[1,1])$}}](2_11){};

\draw(3,3)node[circle, draw, fill=black!100, inner sep=0pt, minimum width=6pt, label={[yshift=0cm, xshift=0.45cm] 90:{\small $(1,[1])$}}](1_1){};
\draw(3,1)node[circle, draw, fill=black!100, inner sep=0pt, minimum width=6pt, label={[yshift=0cm, xshift=0.45cm] 270:{\small $(2,[1])$}}](2_1){};

\draw(-6,2)node[circle, draw, fill=black!100, inner sep=0pt, minimum width=6pt, label={[] 180:{\small $(1,[1,2])$}}](1_12){};
\draw(-6,0)node[circle, draw, fill=black!100, inner sep=0pt, minimum width=6pt, label={[] 180:{\small $(2,[1,2])$}}](2_12){};

\draw(-3,3)node[circle, draw, fill=black!100, inner sep=0pt, minimum width=6pt, label={[yshift=0cm, xshift=-0.45cm] 90:{\small $(1,[2])$}}](1_2){};
\draw(-3,1)node[circle, draw, fill=black!100, inner sep=0pt, minimum width=6pt, label={[yshift=-0.1cm, xshift=-0.45cm] 270:{\small $(2,[2])$}}](2_2){};
\draw(-3,-1)node[circle, draw, fill=black!100, inner sep=0pt, minimum width=6pt, label={[yshift=0cm, xshift=-0.45cm] 270:{\small $(3,[2])$}}](3_2){};

\draw(-6,-4)node[circle, draw, fill=black!100, inner sep=0pt, minimum width=6pt, label={[] 180:{\small $(1,[1,3])$}}](1_13){};
\draw(-6,-6)node[circle, draw, fill=black!100, inner sep=0pt, minimum width=6pt, label={[] 180:{\small $(2,[1,3])$}}](2_13){};

\draw(-3,-3)node[circle, draw, fill=black!100, inner sep=0pt, minimum width=6pt, label={[yshift=-0.15cm, xshift=-0.45cm] 270:{\small $(1,[3])$}}](1_3){};
\draw(-3,-5)node[circle, draw, fill=black!100, inner sep=0pt, minimum width=6pt, label={[yshift=-0.1cm, xshift=-0.45cm] 270:{\small $(2,[3])$}}](2_3){};
\draw(-3,-7)node[circle, draw, fill=black!100, inner sep=0pt, minimum width=6pt, label={[yshift=0cm, xshift=-0.45cm] 270:{\small $(3,[3])$}}](3_3){};

\draw(3,-1)node[circle, draw, fill=black!100, inner sep=0pt, minimum width=6pt, label={[] 0:{\small $(1,[4])$}}](1_4){};
\draw(3,-3)node[circle, draw, fill=black!100, inner sep=0pt, minimum width=6pt, label={[] 0:{\small $(2,[4])$}}](2_4){};

\draw(3,-5)node[circle, draw, fill=black!100, inner sep=0pt, minimum width=6pt, label={[] 0:{\small $(1,[5])$}}](1_5){};
\draw(3,-7)node[circle, draw, fill=black!100, inner sep=0pt, minimum width=6pt, label={[] 0:{\small $(2,[5])$}}](2_5){};

\draw[->, line width=0.3mm, >=latex, shorten <= 0.2cm, shorten >= 0.15cm](1_11)--(1_1);
\draw[->, line width=0.3mm, >=latex, shorten <= 0.2cm, shorten >= 0.15cm](2_11)--(1_1);
\draw[dashed,->, line width=0.3mm, >=latex, shorten <= 0.2cm, shorten >= 0.15cm](2_1)--(1_11);
\draw[dashed,->, line width=0.3mm, >=latex, shorten <= 0.2cm, shorten >= 0.15cm](2_1)--(2_11);

\draw[dashed,->, line width=0.3mm, >=latex, shorten <= 0.2cm, shorten >= 0.15cm](1_3)--(1_13);
\draw[dashed,->, line width=0.3mm, >=latex, shorten <= 0.2cm, shorten >= 0.15cm](1_3) to [out=202, in=85] (2_13);
\draw[dashed,->, line width=0.3mm, >=latex, shorten <= 0.2cm, shorten >= 0.15cm](2_3)--(1_13);
\draw[dashed,->, line width=0.3mm, >=latex, shorten <= 0.2cm, shorten >= 0.15cm](2_3)--(2_13);
\draw[->, line width=0.3mm, >=latex, shorten <= 0.2cm, shorten >= 0.15cm](1_13) to [out=275, in=158] (3_3);
\draw[->, line width=0.3mm, >=latex, shorten <= 0.2cm, shorten >= 0.15cm](2_13)--(3_3);

\draw[dashed,->, line width=0.3mm, >=latex, shorten <= 0.2cm, shorten >= 0.15cm](3_2) to [out=158, in=275] (1_12);
\draw[dashed,->, line width=0.3mm, >=latex, shorten <= 0.2cm, shorten >= 0.15cm](3_2)--(2_12);
\draw[->, line width=0.3mm, >=latex, shorten <= 0.2cm, shorten >= 0.15cm](1_12)--(1_2);
\draw[->, line width=0.3mm, >=latex, shorten <= 0.2cm, shorten >= 0.15cm](2_12)--(1_2);
\draw[->, line width=0.3mm, >=latex, shorten <= 0.2cm, shorten >= 0.15cm](1_12)--(2_2);
\draw[->, line width=0.3mm, >=latex, shorten <= 0.2cm, shorten >= 0.15cm](2_12)--(2_2);

\draw[densely dashdotdotted, ->, line width=0.3mm, >=latex, shorten <= 0.2cm, shorten >= 0.1cm](2_1)--(1_u);
\draw[densely dotted, ->, line width=0.3mm, >=latex, shorten <= 0.2cm, shorten >= 0.1cm](1_1)--(2_u);
\draw[densely dotted, ->, line width=0.3mm, >=latex, shorten <= 0.2cm, shorten >= 0.1cm](1_1)--(3_u);

\draw[->, line width=0.3mm, >=latex, shorten <= 0.2cm, shorten >= 0.1cm](1_2)--(1_u);
\draw[->, line width=0.3mm, >=latex, shorten <= 0.2cm, shorten >= 0.1cm](1_2)--(2_u);
\draw[->, line width=0.3mm, >=latex, shorten <= 0.2cm, shorten >= 0.1cm](2_2)--(3_u);
\draw[densely dashdotdotted, ->, line width=0.3mm, >=latex, shorten <= 0.2cm, shorten >= 0.1cm](3_2)--(2_u);

\draw[->, line width=0.3mm, >=latex, shorten <= 0.2cm, shorten >= 0.1cm](1_3)--(3_u);
\draw[->, line width=0.3mm, >=latex, shorten <= 0.2cm, shorten >= 0.1cm](2_3)--(1_u);
\draw[->, line width=0.3mm, >=latex, shorten <= 0.2cm, shorten >= 0.1cm](2_3)--(2_u);
\draw[densely dotted, ->, line width=0.3mm, >=latex, shorten <= 0.2cm, shorten >= 0.1cm](3_3)--(1_u);
\draw[densely dotted, ->, line width=0.3mm, >=latex, shorten <= 0.2cm, shorten >= 0.1cm](3_3)--(3_u);

\draw[->, line width=0.3mm, >=latex, shorten <= 0.2cm, shorten >= 0.1cm](1_4)--(3_u);
\draw[->, line width=0.3mm, >=latex, shorten <= 0.2cm, shorten >= 0.1cm](2_4)--(3_u);
\draw[->, line width=0.3mm, >=latex, shorten <= 0.2cm, shorten >= 0.1cm](1_5)--(3_u);
\draw[->, line width=0.3mm, >=latex, shorten <= 0.2cm, shorten >= 0.1cm](2_5)--(3_u);
\end{tikzpicture}
\captionsetup{justification=centering}
{\captionof{figure}{Orientation $D_2$ for $\mathcal{H}$ for $s=3$, $A_2=\{[1]\}$, $A_3=\{[2],[3]\}$, $E=\{[4],[5]\}$.}\label{figC6.4.14}}
\end{center}

\begin{center}
\begin{tikzpicture}[thick,scale=0.75]%
\draw(0,3)node[circle, draw, inner sep=0pt, minimum width=3pt](1_u){\scriptsize $(1,\mathtt{c})$};
\draw(0,0)node[circle, draw, inner sep=0pt, minimum width=3pt](2_u){\scriptsize $(2,\mathtt{c})$};
\draw(0,-3)node[circle, draw, inner sep=0pt, minimum width=3pt](3_u){\scriptsize $(3,\mathtt{c})$};
\draw(0,-6)node[circle, draw, inner sep=0pt, minimum width=3pt](4_u){\scriptsize $(4,\mathtt{c})$};
\draw(0,-9)node[circle, draw, inner sep=0pt, minimum width=3pt](5_u){\scriptsize $(5,\mathtt{c})$};

\draw(-6,5)node[circle, draw, fill=black!100, inner sep=0pt, minimum width=6pt, label={[] 180:{\small $(1,[1,1])$}}](1_11){};
\draw(-6,3)node[circle, draw, fill=black!100, inner sep=0pt, minimum width=6pt, label={[] 180:{\small $(2,[1,1])$}}](2_11){};

\draw(-3,5)node[circle, draw, fill=black!100, inner sep=0pt, minimum width=6pt, label={[yshift=0cm, xshift=-0.45cm] 90:{\small $(1,[1])$}}](1_1){};
\draw(-3,3)node[circle, draw, fill=black!100, inner sep=0pt, minimum width=6pt, label={[yshift=0cm, xshift=-0.45cm] 270:{\small $(2,[1])$}}](2_1){};

\draw(6,5)node[circle, draw, fill=black!100, inner sep=0pt, minimum width=6pt, label={[] 0:{\small $(1,[1,2])$}}](1_21){};
\draw(6,3)node[circle, draw, fill=black!100, inner sep=0pt, minimum width=6pt, label={[] 0:{\small $(2,[1,2])$}}](2_21){};

\draw(3,5)node[circle, draw, fill=black!100, inner sep=0pt, minimum width=6pt, label={[yshift=0cm, xshift=0.45cm] 90:{\small $(1,[2])$}}](1_2){};
\draw(3,3)node[circle, draw, fill=black!100, inner sep=0pt, minimum width=6pt, label={[yshift=0cm, xshift=0.45cm] 270:{\small $(2,[2])$}}](2_2){};

\draw(-6,1)node[circle, draw, fill=black!100, inner sep=0pt, minimum width=6pt, label={[] 180:{\small $(1,[1,3])$}}](1_31){};
\draw(-6,-1)node[circle, draw, fill=black!100, inner sep=0pt, minimum width=6pt, label={[] 180:{\small $(2,[1,3])$}}](2_31){};

\draw(-3,1)node[circle, draw, fill=black!100, inner sep=0pt, minimum width=6pt, label={[yshift=0cm, xshift=-0.45cm] 270:{\small $(1,[3])$}}](1_3){};
\draw(-3,-1)node[circle, draw, fill=black!100, inner sep=0pt, minimum width=6pt, label={[yshift=0cm, xshift=-0.45cm] 270:{\small $(2,[3])$}}](2_3){};

\draw(6,1)node[circle, draw, fill=black!100, inner sep=0pt, minimum width=6pt, label={[] 0:{\small $(1,[1,4])$}}](1_41){};
\draw(6,-1)node[circle, draw, fill=black!100, inner sep=0pt, minimum width=6pt, label={[] 0:{\small $(2,[1,4])$}}](2_41){};

\draw(3,1)node[circle, draw, fill=black!100, inner sep=0pt, minimum width=6pt, label={[yshift=0cm, xshift=0.45cm] 270:{\small $(1,[4])$}}](1_4){};
\draw(3,-1)node[circle, draw, fill=black!100, inner sep=0pt, minimum width=6pt, label={[yshift=0cm, xshift=0.45cm] 270:{\small $(2,[4])$}}](2_4){};

\draw(-6,-3)node[circle, draw, fill=black!100, inner sep=0pt, minimum width=6pt, label={[] 180:{\small $(1,[1,5])$}}](1_51){};
\draw(-6,-5)node[circle, draw, fill=black!100, inner sep=0pt, minimum width=6pt, label={[] 180:{\small $(2,[1,5])$}}](2_51){};

\draw(-3,-3)node[circle, draw, fill=black!100, inner sep=0pt, minimum width=6pt, label={[yshift=0cm, xshift=-0.45cm] 270:{\small $(1,[5])$}}](1_5){};
\draw(-3,-5)node[circle, draw, fill=black!100, inner sep=0pt, minimum width=6pt, label={[yshift=0cm, xshift=-0.45cm] 270:{\small $(2,[5])$}}](2_5){};

\draw(6,-3)node[circle, draw, fill=black!100, inner sep=0pt, minimum width=6pt, label={[] 0:{\small $(1,[1,6])$}}](1_61){};
\draw(6,-5)node[circle, draw, fill=black!100, inner sep=0pt, minimum width=6pt, label={[] 0:{\small $(2,[1,6])$}}](2_61){};

\draw(3,-3)node[circle, draw, fill=black!100, inner sep=0pt, minimum width=6pt, label={[yshift=0cm, xshift=0.45cm] 270:{\small $(1,[6])$}}](1_6){};
\draw(3,-5)node[circle, draw, fill=black!100, inner sep=0pt, minimum width=6pt, label={[yshift=0cm, xshift=0.45cm] 270:{\small $(2,[6])$}}](2_6){};

\draw(-6,-8)node[circle, draw, fill=black!100, inner sep=0pt, minimum width=6pt, label={[] 180:{\small $(1,[1,7])$}}](1_71){};
\draw(-6,-10)node[circle, draw, fill=black!100, inner sep=0pt, minimum width=6pt, label={[] 180:{\small $(2,[1,7])$}}](2_71){};

\draw(-3,-7)node[circle, draw, fill=black!100, inner sep=0pt, minimum width=6pt, label={[yshift=-0.15cm, xshift=-0.45cm] 270:{\small $(1,[7])$}}](1_7){};
\draw(-3,-9)node[circle, draw, fill=black!100, inner sep=0pt, minimum width=6pt, label={[yshift=-0.1cm, xshift=-0.45cm] 270:{\small $(2,[7])$}}](2_7){};
\draw(-3,-11)node[circle, draw, fill=black!100, inner sep=0pt, minimum width=6pt, label={[yshift=0cm, xshift=-0.45cm] 270:{\small $(3,[7])$}}](3_7){};

\draw(6,-8)node[circle, draw, fill=black!100, inner sep=0pt, minimum width=6pt, label={[] 0:{\small $(1,[1,8])$}}](1_81){};
\draw(6,-10)node[circle, draw, fill=black!100, inner sep=0pt, minimum width=6pt, label={[] 0:{\small $(2,[1,8])$}}](2_81){};

\draw(3,-7)node[circle, draw, fill=black!100, inner sep=0pt, minimum width=6pt, label={[yshift=-0.15cm, xshift=0.45cm] 270:{\small $(1,[8])$}}](1_8){};
\draw(3,-9)node[circle, draw, fill=black!100, inner sep=0pt, minimum width=6pt, label={[yshift=-0.1cm, xshift=0.45cm] 270:{\small $(2,[8])$}}](2_8){};
\draw(3,-11)node[circle, draw, fill=black!100, inner sep=0pt, minimum width=6pt, label={[yshift=0cm, xshift=0.45cm] 270:{\small $(3,[8])$}}](3_8){};

\draw[->, line width=0.3mm, >=latex, shorten <= 0.2cm, shorten >= 0.15cm](1_11)--(1_1);
\draw[->, line width=0.3mm, >=latex, shorten <= 0.2cm, shorten >= 0.15cm](2_11)--(1_1);
\draw[dashed,->, line width=0.3mm, >=latex, shorten <= 0.2cm, shorten >= 0.15cm](2_1)--(1_11);
\draw[dashed,->, line width=0.3mm, >=latex, shorten <= 0.2cm, shorten >= 0.15cm](2_1)--(2_11);

\draw[->, line width=0.3mm, >=latex, shorten <= 0.2cm, shorten >= 0.15cm](1_21)--(1_2);
\draw[->, line width=0.3mm, >=latex, shorten <= 0.2cm, shorten >= 0.15cm](2_21)--(1_2);
\draw[dashed,->, line width=0.3mm, >=latex, shorten <= 0.2cm, shorten >= 0.15cm](2_2)--(1_21);
\draw[dashed,->, line width=0.3mm, >=latex, shorten <= 0.2cm, shorten >= 0.15cm](2_2)--(2_21);

\draw[->, line width=0.3mm, >=latex, shorten <= 0.2cm, shorten >= 0.15cm](1_31)--(1_3);
\draw[->, line width=0.3mm, >=latex, shorten <= 0.2cm, shorten >= 0.15cm](2_31) to [out=65, in=190] (1_3);
\draw[dashed,->, line width=0.3mm, >=latex, shorten <= 0.2cm, shorten >= 0.15cm](2_3)--(1_31);
\draw[dashed,->, line width=0.3mm, >=latex, shorten <= 0.2cm, shorten >= 0.15cm](2_3)--(2_31);

\draw[->, line width=0.3mm, >=latex, shorten <= 0.2cm, shorten >= 0.15cm](1_41)--(1_4);
\draw[->, line width=0.3mm, >=latex, shorten <= 0.2cm, shorten >= 0.15cm](2_41) to [out=115, in=350] (1_4);
\draw[dashed,->, line width=0.3mm, >=latex, shorten <= 0.2cm, shorten >= 0.15cm](2_4)--(1_41);
\draw[dashed,->, line width=0.3mm, >=latex, shorten <= 0.2cm, shorten >= 0.15cm](2_4)--(2_41);

\draw[->, line width=0.3mm, >=latex, shorten <= 0.2cm, shorten >= 0.15cm](1_51)--(1_5);
\draw[->, line width=0.3mm, >=latex, shorten <= 0.2cm, shorten >= 0.15cm](2_51) to [out=65, in=190] (1_5);
\draw[dashed,->, line width=0.3mm, >=latex, shorten <= 0.2cm, shorten >= 0.15cm](2_5)--(1_51);
\draw[dashed,->, line width=0.3mm, >=latex, shorten <= 0.2cm, shorten >= 0.15cm](2_5)--(2_51);

\draw[->, line width=0.3mm, >=latex, shorten <= 0.2cm, shorten >= 0.15cm](1_61)--(1_6);
\draw[->, line width=0.3mm, >=latex, shorten <= 0.2cm, shorten >= 0.15cm](2_61) to [out=115, in=350] (1_6);
\draw[dashed,->, line width=0.3mm, >=latex, shorten <= 0.2cm, shorten >= 0.15cm](2_6)--(1_61);
\draw[dashed,->, line width=0.3mm, >=latex, shorten <= 0.2cm, shorten >= 0.15cm](2_6)--(2_61);

\draw[->, line width=0.3mm, >=latex, shorten <= 0.2cm, shorten >= 0.15cm](1_71) to [out=275, in=158] (3_7);
\draw[->, line width=0.3mm, >=latex, shorten <= 0.2cm, shorten >= 0.15cm](2_71)--(2_7);
\draw[dashed,->, line width=0.3mm, >=latex, shorten <= 0.2cm, shorten >= 0.15cm](1_7)--(1_71);
\draw[dashed,->, line width=0.3mm, >=latex, shorten <= 0.2cm, shorten >= 0.15cm](2_7)--(1_71);
\draw[dashed,->, line width=0.3mm, >=latex, shorten <= 0.2cm, shorten >= 0.15cm](1_7) to [out=202, in=85] (2_71);
\draw[dashed,->, line width=0.3mm, >=latex, shorten <= 0.2cm, shorten >= 0.15cm](3_7)--(2_71);

\draw[dashed,->, line width=0.3mm, >=latex, shorten <= 0.2cm, shorten >= 0.15cm](1_8)--(1_81);
\draw[dashed,->, line width=0.3mm, >=latex, shorten <= 0.2cm, shorten >= 0.15cm](1_8) to [out=338, in=95] (2_81);
\draw[dashed,->, line width=0.3mm, >=latex, shorten <= 0.2cm, shorten >= 0.15cm](2_8)--(1_81);
\draw[->, line width=0.3mm, >=latex, shorten <= 0.2cm, shorten >= 0.15cm](1_81) to [out=265, in=22] (3_8);
\draw[dashed,->, line width=0.3mm, >=latex, shorten <= 0.2cm, shorten >= 0.15cm](3_8)--(2_81);
\draw[->, line width=0.3mm, >=latex, shorten <= 0.2cm, shorten >= 0.15cm](2_81)--(2_8);
\draw[densely dashdotdotted, ->, line width=0.3mm, >=latex, shorten <= 0.2cm, shorten >= 0.1cm](2_1)--(2_u);
\draw[densely dashdotdotted, ->, line width=0.3mm, >=latex, shorten <= 0.2cm, shorten >= 0.1cm](2_1)--(5_u);
\draw[densely dotted, ->, line width=0.3mm, >=latex, shorten <= 0.2cm, shorten >= 0.1cm](1_1)--(1_u);
\draw[densely dotted, ->, line width=0.3mm, >=latex, shorten <= 0.2cm, shorten >= 0.1cm](1_1)--(2_u);
\draw[densely dotted, ->, line width=0.3mm, >=latex, shorten <= 0.2cm, shorten >= 0.1cm](1_1)--(3_u);

\draw[densely dashdotdotted, ->, line width=0.3mm, >=latex, shorten <= 0.2cm, shorten >= 0.1cm](2_2)--(1_u);
\draw[densely dashdotdotted, ->, line width=0.3mm, >=latex, shorten <= 0.2cm, shorten >= 0.1cm](2_2)--(5_u);
\draw[densely dotted, ->, line width=0.3mm, >=latex, shorten <= 0.2cm, shorten >= 0.1cm](1_2)--(1_u);
\draw[densely dotted, ->, line width=0.3mm, >=latex, shorten <= 0.2cm, shorten >= 0.1cm](1_2)--(2_u);
\draw[densely dotted, ->, line width=0.3mm, >=latex, shorten <= 0.2cm, shorten >= 0.1cm](1_2)--(4_u);

\draw[densely dashdotdotted, ->, line width=0.3mm, >=latex, shorten <= 0.2cm, shorten >= 0.1cm](2_3)--(2_u);
\draw[densely dashdotdotted, ->, line width=0.3mm, >=latex, shorten <= 0.2cm, shorten >= 0.1cm](2_3)--(4_u);
\draw[densely dotted, ->, line width=0.3mm, >=latex, shorten <= 0.2cm, shorten >= 0.1cm](1_3)--(1_u);
\draw[densely dotted, ->, line width=0.3mm, >=latex, shorten <= 0.2cm, shorten >= 0.1cm](1_3)--(2_u);
\draw[densely dotted, ->, line width=0.3mm, >=latex, shorten <= 0.2cm, shorten >= 0.1cm](1_3)--(5_u);

\draw[densely dashdotdotted, ->, line width=0.3mm, >=latex, shorten <= 0.2cm, shorten >= 0.1cm](2_4)--(1_u);
\draw[densely dashdotdotted, ->, line width=0.3mm, >=latex, shorten <= 0.2cm, shorten >= 0.1cm](2_4)--(4_u);
\draw[densely dotted, ->, line width=0.3mm, >=latex, shorten <= 0.2cm, shorten >= 0.1cm](1_4)--(1_u);
\draw[densely dotted, ->, line width=0.3mm, >=latex, shorten <= 0.2cm, shorten >= 0.1cm](1_4)--(3_u);
\draw[densely dotted, ->, line width=0.3mm, >=latex, shorten <= 0.2cm, shorten >= 0.1cm](1_4)--(4_u);

\draw[densely dashdotdotted, ->, line width=0.3mm, >=latex, shorten <= 0.2cm, shorten >= 0.1cm](2_5)--(2_u);
\draw[densely dashdotdotted, ->, line width=0.3mm, >=latex, shorten <= 0.2cm, shorten >= 0.1cm](2_5)--(3_u);
\draw[densely dotted, ->, line width=0.3mm, >=latex, shorten <= 0.2cm, shorten >= 0.1cm](1_5)--(2_u);
\draw[densely dotted, ->, line width=0.3mm, >=latex, shorten <= 0.2cm, shorten >= 0.1cm](1_5)--(3_u);
\draw[densely dotted, ->, line width=0.3mm, >=latex, shorten <= 0.2cm, shorten >= 0.1cm](1_5)--(4_u);

\draw[densely dashdotdotted, ->, line width=0.3mm, >=latex, shorten <= 0.2cm, shorten >= 0.1cm](2_6)--(1_u);
\draw[densely dashdotdotted, ->, line width=0.3mm, >=latex, shorten <= 0.2cm, shorten >= 0.1cm](2_6)--(3_u);
\draw[densely dotted, ->, line width=0.3mm, >=latex, shorten <= 0.2cm, shorten >= 0.1cm](1_6)--(1_u);
\draw[densely dotted, ->, line width=0.3mm, >=latex, shorten <= 0.2cm, shorten >= 0.1cm](1_6)--(3_u);
\draw[densely dotted, ->, line width=0.3mm, >=latex, shorten <= 0.2cm, shorten >= 0.1cm](1_6)--(5_u);

\draw[densely dotted, ->, line width=0.3mm, >=latex, shorten <= 0.2cm, shorten >= 0.1cm](2_7)--(2_u);
\draw[densely dotted, ->, line width=0.3mm, >=latex, shorten <= 0.2cm, shorten >= 0.1cm](2_7)--(3_u);
\draw[densely dotted, ->, line width=0.3mm, >=latex, shorten <= 0.2cm, shorten >= 0.1cm](2_7)--(5_u);
\draw[densely dotted, ->, line width=0.3mm, >=latex, shorten <= 0.2cm, shorten >= 0.1cm](3_7)--(2_u);
\draw[densely dotted, ->, line width=0.3mm, >=latex, shorten <= 0.2cm, shorten >= 0.1cm](3_7)--(3_u);
\draw[densely dotted, ->, line width=0.3mm, >=latex, shorten <= 0.2cm, shorten >= 0.1cm](3_7)--(5_u);

\draw[->, line width=0.3mm, >=latex, shorten <= 0.2cm, shorten >= 0.1cm](1_7)--(1_u);
\draw[->, line width=0.3mm, >=latex, shorten <= 0.2cm, shorten >= 0.1cm](1_7)--(2_u);
\draw[densely dotted, ->, line width=0.3mm, >=latex, shorten <= 0.2cm, shorten >= 0.1cm](2_8)--(1_u);
\draw[densely dotted, ->, line width=0.3mm, >=latex, shorten <= 0.2cm, shorten >= 0.1cm](2_8)--(4_u);
\draw[densely dotted, ->, line width=0.3mm, >=latex, shorten <= 0.2cm, shorten >= 0.1cm](2_8)--(5_u);
\draw[densely dotted, ->, line width=0.3mm, >=latex, shorten <= 0.2cm, shorten >= 0.1cm](3_8)--(1_u);
\draw[densely dotted, ->, line width=0.3mm, >=latex, shorten <= 0.2cm, shorten >= 0.1cm](3_8)--(4_u);
\draw[densely dotted, ->, line width=0.3mm, >=latex, shorten <= 0.2cm, shorten >= 0.1cm](3_8)--(5_u);

\draw[->, line width=0.3mm, >=latex, shorten <= 0.2cm, shorten >= 0.1cm](1_8)--(1_u);
\draw[->, line width=0.3mm, >=latex, shorten <= 0.2cm, shorten >= 0.1cm](1_8)--(2_u);
\end{tikzpicture}
\captionsetup{justification=centering}
{\captionof{figure}{Partial orientation $D_3$ for $\mathcal{H}$ for $s=5$, $2|A_2|+|A_3|= 2{{s}\choose{\lceil{s/2}\rceil}}-1$, $A_2=\{[1],[2],\ldots,[6]\}$, $A_3=\{[7],[8],\ldots,[13]\}$.}\label{figC6.4.15}}
\end{center}

\begin{center}
\begin{tikzpicture}[thick,scale=0.75]%
\draw(0,1)node[circle, draw, inner sep=0pt, minimum width=3pt](1_u){\scriptsize $(1,\mathtt{c})$};
\draw(0,-2)node[circle, draw, inner sep=0pt, minimum width=3pt](2_u){\scriptsize $(2,\mathtt{c})$};
\draw(0,-5)node[circle, draw, inner sep=0pt, minimum width=3pt](3_u){\scriptsize $(3,\mathtt{c})$};
\draw(0,-8)node[circle, draw, inner sep=0pt, minimum width=3pt](4_u){\scriptsize $(4,\mathtt{c})$};
\draw(0,-11)node[circle, draw, inner sep=0pt, minimum width=3pt](5_u){\scriptsize $(5,\mathtt{c})$};

\draw(-6,2)node[circle, draw, fill=black!100, inner sep=0pt, minimum width=6pt, label={[] 180:{\small $(1,[1,9])$}}](1_91){};
\draw(-6,0)node[circle, draw, fill=black!100, inner sep=0pt, minimum width=6pt, label={[] 180:{\small $(2,[1,9])$}}](2_91){};

\draw(-3,3)node[circle, draw, fill=black!100, inner sep=0pt, minimum width=6pt, label={[yshift=0cm, xshift=-0.45cm] 90:{\small $(1,[9])$}}](1_9){};
\draw(-3,1)node[circle, draw, fill=black!100, inner sep=0pt, minimum width=6pt, label={[yshift=-0.1cm, xshift=-0.45cm] 270:{\small $(2,[9])$}}](2_9){};
\draw(-3,-1)node[circle, draw, fill=black!100, inner sep=0pt, minimum width=6pt, label={[yshift=0cm, xshift=-0.45cm] 270:{\small $(3,[9])$}}](3_9){};

\draw(6,-4)node[circle, draw, fill=black!100, inner sep=0pt, minimum width=6pt, label={[] 0:{\small $(1,[1,10])$}}](1_101){};
\draw(6,-6)node[circle, draw, fill=black!100, inner sep=0pt, minimum width=6pt, label={[] 0:{\small $(2,[1,10])$}}](2_101){};

\draw(3,-3)node[circle, draw, fill=black!100, inner sep=0pt, minimum width=6pt, label={[yshift=-0.2cm, xshift=0.45cm] 270:{\small $(1,[10])$}}](1_10){};
\draw(3,-5)node[circle, draw, fill=black!100, inner sep=0pt, minimum width=6pt, label={[yshift=-0.1cm, xshift=0.45cm] 270:{\small $(2,[10])$}}](2_10){};
\draw(3,-7)node[circle, draw, fill=black!100, inner sep=0pt, minimum width=6pt, label={[yshift=0cm, xshift=0.45cm] 270:{\small $(3,[10])$}}](3_10){};

\draw(-6,-4)node[circle, draw, fill=black!100, inner sep=0pt, minimum width=6pt, label={[] 180:{\small $(1,[1,11])$}}](1_111){};
\draw(-6,-6)node[circle, draw, fill=black!100, inner sep=0pt, minimum width=6pt, label={[] 180:{\small $(2,[1,11])$}}](2_111){};

\draw(-3,-3)node[circle, draw, fill=black!100, inner sep=0pt, minimum width=6pt, label={[yshift=-0.2cm, xshift=-0.45cm] 270:{\small $(1,[11])$}}](1_11){};
\draw(-3,-5)node[circle, draw, fill=black!100, inner sep=0pt, minimum width=6pt, label={[yshift=-0.1cm, xshift=-0.45cm] 270:{\small $(2,[11])$}}](2_11){};
\draw(-3,-7)node[circle, draw, fill=black!100, inner sep=0pt, minimum width=6pt, label={[yshift=0cm, xshift=-0.45cm] 270:{\small $(3,[11])$}}](3_11){};

\draw(6,-10)node[circle, draw, fill=black!100, inner sep=0pt, minimum width=6pt, label={[] 0:{\small $(1,[1,12])$}}](1_121){};
\draw(6,-12)node[circle, draw, fill=black!100, inner sep=0pt, minimum width=6pt, label={[] 0:{\small $(2,[1,12])$}}](2_121){};

\draw(3,-9)node[circle, draw, fill=black!100, inner sep=0pt, minimum width=6pt, label={[yshift=-0.2cm, xshift=0.45cm] 270:{\small $(1,[12])$}}](1_12){};
\draw(3,-11)node[circle, draw, fill=black!100, inner sep=0pt, minimum width=6pt, label={[yshift=-0.1cm, xshift=0.45cm] 270:{\small $(2,[12])$}}](2_12){};
\draw(3,-13)node[circle, draw, fill=black!100, inner sep=0pt, minimum width=6pt, label={[yshift=0cm, xshift=0.45cm] 270:{\small $(3,[12])$}}](3_12){};

\draw(-6,-10)node[circle, draw, fill=black!100, inner sep=0pt, minimum width=6pt, label={[] 180:{\small $(1,[1,13])$}}](1_131){};
\draw(-6,-12)node[circle, draw, fill=black!100, inner sep=0pt, minimum width=6pt, label={[] 180:{\small $(2,[1,13])$}}](2_131){};

\draw(-3,-9)node[circle, draw, fill=black!100, inner sep=0pt, minimum width=6pt, label={[yshift=-0.2cm, xshift=-0.45cm] 270:{\small $(1,[13])$}}](1_13){};
\draw(-3,-11)node[circle, draw, fill=black!100, inner sep=0pt, minimum width=6pt, label={[yshift=-0.1cm, xshift=-0.45cm] 270:{\small $(2,[13])$}}](2_13){};
\draw(-3,-13)node[circle, draw, fill=black!100, inner sep=0pt, minimum width=6pt, label={[yshift=0cm, xshift=-0.45cm] 270:{\small $(3,[13])$}}](3_13){};

\draw(3,5)node[circle, draw, fill=black!100, inner sep=0pt, minimum width=6pt, label={[] 0:{\small $(1,[14])$}}](1_14){};
\draw(3,3)node[circle, draw, fill=black!100, inner sep=0pt, minimum width=6pt, label={[] 0:{\small $(2,[14])$}}](2_14){};
\draw(3,1)node[circle, draw, fill=black!100, inner sep=0pt, minimum width=6pt, label={[] 0:{\small $(1,[15])$}}](1_15){};
\draw(3,-1)node[circle, draw, fill=black!100, inner sep=0pt, minimum width=6pt, label={[] 0:{\small $(2,[15])$}}](2_15){};

\draw[dashed,->, line width=0.3mm, >=latex, shorten <= 0.2cm, shorten >= 0.15cm](1_9)--(1_91);
\draw[dashed,->, line width=0.3mm, >=latex, shorten <= 0.2cm, shorten >= 0.15cm](1_9)--(2_91);
\draw[dashed,->, line width=0.3mm, >=latex, shorten <= 0.2cm, shorten >= 0.15cm](2_9)--(1_91);
\draw[->, line width=0.3mm, >=latex, shorten <= 0.2cm, shorten >= 0.15cm](1_91) to [out=275, in=158] (3_9);
\draw[dashed,->, line width=0.3mm, >=latex, shorten <= 0.2cm, shorten >= 0.15cm](3_9)--(2_91);
\draw[->, line width=0.3mm, >=latex, shorten <= 0.2cm, shorten >= 0.15cm](2_91)--(2_9);

\draw[->, line width=0.3mm, >=latex, shorten <= 0.2cm, shorten >= 0.15cm](1_101) to [out=265, in=22] (3_10);
\draw[->, line width=0.3mm, >=latex, shorten <= 0.2cm, shorten >= 0.15cm](2_101)--(2_10);
\draw[dashed,->, line width=0.3mm, >=latex, shorten <= 0.2cm, shorten >= 0.15cm](1_10)--(1_101);
\draw[dashed,->, line width=0.3mm, >=latex, shorten <= 0.2cm, shorten >= 0.15cm](2_10)--(1_101);
\draw[dashed,->, line width=0.3mm, >=latex, shorten <= 0.2cm, shorten >= 0.15cm](1_10) to [out=338, in=95] (2_101);
\draw[dashed,->, line width=0.3mm, >=latex, shorten <= 0.2cm, shorten >= 0.15cm](3_10)--(2_101);

\draw[->, line width=0.3mm, >=latex, shorten <= 0.2cm, shorten >= 0.15cm](1_111)--(1_11);
\draw[->, line width=0.3mm, >=latex, shorten <= 0.2cm, shorten >= 0.15cm](2_111) to [out=85, in=202] (1_11);
\draw[->, line width=0.3mm, >=latex, shorten <= 0.2cm, shorten >= 0.15cm](1_111)--(2_11);
\draw[dashed,->, line width=0.3mm, >=latex, shorten <= 0.2cm, shorten >= 0.15cm](2_11)--(2_111);
\draw[->, line width=0.3mm, >=latex, shorten <= 0.2cm, shorten >= 0.15cm](2_111)--(3_11);
\draw[dashed,->, line width=0.3mm, >=latex, shorten <= 0.2cm, shorten >= 0.15cm](3_11) to [out=158, in=275] (1_111);

\draw[->, line width=0.3mm, >=latex, shorten <= 0.2cm, shorten >= 0.15cm](1_121)--(1_12);
\draw[->, line width=0.3mm, >=latex, shorten <= 0.2cm, shorten >= 0.15cm](2_121) to [out=95, in=338] (1_12);
\draw[->, line width=0.3mm, >=latex, shorten <= 0.2cm, shorten >= 0.15cm](1_121)--(2_12);
\draw[dashed,->, line width=0.3mm, >=latex, shorten <= 0.2cm, shorten >= 0.15cm](2_12)--(2_121);
\draw[->, line width=0.3mm, >=latex, shorten <= 0.2cm, shorten >= 0.15cm](2_121)--(3_12);
\draw[dashed,->, line width=0.3mm, >=latex, shorten <= 0.2cm, shorten >= 0.15cm](3_12) to [out=22, in=265] (1_121);

\draw[->, line width=0.3mm, >=latex, shorten <= 0.2cm, shorten >= 0.15cm](1_131)--(1_13);
\draw[->, line width=0.3mm, >=latex, shorten <= 0.2cm, shorten >= 0.15cm](2_131) to [out=85, in=202] (1_13);
\draw[->, line width=0.3mm, >=latex, shorten <= 0.2cm, shorten >= 0.15cm](1_131)--(2_13);
\draw[dashed,->, line width=0.3mm, >=latex, shorten <= 0.2cm, shorten >= 0.15cm](2_13)--(2_131);
\draw[->, line width=0.3mm, >=latex, shorten <= 0.2cm, shorten >= 0.15cm](2_131)--(3_13);
\draw[dashed,->, line width=0.3mm, >=latex, shorten <= 0.2cm, shorten >= 0.15cm](3_13) to [out=158, in=275] (1_131);

\draw[densely dotted, ->, line width=0.3mm, >=latex, shorten <= 0.2cm, shorten >= 0.1cm](2_9)--(2_u);
\draw[densely dotted, ->, line width=0.3mm, >=latex, shorten <= 0.2cm, shorten >= 0.1cm](2_9)--(4_u);
\draw[densely dotted, ->, line width=0.3mm, >=latex, shorten <= 0.2cm, shorten >= 0.1cm](2_9)--(5_u);
\draw[densely dotted, ->, line width=0.3mm, >=latex, shorten <= 0.2cm, shorten >= 0.1cm](3_9)--(2_u);
\draw[densely dotted, ->, line width=0.3mm, >=latex, shorten <= 0.2cm, shorten >= 0.1cm](3_9)--(4_u);
\draw[densely dotted, ->, line width=0.3mm, >=latex, shorten <= 0.2cm, shorten >= 0.1cm](3_9)--(5_u);

\draw[->, line width=0.3mm, >=latex, shorten <= 0.2cm, shorten >= 0.1cm](1_9)--(1_u);
\draw[->, line width=0.3mm, >=latex, shorten <= 0.2cm, shorten >= 0.1cm](1_9)--(2_u);

\draw[densely dotted, ->, line width=0.3mm, >=latex, shorten <= 0.2cm, shorten >= 0.1cm](2_10)--(3_u);
\draw[densely dotted, ->, line width=0.3mm, >=latex, shorten <= 0.2cm, shorten >= 0.1cm](2_10)--(4_u);
\draw[densely dotted, ->, line width=0.3mm, >=latex, shorten <= 0.2cm, shorten >= 0.1cm](2_10)--(5_u);
\draw[densely dotted, ->, line width=0.3mm, >=latex, shorten <= 0.2cm, shorten >= 0.1cm](3_10)--(3_u);
\draw[densely dotted, ->, line width=0.3mm, >=latex, shorten <= 0.2cm, shorten >= 0.1cm](3_10)--(4_u);
\draw[densely dotted, ->, line width=0.3mm, >=latex, shorten <= 0.2cm, shorten >= 0.1cm](3_10)--(5_u);

\draw[->, line width=0.3mm, >=latex, shorten <= 0.2cm, shorten >= 0.1cm](1_10)--(1_u);
\draw[->, line width=0.3mm, >=latex, shorten <= 0.2cm, shorten >= 0.1cm](1_10)--(2_u);

\draw[densely dashdotdotted, ->, line width=0.3mm, >=latex, shorten <= 0.2cm, shorten >= 0.1cm](2_11)--(4_u);
\draw[densely dashdotdotted, ->, line width=0.3mm, >=latex, shorten <= 0.2cm, shorten >= 0.1cm](2_11)--(5_u);
\draw[densely dashdotdotted, ->, line width=0.3mm, >=latex, shorten <= 0.2cm, shorten >= 0.1cm](3_11)--(4_u);
\draw[densely dashdotdotted, ->, line width=0.3mm, >=latex, shorten <= 0.2cm, shorten >= 0.1cm](3_11)--(5_u);

\draw[->, line width=0.3mm, >=latex, shorten <= 0.2cm, shorten >= 0.1cm](1_11)--(1_u);
\draw[->, line width=0.3mm, >=latex, shorten <= 0.2cm, shorten >= 0.1cm](1_11)--(2_u);

\draw[densely dashdotdotted, ->, line width=0.3mm, >=latex, shorten <= 0.2cm, shorten >= 0.1cm](2_12)--(3_u);
\draw[densely dashdotdotted, ->, line width=0.3mm, >=latex, shorten <= 0.2cm, shorten >= 0.1cm](2_12)--(5_u);
\draw[densely dashdotdotted, ->, line width=0.3mm, >=latex, shorten <= 0.2cm, shorten >= 0.1cm](3_12)--(3_u);
\draw[densely dashdotdotted, ->, line width=0.3mm, >=latex, shorten <= 0.2cm, shorten >= 0.1cm](3_12)--(5_u);

\draw[->, line width=0.3mm, >=latex, shorten <= 0.2cm, shorten >= 0.1cm](1_12)--(1_u);
\draw[->, line width=0.3mm, >=latex, shorten <= 0.2cm, shorten >= 0.1cm](1_12)--(2_u);

\draw[densely dashdotdotted, ->, line width=0.3mm, >=latex, shorten <= 0.2cm, shorten >= 0.1cm](2_13)--(3_u);
\draw[densely dashdotdotted, ->, line width=0.3mm, >=latex, shorten <= 0.2cm, shorten >= 0.1cm](2_13)--(4_u);
\draw[densely dashdotdotted, ->, line width=0.3mm, >=latex, shorten <= 0.2cm, shorten >= 0.1cm](3_13)--(3_u);
\draw[densely dashdotdotted, ->, line width=0.3mm, >=latex, shorten <= 0.2cm, shorten >= 0.1cm](3_13)--(4_u);

\draw[->, line width=0.3mm, >=latex, shorten <= 0.2cm, shorten >= 0.1cm](1_13)--(1_u);
\draw[->, line width=0.3mm, >=latex, shorten <= 0.2cm, shorten >= 0.1cm](1_13)--(2_u);

\draw[->, line width=0.3mm, >=latex, shorten <= 0.2cm, shorten >= 0.1cm](1_14)--(1_u);
\draw[->, line width=0.3mm, >=latex, shorten <= 0.2cm, shorten >= 0.1cm](1_14)--(2_u);
\draw[->, line width=0.3mm, >=latex, shorten <= 0.2cm, shorten >= 0.1cm](2_14)--(1_u);
\draw[->, line width=0.3mm, >=latex, shorten <= 0.2cm, shorten >= 0.1cm](2_14)--(2_u);

\draw[->, line width=0.3mm, >=latex, shorten <= 0.2cm, shorten >= 0.1cm](1_15)--(1_u);
\draw[->, line width=0.3mm, >=latex, shorten <= 0.2cm, shorten >= 0.1cm](1_15)--(2_u);
\draw[->, line width=0.3mm, >=latex, shorten <= 0.2cm, shorten >= 0.1cm](2_15)--(1_u);
\draw[->, line width=0.3mm, >=latex, shorten <= 0.2cm, shorten >= 0.1cm](2_15)--(2_u);
\end{tikzpicture}
\captionsetup{justification=centering}
{\captionof{figure}{Partial orientation $D_3$ for $\mathcal{H}$ for $s=5$, $2|A_2|+|A_3|= 2{{s}\choose{\lceil{s/2}\rceil}}-1$, $A_2=\{[1],[2],\ldots,[6]\}$, $A_3=\{[7],[8],\ldots,[13]\}$, $E=\{[14],[15]\}$.}\label{figC6.4.16}}
\end{center}

\begin{cor}\label{corC6.4.10}
Suppose $s\ge 5$ is odd, $A_{\ge 4}=\emptyset$, $|A_2|\ge \lceil\frac{s}{2}\rceil \lfloor\frac{s}{2}\rfloor$, $|A_3|\ge 2$, and $2|A_2|+|A_3|=2{{s}\choose{\lceil{s/2}\rceil}}-1$ for a $\mathcal{T}$. If $D$ is an optimal orientation of $\mathcal{T}$, then either $D$ or $\tilde{D}$ fulfills the following, after a suitable relabelling of vertices and with $A_3$ partitioned into $A^O_3$  and $A^I_3$.
\\(I) $|O^\mathtt{c}((1,[i]))|=|I^\mathtt{c}((2,[i]))|=\lceil\frac{s}{2}\rceil$ for all $[i]\in A_2$,
\\(II) $|O^\mathtt{c}((3,[j]))|=|I^\mathtt{c}((3,[k]))|=\lceil\frac{s}{2}\rceil$ for all $[j]\in A^O_3$, $[k]\in A^I_3$, and
\\(III) $(2,[i])\rightarrow (p,[\alpha,i])\rightarrow (1,[i])$ for all $[i]\in A_2$, all $1\le\alpha\le deg_T([i])-1$ and $1\le p\le s_{[\alpha,i]}$.
\\(1A) $|A_2\cup A^O_3|= {{s}\choose{\lceil{s/2}\rceil}}$ and $|A_2\cup A^I_3|= {{s}\choose{\lceil{s/2}\rceil}}-1$;
\\(2A) $|O^\mathtt{c}((1,[i]))|=\lfloor\frac{s}{2}\rfloor$ and $|O^\mathtt{c}((2,[i]))|=\lceil\frac{s}{2}\rceil$ for all $[i]\in A^O_3$;
\\(3A) $O^\mathtt{c}((1,[i]))=O^\mathtt{c}((1,[j]))$ for all $[i],[j]\in A^O_3$;
\\(4A) $O^\mathtt{c}((2,[i]))=O^\mathtt{c}((3,[i]))$ for all $[i]\in A^O_3$;
\\(5A) $|O^\mathtt{c}((1,[j]))|=|O^\mathtt{c}((2,[j]))|=\lfloor\frac{s}{2}\rfloor$ for all $[j]\in A^I_3$;
\\(6A) $O^\mathtt{c}((1,[j]))= O^\mathtt{c}((1,[i]))$ and $O^\mathtt{c}((2,[j]))=O^\mathtt{c}((3,[j]))$ for all $[i]\in A^O_3$, $[j]\in A^I_3$.
\end{cor}
\noindent\textit{Proof}: Since $d(D)=d(\tilde{D})=4$ and $2|A_2|+|A_3|=2{{s}\choose{\lfloor{s/2}\rfloor}}-1$, it follows from the proof of Proposition \ref{ppnC6.4.3} that only Subcases 2.3.6.3 and 2.3.6.4 are possible. (I) and (II) follow from the assumption of Subcase 2.3.6, $B^O_2\cup B^O_3\cup B^I_2\cup B^I_3\subseteq {{(\mathbb{N}_s,\mathtt{c})}\choose{\lceil{s/2}\rceil}}$. Also, (1A)-(6A) follow from (\ref{eqC6.4.11}), Claims 2A-7A  in Subcase 2.3.6.3 and their analogues in Subcase 2.3.6.4.
\indent\par To prove (III), suppose for some $[i]\in A_2$ that there exists some $1\le p,q\le s_{[\alpha,i]}$, $p\neq q$, such that $(p,[\alpha,i])\rightarrow (2,[i])$ and $(q,[\alpha,i])\rightarrow (1,[i])$ for a contradiction. By Lemma \ref{lemC6.2.13}(a), $(1,[i])\rightarrow (p,[\alpha,i])\rightarrow (2,[i])$ and $(2,[i])\rightarrow (q,[\alpha,i])\rightarrow (1,[i])$. This implies $O^\mathtt{c}((1,[i])), I^\mathtt{c}((2,[i]))\in B^I_2\cup B^I_3$. Therefore, $|O^\mathtt{c}((1,[i]))|=|I^\mathtt{c}((2,[i]))|=\lceil\frac{s}{2}\rceil$ implies $|I^\mathtt{c}((1,[i]))|=|O^\mathtt{c}((2,[i]))|=\lfloor\frac{s}{2}\rfloor$. By Lemmas \ref{lemC6.2.15} and \ref{lemC6.2.16}, $\{O^\mathtt{c}((2,[i]))\}\cup B^O_2\cup B^O_3-\{O^\mathtt{c}((1,[i]))\}$ and $\{I^\mathtt{c}((1,[i]))\}\cup B^I_2\cup B^I_3-\{I^\mathtt{c}((2,[i]))\}$ are antichains, each containing contains elements of different sizes. Hence, $|A_2|+|A^O_3|=|B^O_2|+|B^O_3|=|\{O^\mathtt{c}((1,[i]))\}\cup B^O_2\cup B^O_3-\{O^\mathtt{c}((2,[i]))\}|\le {{s}\choose{\lfloor{s/2}\rfloor}}-1$ by Sperner's theorem. Similarly, $|A_2|+|A^I_3|=|B^I_2|+|B^I_3|=|\{I^\mathtt{c}((1,[i]))\}\cup B^I_2\cup B^I_3-\{I^\mathtt{c}((2,[i]))\}|\le {{s}\choose{\lfloor{s/2}\rfloor}}-1$. Hence, $2|A_2|+|A_3|=(|A_2|+|A^O_3|)+(|A_2|+|A^I_3|)\le 2[{{s}\choose{\lfloor{s/2}\rfloor}}-1]=2{{s}\choose{\lceil{s/2}\rceil}}-2$, a contradiction.
\qed

\indent\par The optimal orientation(s) described in Corollary \ref{corC6.4.10} was extended to the construction $D_3$ in Case 3 of the proof of Proposition \ref{ppnC6.4.3}.

\begin{ppn}\label{ppnC6.4.11}
Suppose $s\ge 3$ is odd, $A_2\neq \emptyset$, $A_3=\emptyset$ and $A_{\ge 4}\neq \emptyset$ for a $\mathcal{T}$. Then, $\mathcal{T}\in \mathscr{C}_0$ if and only if $|A_2|\le{{s}\choose{\lceil{s/2}\rceil}}-1$.
\end{ppn}
\noindent\textit{Proof}: $(\Rightarrow)$ Since $\mathcal{T}\in \mathscr{C}_0$, there exists an orientation $D$ of $\mathcal{T}$, where $d(D)=4$. As $A_2\neq\emptyset$, we assume (\ref{eqC6.3.1})-(\ref{eqC6.3.2}) here. Let $[j]\in A_{\ge 4}$. If $|O^\mathtt{c}((1,[j]))|\ge\lceil\frac{s}{2}\rceil$, then for any $[i]\in A_2$, $d_D((1,[1,i]),(1,[j]))=3$ implies $O^\mathtt{c}((1,[i]))\cap I^\mathtt{c}((1,[j]))\neq \emptyset$. Hence, by Lih's theorem, $|A_2|=|B^O_2|\le {{s}\choose{\lceil{s/2}\rceil}}-{{s-|I^\mathtt{c}((1,[i]))|}\choose{\lceil{s/2}\rceil}}\le {{s}\choose{\lceil{s/2}\rceil}}-1$. Suppose $|O^\mathtt{c}((1,[j]))|\le\lfloor\frac{s}{2}\rfloor$. Equivalently, $|I^\mathtt{c}((1,[j]))|\ge\lceil\frac{s}{2}\rceil$. So, this case follows from the previous case by the Duality Lemma.
\\
\\$(\Leftarrow)$ If $|A_2|+|A_{\ge 4}|\le {{s}\choose{\lceil{s/2}\rceil}}-1$, then by Corollary \ref{corC6.3.8}(i), $\mathcal{T}\in\mathscr{C}_0$. Hence, we assume $|A_2|+|A_{\ge 4}|\ge {{s}\choose{\lceil{s/2}\rceil}}$ hereafter, on top of the hypothesis that $|A_2|\le{{s}\choose{\lceil{s/2}\rceil}}-1$. If $|A_2|\ge s-1$, define $A^{\diamond}_2=A_2$. Otherwise, $A^{\diamond}_2=A_2\cup A^*$, where $A^*$ is an arbitrary subset of $A_{\ge 4}$ such that $|A^{\diamond}_2|=s-1$. Then, let $A^{\diamond}_4=A_2\cup A_{\ge 4}-A^{\diamond}_2$. Furthermore, assume without loss of generality that $A^{\diamond}_2=\{[i]\mid i\in\mathbb{N}_{|A^{\diamond}_2|}\}$ and $A^{\diamond}_4=\{[i]\mid i\in\mathbb{N}_{|A^{\diamond}_2|+|A^{\diamond}_4|}-\mathbb{N}_{|A^{\diamond}_2|}\}$.
\indent\par Let $\mathcal{H}=T(t_1,t_2,\ldots, t_n)$ be the subgraph of $\mathcal{T}$, where $t_\mathtt{c}=s$, $t_{[i]}=4$ for all $[i]\in \mathcal{T}(A^{\diamond}_4)$ and $t_v=2$ otherwise. We will use $A_j$ for $\mathcal{H}(A_j)$ for the remainder of this proof. Define an orientation $D$ of $\mathcal{H}$ as follows.
\begin{align}
& (2,[i])\rightarrow \{(1,[\alpha,i]),(2,[\alpha,i])\}\rightarrow (1,[i]),\text{ and}\label{eqC6.4.38}\\
&\bar{\lambda}_{i+1}\rightarrow (1,[i])\rightarrow \lambda_{i+1}\rightarrow (2,[i]) \rightarrow \bar{\lambda}_{i+1}\label{eqC6.4.39}
\end{align}
for all $1\le i\le |A_2|$ and all $1\le \alpha\le deg_T([i])-1$, i.e., the $\lceil\frac{s}{2}\rceil$-sets $\lambda_2,\lambda_3,\ldots,\lambda_{|A_2|+1}$ are used as `in-sets' (`out-sets' resp.) to construct $B^I_2$ ($B^O_2$ resp.).
\begin{align}
&(2,[\beta,j])\rightarrow \{(2,[j]),(4,[j])\}\rightarrow (1,[\beta,j])\rightarrow \{(1,[j]), (3,[j])\}\rightarrow (2,[\beta,j]),\label{eqC6.4.40}\\
&\text{and }\bar{\lambda}_1 \rightarrow \{(1,[j]), (4,[j])\}\rightarrow \lambda_1 \rightarrow\{(2,[j]), (3,[j])\}\rightarrow \bar{\lambda}_1 \label{eqC6.4.41}
\end{align}
for all $[j]\in A_4$ and all $1\le \beta\le deg_T([j])-1$. 
\begin{align}
\lambda_1 \rightarrow \{(1,[k]), (2,[k])\} \rightarrow \bar{\lambda}_1 \label{eqC6.4.42}
\end{align}
for all $[k]\in E$. (See Figure \ref{figC6.4.17} for $D$ when $s=3$.)
\\
\\Claim: $d_{D}(v,w)\le 4$ for all $v,w\in V(D)$.
\\
\\Case 1.1. $v,w \in \{(1,[\alpha,i]),(2,[\alpha,i]),(1,[i]),(2,[i])\}$ for each $[i]\in A_2$ and $1\le\alpha\le deg_T([i])-1$.
\indent\par Since the orientation defined for $A_2$ (see (\ref{eqC6.4.38})-(\ref{eqC6.4.39})) is similar to $D_1$ in Proposition \ref{ppnC6.4.3} (see (\ref{eqC6.4.12})-(\ref{eqC6.4.13})), this case follows from Case 1.1.1 of Proposition \ref{ppnC6.4.3}.
\\
\\Case 1.2. $v,w \in \{(1,[\alpha,i]),(2,[\alpha,i]),(1,[i]),(2,[i]), (3,[i]), (4,[i])\}$ for each $[i]\in A_4$ and $1\le\alpha\le deg_T([i])-1$.
\indent\par Since the orientation defined for $A_4$ (see (\ref{eqC6.4.40})-(\ref{eqC6.4.41})) is similar to that in Proposition \ref{ppnC6.4.1} (see (\ref{eqC6.4.7})-(\ref{eqC6.4.8})), this case follows from Case 1.3 of Proposition \ref{ppnC6.4.1}.
\\
\\Case 2. For each $[i],[j]\in A_2$, $i \neq j$, each $1\le\alpha\le deg_T([i])-1$, and each $1\le\beta\le deg_T([j])-1$,
\\(i) $v=(p,[\alpha,i]), w=(q,[j])$ for $p,q=1,2$.
\\(ii) $v=(p,[i]), w=(q,[\beta,j])$ for $p,q=1,2$.
\\(iii) $v=(p,[\alpha,i]), w=(q,[\beta,j])$ for each $p,q=1,2$.
\indent\par Since the orientation defined for $A_2$ (see (\ref{eqC6.4.38})-(\ref{eqC6.4.39})) is similar to $D_1$ in Proposition \ref{ppnC6.4.3} (see (\ref{eqC6.4.12})-(\ref{eqC6.4.13})), this case follows from Case 1.2 of Proposition \ref{ppnC6.4.3}.
\\
\\Case 3. For each $[i],[j]\in A_4$, $i \neq j$, each $1\le\alpha\le deg_T([i])-1$, and each $1\le\beta\le deg_T([j])-1$,
\\(i) $v=(p,[\alpha,i]), w=(q,[j])$ for $p=1,2$, and $q=1,2,3,4$.
\\(ii) $v=(p,[i]), w=(q,[\beta,j])$ for for $p=1,2,3,4$, and $q=1,2$.
\\(iii) $v=(p,[\alpha,i]), w=(q,[\beta,j])$ for $p,q=1,2$.
\indent\par Since the orientation defined for $A_4$ (see (\ref{eqC6.4.40})-(\ref{eqC6.4.41})) is similar to that in Proposition \ref{ppnC6.4.1} (see (\ref{eqC6.4.7})-(\ref{eqC6.4.8})), this case follows from Case 4 of Proposition \ref{ppnC6.4.1}.
\\
\\Case 4.  For each $[i]\in A_2$, each $[j]\in A_4$, each $1\le\alpha\le deg_T([i])-1$, and each $1\le\beta\le deg_T([j])-1$,
\\(i) $v=(p,[\alpha,i]), w=(q,[j])$ for each $p=1,2$ and $q=1,2,3,4$.
\indent\par By (\ref{eqC6.4.38})-(\ref{eqC6.4.39}) and (\ref{eqC6.4.41}), $\{(1,[\alpha,i]),(2,[\alpha,i])\}\rightarrow (1,[i])\rightarrow \lambda_{i+1}$, $\lambda_1\rightarrow \{(2,[j]), (3,[j])\}$, and $\bar{\lambda}_1\rightarrow \{(1,[j]), (4,[j])\}$. Since $\lambda_{i+1}\in {{(\mathbb{N}_s,\mathtt{c})}\choose{\lceil s/2\rceil}}-\{\lambda_1\}$, there exist some $x_i\in\lambda_{i+1}\cap\lambda_1$ and $y_i\in\lambda_{i+1}\cap\bar{\lambda}_1$ such that $(1,[i])\rightarrow \{x_i, y_i\}$.
\\
\\(ii) $v=(q,[j]), w=(p,[\alpha,i])$ for each $p=1,2$ and $q=1,2,3,4$.
\indent\par By (\ref{eqC6.4.38})-(\ref{eqC6.4.39}) and (\ref{eqC6.4.41}), $\{(1,[j]), (4,[j])\}\rightarrow \lambda_1$, $\{(2,[j]), (3,[j])\}\rightarrow \bar{\lambda}_1$, and $\lambda_{i+1}\rightarrow (2,[i])\rightarrow \{(1,[\alpha,i]),(2,[\alpha,i])\}$. Since $\lambda_{i+1}\in {{(\mathbb{N}_s,\mathtt{c})}\choose{\lceil s/2\rceil}}-\{\lambda_1\}$, there exist some $x_{i}\in \lambda_{i+1}\cap\lambda_1$ and $y_{i}\in\lambda_{i+1}\cap\bar{\lambda}_{1}$ such that $\{x_{i}, y_{i}\}\rightarrow (2,[i])$.
\\
\\(iii) $v=(p,[i]), w=(q,[\beta,j])$ for each $p,q=1,2$.
\indent\par From (i), $d_D((1,[i]), (r,[j]))=2$ for $r=2,3$. Furthermore by (\ref{eqC6.4.39}) and (\ref{eqC6.4.41}), since $(2,[i])\rightarrow \bar{\lambda}_{i+1}$ and $i+1\neq 1$, there exists some $x_i\in\bar{\lambda}_{i+1}\cap\lambda_1$ such that $(2,[i])\rightarrow x_i\rightarrow \{(2,[j]),(3,[j])\}$. Since $(2,[j])\rightarrow (1,[\beta,j])$ and $(3,[j])\rightarrow (2,[\beta,j])$ by (\ref{eqC6.4.40}), this subcase follows.
\\
\\(iv) $v=(q,[\beta,j]), w=(p,[i])$ for each $p,q=1,2$.
\indent\par By (\ref{eqC6.4.39})-(\ref{eqC6.4.41}), $(2,[\beta,j])\rightarrow \{(2,[j]),(4,[j])\}$, $(1,[\beta,j])\rightarrow \{(1,[j]), (3,[j])\}$, $\{(1,[j]),$ $(4,[j])\}$ $\rightarrow \lambda_1$, $\{(2,[j]), (3,[j])\}\rightarrow \bar{\lambda}_1$, and $|I^\mathtt{c}((p,[i]))|>0$.
\\
\\(v) $v=(p,[\alpha,i]), w=(q,[\beta,j])$ for each $p,q=1,2$.
\indent\par From (i), $d_D((p,[\alpha,i]),(r,[j]))=3$ for $r=2,3$. Since $(2,[j])\rightarrow (1,[\beta,j])$ and $(3,[j])\rightarrow (2,[\beta,j])$ by (\ref{eqC6.4.40}), this subcase follows.
\\
\\(vi) $v=(q,[\beta,j]), w=(p,[\alpha,i])$ for each $p,q=1,2$.
\indent\par From (iv), $d_D((q,[\beta,j]),(2,[i]))=3$. Since $(2,[i])\rightarrow \{(1,[\alpha,i]), (2,[\alpha,i])\}$ by (\ref{eqC6.4.38}), this subcase follows.
\\
\\Case 5. For each $[i]\in A_2$, each $1\le\alpha\le deg_T([i])-1$, and each $[j]\in E$,
\\(i) $v=(p,[\alpha,i]), w=(q,[j])$ for each $p=1,2$ and $q=1,2$.
\\(ii) $v=(q,[j]), w=(p,[\alpha,i])$ for each $p=1,2$ and $q=1,2$.
\indent\par Let $[k]\in A_4$. Since $\lambda_1\rightarrow \{(q,[j]),(2,[k])\}\rightarrow \bar{\lambda}_1$ by (\ref{eqC6.4.41})-(\ref{eqC6.4.42}), this case follows from Cases 4(i)-(ii).
\\
\\Case 6. For each $[i]\in A_4$, each $1\le\alpha\le deg_T([i])-1$, and each $[j]\in E$,
\\(i) $v=(p,[\alpha,i]), w=(q,[j])$ for each $p=1,2$ and $q=1,2$.
\\(ii) $v=(q,[j]), w=(p,[\alpha,i])$ for each $p=1,2$ and $q=1,2$.
\indent\par Let $[k]\in A_4$, where $k\neq i$. Since $\lambda_1\rightarrow \{(q,[j]),(2,[k])\}\rightarrow \bar{\lambda}_1$ by (\ref{eqC6.4.41})-(\ref{eqC6.4.42}), this case follows from Cases 3(i)-(ii).
\\
\\Case 7. $v=(r_1,\mathtt{c})$ and $w=(r_2,\mathtt{c})$ for $r_1\neq r_2$ and $1\le r_1, r_2\le s$.
\indent\par Here, we want to prove a stronger claim, $d_D((r_1,\mathtt{c}), (r_2,\mathtt{c}))=2$. Let $x_1=(2,[|A_2|+1])$, $x_{k+1}=(2,[k])$ for $1\le k\le s-1$. Observe from (\ref{eqC6.4.39}) and (\ref{eqC6.4.41}) that $\lambda_k\rightarrow x_k\rightarrow \bar{\lambda}_k$ for all $1\le k\le s$ and the subgraph induced by $V_1=(\mathbb{N}_s,\mathtt{c})$ and $V_2=\{x_i\mid 1\le i \le s\}$ is a complete bipartite graph $K(V_1,V_2)$. By Lemma \ref{lemC6.2.19}, $d_{D}((r_1,\mathtt{c}), (r_2,\mathtt{c}))=2$.
\\
\\Case 8. $v\in \{(1,[i]), (2,[i]), (3,[i]), (4,[i]), (1,[\alpha,i]), (2,[\alpha,i])\}$ for each $1\le i\le deg_T(\mathtt{c})$ and $1\le\alpha\le deg_T([i])-1$, and $w=(r,\mathtt{c})$ for $1\le r\le s$.
\indent\par Note that there exists some $1\le k\le s$ such that $d_D(v,(k,\mathtt{c}))\le 2$, and $d_D((k,\mathtt{c}),w)\le 2$ by Case 7. Hence, it follows that $d_D(v,w)\le d_D(v,(k,\mathtt{c}))+d_D((k,\mathtt{c}),w)\le 4$.
\\
\\Case 9. $v=(r,\mathtt{c})$ for $r=1,2,\ldots, s$ and $w\in \{(1,[i]), (2,[i]), (3,[i]), (4,[i]), (1,[\alpha,i]), (2,[\alpha,i])\}$ for each $1\le i\le deg_T(\mathtt{c})$ and $1\le\alpha\le deg_T([i])-1$.
\indent\par Note that there exists some $1\le k\le s$ such that $d_D((k,\mathtt{c}), w)\le 2$, and $d_D(v,(k,\mathtt{c}))\le 2$ by Case 7. Hence, it follows that $d_D(v,w)\le d_D(v,(k,\mathtt{c}))+d_D((k,\mathtt{c}),w)\le 4$.
\\
\\Case 10. $v=(p,[i])$ and $w=(q, [j])$, where $1\le p,q\le 4$ and $1\le i,j\le deg_T(\mathtt{c})$.
\noindent\par This follows from the fact that $|O^\mathtt{c}((p,[i]))|>0$, $|I^\mathtt{c}((q,[j]))|>0$, and $d_{D}((r_1,\mathtt{c}), (r_2,\mathtt{c}))$ $=2$ for any $r_1\neq r_2$ and $1\le r_1, r_2\le s$.
\\
\indent\par Therefore, the claim follows. Since every vertex lies in a directed $C_4$ for $D$ and $d(D)=4$, $\bar{d}(\mathcal{T})\le \max \{4, d(D)\}$ by Lemma \ref{lemC6.1.3}, and thus $\bar{d}(\mathcal{T})= 4$.
\qed

\begin{center}
\begin{tikzpicture}[thick,scale=0.75]%
\draw(0,1)node[circle, draw, inner sep=0pt, minimum width=3pt](1_u){\scriptsize $(1,\mathtt{c})$};
\draw(0,-3)node[circle, draw, inner sep=0pt, minimum width=3pt](2_u){\scriptsize $(2,\mathtt{c})$};
\draw(0,-7)node[circle, draw, inner sep=0pt, minimum width=3pt](3_u){\scriptsize $(3,\mathtt{c})$};

\draw(-6,4)node[circle, draw, fill=black!100, inner sep=0pt, minimum width=6pt, label={[] 180:{\small $(1,[1,1])$}}](1_11){};
\draw(-6,2)node[circle, draw, fill=black!100, inner sep=0pt, minimum width=6pt, label={[] 180:{\small $(2,[1,1])$}}](2_11){};

\draw(-3,4)node[circle, draw, fill=black!100, inner sep=0pt, minimum width=6pt, label={[yshift=0cm, xshift=-0.45cm] 90:{\small $(1,[1])$}}](1_1){};
\draw(-3,2)node[circle, draw, fill=black!100, inner sep=0pt, minimum width=6pt, label={[yshift=0cm, xshift=-0.45cm] 270:{\small $(2,[1])$}}](2_1){};

\draw(6,4)node[circle, draw, fill=black!100, inner sep=0pt, minimum width=6pt, label={[] 0:{\small $(1,[1,2])$}}](1_21){};
\draw(6,2)node[circle, draw, fill=black!100, inner sep=0pt, minimum width=6pt, label={[] 0:{\small $(2,[1,2])$}}](2_21){};

\draw(3,4)node[circle, draw, fill=black!100, inner sep=0pt, minimum width=6pt, label={[yshift=0cm, xshift=0.45cm] 90:{\small $(1,[2])$}}](1_2){};
\draw(3,2)node[circle, draw, fill=black!100, inner sep=0pt, minimum width=6pt, label={[yshift=0cm, xshift=0.45cm] 270:{\small $(2,[2])$}}](2_2){};

\draw(-6,-2)node[circle, draw, fill=black!100, inner sep=0pt, minimum width=6pt, label={[] 180:{\small $(1,[1,3])$}}](1_31){};
\draw(-6,-4)node[circle, draw, fill=black!100, inner sep=0pt, minimum width=6pt, label={[] 180:{\small $(2,[1,3])$}}](2_31){};

\draw(-3,0)node[circle, draw, fill=black!100, inner sep=0pt, minimum width=6pt, label={[yshift=-0.15cm, xshift=-0.45cm] 270:{\small $(1,[3])$}}](1_3){};
\draw(-3,-2)node[circle, draw, fill=black!100, inner sep=0pt, minimum width=6pt, label={[yshift=0cm, xshift=-0.45cm] 270:{\small $(2,[3])$}}](2_3){};
\draw(-3,-4)node[circle, draw, fill=black!100, inner sep=0pt, minimum width=6pt, label={[yshift=0.05cm, xshift=-0.45cm] 270:{\small $(3,[3])$}}](3_3){};
\draw(-3,-6)node[circle, draw, fill=black!100, inner sep=0pt, minimum width=6pt, label={[yshift=0cm, xshift=-0.45cm] 270:{\small $(4,[3])$}}](4_3){};

\draw(6,-2)node[circle, draw, fill=black!100, inner sep=0pt, minimum width=6pt, label={[] 0:{\small $(1,[1,4])$}}](1_41){};
\draw(6,-4)node[circle, draw, fill=black!100, inner sep=0pt, minimum width=6pt, label={[] 0:{\small $(2,[1,4])$}}](2_41){};

\draw(3,0)node[circle, draw, fill=black!100, inner sep=0pt, minimum width=6pt, label={[yshift=-0.15cm, xshift=0.45cm] 270:{\small $(1,[4])$}}](1_4){};
\draw(3,-2)node[circle, draw, fill=black!100, inner sep=0pt, minimum width=6pt, label={[yshift=0cm, xshift=0.45cm] 270:{\small $(2,[4])$}}](2_4){};
\draw(3,-4)node[circle, draw, fill=black!100, inner sep=0pt, minimum width=6pt, label={[yshift=0.05cm, xshift=0.45cm] 270:{\small $(3,[4])$}}](3_4){};
\draw(3,-6)node[circle, draw, fill=black!100, inner sep=0pt, minimum width=6pt, label={[yshift=0cm, xshift=0.45cm] 270:{\small $(4,[4])$}}](4_4){};

\draw(-3,-8)node[circle, draw, fill=black!100, inner sep=0pt, minimum width=6pt, label={[] 180:{\small $(1,[5])$}}](1_5){};
\draw(-3,-10)node[circle, draw, fill=black!100, inner sep=0pt, minimum width=6pt, label={[] 180:{\small $(2,[5])$}}](2_5){};

\draw(3,-8)node[circle, draw, fill=black!100, inner sep=0pt, minimum width=6pt, label={[] 0:{\small $(1,[6])$}}](1_6){};
\draw(3,-10)node[circle, draw, fill=black!100, inner sep=0pt, minimum width=6pt, label={[] 0:{\small $(2,[6])$}}](2_6){};

\draw[->, line width=0.3mm, >=latex, shorten <= 0.2cm, shorten >= 0.15cm](1_11)--(1_1);
\draw[->, line width=0.3mm, >=latex, shorten <= 0.2cm, shorten >= 0.15cm](2_11)--(1_1);
\draw[dashed,->, line width=0.3mm, >=latex, shorten <= 0.2cm, shorten >= 0.15cm](2_1)--(1_11);
\draw[dashed,->, line width=0.3mm, >=latex, shorten <= 0.2cm, shorten >= 0.15cm](2_1)--(2_11);

\draw[dashed,->, line width=0.3mm, >=latex, shorten <= 0.2cm, shorten >= 0.15cm](2_2)--(1_21);
\draw[dashed,->, line width=0.3mm, >=latex, shorten <= 0.2cm, shorten >= 0.15cm](2_2)--(2_21);
\draw[->, line width=0.3mm, >=latex, shorten <= 0.2cm, shorten >= 0.15cm](1_21)--(1_2);
\draw[->, line width=0.3mm, >=latex, shorten <= 0.2cm, shorten >= 0.15cm](2_21)--(1_2);

\draw[dashed,->, line width=0.3mm, >=latex, shorten <= 0.2cm, shorten >= 0.15cm](2_4)--(1_41);
\draw[dashed,->, line width=0.3mm, >=latex, shorten <= 0.2cm, shorten >= 0.15cm](1_4) to [out=345, in=100] (2_41);
\draw[dashed,->, line width=0.3mm, >=latex, shorten <= 0.2cm, shorten >= 0.15cm](4_4) to [out=40, in=250] (1_41);
\draw[dashed,->, line width=0.3mm, >=latex, shorten <= 0.2cm, shorten >= 0.15cm](3_4)--(2_41);
\draw[->, line width=0.3mm, >=latex, shorten <= 0.2cm, shorten >= 0.15cm](1_41) to [out=115, in=0] (1_4);
\draw[->, line width=0.3mm, >=latex, shorten <= 0.2cm, shorten >= 0.15cm](2_41) to [out=115, in=350] (2_4);
\draw[->, line width=0.3mm, >=latex, shorten <= 0.2cm, shorten >= 0.15cm](1_41)--(3_4);
\draw[->, line width=0.3mm, >=latex, shorten <= 0.2cm, shorten >= 0.15cm](2_41)--(4_4);

\draw[dashed,->, line width=0.3mm, >=latex, shorten <= 0.2cm, shorten >= 0.15cm](2_3)--(1_31);
\draw[dashed,->, line width=0.3mm, >=latex, shorten <= 0.2cm, shorten >= 0.15cm](1_3) to [out=195, in=80] (2_31);
\draw[dashed,->, line width=0.3mm, >=latex, shorten <= 0.2cm, shorten >= 0.15cm](4_3) to [out=140, in=290] (1_31);
\draw[dashed,->, line width=0.3mm, >=latex, shorten <= 0.2cm, shorten >= 0.15cm](3_3)--(2_31);
\draw[->, line width=0.3mm, >=latex, shorten <= 0.2cm, shorten >= 0.15cm](1_31) to [out=65, in=180] (1_3);
\draw[->, line width=0.3mm, >=latex, shorten <= 0.2cm, shorten >= 0.15cm](2_31) to [out=65, in=190] (2_3);
\draw[->, line width=0.3mm, >=latex, shorten <= 0.2cm, shorten >= 0.15cm](1_31)--(3_3);
\draw[->, line width=0.3mm, >=latex, shorten <= 0.2cm, shorten >= 0.15cm](2_31)--(4_3);

\draw[densely dashdotdotted, ->, line width=0.3mm, >=latex, shorten <= 0.2cm, shorten >= 0.1cm](2_1)--(1_u);
\draw[densely dotted, ->, line width=0.3mm, >=latex, shorten <= 0.2cm, shorten >= 0.1cm](1_1)--(2_u);
\draw[densely dotted, ->, line width=0.3mm, >=latex, shorten <= 0.2cm, shorten >= 0.1cm](1_1)--(3_u);

\draw[densely dashdotdotted, ->, line width=0.3mm, >=latex, shorten <= 0.2cm, shorten >= 0.1cm](2_2)--(2_u);
\draw[densely dotted, ->, line width=0.3mm, >=latex, shorten <= 0.2cm, shorten >= 0.1cm](1_2)--(1_u);
\draw[densely dotted, ->, line width=0.3mm, >=latex, shorten <= 0.2cm, shorten >= 0.1cm](1_2)--(3_u);

\draw[->, line width=0.3mm, >=latex, shorten <= 0.2cm, shorten >= 0.1cm](2_3)--(3_u);
\draw[->, line width=0.3mm, >=latex, shorten <= 0.2cm, shorten >= 0.1cm](3_3)--(3_u);
\draw[->, line width=0.3mm, >=latex, shorten <= 0.2cm, shorten >= 0.1cm](1_3)--(1_u);
\draw[->, line width=0.3mm, >=latex, shorten <= 0.2cm, shorten >= 0.1cm](1_3)--(2_u);
\draw[->, line width=0.3mm, >=latex, shorten <= 0.2cm, shorten >= 0.1cm](4_3)--(1_u);
\draw[->, line width=0.3mm, >=latex, shorten <= 0.2cm, shorten >= 0.1cm](4_3)--(2_u);

\draw[->, line width=0.3mm, >=latex, shorten <= 0.2cm, shorten >= 0.1cm](2_4)--(3_u);
\draw[->, line width=0.3mm, >=latex, shorten <= 0.2cm, shorten >= 0.1cm](3_4)--(3_u);
\draw[->, line width=0.3mm, >=latex, shorten <= 0.2cm, shorten >= 0.1cm](1_4)--(1_u);
\draw[->, line width=0.3mm, >=latex, shorten <= 0.2cm, shorten >= 0.1cm](1_4)--(2_u);
\draw[->, line width=0.3mm, >=latex, shorten <= 0.2cm, shorten >= 0.1cm](4_4)--(1_u);
\draw[->, line width=0.3mm, >=latex, shorten <= 0.2cm, shorten >= 0.1cm](4_4)--(2_u);
\draw[->, line width=0.3mm, >=latex, shorten <= 0.2cm, shorten >= 0.1cm](1_5)--(3_u);
\draw[->, line width=0.3mm, >=latex, shorten <= 0.2cm, shorten >= 0.1cm](2_5)--(3_u);
\draw[->, line width=0.3mm, >=latex, shorten <= 0.2cm, shorten >= 0.1cm](1_6)--(3_u);
\draw[->, line width=0.3mm, >=latex, shorten <= 0.2cm, shorten >= 0.1cm](2_6)--(3_u);

\end{tikzpicture}
\captionsetup{justification=centering}
{
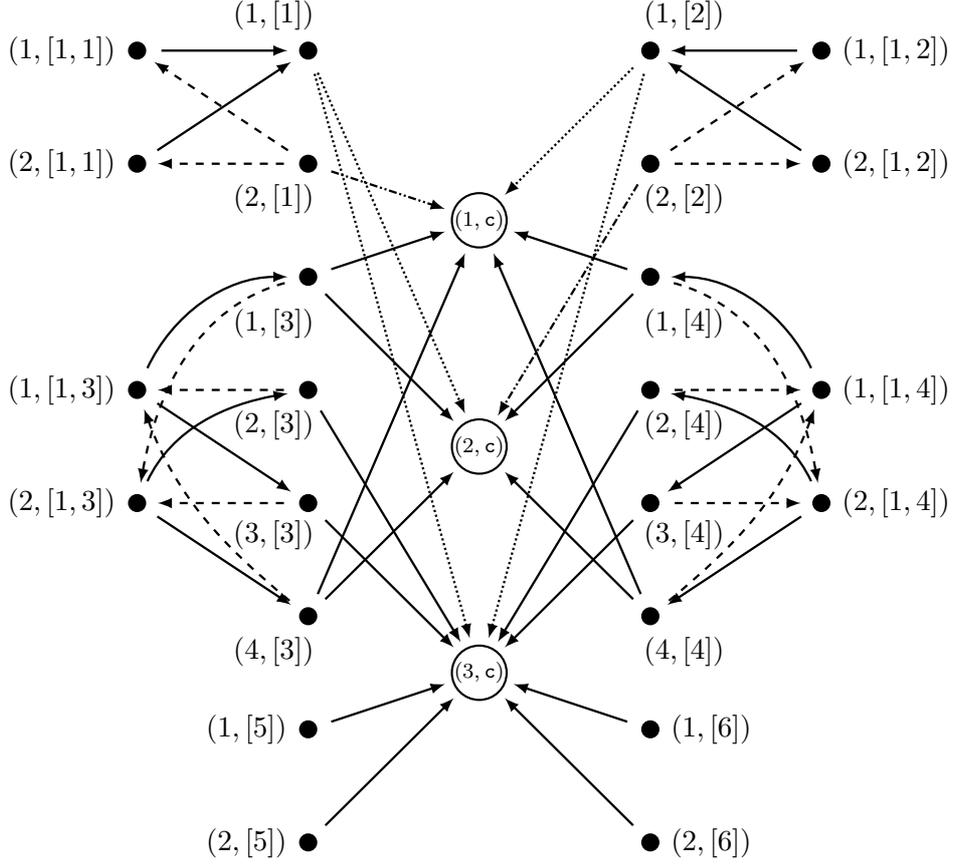
\captionof{figure}{Orientation $D$ for $\mathcal{H}$ for $s=3$, $A_2=\{[1],[2]\}$, $A_4=\{[3],[4]\}$, $E=\{[5],[6]\}$.}\label{figC6.4.17}}
\end{center}

\begin{cor}\label{corC6.4.12}
If $s\ge 3$ is odd and $|A_{\le 3}|\le{{s}\choose{\lceil{s/2}\rceil}}-1$ for a $\mathcal{T}$, then $\mathcal{T}\in \mathscr{C}_0$.
\end{cor}
\noindent\textit{Proof}: Note that every vertex lies in a directed $C_4$ for the orientation $D$ defined in Proposition \ref{ppnC6.4.11}, $\bar{d}(\mathcal{T})\le \max \{4, d(D)\}$ by Lemma \ref{lemC6.1.3}. This implies $\bar{d}(\mathcal{T})=4$.
\qed

\indent\par Corollary \ref{corC6.4.12} is an improvement of Corollary \ref{corC6.3.8}(i) for odd $s$ as the count $|A_{\le 3}|$ now excludes $|A_{\ge 4}|$, compared to the previous $|A_{\ge 2}|$.

\begin{ppn}\label{ppnC6.4.13}
Suppose $s\ge 3$ is odd, $A_2\neq\emptyset$, $A_3\neq\emptyset$ and $A_{\ge 4}\neq\emptyset$ for a $\mathcal{T}$. Then,
\begin{equation}
\mathcal{T}\in \mathscr{C}_0\iff \left\{
  \begin{array}{@{}ll@{}}
    |A_2|\le{{s}\choose{\lceil{s/2}\rceil}}-2, & \text{if}\ |A_3|=1, \nonumber\\
    2|A_2|+|A_3|\le 2{{s}\choose{\lceil{s/2}\rceil}}-2, & \text{if}\ |A_3|\ge 2. \nonumber
  \end{array}\right.
\end{equation}
\end{ppn}
\noindent\textit{Proof}: $(\Rightarrow)$ Since $\mathcal{T}\in \mathscr{C}_0$, there exists an orientation $D$ of $\mathcal{T}$, where $d(D)=4$. As $A_2\neq\emptyset$ and $A_3\neq\emptyset$, we assume (\ref{eqC6.3.1})-(\ref{eqC6.3.5}) here.
\\
\\Case 1. $|A_3|=1$.
\indent\par Partition $A_{\ge 4}$ into $A^O_{\ge 4}$, $A^I_{\ge 4}$, and $A^N_{\ge 4}$ as follows. Let $A^N_{\ge 4}=\{[i]\in A_{\ge 4} \mid |O((p,[\alpha,i]))|\ge 2$ and $|I((p,[\alpha,i]))|\ge 2$ for all $1\le \alpha \le deg_T([i])-1$ and all $1\le p\le s_{[\alpha,i]}\}$, and $A^O_{\ge 4}=\{[i]\in A_{\ge 4}\mid |O((p,[\alpha,i]))|=1 \text{ for some } 1\le \alpha \le deg_T([i])-1 \text{ and some }1\le p\le s_{[\alpha,i]} \}$. Furthermore, let $A^I_{\ge 4}=A_{\ge 4}-A^O_{\ge 4}\cup A^N_{\ge 4}$, i.e., for every $[i]\in A^I_{\ge 4}$, there exist some $1\le \alpha \le deg_T([i])-1$ and some $1\le p\le s_{[\alpha,i]}$ such that $|I((p,[\alpha,i]))|=1$.
\indent\par Without loss of generality, we assume 
\begin{align*}
(\mathbb{N}_{s_{[i]}},[i])-\{(4,[i])\}\rightarrow &(1,[1,i])\rightarrow (4,[i])\text{ if }[i]\in A^O_{\ge 4},\\
(4,[i])\rightarrow &(1,[1,i])\rightarrow (\mathbb{N}_{s_{[i]}},[i])-\{(4,[i])\}\text{ if }[i]\in A^I_{\ge 4},\\
\text{and }\{(1,[i]), (2,[i])\}\rightarrow &(1,[1,i])\rightarrow \{(3,[i]),(4,[i])\}\text{ if }[i]\in A^N_{\ge 4}.
\end{align*}
Also, we let
\begin{align*} 
B^O_{\ge 4}=\{O^\mathtt{c}((4,[i]))\mid [i]\in A^O_{\ge 4}\} \text{ and }B^I_{\ge 4}=\{I^\mathtt{c}((4,[i]))\mid [i]\in A^I_{\ge 4}\}.
\end{align*}
Note that each of $B^O_2\cup B^O_3\cup B^O_{\ge 4}$ and $B^I_2\cup B^I_3\cup B^I_{\ge 4}$ is an antichain by Lemmas \ref{lemC6.2.15} and \ref{lemC6.2.16} respectively. Hence, $|A_2|+|A^O_3|+|A^O_{\ge 4}|=|B^O_2\cup B^O_3\cup B^O_{\ge 4}|\le {{s}\choose{\lfloor{s/2}\rfloor}}$ and $|A_2|+|A^I_3|+|A^I_{\ge 4}|=|B^I_2\cup B^I_3\cup B^I_{\ge 4}|\le {{s}\choose{\lfloor{s/2}\rfloor}}$ by Sperner's theorem.
\\
\\Subcase 1.1. $|A^O_3|=1$. Note that $|A_3|=1$ implies $|A^I_3|=|B^I_3|=0$.
\indent\par If $|A^O_{\ge 4}|>0$, then $|A_2|+|A^O_3|+|A^O_{\ge 4}|\le {{s}\choose{\lfloor{s/2}\rfloor}}$ implies $|A_2|\le {{s}\choose{\lfloor{s/2}\rfloor}}-2$. Hence, suppose $|A^O_{\ge 4}|=0$. 
\\
\\Subcase 1.1.1. $|O^\mathtt{c}((q,[j]))|\ge\lceil\frac{s}{2}\rceil$ for some $[j]\in A^I_{\ge 4}\cup A^N_{\ge 4}$ and some $1\le q\le s_{[j]}$. 
\noindent\par For any $[i]\in A_2\cup A^O_3$, $d_D((1,[1,i]),(q,[j]))=3$ implies $X\cap I^\mathtt{c}((q,[j]))\neq \emptyset$ for all $X\in B^O_2\cup B^O_3$. Hence, by Lih's theorem, $|A_2|+|A^O_3|=|B^O_2\cup B^O_3|\le {{s}\choose{\lceil{s/2}\rceil}}-{{s-|I^\mathtt{c}((1,[i]))|}\choose{\lceil{s/2}\rceil}}\le {{s}\choose{\lceil{s/2}\rceil}}-{{\lceil{s/2}\rceil}\choose{\lceil{s/2}\rceil}}\le {{s}\choose{\lceil{s/2}\rceil}}-1$. Since $|A^O_3|=1$, it follows that $|A_2|\le{{s}\choose{\lfloor{s/2}\rfloor}}-2$.
\\
\\Subcase 1.1.2. $|O^\mathtt{c}((q,[j]))|<\lfloor\frac{s}{2}\rfloor$ for some $[j]\in A^I_{\ge 4}\cup A^N_{\ge 4}$ and some $1\le q\le s_{[j]}$. 
\noindent\par For any $[i]\in A_2$, $d_D((q,[j]),(1,[1,i]))=3$ implies $O^\mathtt{c}((q,[j]))\cap X \neq \emptyset$ for all $X\in B^I_2$. Hence, by Lih's theorem, $|A_2|=|B^I_2|\le {{s}\choose{\lceil{s/2}\rceil}}-{{s-|O^\mathtt{c}((q,[j]))|}\choose{\lceil{s/2}\rceil}}\le {{s}\choose{\lceil{s/2}\rceil}}-{{\lceil{s/2}\rceil+1}\choose{\lceil{s/2}\rceil}}={{s}\choose{\lceil{s/2}\rceil}}-(\lceil{\frac{s}{2}}\rceil+1)\le{{s}\choose{\lceil{s/2}\rceil}}-3$.
\\
\\Subcase 1.1.3. $|O^\mathtt{c}((q,[j]))|=\lfloor\frac{s}{2}\rfloor$ for all $[j]\in A^I_{\ge 4}\cup A^N_{\ge 4}$ and all $1\le q\le s_{[j]}$. 
\indent\par For any $[i]\in A_2$, $d_D((q,[j]),(1,[1,i])=3$ implies $X\cap O^\mathtt{c}((q,[j])) \neq \emptyset$ for all $X\in B^I_2$. Hence, by Lih's theorem, $|A_2|=|B^I_2|\le {{s}\choose{\lceil{s/2}\rceil}}-{{s-|O^\mathtt{c}((q,[j]))|}\choose{\lceil{s/2}\rceil}}={{s}\choose{\lceil{s/2}\rceil}}-{{\lceil{s/2}\rceil}\choose{\lceil{s/2}\rceil}}={{s}\choose{\lceil{s/2}\rceil}}-1$. By Griggs' theorem, $|A_2|={{s}\choose{\lceil{s/2}\rceil}}-1$ if and only if $B^I_2$ consists of only $\lceil\frac{s}{2}\rceil$-sets.
\indent\par If $|A_2|={{s}\choose{\lceil{s/2}\rceil}}-1$, then it must follow that $O^\mathtt{c}((p,[j]))=O^\mathtt{c}((q,[j]))$ for all $[j]\in A^I_{\ge 4}\cup A^N_{\ge 4}$ and all $1\le p,q\le s_{[j]}$. Otherwise, $d_D((r,[j]),(1,[1,i]))=3$ for all $1\le r\le s_{[j]}$ and all $[i]\in A_2$ implies $\bar{X}\neq O^\mathtt{c}((r,[j]))$ for all $X\in B^I_2$, i.e., $B^I_2\subseteq {{(\mathbb{N}_s,\mathtt{c})}\choose{\lceil{s/2}\rceil}}-\{O^\mathtt{c}((p,[j])), O^\mathtt{c}((q,[j]))\}$. Hence, $|A_2|=|B^I_2|\le{{s}\choose{\lceil{s/2}\rceil}}-2$, a contradiction.
\indent\par Now, $B^O_2\cup B^O_3\cup \{O^\mathtt{c}((3,[j]))\mid [j]\in A^I_{\ge 4}\cup A^N_{\ge 4}\}$ is an antichain by Lemma \ref{lemC6.2.15}. So, $|A_2|+|A_3|+|A_4|=|B^O_2|+|B^O_3|+|A^I_{\ge 4}|+|A^N_{\ge 4}|\le {{s}\choose{\lceil{s/2}\rceil}}$. Since $|A_3|=1$ and $|A_4|\ge 1$, it follows that $|A_2|\le {{s}\choose{\lceil{s/2}\rceil}}-2$.
\\
\\Subcase 1.2. $|A^I_3|=1$. Note that $|A_3|=1$ implies $|A^O_3|=|B^O_3|=0$.
\indent\par This follows from an argument similar to Subcase 1.1.
\\
\\Case 2. $|A_3|\ge2$.
\indent\par For any $[i],[j] \in A_2\cup A_3$, $i\neq j$, $[k]\in A_{\ge 4}$, $1\le\alpha\le deg_T([i])-1$, $1\le\gamma\le deg_T([k])-1$, $1\le x\le s_{[\alpha,i]}$, and $1\le z\le s_{[\gamma,k]}$, $1\le q\le 3$ (where applicable), $1\le r\le 4$, observe in $D$ that the vertices $(r, [k])$ and $(z,[\gamma,k])$ do not lie on any shortest path between $(x,[\alpha,i])$ and $(q,[j])$. By the proof in Proposition \ref{ppnC6.4.3}, we have $2|A_2|+|A_3|\le 2{{s}\choose{\lceil{s/2}\rceil}}-1$, where equality is possible only as in Subcases 2.3.6.3 and 2.3.6.4.
\indent\par Suppose $2|A_2|+|A_3|=2{{s}\choose{\lceil{s/2}\rceil}}-1$ holds as in Subcase 2.3.6.3 with $[i^*]$ as given, and recall (\ref{eqC6.4.11}). In particular, $B^I_2\cup B^I_3={{(\mathbb{N}_s,\mathtt{c})}\choose{\lceil{s/2}\rceil}}-\{I^\mathtt{c}((1,[i^*]))\}$ by Claims 6A and 7A.
\indent\par Let $[i]\in A_{\ge 4}$. If there exists some $1\le p\le s_{[i]}$ such that $|O^\mathtt{c}((p,[i]))|\ge \lceil\frac{s}{2}\rceil$, then $X\subseteq O^\mathtt{c}((p,[i]))$ for some $X\in B^O_2\cup B^O_3$. This implies that $d_D((1,[1,j]),(p,[i]))>4$ for some $[j]\in A_2\cup A^O_3$, a contradiction. If there exists some $1\le p\le s_{[i]}$ such that $|O^\mathtt{c}((p,[i]))|<\lfloor\frac{s}{2}\rfloor$, or $|O^\mathtt{c}((p,[i]))|=\lfloor\frac{s}{2}\rfloor$ and $O^\mathtt{c}((p,[i]))\neq O^\mathtt{c}((1,[i^*]))$, then $O^\mathtt{c}((p,[i]))\subseteq \bar{X}$ for some $X\in B^I_2\cup B^I_3$. It follows that $d_D((p,[i]),(1,[1,j]))>4$ for some $[j]\in A_2\cup A^I_3$, a contradiction. Thus, it remains that $O^\mathtt{c}((p,[i]))=O^\mathtt{c}((1,[i^*]))$ for all $[i]\in A_{\ge 4}$ and all $1\le p\le s_{[i]}$. By Lemma \ref{lemC6.2.15}, $Q=\{O^\mathtt{c}((4,[i]))\mid [i]\in A^O_{\ge 4}\}\cup \{O^\mathtt{c}((3,[i]))\mid [i]\in A^I_{\ge 4}\cup A^N_{\ge 4}\}\cup B^O_2\cup B^O_3$ is an antichain. However, this contradicts Sperner's theorem as $|Q|>|B^O_2|+|B^O_3|={{s}\choose{\lceil{s/2}\rceil}}$.
\indent\par Suppose $2|A_2|+|A_3|=2{{s}\choose{\lceil{s/2}\rceil}}-1$ holds as in Subcase 2.3.6.4 with $[i^*]$ as given. Analogous to (\ref{eqC6.4.11}) we have $|B^I_2|+|B^I_3|= {{s}\choose{\lceil{s/2}\rceil}}$ and $|B^O_2|+|B^O_3|= {{s}\choose{\lceil{s/2}\rceil}}-1$. In particular, $B^O_2\cup B^O_3={{(\mathbb{N}_s,\mathtt{c})}\choose{\lceil{s/2}\rceil}}-\{O^\mathtt{c}((1,[i^*]))\}$ by Claims 6B and 7B.
\indent\par Let $[i]\in A_{\ge 4}$. If there exists some $1\le p\le s_{[i]}$ such that $|I^\mathtt{c}((p,[i]))|\ge \lceil\frac{s}{2}\rceil$, then $X\subseteq I^\mathtt{c}((p,[i]))$ for some $X\in B^I_2\cup B^I_3$. This implies $d_D((p,[i]),(1,[1,j]))>4$ for some $[j]\in A_2\cup A^I_3$. If there exists some $1\le p\le s_{[i]}$ such that $|I^\mathtt{c}((p,[i]))|<\lfloor\frac{s}{2}\rfloor$, or $|I^\mathtt{c}((p,[i]))|=\lfloor\frac{s}{2}\rfloor$ and $I^\mathtt{c}((p,[i]))\neq I^\mathtt{c}((1,[i^*]))$, then $I^\mathtt{c}((p,[i]))\subseteq \bar{X}$ for some $X\in B^O_2\cup B^O_3$. It follows that $d_D((1,[1,j]),(p,[i]))>4$ for some $[j]\in A_2\cup A^O_3$, a contradiction. Thus, it remains that $I^\mathtt{c}((p,[i]))=I^\mathtt{c}((1,[i^*]))$ for all $[i]\in A_{\ge 4}$ and all $1\le p\le s_{[i]}$. By Lemma \ref{lemC6.2.16}, $Q=\{I^\mathtt{c}((4,[i]))\mid [i]\in A^I_{\ge 4}\}\cup \{I^\mathtt{c}((1,[i]))\mid [i]\in A^O_{\ge 4}\cup A^N_{\ge 4}\}\cup B^I_2\cup B^I_3$ is an antichain. However, this contradicts Sperner's theorem as $|Q|>|B^I_2|+|B^I_3|={{s}\choose{\lceil{s/2}\rceil}}$.
\noindent\par So, we conclude that $2|A_2|+|A_3|=2{{s}\choose{\lceil{s/2}\rceil}}-1$ is impossible in both subcases and $2|A_2|+|A_3|\le 2{{s}\choose{\lceil{s/2}\rceil}}-2$.
\\
\\$(\Leftarrow)$ By Corollary \ref{corC6.4.12}, $\mathcal{T}\in \mathscr{C}_0$ if $|A_2|\le{{s}\choose{\lceil{s/2}\rceil}}-2$ and $|A_3|=1$, or $|A_2|+|A_3|\le {{s}\choose{\lceil{s/2}\rceil}}-1$ and $|A_3|\ge 2$. Hence, we assume $|A_2|+|A_3|\ge{{s}\choose{\lceil{s/2}\rceil}}$, on top of the hypothesis that $2|A_2|+|A_3|\le 2{{s}\choose{\lceil{s/2}\rceil}}-2$ and $|A_3|\ge 2$. Furthermore, assume without loss of generality that $A_2=\{[i]\mid i\in\mathbb{N}_{|A_2|}\}$, and $A_3=\{[i]\mid i\in\mathbb{N}_{|A_2|+|A_3|}-\mathbb{N}_{|A_2|}\}$.
\indent\par Let $\mathcal{H}=T(t_1,t_2,\ldots, t_n)$ be the subgraph of $\mathcal{T}$, where $t_\mathtt{c}=s$, $t_{[i]}=3$ for all $[i]\in \mathcal{T}(A_3)$, $t_{[j]}=4$ for all $[j]\in \mathcal{T}(A_{\ge 4})$ and $t_v=2$ otherwise. We will use $A_j$ for $\mathcal{H}(A_j)$ for the remainder of this proof. Define an orientation $D$ of $\mathcal{H}$ as follows.
\begin{align}
& (2,[i])\rightarrow \{(1,[\alpha,i]),(2,[\alpha,i])\}\rightarrow (1,[i]),\text{ and}\label{eqC6.4.43}\\
& \bar{\lambda}_{i+1}\rightarrow (1,[i])\rightarrow \lambda_{i+1}\rightarrow (2,[i]) \rightarrow \bar{\lambda}_{i+1}\label{eqC6.4.44}
\end{align}
for all $1\le i\le |A_2|$ and all $1\le \alpha\le deg_T([i])-1$, i.e., the $\lceil\frac{s}{2}\rceil$-sets $\lambda_2,\lambda_3,\ldots,\lambda_{|A_2|+1}$ are used as `in-sets' (`out-sets' resp.) to construct $B^I_2$ ($B^O_2$ resp.).
\begin{align}
& (3,[j])\rightarrow \{(1,[\beta,j]),(2,[\beta,j])\}\rightarrow \{(1,[j]),(2,[j])\},\label{eqC6.4.45}\\
& \bar{\lambda}_1 \rightarrow (1,[j])\rightarrow \lambda_1\rightarrow (2,[j])\rightarrow \bar{\lambda}_1,\text{ and}\label{eqC6.4.46}\\
&\lambda_{j+1}\rightarrow (3,[j]) \rightarrow \bar{\lambda}_{j+1}\label{eqC6.4.47}
\end{align}
for all $|A_2|+1\le j\le {{s}\choose{\lceil{s/2}\rceil}}-1$ and all $1\le \beta\le deg_T([j])-1$, i.e., the $\lceil\frac{s}{2}\rceil$-sets $\lambda_{|A_2|+2},\lambda_{|A_2|+3},\ldots,\lambda_{{{s}\choose{\lceil{s/2}\rceil}}}$ are used as `in-sets' to construct $B^I_3$.
\begin{align}
& \{(1,[k]),(2,[k])\}\rightarrow \{(1,[\gamma,k]),(2,[\gamma,k])\}\rightarrow (3,[k]),\label{eqC6.4.48}\\
& \lambda_1 \rightarrow (1,[k])\rightarrow \bar{\lambda}_1\rightarrow (2,[k])\rightarrow \lambda_1,\text{ and}\label{eqC6.4.49}\\
&\bar{\lambda}_{k-{{s}\choose{\lceil{s/2}\rceil}}+|A_2|+2}\rightarrow (3,[k]) \rightarrow \lambda_{k-{{s}\choose{\lceil{s/2}\rceil}}+|A_2|+2}\label{eqC6.4.50}
\end{align}
for all ${{s}\choose{\lceil{s/2}\rceil}}\le k\le |A_2|+|A_3|$ and all $1\le\gamma\le deg_T([k])-1$, i.e., the $\lceil\frac{s}{2}\rceil$-sets $\lambda_{|A_2|+2},\lambda_{|A_2|+3},\ldots,\lambda_{2|A_2|+|A_3|+2-{{s}\choose{\lceil{s/2}\rceil}}}$ are used as `out-sets' to construct $B^O_3$.
\begin{align}
&(2,[\tau,l])\rightarrow \{(2,[l]),(4,[l])\}\rightarrow (1,[\tau,l])\rightarrow \{(1,[l]), (3,[l])\}\rightarrow (2,[\tau,l]),\text{ and}\label{eqC6.4.51}\\
&\bar{\lambda}_1 \rightarrow \{(1,[l]), (4,[l])\}\rightarrow \lambda_1 \rightarrow\{(2,[l]), (3,[l])\}\rightarrow \bar{\lambda}_1\label{eqC6.4.52}
\end{align}
for all $[l]\in A_4$ and all $1\le\tau\le deg_T([l])-1$.
\begin{align}
\lambda_1 \rightarrow \{(1,[m]), (2,[m])\} \rightarrow \bar{\lambda}_1 \label{eqC6.4.53}
\end{align}
for all $[m]\in E$. (See Figure \ref{figC6.4.18} for $D$ when $s=3$.)
\\
\\Claim: $d_{D}(v,w)\le 4$ for all $v,w\in V(D)$.
\\
\\Case 1. $v \in \{(1,[\alpha,i]),(2,[\alpha,i]),(1,[i]),(2,[i])\}$ and $w\in \{(1,[\beta,j]),(2,[\beta,j]),(1,[j]),$ $(2,[j])\}$ for each $[i],[j]\in A_2$, each $1\le\alpha\le deg_T([i])-1$, and each $1\le\beta\le deg_T([j])-1$.
\noindent\par Since the orientation defined for $A_2$ (see (\ref{eqC6.4.43})-(\ref{eqC6.4.44})) is similar to that in Proposition \ref{ppnC6.4.11} (see (\ref{eqC6.4.38})-(\ref{eqC6.4.39})), this case follows from Cases 1.1 and 2 of Proposition \ref{ppnC6.4.11}, except possibly for $v=(p,[i])$ and $w=(q, [j])$, which will be covered in Case 13 later.
\\
\\Case 2. $v \in \{(1,[\alpha,i]),(2,[\alpha,i]),(1,[i]),(2,[i]),(3,[i])\}$ and $w\in \{(1,[\beta,j]),(2,[\beta,j]),$ $(1,[j]),(2,[j]),(3,[j])\}$ for each $[i],[j]\in A_3$, each $1\le\alpha\le deg_T([i])-1$, and each $1\le\beta\le deg_T([j])-1$.
\indent\par Since the orientation defined for $A_3$ (see (\ref{eqC6.4.45})-(\ref{eqC6.4.50})) is similar to that in Proposition \ref{ppnC6.4.1} (see (\ref{eqC6.4.1})-(\ref{eqC6.4.6})), this case follows from Cases 1.1-1.2, 2-3 and 5 of Proposition \ref{ppnC6.4.1}, except possibly for $v=(p,[i])$ and $w=(q, [j])$, which will be covered in Case 13 later.
\\
\\Case 3. $v \in \{(1,[\alpha,i]),(2,[\alpha,i]),(1,[i]),(2,[i]),(3,[i]),(4,[i])\}$ and $w\in \{(1,[\beta,j]),(2,[\beta,j]),$ $(1,[j]),(2,[j]),(3,[j]),(4,[j])\}$ for each $[i],[j]\in A_4$, each $1\le\alpha\le deg_T([i])-1$, and each $1\le\beta\le deg_T([j])-1$.
\indent\par Since the orientation defined for $A_4$ (see (\ref{eqC6.4.51})-(\ref{eqC6.4.52})) is similar to that in Proposition \ref{ppnC6.4.11} (see (\ref{eqC6.4.40})-(\ref{eqC6.4.41})), this case follows from Cases 1.2 and 3 of Proposition \ref{ppnC6.4.11}, except possibly for $v=(p,[i])$ and $w=(q, [j])$, which will be covered in Case 13 later.
\\
\\Case 4. $v \in \{(1,[\alpha,i]),(2,[\alpha,i]),(1,[i]),(2,[i])\}$ and $w\in \{(1,[\beta,j]),(2,[\beta,j]),(1,[j]),$ $(2,[j]),(3,[j])\}$ for each $[i]\in A_2$, each $[j]\in A_3$, each $1\le\alpha\le deg_T([i])-1$, and each $1\le\beta\le deg_T([j])-1$.
\indent\par Since the orientation defined for $A_2$ and $A_3$ (see (\ref{eqC6.4.43})-(\ref{eqC6.4.50})) are similar to those of $D_2$ in Proposition \ref{ppnC6.4.3} (see (\ref{eqC6.4.18})-(\ref{eqC6.4.25})), this case follows from Cases 2.4 and 2.5 of Proposition \ref{ppnC6.4.3}, except possibly for $v=(p,[i])$ and $w=(q, [j])$, which will be covered in Case 13 later.
\\
\\Case 5. $v \in \{(1,[\alpha,i]),(2,[\alpha,i]),(1,[i]),(2,[i])\}$ and $w\in \{(1,[\beta,j]),(2,[\beta,j]),(1,[j]),$ $(2,[j]),(3,[j]),(4,[j])\}$ for each $[i]\in A_2$, each $[j]\in A_4$, each $1\le\alpha\le deg_T([i])-1$, and each $1\le\beta\le deg_T([j])-1$.
\indent\par Since the orientation defined for $A_2$ and $A_4$ (see (\ref{eqC6.4.43})-(\ref{eqC6.4.44}) and (\ref{eqC6.4.51})-(\ref{eqC6.4.52})) are similar to those in Proposition \ref{ppnC6.4.11} (see (\ref{eqC6.4.38})-(\ref{eqC6.4.41})), this case follows from Case 4 of Proposition \ref{ppnC6.4.11}, except possibly for $v=(p,[i])$ and $w=(q, [j])$, which will be covered in Case 13 later.
\\
\\Case 6. $v \in \{(1,[\alpha,i]),(2,[\alpha,i]),(1,[i]),(2,[i]),(3,[i])\}$ and $w\in \{(1,[\beta,j]),(2,[\beta,j]),(1,[j]),$ $(2,[j]),(3,[j]),(4,[j])\}$ for each $[i]\in A_3$, each $[j]\in A_4$, each $1\le\alpha\le deg_T([i])-1$, and each $1\le\beta\le deg_T([j])-1$.
\indent\par Since the orientation defined for $A_3$ and $A_4$ (see (\ref{eqC6.4.45})-(\ref{eqC6.4.52})) are similar to those in Proposition \ref{ppnC6.4.1} (see (\ref{eqC6.4.1})-(\ref{eqC6.4.8})), this case follows from Cases 6 and 7 of Proposition \ref{ppnC6.4.1}, except possibly for $v=(p,[i])$ and $w=(q, [j])$, which will be covered in Case 13 later.
\\
\\Case 7. For each $[i]\in A_2$, each $1\le\alpha\le deg_T([i])-1$, and each $[j]\in E$,
\\(i) $v=(p,[\alpha,i]), w=(q,[j])$ for each $p=1,2$ and $q=1,2$.
\\(ii) $v=(q,[j]), w=(p,[\alpha,i])$ for each $p=1,2$ and $q=1,2$.
\indent\par Let $[k]\in A_4$. Since $\lambda_1\rightarrow \{(q,[j]),(2,[k])\}\rightarrow \bar{\lambda}_1$ by (\ref{eqC6.4.52})-(\ref{eqC6.4.53}), this case follows from Case 5.
\\
\\Case 8. For each $[i]\in A_3$, each $1\le\alpha\le deg_T([i])-1$, and each $[j]\in E$,
\\(i) $v=(p,[\alpha,i]), w=(q,[j])$ for each $p=1,2$ and $q=1,2$.
\\(ii) $v=(q,[j]), w=(p,[\alpha,i])$ for each $p=1,2$ and $q=1,2$.
\indent\par Let $[k]\in A_4$. Since $\lambda_1\rightarrow \{(q,[j]),(2,[k])\}\rightarrow \bar{\lambda}_1$ by (\ref{eqC6.4.52})-(\ref{eqC6.4.53}), this case follows from  Case 6.
\\
\\Case 9. For each $[i]\in A_4$, each $1\le\alpha\le deg_T([i])-1$, and each $[j]\in E$,
\\(i) $v=(p,[\alpha,i]), w=(q,[j])$ for each $p=1,2$ and $q=1,2$.
\\(ii) $v=(q,[j]), w=(p,[\alpha,i])$ for each $p=1,2$ and $q=1,2$.
\indent\par Since $\lambda_1\rightarrow \{(q,[j]),(2,[|A_2|+1])\}\rightarrow \bar{\lambda}_1$ by (\ref{eqC6.4.46}) and (\ref{eqC6.4.53}), this case follows from Case 6.
\\
\\Case 10. $v=(r_1,\mathtt{c})$ and $w=(r_2,\mathtt{c})$ for $r_1\neq r_2$ and $1\le r_1, r_2\le s$.
\indent\par Here, we want to prove a stronger claim, $d_{D}((r_1,\mathtt{c}), (r_2,\mathtt{c}))=2$. Let $x_1=(2,[l])$ for some $[l] \in A_4$, $x_{k+1}=(2,[k])$ for $1\le k\le |A_2|$, and $x_{k+1}=(3,[k])$ for $|A_2|+1\le k\le {{s}\choose{\lceil{s/2}\rceil}}-1$. Observe from (\ref{eqC6.4.52}), (\ref{eqC6.4.44}) and (\ref{eqC6.4.47}) that $\lambda_k\rightarrow x_k\rightarrow \bar{\lambda}_k$ for all $1\le k\le s$ and the subgraph induced by $V_1=(\mathbb{N}_s,\mathtt{c})$ and $V_2=\{x_i\mid 1\le i \le s\}$ is a complete bipartite graph $K(V_1,V_2)$. By Lemma \ref{lemC6.2.19}, $d_{D}((r_1,\mathtt{c}), (r_2,\mathtt{c}))=2$.
\\
\\Case 11. $v\in \{(1,[i]), (2,[i]), (3,[i]), (4,[i]), (1,[\alpha,i]), (2,[\alpha,i])\}$ for each $1\le i\le deg_T(\mathtt{c})$ and $1\le\alpha\le deg_T([i])-1$, and $w=(r,\mathtt{c})$ for $1\le r\le s$.
\indent\par Note that there exists some $1\le k\le s$ such that $d_D(v,(k,\mathtt{c}))\le 2$, and $d_D((k,\mathtt{c}),w)\le 2$ by Case 10. Hence, it follows that $d_D(v,w)\le d_D(v,(k,\mathtt{c}))+d_D((k,\mathtt{c}),w)\le 4$.
\\
\\Case 12. $v=(r,\mathtt{c})$ for $r=1,2,\ldots, s$ and $w\in \{(1,[i]), (2,[i]), (3,[i]), (4,[i]), (1,[\alpha,i]), (2,[\alpha,i])\}$ for each $1\le i\le deg_T(\mathtt{c})$ and $1\le\alpha\le deg_T([i])-1$.
\indent\par Note that there exists some $1\le k\le s$ such that $d_D((k,\mathtt{c}), w)\le 2$, and $d_D(v,(k,\mathtt{c}))\le 2$ by Case 10. Hence, it follows that $d_D(v,w)\le d_D(v,(k,\mathtt{c}))+d_D((k,\mathtt{c}),w)\le 4$.
\\
\\Case 13. $v=(p,[i])$ and $w=(q, [j])$, where $1\le p,q\le 4$ and $1\le i,j\le deg_T(\mathtt{c})$.
\noindent\par This follows from the fact that $|O^\mathtt{c}((p,[i]))|>0$, $|I^\mathtt{c}((q,[j]))|>0$, and $d_{D}((r_1,\mathtt{c}), (r_2,\mathtt{c}))$ $=2$ for any $r_1\neq r_2$ and $1\le r_1, r_2\le s$.
\\
\indent\par Therefore, the claim follows. Since every vertex lies in a directed $C_4$ for $D$ and $d(D)=4$, $\bar{d}(\mathcal{T})\le \max \{4, d(D)\}$ by Lemma \ref{lemC6.1.3}, and thus $\bar{d}(\mathcal{T})=4$.
\qed
\begin{center}
\begin{tikzpicture}[thick,scale=0.7]%
\draw(0,2)node[circle, draw, inner sep=0pt, minimum width=3pt](1_u){\scriptsize $(1,\mathtt{c})$};
\draw(0,-2)node[circle, draw, inner sep=0pt, minimum width=3pt](2_u){\scriptsize $(2,\mathtt{c})$};
\draw(0,-6)node[circle, draw, inner sep=0pt, minimum width=3pt](3_u){\scriptsize $(3,\mathtt{c})$};

\draw(-6,6)node[circle, draw, fill=black!100, inner sep=0pt, minimum width=6pt, label={[] 180:{\small $(1,[1,1])$}}](1_11){};
\draw(-6,4)node[circle, draw, fill=black!100, inner sep=0pt, minimum width=6pt, label={[] 180:{\small $(2,[1,1])$}}](2_11){};

\draw(-3,6)node[circle, draw, fill=black!100, inner sep=0pt, minimum width=6pt, label={[yshift=0cm, xshift=-0.45cm] 90:{\small $(1,[1])$}}](1_1){};
\draw(-3,4)node[circle, draw, fill=black!100, inner sep=0pt, minimum width=6pt, label={[yshift=0cm, xshift=-0.45cm] 270:{\small $(2,[1])$}}](2_1){};

\draw(-6,1)node[circle, draw, fill=black!100, inner sep=0pt, minimum width=6pt, label={[] 180:{\small $(1,[1,2])$}}](1_21){};
\draw(-6,-1)node[circle, draw, fill=black!100, inner sep=0pt, minimum width=6pt, label={[] 180:{\small $(2,[1,2])$}}](2_21){};

\draw(-3,2)node[circle, draw, fill=black!100, inner sep=0pt, minimum width=6pt, label={[yshift=-0.15cm, xshift=-0.45cm] 270:{\small $(1,[2])$}}](1_2){};
\draw(-3,0)node[circle, draw, fill=black!100, inner sep=0pt, minimum width=6pt, label={[yshift=-0.1cm, xshift=-0.45cm] 270:{\small $(2,[2])$}}](2_2){};
\draw(-3,-2)node[circle, draw, fill=black!100, inner sep=0pt, minimum width=6pt, label={[yshift=0cm, xshift=-0.45cm] 270:{\small $(3,[2])$}}](3_2){};

\draw(6,1)node[circle, draw, fill=black!100, inner sep=0pt, minimum width=6pt, label={[] 0:{\small $(1,[1,3])$}}](1_31){};
\draw(6,-1)node[circle, draw, fill=black!100, inner sep=0pt, minimum width=6pt, label={[] 0:{\small $(2,[1,3])$}}](2_31){};

\draw(3,2)node[circle, draw, fill=black!100, inner sep=0pt, minimum width=6pt, label={[yshift=-0.15cm, xshift=0.45cm] 270:{\small $(1,[3])$}}](1_3){};
\draw(3,0)node[circle, draw, fill=black!100, inner sep=0pt, minimum width=6pt, label={[yshift=-0.1cm, xshift=0.45cm] 270:{\small $(2,[3])$}}](2_3){};
\draw(3,-2)node[circle, draw, fill=black!100, inner sep=0pt, minimum width=6pt, label={[yshift=0cm, xshift=0.45cm] 270:{\small $(3,[3])$}}](3_3){};

\draw(-6,-6)node[circle, draw, fill=black!100, inner sep=0pt, minimum width=6pt, label={[] 180:{\small $(1,[1,4])$}}](1_41){};
\draw(-6,-8)node[circle, draw, fill=black!100, inner sep=0pt, minimum width=6pt, label={[] 180:{\small $(2,[1,4])$}}](2_41){};

\draw(-3,-4)node[circle, draw, fill=black!100, inner sep=0pt, minimum width=6pt, label={[yshift=-0.1cm, xshift=-0.45cm] 270:{\small $(1,[4])$}}](1_4){};
\draw(-3,-6)node[circle, draw, fill=black!100, inner sep=0pt, minimum width=6pt, label={[yshift=0cm, xshift=-0.45cm] 270:{\small $(2,[4])$}}](2_4){};
\draw(-3,-8)node[circle, draw, fill=black!100, inner sep=0pt, minimum width=6pt, label={[yshift=0.05cm, xshift=-0.45cm] 270:{\small $(3,[4])$}}](3_4){};
\draw(-3,-10)node[circle, draw, fill=black!100, inner sep=0pt, minimum width=6pt, label={[yshift=0cm, xshift=-0.45cm] 270:{\small $(4,[4])$}}](4_4){};

\draw(6,-6)node[circle, draw, fill=black!100, inner sep=0pt, minimum width=6pt, label={[] 0:{\small $(1,[1,5])$}}](1_51){};
\draw(6,-8)node[circle, draw, fill=black!100, inner sep=0pt, minimum width=6pt, label={[] 0:{\small $(2,[1,5])$}}](2_51){};

\draw(3,-4)node[circle, draw, fill=black!100, inner sep=0pt, minimum width=6pt, label={[yshift=-0.1cm, xshift=0.45cm] 270:{\small $(1,[5])$}}](1_5){};
\draw(3,-6)node[circle, draw, fill=black!100, inner sep=0pt, minimum width=6pt, label={[yshift=0cm, xshift=0.45cm] 270:{\small $(2,[5])$}}](2_5){};
\draw(3,-8)node[circle, draw, fill=black!100, inner sep=0pt, minimum width=6pt, label={[yshift=0.05cm, xshift=0.45cm] 270:{\small $(3,[5])$}}](3_5){};
\draw(3,-10)node[circle, draw, fill=black!100, inner sep=0pt, minimum width=6pt, label={[yshift=0cm, xshift=0.45cm] 270:{\small $(4,[5])$}}](4_5){};

\draw(-3,-12)node[circle, draw, fill=black!100, inner sep=0pt, minimum width=6pt, label={[] 180:{\small $(1,[7])$}}](1_7){};
\draw(-3,-14)node[circle, draw, fill=black!100, inner sep=0pt, minimum width=6pt, label={[] 180:{\small $(2,[7])$}}](2_7){};

\draw(3,-12)node[circle, draw, fill=black!100, inner sep=0pt, minimum width=6pt, label={[] 0:{\small $(1,[6])$}}](1_6){};
\draw(3,-14)node[circle, draw, fill=black!100, inner sep=0pt, minimum width=6pt, label={[] 0:{\small $(2,[6])$}}](2_6){};

\draw[->, line width=0.3mm, >=latex, shorten <= 0.2cm, shorten >= 0.15cm](1_11)--(1_1);
\draw[->, line width=0.3mm, >=latex, shorten <= 0.2cm, shorten >= 0.15cm](2_11)--(1_1);
\draw[dashed,->, line width=0.3mm, >=latex, shorten <= 0.2cm, shorten >= 0.15cm](2_1)--(1_11);
\draw[dashed,->, line width=0.3mm, >=latex, shorten <= 0.2cm, shorten >= 0.15cm](2_1)--(2_11);

\draw[dashed,->, line width=0.3mm, >=latex, shorten <= 0.2cm, shorten >= 0.15cm](3_2) to [out=158, in=275] (1_21);
\draw[dashed,->, line width=0.3mm, >=latex, shorten <= 0.2cm, shorten >= 0.15cm](3_2)--(2_21);
\draw[->, line width=0.3mm, >=latex, shorten <= 0.2cm, shorten >= 0.15cm](1_21)--(1_2);
\draw[->, line width=0.3mm, >=latex, shorten <= 0.2cm, shorten >= 0.15cm](2_21) to [out=85, in=202] (1_2);
\draw[->, line width=0.3mm, >=latex, shorten <= 0.2cm, shorten >= 0.15cm](1_21)--(2_2);
\draw[->, line width=0.3mm, >=latex, shorten <= 0.2cm, shorten >= 0.15cm](2_21)--(2_2);

\draw[dashed,->, line width=0.3mm, >=latex, shorten <= 0.2cm, shorten >= 0.15cm](1_3)--(1_31);
\draw[dashed,->, line width=0.3mm, >=latex, shorten <= 0.2cm, shorten >= 0.15cm](1_3) to [out=338, in=95] (2_31);
\draw[dashed,->, line width=0.3mm, >=latex, shorten <= 0.2cm, shorten >= 0.15cm](2_3)--(1_31);
\draw[dashed,->, line width=0.3mm, >=latex, shorten <= 0.2cm, shorten >= 0.15cm](2_3)--(2_31);
\draw[->, line width=0.3mm, >=latex, shorten <= 0.2cm, shorten >= 0.15cm](1_31) to [out=265, in=22] (3_3);
\draw[->, line width=0.3mm, >=latex, shorten <= 0.2cm, shorten >= 0.15cm](2_31)--(3_3);

\draw[dashed,->, line width=0.3mm, >=latex, shorten <= 0.2cm, shorten >= 0.15cm](2_5)--(1_51);
\draw[dashed,->, line width=0.3mm, >=latex, shorten <= 0.2cm, shorten >= 0.15cm](1_5) to [out=345, in=100] (2_51);
\draw[dashed,->, line width=0.3mm, >=latex, shorten <= 0.2cm, shorten >= 0.15cm](4_5) to [out=40, in=250] (1_51);
\draw[dashed,->, line width=0.3mm, >=latex, shorten <= 0.2cm, shorten >= 0.15cm](3_5)--(2_51);
\draw[->, line width=0.3mm, >=latex, shorten <= 0.2cm, shorten >= 0.15cm](1_51) to [out=115, in=0] (1_5);
\draw[->, line width=0.3mm, >=latex, shorten <= 0.2cm, shorten >= 0.15cm](2_51) to [out=115, in=350] (2_5);
\draw[->, line width=0.3mm, >=latex, shorten <= 0.2cm, shorten >= 0.15cm](1_51)--(3_5);
\draw[->, line width=0.3mm, >=latex, shorten <= 0.2cm, shorten >= 0.15cm](2_51)--(4_5);

\draw[dashed,->, line width=0.3mm, >=latex, shorten <= 0.2cm, shorten >= 0.15cm](2_4)--(1_41);
\draw[dashed,->, line width=0.3mm, >=latex, shorten <= 0.2cm, shorten >= 0.15cm](1_4) to [out=195, in=80] (2_41);
\draw[dashed,->, line width=0.3mm, >=latex, shorten <= 0.2cm, shorten >= 0.15cm](4_4) to [out=140, in=290] (1_41);
\draw[dashed,->, line width=0.3mm, >=latex, shorten <= 0.2cm, shorten >= 0.15cm](3_4)--(2_41);
\draw[->, line width=0.3mm, >=latex, shorten <= 0.2cm, shorten >= 0.15cm](1_41) to [out=65, in=180] (1_4);
\draw[->, line width=0.3mm, >=latex, shorten <= 0.2cm, shorten >= 0.15cm](2_41) to [out=65, in=190] (2_4);
\draw[->, line width=0.3mm, >=latex, shorten <= 0.2cm, shorten >= 0.15cm](1_41)--(3_4);
\draw[->, line width=0.3mm, >=latex, shorten <= 0.2cm, shorten >= 0.15cm](2_41)--(4_4);

\draw[densely dashdotdotted, ->, line width=0.3mm, >=latex, shorten <= 0.2cm, shorten >= 0.1cm](2_1)--(2_u);
\draw[densely dotted, ->, line width=0.3mm, >=latex, shorten <= 0.2cm, shorten >= 0.1cm](1_1)--(1_u);
\draw[densely dotted, ->, line width=0.3mm, >=latex, shorten <= 0.2cm, shorten >= 0.1cm](1_1)--(3_u);

\draw[->, line width=0.3mm, >=latex, shorten <= 0.2cm, shorten >= 0.1cm](1_3)--(3_u);
\draw[->, line width=0.3mm, >=latex, shorten <= 0.2cm, shorten >= 0.1cm](2_3)--(1_u);
\draw[->, line width=0.3mm, >=latex, shorten <= 0.2cm, shorten >= 0.1cm](2_3)--(2_u);
\draw[densely dotted, ->, line width=0.3mm, >=latex, shorten <= 0.2cm, shorten >= 0.1cm](3_3)--(2_u);
\draw[densely dotted, ->, line width=0.3mm, >=latex, shorten <= 0.2cm, shorten >= 0.1cm](3_3)--(3_u);

\draw[->, line width=0.3mm, >=latex, shorten <= 0.2cm, shorten >= 0.1cm](1_2)--(1_u);
\draw[->, line width=0.3mm, >=latex, shorten <= 0.2cm, shorten >= 0.1cm](1_2)--(2_u);
\draw[->, line width=0.3mm, >=latex, shorten <= 0.2cm, shorten >= 0.1cm](2_2)--(3_u);
\draw[densely dashdotdotted, ->, line width=0.3mm, >=latex, shorten <= 0.2cm, shorten >= 0.1cm](3_2)--(1_u);

\draw[->, line width=0.3mm, >=latex, shorten <= 0.2cm, shorten >= 0.1cm](2_4)--(3_u);
\draw[->, line width=0.3mm, >=latex, shorten <= 0.2cm, shorten >= 0.1cm](3_4)--(3_u);
\draw[->, line width=0.3mm, >=latex, shorten <= 0.2cm, shorten >= 0.1cm](1_4)--(1_u);
\draw[->, line width=0.3mm, >=latex, shorten <= 0.2cm, shorten >= 0.1cm](1_4)--(2_u);
\draw[->, line width=0.3mm, >=latex, shorten <= 0.2cm, shorten >= 0.1cm](4_4)--(1_u);
\draw[->, line width=0.3mm, >=latex, shorten <= 0.2cm, shorten >= 0.1cm](4_4)--(2_u);

\draw[->, line width=0.3mm, >=latex, shorten <= 0.2cm, shorten >= 0.1cm](2_5)--(3_u);
\draw[->, line width=0.3mm, >=latex, shorten <= 0.2cm, shorten >= 0.1cm](3_5)--(3_u);
\draw[->, line width=0.3mm, >=latex, shorten <= 0.2cm, shorten >= 0.1cm](1_5)--(1_u);
\draw[->, line width=0.3mm, >=latex, shorten <= 0.2cm, shorten >= 0.1cm](1_5)--(2_u);
\draw[->, line width=0.3mm, >=latex, shorten <= 0.2cm, shorten >= 0.1cm](4_5)--(1_u);
\draw[->, line width=0.3mm, >=latex, shorten <= 0.2cm, shorten >= 0.1cm](4_5)--(2_u);
\draw[->, line width=0.3mm, >=latex, shorten <= 0.2cm, shorten >= 0.1cm](1_6)--(3_u);
\draw[->, line width=0.3mm, >=latex, shorten <= 0.2cm, shorten >= 0.1cm](2_6)--(3_u);
\draw[->, line width=0.3mm, >=latex, shorten <= 0.2cm, shorten >= 0.1cm](1_7)--(3_u);
\draw[->, line width=0.3mm, >=latex, shorten <= 0.2cm, shorten >= 0.1cm](2_7)--(3_u);
\end{tikzpicture}
\captionsetup{justification=centering}
{
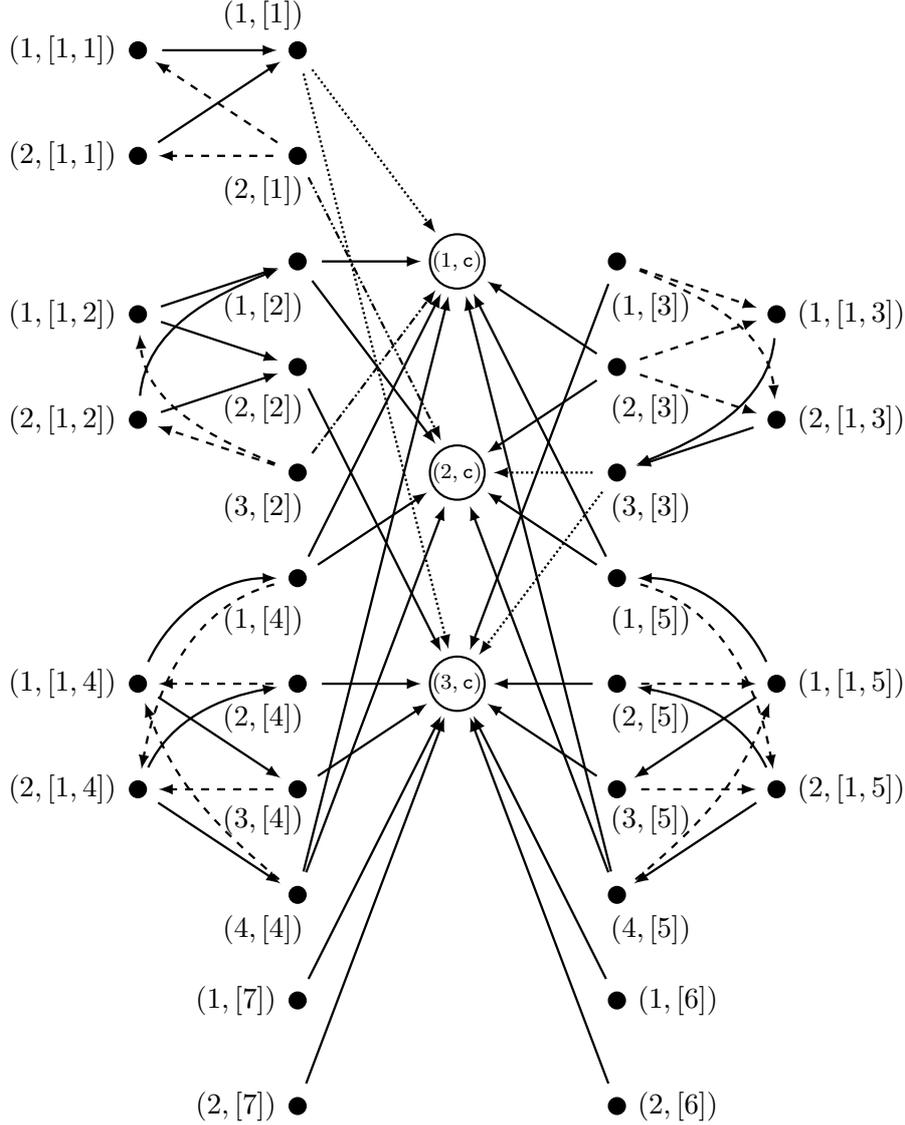
\captionof{figure}{Orientation $D$ for $\mathcal{H}$ for $s=3$,
\\$A_2=\{[1]\}$, $A_3=\{[2],[3]\}$, $A_4=\{[4],[5]\}$, $E=\{[6],[7]\}$.}\label{figC6.4.18}}
\end{center}
\indent\par This concludes the proof of Theorem \ref{thmC6.1.8}. 
\section{Conclusion}
In this paper, we almost completely characterise the case of even $s$ and give a complete characterisation for the case of odd $s\ge 3$. With the current approach of searching for optimal orientation(s) in tree vertex-multiplications, the complexity and quantity of the subcases increase sharply when the subsets $O^\mathtt{c}((p,[i]))$ are of `middle' size ($\lfloor\frac{s}{2}\rfloor$ or $\lceil\frac{s}{2}\rceil$). For instance in Proposition \ref{ppnC6.4.3}, it is relatively easy to settle Subcases 2.1 and 2.2 but Subcase 2.3 is rather involved. Furthermore, the even case (see Proposition \ref{ppnC6.3.12}) illustrates a similar yet more complicated situation. It seems that a new approach may be needed to cut through this entanglement. Since this paper focuses on trees of diameter 4, we end off by proposing the following problem.
\begin{prob}
For trees $T$ with $d(T)=3$, characterise the tree vertex-multiplications $T(s_1,s_2,\ldots,s_n)$ that belong to $\mathscr{C}_0$.  
\end{prob}

\section*{Acknowledgements}
The first author would like to thank the National Institute of Education, Nanyang Technological University of Singapore, for the generous support of the Nanyang Technological University Research Scholarship.

\end{document}